\documentclass[10pt]{amsart}

\usepackage{xr-hyper}

\usepackage[utf8]{inputenc}         %
\usepackage[T1]{fontenc}            %

\usepackage{a4wide}                 %
\usepackage{amsmath, amssymb}       %
\usepackage{amsthm}                 %
\usepackage{appendix}               %

\usepackage{tikz}                   %
\usepackage{tikz-cd}                %
\tikzcdset{ %
arrow style=tikz
}
\usetikzlibrary{arrows}             %

\usepackage{graphicx}               %
\usepackage{stackrel}               %

\usepackage{color}                  %
\usepackage[normalem]{ulem}         %
\usepackage{enumitem}               %
\setlist[itemize]{leftmargin=*}
\setlist[enumerate]{leftmargin=*,label=\roman*),ref=\roman*)}
\newlist{subenumerate}{enumerate}{2}
\setlist[subenumerate]{leftmargin=*,label=\alph*),ref=\alph*)}

\usepackage{blkarray}               %

\definecolor{darkblue}{rgb}{0,0,0.6} %
\usepackage[ocgcolorlinks,colorlinks=true, citecolor=darkblue, filecolor=darkblue, linkcolor=darkblue, urlcolor=darkblue]{hyperref}         %

\usepackage[alphabetic]{amsrefs}

\usepackage{stix}                   %
\usepackage{bbm}                    %
\usepackage{eucal}                  %

\DeclareFontFamily{T1}{cbgreek}{}
\DeclareFontShape{T1}{cbgreek}{m}{n}{<-6>  grmn0500 <6-7> grmn0600 <7-8> grmn0700 <8-9> grmn0800 <9-10> grmn0900 <10-12> grmn1000 <12-17> grmn1200 <17-> grmn1728}{}
\DeclareSymbolFont{quadratics}{T1}{cbgreek}{m}{n}
\DeclareMathSymbol{\qoppa}{\mathord}{quadratics}{19}
\DeclareMathSymbol{\Qoppa}{\mathord}{quadratics}{21}

\makeatletter
\let\old@font@info\@font@info
\def\@font@info#1{%
\expandafter\ifx\csname\detokenize{#1}\endcsname\relax
  \old@font@info{#1}%
\fi
\expandafter\xdef\csname\detokenize{#1}\endcsname{}%
}
\makeatother

\usepackage{comment}                %
\usepackage{xspace}                 %

\newtheoremstyle{thms}
	{}{}{\itshape}{}{\bfseries }{.}{ }
	{\thmname{#1} \thmnumber{#2}. \thmnote{\bfseries{[#3]}}}
\newtheoremstyle{thms2}
	{}{}{\itshape}{}{\bfseries }{.}{ }
	{}
\newtheoremstyle{ithreethms}
	{}{}{\itshape}{}{\bfseries }{}{ }
	{\thmname{#1} \thmnumber{#2}. \thmnote{\bfseries{[#3]}}}

\newtheoremstyle{name}
	{}{}{\itshape}{}{\bfseries }{.}{ }
	{\thmname{#1}\thmnumber{#2}\thmnote{\bfseries{[#3]}}}

\newtheoremstyle{defs}
	{}{}{\normalfont}{}{\bfseries }{.}{ }
	{\thmname{#1} \thmnumber{#2}. \thmnote{\bfseries{(#3)}}}
\newtheoremstyle{defs2}
	{}{12pt}{\normalfont}{}{\bfseries }{.}{ }
	{\thmname{#1}\thmnumber{#2}. \thmnote{\bfseries{(#3)}}}
\newtheoremstyle{ithreedefs}
	{}{}{\normalfont}{}{\bfseries }{}{ }
	{\thmname{#1} \thmnumber{#2}. \thmnote{\bfseries{(#3)}}}

\newtheoremstyle{rmk}
	{}{}{\normalfont}{}{\itshape }{.}{ }
        {}

\newtheoremstyle{claim}
	{}{}{\normalfont}{}{\itshape}{.}{ }
        {\thmname{#1} \thmnumber{#2}. \thmnote{#3}}

\theoremstyle{thms2}

\theoremstyle{thms2}

\newtheorem*{namedthm}{\namedthmname}
\newcounter{namedthm}

\makeatletter
\newenvironment{named}[1]
{\def\namedthmname{#1}%
\refstepcounter{namedthm}%
\namedthm\def\@currentlabel{#1}}%
{\endnamedthm}
\makeatother

\theoremstyle{rmk}

\theoremstyle{ithreethms}

\theoremstyle{ithreedefs}

\theoremstyle{thms2}

\newtheorem{theoremintrofour}{Theorem}

\newtheorem{corollaryintrofour}[theoremintrofour]{Corollary}
\newtheorem{propositionintrofour}[theoremintrofour]{Proposition}

\swapnumbers         %

\newcounter{rthree}
\theoremstyle{defs}
\newtheorem{definition-r-three}{Definition}[rthree]
\newtheorem{notation-r-three}[definition-r-three]{Notation}
\newtheorem{remark-r-three}[definition-r-three]{Remark}
\newtheorem{example-r-three}[definition-r-three]{Example}

\theoremstyle{thms}
\newtheorem{proposition-r-three}[definition-r-three]{Proposition}
\newtheorem{corollary-r-three}[definition-r-three]{Corollary}
\newtheorem{theorem-r-three}[definition-r-three]{Theorem}

\newcounter{rfour}
\theoremstyle{defs}
\newtheorem{definition-r-four}{Definition}[rthree]
\newtheorem{notation-r-four}[definition-r-four]{Notation}
\newtheorem{remark-r-four}[definition-r-four]{Remark}
\newtheorem{example-r-four}[definition-r-four]{Example}

\theoremstyle{thms}
\newtheorem{proposition-r-four}[definition-r-four]{Proposition}
\newtheorem{corollary-r-four}[definition-r-four]{Corollary}
\newtheorem{theorem-r-four}[definition-r-four]{Theorem}

\theoremstyle{thms}
\newtheorem{proposition}{Proposition}[subsection]
\newtheorem{theorem}[proposition]{Theorem}
\newtheorem{lemma}[proposition]{Lemma}
\newtheorem{corollary}[proposition]{Corollary}

\theoremstyle{defs}
\newtheorem{definition}[proposition]{Definition}%
\newtheorem{notation}[proposition]{Notation}
\newtheorem{variant}[proposition]{Variant}
\newtheorem{construction}[proposition]{Construction}
\newtheorem{example}[proposition]{Example}
\newtheorem{examples}[proposition]{Examples}
\newtheorem{remark}[proposition]{Remark}
\newtheorem{remarks}[proposition]{Remarks}

\newtheorem{assumption}[proposition]{Assumption}

\theoremstyle{defs2}

\theoremstyle{rmk}

\theoremstyle{claim}

\newcommand{\defi}[1]{\emph{#1}}                 %

\newcommand{\spacefont}{\mathcal}                %

\newcommand{\hrar}{\hookrightarrow}

\newcommand{\wtl}{\widetilde}

\newcommand{\ovl}{\overline}

\newcommand{\adj}{\mathbin{%
\begin{tikzpicture}[baseline,thick] 
\coordinate (source) at (0ex,.5ex);
\coordinate (target) at (3ex,.5ex);
\draw[->] ([yshift=1ex]source) -- ([yshift=1ex]target); 
\draw[->] ([yshift=-.5ex]target) -- ([yshift=-.5ex]source);
\node at (1.5ex,.8ex) {$\scriptscriptstyle \perp$};
\end{tikzpicture}%
}}

\newcommand{\locadj}{\mathbin{%
    \begin{tikzpicture}[baseline=-0.8ex]
	\coordinate (source) at (0ex,0ex);
	\coordinate (target) at (3ex,0ex);
	\draw[-stealth] ([yshift=0.7ex]source) to ([yshift=0.7ex]target);
    \draw[left hook-stealth] ([yshift=-0.85 ex]target) to ([yshift=-0.85 ex]source);
	\node at (1.5ex,0  ex) {$\scriptscriptstyle \perp$};
    \end{tikzpicture}
}}

\newcommand{\locadjL}[1]{\mathbin{%
    \begin{tikzpicture}[baseline=-0.8ex]
	\coordinate (source) at (0ex,0ex);
	\coordinate (target) at (3ex,0ex);
	\draw[-stealth] ([yshift=0.7ex]source) to node[above]{\scriptsize $#1$}([yshift=0.7ex]target);
    \draw[left hook-stealth] ([yshift=-0.85 ex]target) to ([yshift=-0.85 ex]source);
	\node at (1.5ex,0  ex) {$\scriptscriptstyle \perp$};
    \end{tikzpicture}
}}

\newcommand{\cocolon}{\nobreak \mskip6mu plus1mu \mathpunct{}\nonscript\mkern-\thinmuskip {:}\mskip2mu \relax}

\makeatletter
\newcommand{\dovl}[1]{\overline{\dbl@overline{#1}}}
\newcommand{\dbl@overline}[1]{\mathpalette\dbl@@overline{#1}}
\newcommand{\dbl@@overline}[2]{%
  \begingroup
  \sbox\z@{$\m@th#1\overline{#2}$}%
  \ht\z@=\dimexpr\ht\z@-2\dbl@adjust{#1}\relax
  \box\z@
  \ifx#1\scriptstyle\kern-\scriptspace\else
  \ifx#1\scriptscriptstyle\kern-\scriptspace\fi\fi
  \endgroup
}
\newcommand{\dbl@adjust}[1]{%
  \fontdimen8
  \ifx#1\displaystyle\textfont\else
  \ifx#1\textstyle\textfont\else
  \ifx#1\scriptstyle\scriptfont\else
  \scriptscriptfont\fi\fi\fi 3
}
\makeatother

\newcommand{\alp}{\alpha}
\newcommand{\bet}{\beta}
\newcommand{\gam}{\gamma}
\newcommand{\Del}{\Delta}

\newcommand{\Lam}{\Lambda}
\newcommand{\Om}{\Omega}

\newcommand{\Sig}{\Sigma}
\newcommand{\sig}{\sigma}
\newcommand{\vphi}{\varphi}
\newcommand{\eps}{\epsilon}

\newcommand{\NN}{\mathbb{N}}               %
\newcommand{\RR}{\mathbb{R}}               %
\newcommand{\CC}{\mathbb{C}}               %
\newcommand{\ZZ}{\mathbb{Z}}               %
\newcommand{\QQ}{\mathbb{Q}}               %
\newcommand{\FF}{\mathbb{F}}               %

\newcommand{\im}{\mathrm{im}}              %
\newcommand{\coker}{\mathrm{coker}}        %

\newcommand{\fp}{\mathfrak{p}}             %

\newcommand{\Hom}{\operatorname{Hom}}            %
\newcommand{\Aut}{\operatorname{Aut}}            %

\newcommand{\fib}{\operatorname{fib}}            %
\newcommand{\cof}{\operatorname{cof}}            %
\newcommand{\coeq}{\operatorname{coeq}}		 %

\newcommand{\colim}{\mathop{\mathrm{colim}}}     %

\newcommand{\Map}{\operatorname{Map}}            %

\newcommand{\grpcr}{\iota}                       %
\newcommand{\core}{\mathrm{Cr}}                  %
\newcommand{\dec}{\mathrm{dec}}                  %

\newcommand{\Pro}{\operatorname{Pro}}            %
\newcommand{\Ind}{\operatorname{Ind}}            %
\newcommand{\Tate}{\operatorname{Tate}}          %
\newcommand{\Latt}{\operatorname{Latt}}          %

\newcommand{\op}{^\mathrm{op}}                   %
\newcommand{\mop}{\mathrm{op}}                   %

\newcommand{\Ar}{\operatorname{Ar}}              %
\newcommand{\Twar}{\operatorname{TwAr}^{\mathrm{r}}}          %

\newcommand{\cp}{\omega}                         %
\newcommand{\perf}{\mathrm{p}}                   %

\newcommand{\Mon}{\operatorname{Mon}}            %
\newcommand{\Grp}{\operatorname{Grp}}            %
\newcommand{\grp}{\mathrm{grp}}                  %

\newcommand{\pr}{\mathrm{pr}}                %
\newcommand{\id}{\mathrm{id}}                %
\newcommand{\const}{\mathrm{const}}          %

\newcommand{\Spa}{{\mathcal Sp}}             %
\newcommand{\Spaf}{\Spa^{\free}}             %
\newcommand{\Sps}{\mathcal S}                %
\newcommand{\cwedge}{{\scriptscriptstyle\wedge}} %
\renewcommand{\hom}{\operatorname{hom}}        %
\newcommand{\map}{\operatorname{hom}}          %

\renewcommand{\SS}{\mathbb{S}}                 %

\newcommand{\h}{\mathrm{h}}                    %
\newcommand{\Ct}{\mathrm{C_2}}                 %
\newcommand{\hC}{{\h\Ct}}                      %
\newcommand{\tC}{{\mathrm{t}\Ct}}              %
\newcommand{\geofix}{{\varphi\Ct}}             %
\newcommand{\gC}{{\mathrm{g}\Ct}}              %
\newcommand{\gCt}{\gC}                         %

\newcommand{\Th}{\mathrm{Th}}		       %

\newcommand{\N}{\mathrm{N}}                    %

\newcommand{\Ab}{\mathcal{A}b}                   %
\newcommand{\Mack}{\mathcal{M}a_{\Ct}}                   %

\newcommand{\one}{\mathbf{1}} 			 %

\newcommand{\Mod}{\operatorname{Mod}}            %

\newcommand{\fg}{{\mathrm{f}}}                   %

\newcommand{\Proj}{\mathrm{Proj}}                %

\newcommand{\Modcp}{\Mod^{\cp}}                  %

\newcommand{\Modp}[1]{\Modcp({#1})}              %
\newcommand{\Modf}[1]{\mathrm{Mod}^\fg({#1})}    %

\newcommand{\gl}{\mathrm{gl}}                    %

\newcommand{\Ering}{\mathrm{E}}			 %
\newcommand{\Ek}[1]{{\Ering_{#1}}}               %
\newcommand{\Einf}{\Ek{\infty}}          	 %
\newcommand{\Eone}{\Ek{1}}               	 %

\newcommand{\Alg}{\mathrm{Alg}}                  %

\newcommand{\GEM}{\mathrm{H}}                    %

\newcommand{\Dfree}{{\mathcal D}^\fg}            %
\newcommand{\Dperf}{{\mathcal D}^\perf}          %
\newcommand{\Perf}{\Dperf}                       %

\newcommand{\free}{\mathrm{f}}                   %

\newcommand{\Fr}{\mathrm{Fr}}                    %

\newcommand{\Com}{\mathrm{Com}}       %
\newcommand{\Day}{\mathrm{Day}}       %
\newcommand{\Sm}{\mathrm{Sm}}			%
\newcommand{\MCom}{\mathrm{MCom}}     %
\newcommand{\CAlg}{\mathrm{CAlg}}     %
\newcommand{\Op}{\mathrm{Op}_\infty}  %

\newcommand{\Mul}{\mathrm{Mul}}       %

\newcommand{\<}{\langle}			%
\renewcommand{\>}{\rangle}			%

\newcommand{\act}{\mathrm{act}}			%
\newcommand{\mon}{\mathrm{mon}}			%

\newcommand{\Cat}{\mathrm{Cat}_\infty}                  %
\newcommand{\PrL}{\mathrm{Pr}^{\mathrm{L}}}  %

\newcommand{\Catx}{\Cat^\ex}     %
\newcommand{\Catxi}{\mathrm{Cat}^{\ex,\natural}_{\infty}}                 %

\newcommand{\Cath}{\mathrm{Cat}^\mathrm{h}_\infty}      %
\newcommand{\Catps}{\mathrm{Cat}^\mathrm{ps}_\infty}    %
\newcommand{\Catp}{\mathrm{Cat}^{\mathrm p}_\infty}     %
\newcommand{\Catpi}{\mathrm{Cat}^{\mathrm p,\natural}_{\infty}}                 %

\newcommand{\Hlgy}{\operatorname{Hlgy}}          %
\newcommand{\Cob}{\mathrm{Cob}}                  %

\newcommand{\Fun}{\operatorname{Fun}}            %
\newcommand{\Nat}{\operatorname{Nat}}            %

\newcommand{\sm}{\mathrm{sm}}			 %

\newcommand{\ex}{\mathrm{ex}}                    %
\newcommand{\fin}{\mathrm{fin}}                  %
\newcommand{\red}{\mathrm{red}}                  %

\newcommand{\add}{\mathrm{add}}                  %
\newcommand{\vloc}{\mathrm{vloc}}                %
\newcommand{\kloc}{\mathrm{kloc}}                %
\newcommand{\brd}{\bord}                %
\newcommand{\binv}{\mathrm{binv}}

\newcommand{\Funadd}{\Fun^\add}                  %

\newcommand{\Funvloc}{\Fun^\vloc}                  %
\newcommand{\Funkloc}{\Fun^\kloc}                  %
\newcommand{\Funbord}{\Fun^\brd}                  %

\newcommand{\dM}{\mathcal{M}}                    %

\newcommand{\Funx}{\operatorname{Fun^{ex}}}      %
\newcommand{\Funq}{\operatorname{Fun^q}}         %
\newcommand{\Funs}{\operatorname{Fun^s}}         %
\newcommand{\Funh}{\operatorname{Fun^h}}         %
\newcommand{\Pairings}{\operatorname{Pair}}      %

\newcommand{\Motp}[1][]{{\if\relax\detokenize{#1}\mathrm{Mot^p}\relax\else\mathrm{Mot}_{#1}^{\mathrm{p}}\fi}}  %
\newcommand{\Motpun}[1][]{\mathrm{Mot^p_{un\if\relax\detokenize{#1}\relax\else{,}#1\fi}}}           %
\newcommand{\Motkar}[1][]{{\if\relax\detokenize{#1}\mathrm{Mot^{k}}\relax\else\mathrm{Mot^k_{#1}}\fi}}  %
\newcommand{\Motpbord}[1][]{{\if\relax\detokenize{#1}\mathrm{Mot^{pw}}\relax\else\mathrm{Mot}_{#1}^{\mathrm{pw}}\fi}} %
\newcommand{\Mot}[1][]{{\if\relax\detokenize{#1}\mathrm{Mot}\relax\else\mathrm{Mot_{#1}}\fi}}            %

\newcommand{\MotpE}[2]{\mathrm{Mot^p_{{#2}\if\relax\detokenize{#1}\relax\else{,}#1\fi}}}           %

\newcommand{\MotkwE}[2]{\mathrm{Mot^{kw}_{{#2}\if\relax\detokenize{#1}\relax\else{,}#1\fi}}}           %

\newcommand{\Motwun}[1][]{\mathrm{Mot^w_{un\if\relax\detokenize{#1}\relax\else{,}#1\fi}}}           %
\newcommand{\MotwE}[2]{\mathrm{Mot^w_{{#2}\if\relax\detokenize{#1}\relax\else{,}#1\fi}}}           %

\newcommand{\MotpwE}[2]{\mathrm{Mot^{pw}_{{#2}\if\relax\detokenize{#1}\relax\else{,}#1\fi}}}           %

\newcommand{\Motw}[1][]{{\if\relax\detokenize{#1}\mathrm{Mot^w}\relax\else\mathrm{Mot^w_{#1}}\fi}}  %

\newcommand{\Motpw}[1][]{{\if\relax\detokenize{#1}\mathrm{Mot^{pw}}\relax\else\mathrm{Mot^{pw}_{#1}}\fi}}  %

\newcommand{\Motkw}[1][]{{\if\relax\detokenize{#1}\mathrm{Mot^{kw}}\relax\else\mathrm{Mot^{kw}_{#1}}\fi}}  %

\newcommand{\MotkarE}[2]{\mathrm{Mot^k_{{#2}\if\relax\detokenize{#1}\relax\else{,}#1\fi}}}           %

\newcommand{\Sp}{\mathrm{Sp}}                    %
\newcommand{\ku}{\mathrm{ku}}
\newcommand{\KU}{\mathrm{KU}}
\newcommand{\SU}{\mathrm{SU}}			%

\newcommand{\MSO}{\mathrm{MSO}}                  %

\newcommand{\LA}{\mathrm{LA}}                    %

\newcommand{\Pic}{\operatorname{Pic}}            %

\newcommand{\Kspace}{\operatorname{\spacefont{K}}} %
\newcommand{\ko}{\mathrm{ko}}                    %
\newcommand{\K}{\operatorname K}                 %
\newcommand{\KK}{\mathbb K}                      %

\newcommand{\topo}{\mathrm{top}}                 %

\newcommand{\KR}{\operatorname{KR}}              %
\newcommand{\KKR}{\mathbb{KR}}			 %

\newcommand{\KO}{\mathrm{KO}}                    %

\newcommand{\Lspace}{\operatorname{\spacefont{L}}}%
\renewcommand{\L}{\operatorname L}               %
\newcommand{\NL}{{\nu\!\L}}         %
\newcommand{\LL}{\mathbb{L}}                     %

\newcommand{\KL}{\mathbb{L}}                     %
\newcommand{\W}{\operatorname W}                 %

\newcommand{\GWspace}{\operatorname{\spacefont{GW}}} %
\newcommand{\GW}{\operatorname{GW}}              %
\newcommand{\KGW}{\mathbb{GW}}                   %

\newcommand{\bord}{\mathrm{bord}}                %
\newcommand{\rbord}{{\mathrm{c}\bord}}		 	 %
\newcommand{\fgt}{\mathrm{fgt}}                  %

\newcommand{\Q}{\operatorname{Q}}                %
\newcommand{\dualQ}{\mathrm d\Q}                 %
\newcommand{\cyl}{\mathrm{cyl}}                  %
\newcommand{\dMet}{\mathrm{d}\hMet}              %
\newcommand{\hMet}{\mathrm{Null}}                %

\newcommand{\vsign}{\sigma^{\vis}}         	 %
\newcommand{\tr}{\mathrm{tr}}   		 %

\newcommand{\fpm}{\lambda}			 %

\newcommand{\FPM}{\mathrm {FPM}}		%
\newcommand{\FPMp}{\FPM^{\mathrm p}}		%
\newcommand{\utens}{\otimes_\ZZ}  		%
\newcommand{\boxutens}{\boxtimes_\ZZ}  		%
\newcommand{\Ringp}{\mathrm{Ring}^{\mathrm{p}}} %
\newcommand{\mack}{\mathrm{mack}}		%

\newcommand{\Nmack}{\mathrm{N}^\mack}		%

\newcommand{\B}{\mathrm{B}}                      %

\newcommand{\Bil}{\mathrm{B}}               %
\newcommand{\Lin}{\Lambda}               %

\newcommand{\Dual}{\mathrm{D}}              %

\newcommand{\ev}{\mathrm{ev}}                    %

\newcommand{\qshift}[1]{^{[#1]}}                 %

\newcommand{\Hyp}{\operatorname{Hyp}}            %
\newcommand{\ghyp}{{\mathrm{g}\hyp}}             %
\newcommand{\Met}{\operatorname{Met}}            %
\newcommand{\met}{\mathrm{met}}                  %
\newcommand{\cconorm}{N^c}                      %
\newcommand{\ccconorm}{N^{cc}}                  %

\newcommand{\hyp}{\mathrm{hyp}}                  %

\newcommand{\cl}{\mathrm{cl}}                    %

\newcommand{\sym}{\mathrm{s}}                    %
\newcommand{\s}{\mathrm{s}}                      %
\newcommand{\vis}{\mathrm{v}}                    %
\newcommand{\qdr}{\mathrm{q}}                     %
\newcommand{\f}{\mathrm{f}}                      %
\newcommand{\g}{\mathrm{g}}                      %
\newcommand{\gs}{\mathrm{gs}}                    %
\newcommand{\gev}{\mathrm{ge}}                   %
\newcommand{\tate}{\mathrm{t}}                   %
\newcommand{\fr}{\mathrm{f}}					  %

\newcommand{\burn}{\mathrm{b}}                   %

\newcommand{\uni}{\mathrm{u}}                    %

\newcommand{\Poinc}{\mathrm{Pn}}                 %
\newcommand{\catforms}{\mathrm{He}}              %

\newcommand{\spsforms}{\mathrm{Fm}}              %
\newcommand{\R}{\mathcal{R}}			%
\newcommand{\alg}{\mathrm{alg}} %

\newcommand{\Exc}{\mathrm{Exc}} %

\newcommand{\md}{\mathrm{mod}}             %

\newcommand{\cN}{\mathcal{N}}      %
\newcommand{\equ}{\mathrm{eq}}     %
\newcommand{\ass}{\mathrm{ass}}    %

\newcommand{\TotFib}{\mathrm{TotFib}}    %
 
\newcommand{\cotower}[1]{\langle #1 \rangle}   %

\newcommand{\cR}{\mathcal{R}}      %

\newcommand{\cP}{\mathcal{P}}              %
\newcommand{\T}{\mathcal{T}}               %

\newcommand{\I}{\mathcal{I}}               %
\newcommand{\J}{\mathcal{J}}               %
\newcommand{\cK}{\mathcal{K}}              %

\newcommand{\A}{\mathcal{A}}               %
\newcommand{\cB}{\mathcal{B}}              %
\newcommand{\E}{\mathcal{E}}               %

\newcommand{\M}{\mathcal{M}}               %

\newcommand{\cO}{\mathcal{O}}              %

\newcommand{\C}{\mathcal C}                %
\newcommand{\D}{\mathcal{D}}               %
\newcommand{\Dtwo}{\mathcal{D}'}           %

\newcommand{\cA}{\mathcal{A}}              %

\newcommand{\F}{\mathcal{F}}               %
\newcommand{\G}{\mathcal{G}}               %

\newcommand{\x}{x}                         %
\newcommand{\y}{y} 			   			   %
\newcommand{\z}{z}                         %
\newcommand{\cob}{w}                       %

\newcommand{\QF}{\Qoppa}                   %
\newcommand{\QFtwo}{{\QF'}}                %
\newcommand{\QFD}{\Phi}                    %
\newcommand{\QFDtwo}{{\QFD'}}              %
\newcommand{\QFE}{\Psi}                    %

\newcommand{\GG}{\mathbb{G}}               %

\newcommand{\Lag}{\mathcal{L}}             %

\makeatletter

\renewcommand{\tocsection}[3]{%
\indentlabel{\@ifnotempty{#2}{\parbox[b]{3ex}{\bfseries\ignorespaces#1 #2}}}\bfseries#3} 

\renewcommand{\tocsubsection}[3]{%
\indentlabel{\@ifnotempty{#2}{\hspace{1.6em}\parbox[b]{5ex}{\ignorespaces#1 #2}}}#3}

\renewcommand{\tocsubsubsection}[3]{%
\indentlabel{\@ifnotempty{#2}{\hspace{3.9em}\parbox[b]{5ex}{\ignorespaces#1 #2}}}#3}

\makeatother

\DeclareRobustCommand{\SkipTocEntry}[5]{} 

\newcommand{\introsubsection}[1]{\addtocontents{toc}{\SkipTocEntry}\subsection*{#1}}
\newcommand{\bibsubsubsection}[1]{\addtocontents{toc}{\SkipTocEntry}\subsubsection*{#1}}

\BibSpec{collection.article}{%
    +{}  {\PrintAuthors}                {author}
    +{,} { \textit}                     {title}
    +{.} { }                            {part}
    +{:} { \textit}                     {subtitle}
    +{,} { \PrintContributions}         {contribution}
    +{,} { \PrintConference}            {conference}
    +{}  {\PrintBook}                   {book}
    +{,} { }                            {booktitle}
    +{,} { }                            {series}
    +{,} { \voltext}                    {volume}
    +{,} { }                            {publisher}
    +{,} { }                            {organization}
    +{,} { }                            {address}
    +{,} { \PrintDateB}                 {date}
    +{,} { \eprintpages}                {pages}
    +{,} { }                            {status}
    +{,} { available at \eprint}        {eprint}
    +{}  { \parenthesize}               {language}
    +{}  { \PrintTranslation}           {translation}
    +{;} { \PrintReprint}               {reprint}
    +{.} { }                            {note}
    +{.} {}                             {transition}
    +{}  {\SentenceSpace \PrintReviews} {review}
}

\BibSpec{article}{%
    +{}  {\PrintAuthors}                {author}
    +{.} { \textit}                     {title}
    +{.} { }                            {part}
    +{:} { \textit}                     {subtitle}
    +{,} { \PrintContributions}         {contribution}
    +{.} { \PrintPartials}              {partial}
    +{,} { }                            {journal}
    +{}  { \textbf}                     {volume}
    +{} { \PrintDate}                  {date}
    +{,} { \issuetext}                  {number}
    +{,} { \eprintpages}                {pages}
    +{,} { }                            {status}
    +{,} { available at \eprint}        {eprint}
    +{}  { \parenthesize}               {language}
    +{}  { \PrintTranslation}           {translation}
    +{;} { \PrintReprint}               {reprint}
    +{.} { }                            {note}
    +{.} {}                             {transition}
    +{}  {\SentenceSpace \PrintReviews} {review}
}

\newcommand{\reffour}[1]{\ref{#1}}
\newcommand{\eqreffour}[1]{\eqref{#1}}
\newcommand{\reffouritem}[1]{\ref{#1}}

\newcommand{\refone}[1]{\cite{Part-one}.\ref{I-#1}}

\newcommand{\refoneitem}[1]{\ref{I-#1}}

\newcommand{\reftwo}[1]{\cite{Part-two}.\ref{II-#1}}

\newcommand{\reftwoitem}[1]{\ref{II-#1}}

\newcommand{\refthree}[1]{\cite{Part-three}.\ref{III-#1}}

\newcommand{\refthreeitem}[1]{\ref{III-#1}}

\newcommand{\reffive}[1]{\cite{Part-five}.\ref{V-#1}}

\newcommand{\paperone}{Paper~\cite{Part-one}\xspace}
\newcommand{\papertwo}{Paper~\cite{Part-two}\xspace}
\newcommand{\paperthree}{Paper~\cite{Part-three}\xspace}

\newcommand{\paperfive}{Paper~\cite{Part-five}\xspace}

\newcommand{\refundefined}[1]{\textcolor{red}{Undefined ref}}

\title[Hermitian K-theory for stable $\infty$-categories IV: Karoubi-Grothendieck-Witt theory]{Hermitian K-theory for stable $\infty$-categories IV:\\
Karoubi-Grothendieck-Witt theory}

\author[Calmès]{Baptiste Calmès}
\address{Université d'Artois, Laboratoire de Mathématiques de Lens, Lens, France}
\email{baptiste.calmes@univ-artois.fr}

\author[Dotto]{Emanuele Dotto}
\address{University of Warwick, Mathematics Institute, Coventry, United Kingdom}
\email{emanuele.dotto@warwick.ac.uk}

\author[Harpaz]{Yonatan Harpaz}
\address{Université Paris Cité, Sorbonne Université, Paris, France}
\email{harpaz@imj-prg.fr}

\author[Hebestreit]{Fabian Hebestreit}
\address{Universität Bielefeld, Fakultät für Mathematik, Bielefeld, Germany}
\email{hebestreit@math.uni-bielefeld.de}

\author[Land]{Markus Land}
\address{Johannes Gutenberg University Mainz, Institute of Mathematics, Mainz, Germany}
\email{mland@uni-mainz.de}

\author[Moi]{Kristian Moi}
\address{KTH, Institutionen för matematik, Stockholm, Sweden}
\email{kristian.moi@gmail.com}

\author[Nardin]{Denis Nardin}
\address{Universität Regensburg, Mathematisches Institut, Regensburg, Germany}

\author[Nikolaus]{Thomas Nikolaus}
\address{WWU Münster, Mathematisches Institut, Münster, Germany}
\email{nikolaus@uni-muenster.de}

\author[Steimle]{Wolfgang Steimle}
\address{Universität Augsburg, Institut für Mathematik, Augsburg, Germany}
\email{wolfgang.steimle@math.uni-augsburg.de}

\dedicatory{To Andrew Ranicki.}

\date{\today}

\begin{document}

\begin{abstract}
In analogy with the passage from connective to non-connective algebraic $\K$-theory, Schlichting put forward the Karoubi-Grothendieck-Witt spectrum as an approximation to Grothendieck-Witt theory by a functor invariant under idempotent completion. In this paper, we extend his construction to the setting of Poincaré $\infty$-categories and characterise the resulting functor by the expected universal property.

Our main result places these Karoubi-Grothendieck-Witt-spectra in an exact sequence involving the ultimately decorated $\L$-spectrum and homotopy orbits of non-connective $\K$-theory. This confirms a conjecture of Williams which had remained open even for rings in which $2$ is invertible, and allows us to deduce universal characterisations for ultimately decorated $\L$-theory. We obtain these results by unifying the Bass-Heller-Swan and Ranicki-Shaneson splittings into a general decomposition theorem that yields, in particular, formulas for Grothendieck-Witt groups of Laurent polynomial rings.

Finally, we construct multiplicative structures on both versions of our Grothendieck-Witt and $\L$-functors and refine their universal properties to multiplicative ones. Along the way, we show that the original versions are even idempotent with respect to Day convolution.

\end{abstract}

\maketitle
\tableofcontents

\section*{Introduction}
This paper is the fourth in a series on the foundations of hermitian $\K$-theory using stable $\infty$-categories. In this instalment, we introduce Karoubi-localising versions of hermitian $\K$-functors and study multiplicative structures on both the original and the Karoubi-localising version.

We work with Lurie's notion of Poincaré $\infty$-categories, that is, pairs $(\C,\QF)$ where $\C$ is a stable $\infty$-category and $\QF\colon \C\op \to \Spa$ is a Poincaré structure on $\C$, that is, 
a perfect quadratic functor encoding a notion of unimodular (or Poincaré) forms in $\C$, the space of which we denote by $\Poinc(\C,\QF)$.
The fundamental invariant we consider is %
the Grothendieck-Witt functor $\GW\colon \Catp \to \Sp$, introduced in \papertwo, where the source is the $\infty$-category of Poincaré $\infty$-categories and the target the $\infty$-category of spectra.
For a ring $R$ equipped with a form parameter we in particular produced a Poincar\'e structure on the perfect derived category $\Dperf(R)$ whose Grothendieck-Witt spectrum recovers the classical Grothendieck-Witt space defined by Karoubi in \cites{karoubi-quillen,Karoubi-Le-theoreme-fondamental} as its infinite loop space, see \cite{comparison}. 

Applications of our theory given in \papertwo and \paperthree include a proof of Karoubi's periodicity conjecture, the finite generation of Grothendieck-Witt-groups for rings of integers in number fields, and explicit computations for $\ZZ$.
These computations are also a key input for forthcoming work on the stable cohomology of orthogonal and symplectic groups over the integers \cite{HLN-cohomology2,HLN-cohomologyp}. 
To obtain these applications,
we established in \papertwo what we termed the fundamental fibre sequence
\[
\K(\C,\QF)_\hC \xrightarrow{\hyp} \GW(\C,\QF) \xrightarrow{\bord} \L(\C,\QF),
\]
expressing the Grothendieck-Witt functor as an extension of the simpler algebraic $\K$- and $\L$-spectra, where the $\Ct$-action on $\K(\C)$ is induced by the duality associated to $\QF$. Since the functors $\K$ and $\L$ are Verdier-localising the same then follows for $\GW$, that is, all three take exact sequences in $\Catp$ %
to exact sequences of spectra. In fact, we showed in \papertwo that $\GW$ is the initial Verdier-localising functor under 
$\SS[\Poinc]$, %
while $\L$ is the initial such which is furthermore bordism invariant, that is, vanishes on hyperbolic Poincaré $\infty$-categories.\\

For many applications to rings and schemes, being Verdier-localising is, however, not quite sufficient. For example, on the level of perfect derived $\infty$-categories, localisations of rings yield Verdier quotients, and hence extend to exact sequences,
only up to idempotent completion. As $\K,\GW$, and $\L$ are not invariant under idempotent completion they consequently fail to satisfy Zariski descent. 
A similar problem occurs 
in the Bass-Heller-Swan and Ranicki-Shaneson decompositions concerning the $\K$- and $\L$-groups of Laurent polynomial rings, or in excision sequences for quotients by ideals.

In algebraic $\K$-theory these points were resolved with the introduction of non-connective $\K$-theory $\KK$, which is both Verdier-localising and invariant under idempotent completion, a property we call Karoubi-localising. 
Historically, there are two approaches to this extension: One, originally due to Bass and in spectral form to Thomason and Trobaugh \cites{Bass,thomason-trobaugh}, is based on the behaviour of $\K$-spectra of rings of Laurent polynomials. The other one, based on suspension rings, was suggested by Karoubi and upgraded to the spectrum level by Gersten and Wagoner in \cite{karoubi-foncteurs,gersten-spectrum, wagoner-delooping}. 
The suspension ring approach was extended by Karoubi to Grothendieck-Witt spectra, yielding what is nowadays called the Karoubi-Grothendieck-Witt functor $\KGW$,
and which indeed satisfies Zariski descent on schemes by work of Schlichting \cites{schlichting-mv,schlichting-derived}. On the other hand, Bass' construction was adapted to $\L$-spectra by Ranicki, leading to ultimately decorated $\L$-theory $\L^{\cotower{-\infty}}$, which admits well-behaved splittings for Laurent polynomial rings and %
features in %
surgery theoretic approaches to Borel's rigidity conjecture for aspherical manifolds, see 
\cite{Lueck-FJbook}. 

Unlike the case of $\KK$ and $\KGW$, which agree with $\K$ and $\GW$ respectively in positive degrees, %
the functor $\L^{\cotower{-\infty}}$ satisfies the same periodicity properties as $\L$
and 
differs from it in arbitrary degrees. In addition, while $\KK$ was recognised as universal inside Karoubi-localising functors 
on stable $\infty$-categories in~\cite{BGT}, no universal property was ever established for $\KGW$ or $\L^{\cotower{-\infty}}$.
Despite their differences, Williams conjectured in \cite{williams-quadratic} that for group rings in which 2 is invertible, $\KK$, $\KGW$, and $\L^{\cotower{-\infty}}$ should be related by an exact sequence analogous to the fundamental fibre sequence connecting the $\K$-, $\GW$- and $\L$-spectra above.

In the present paper, we extend the constructions of $\KGW$ and $\L^{\cotower{-\infty}}$ to Poincaré $\infty$-categories, characterise both by universal properties and establish Williams' conjecture in the full generality of Poincar\'e $\infty$-categories:

\begin{named}{Main theorem}
\label{theorem:pfourmaintheorem}%
For a Poincaré $\infty$-category $(\C,\QF)$ there is a natural exact sequence
\[
\KK(\C,\QF)_\hC \xrightarrow{\hyp} \KGW(\C,\QF) \xrightarrow{\bord} \L^{\cotower{-\infty}}(\C,\QF),
\]
and the functors $\KGW$ and $\L^{\cotower{-\infty}}$ enjoy the following universal properties:
\begin{enumerate}
\item
\label{item:main-theorem-i}%
they are the initial Karoubi-localising (and bordism invariant) functors $\Catp\to \Sp$ under $\SS[\Poinc]$,
respectively and
\item
\label{item:main-theorem-ii}%
they underlie the initial lax symmetric monoidal Karoubi-localising (and bordism invariant) functors $\Catp\to \Sp$, respectively.
\end{enumerate}
\end{named}

By the general results of \papertwo, 
the fibre sequence of the \reffour{theorem:pfourmaintheorem} also implies that $\L^{\langle-\infty\rangle}$ is the bordification of $\KGW$.
Moreover, given the universal properties of $\GW$ and $\L$ established in \papertwo and the general relation between bordism invariance and Karoubi-localisation, point \reffouritem{item:main-theorem-i} above shows that $\KGW$ and $\L^{\cotower{-\infty}}$ are also the initial Karoubi-localising functors under $\GW$ and $\L$, respectively. The above universal properties of $\KGW$ are thus the hermitian analogs of the universal properties of $\KK$ established in the celebrated work of Blumberg--Gepner--Tabuada \cite{BGT}.

Our \reffour{theorem:pfourmaintheorem} can be used to infer new information on both $\KGW$ and $\L^{\cotower{-\infty}}$:
For example, one can deduce that $\KGW$ satisfies the Farrell-Jones conjecture for a given group whenever its $\KK$- and $\L^{\cotower{-\infty}}$-theoretic versions have been established for that group. This has been shown to be the case for a large class of groups by Bartels, Lück, Reich and others \cite{Lueck-FJbook}. 
At the same time, the new universal properties of $\L^{\cotower{-\infty}}$ can be used to construct a trace map from $\L^{\cotower{-\infty}}$ to the geometric fixed points of real topological cyclic homology. We again relegate the development and applications of such hermitian trace methods to forthcoming work \cite{Harpaz-Nikolaus-Shah}, building on \cites{HNS-algtrace, HNRS-algtrace}. 

We define the spectra $\L^{\cotower{-\infty}}(\C,\QF)$ by means of an inductive procedure, extending the construction for rings. Just as Bass' construction of negative $\K$-groups is based on the Bass-Heller-Swan decomposition computing the $\K$-groups of Laurent polynomial rings, the construction of $\L^{\cotower{-\infty}}$ is based on the analogous Ranicki-Shaneson decomposition. As an important ingredient in the above main theorem, we generalise such decomposition theorems to arbitrary Verdier-localising functors and arbitrary Poincar\'e $\infty$-categories, subsuming, extending and unifying these previously analogous, but not directly connected results.

Let us now describe in more detail the contents of this paper.

\introsubsection{Karoubi-localising invariants on Poincaré \texorpdfstring{$\infty$-}{}categories}

Let us begin with our general results about Karoubi-localising approximations. Recall to that end that the functors $\GW$ and $\L$ are Verdier-localising but not Karoubi-localising.
Denoting with superscripts the full subcategories spanned by Verdier and Karoubi-localising functors, we show: 
\begin{theoremintrofour}\
\label{theorem:intro-B}%
In the solid square
\[
\begin{tikzcd}
\Fun^{\kloc}(\Catp,\Sp)\ar[r,shift right=2pt]\ar[d,"\Omega^\infty"'] & \Fun^{\vloc}(\Catp,\Sp)\ar[d,"\Omega^\infty"] \ar[l,dashed, shift right=2pt,"(-)^\mathrm{kar}"'] \\
\Fun^{\kloc}(\Catp,\Grp_{\Einf})\ar[r,shift right=2pt] & \Fun^{\vloc}(\Catp,\Grp_{\Einf}) \ar[l,dashed,shift right=2pt,"(-)^\mathrm{kar}"'] \ ,
\end{tikzcd}
\]
\begin{enumerate}
\item
\label{item:delooping}%
the left vertical functor is an equivalence, and
\item
\label{item:adjoint}%
the horizontal functors admit left adjoints $(-)^\mathrm{kar}$ making the square left adjointable.
\end{enumerate}
In particular, the canonical map $\Omega^\infty \F(-^\natural) \Rightarrow \Omega^\infty \F^{\mathrm{kar}}(-)$ is an equivalence if and only if its source is Karoubi-localising.

\end{theoremintrofour}

Theorem~\reffour{theorem:intro-B} holds equally for $\Catx$ in place of $\Catp$. In this version, \reffouritem{item:delooping} appeared previously in work of Cisinski and Khan \cite{CK}. The proof of Theorem \reffour{theorem:intro-B} also describes the horizontal adjoints and left vertical inverse explicitly: Briefly, both the successive deloopings of a Karoubi-localising $\Einf$-group-valued functor and the adjoint evaluated on a Verdier-localising spectrum-valued functor are described using (a countable version of) the Lattice-Tate sequence studied in \papertwo; it serves as a hermitian substitute of the Ind-Calkin sequence (used, for example, to produce non-connective K-theory in \cite{BGT} and continuous K-theory in \cite{Efimov}), which can itself be thought of as a categorical form of the cone-suspension ring approach of Karoubi and Villamayor. %
We also show that the left adjoint admits an alternative description using decorations, see the discussion after Theorem \reffour{theorem:intro-E} below, in our solution of Williams' conjecture.

Applying the left adjoint in Theorem \reffour{theorem:intro-B}\;\reffouritem{item:adjoint} to the spectrum-valued Verdier-localising functors $\GW$ and $\L$ we obtain Karoubi-localising approximations 
\[
\GW \Rightarrow \KGW \quad \text{ and } \quad \L \Rightarrow \KL
\]
that is, initial Karoubi-localising functors under $\GW$ and $\L$, respectively, and from the general interaction of Karoubification with bordification it follows that $\KGW^\bord \simeq \KL$. The results from \papertwo then show that $\KL$ obeys similar periodicity properties as $\L$ and yield the fundamental exact square
\[
\begin{tikzcd}
\KGW(\C,\QF) \ar[r,"\bord"] \ar[d,"\fgt"] & \KL(\C,\QF) \ar[d,"\Xi"] \\
\KK(\C,\QF)^\hC \ar[r,"\mathrm{can}"] & \KK(\C,\QF)^\tC,
\end{tikzcd}
\]
which is also the isotropy separation square of a genuine $\Ct$-equivariant non-connective real algebraic $\K$-spectrum $\KKR(\C,\QF)$.

To leverage the criterion from Theorem~\reffour{theorem:intro-B}, we consider the behaviour of $\GW$- and $\L$-spectra under Karoubi equivalences, that is, Poincar\'e functors $(f,\eta) \colon (\C,\QF) \rightarrow (\D,\QFD)$ where $f$ is fully faithful and induces an equivalence upon idempotent completion and $\eta\colon \QF \to f^*\QFD$ is an equivalence. 
As a straightforward generalisation of a result of Schlichting's in the setting of differential graded categories on which $2$ acts invertibly \cite{schlichting-derived} there are
exact sequences
\[
\GW(\C,\QF)\longrightarrow \GW(\D,\QFD) \xrightarrow{\fgt} [\K_0(\D)/\K_0(\C)]^\hC \quad \text{and} \quad \L(\C,\QF)\longrightarrow \L(\D,\QFD) \xrightarrow{\Xi} [\K_0(\D)/\K_0(\C)]^\tC.
\]
The left hand one in particular implies that the map $\GW(\C,\QF) \rightarrow \GW(\D,\QFD)$ is an isomorphism on positive homotopy groups and an injection on components, just as for $\K(\C) \rightarrow \K(\D)$; we refer to this result as the hermitian cofinality theorem. The right hand sequence is also familiar in surgery theory, entirely independent from Grothendieck-Witt theory: The associated long exact sequence is the Ranicki-Rothenberg sequence from \cite{Ranicki-lower} in the case of rings or additive categories with duality and symmetric or quadratic Poincar\'e structures.

As a consequence of the hermitian cofinality theorem, we find that $\GWspace(-^\natural)$ is Karoubi-localising, so that Theorem~\reffour{theorem:intro-B} implies that the map $\GWspace(-^\natural) \to \Omega^\infty\KGW(-)$ is an equivalence, but as a consequence of the Ranicki-Rothenberg sequences, the same is not true in $\L$-theory. To control the cofibre of the map $\L \to \KL$ (and of $\GW \rightarrow \KGW$ at the spectrum level), we finally extend the full statement of Schlichting's cofinality theorem \cite{schlichting-derived} to our setup, again removing the assumption that $2$ be invertible in the input:

\begin{propositionintrofour}
\label{proposition:intro-D}%
For any Poincaré $\infty$-category $(\C,\QF)$ the following two squares are cartesian.
\[
\begin{tikzcd}
\GW(\C,\QF) \ar[r] \ar[d,"\fgt"] & \KGW(\C,\QF) \ar[d,"\fgt"] && \L(\C,\QF) \ar[r] \ar[d,"\Xi"] & \KL(\C,\QF) \ar[d,"\Xi"]\\
\K(\C,\QF)^\hC \ar[r] & \KK(\C,\QF)^\hC &&\K(\C,\QF)^\tC \ar[r] & \KK(\C,\QF)^\tC
\end{tikzcd}
\]
\end{propositionintrofour}

In particular, if $\C$ is idempotent complete and $\KK(\C)$ happens to be connective, e.g.\ for $\C = \Dperf(R)$ with $R$ regular noetherian, we find $\GW(\C,\QF) \simeq \KGW(\C,\QF)$ and $\L(\C,\QF) \simeq \KL(\C,\QF)$. Moreover, we generally have
\[
\pi_n\cof\big(\GW(\C,\QF) \to \KGW(\C,\QF)\big)[\tfrac{1}{2}] \cong \left((\KK_n(\C)/\K_n(\C))[\tfrac 1 2]\right)^{\Ct}
\]
for all $n \in \mathbb Z$ with the right hand side being $0$ for $n>0$ and $\KK_n(\C)[\tfrac{1}{2}]^{\Ct}$ for $n<0$, while the map $\L(\C,\QF)[\frac{1}{2}] \to \KL(\C,\QF)[\frac{1}{2}]$ is an equivalence.

The above results have been used to develop Karoubi-Grothendieck-Witt-spectra in the context of schemes, where they enjoy many of the same advantages over Grothendieck-Witt spectra that non-connective $\K$-spectra have over the connective ones (and often the same is then true for $\KL$-spectra). This thread was taken up in \cite{motives} and \cite{HL}, where Zariski- and Nisnevich-descent for these spectra is established, generalising Schlichting's celebrated Mayer-Vietoris principle from \cite{schlichting-mv}, and ultimately motivic spectra (non-$\mathbb A^1$-invariant in the latter case) representing homotopical and genuine symmetric Grothendieck-Witt-theory are constructed, respectively.

\introsubsection{The decoration filtration and Williams' conjecture}

Let us now turn our attention to ultimately decorated $\L$-theory and Williams's conjecture. To explain its context, recall first the Bass-Heller-Swan %
decomposition %
\[
\KK(R[\ZZ]) \simeq \KK(R) \oplus \Sigma \KK(R) \oplus \N\KK(R)  \oplus \N\KK(R),
\]
where $R[\ZZ]$ is the group $R$-algebra of the integers and the first two summands are embedded via the assembly map $\B\ZZ \otimes \KK(R) \to \KK(R[\mathbb Z])$; the last two summands are known as the nil-terms, defined for example as the complimentary summand of $\KK(R)$ in $\KK(R[\mathbb N])$. 

Historically, the fact that in the analogous splitting of $\K(R[\mathbb Z])$ the second summand agrees with $\Sigma \K(R)$ only on connected covers, served as one of the primary motivations for introducing negative $\K$-groups. More precisely, Bass defined the spectrum $\K^{\langle -1 \rangle}(R)$ as the desuspension of this second summand and consequently $\K_{-1}(R)$ as $\pi_{-1}\K^{\langle -1 \rangle}(R)$, with the assembly map $\mathrm B\mathbb Z \otimes \K(R) \rightarrow \K(R[\mathbb Z])$ providing a comparison map $\K(R) \rightarrow \K^{\langle -1 \rangle}(R)$ that is an isomorphism in non-negative degrees. He then iterated this procedure producing spectra $\K^{\langle -n \rangle}(R)$, and $\KK(R) = \colim_{n \in \mathbb N} \K^{\langle -n \rangle}(R)$ was then the first definition of non-connective $\K$-spectra (so that, a-posteriori, $\K^{\langle -n \rangle}(R)$ coincides with $\tau_{\geq -n} \KK(R)$). 

This construction was transported to $\L$-theory to address analogous shortcomings of the Ranicki-Shaneson splitting of $\L(R[\mathbb Z])$; in fact, it simplifies slightly due to the non-occurence of nil-terms. More precisely, if $R$ is equipped with an invertible module with involution $M$ then we can equip the group ring $R[\ZZ]$ with an induced module with involution $M[\ZZ] = \oplus_{n \in \ZZ} M$ extending the involution by swapping the $n$ and $-n$ summands, and writing $\L$ for either symmetric or quadratic $\L$-spectra one puts
\[
\L^{\cotower{-1}}(R,M) \coloneq \SS^{-1} \otimes [\L(R[\ZZ],M[\ZZ])/\L(R,M)].
\]
Iterating this process again one obtains what is traditionally called the decoration sequence
\[
\L(R,M) = \L^{\langle 0\rangle}(R,M) \longrightarrow \L^{\cotower{-1}}(R,M) \longrightarrow \dots \longrightarrow \L^{\langle -n\rangle}(R,M) \longrightarrow \L^{\langle -n-1 \rangle}(R,M) \longrightarrow \dots
\]
whose colimit 
$
\L^{\cotower{-\infty}}(R,M) \coloneq \colim_{n \in \mathbb N} \L^{\cotower{-n}}(R,M)
$
is known as ultimately decorated $\L$-theory. In distinction with the case of algebraic $\K$-theory, the connecting maps in this sequence generally enjoy no particular connectivity properties. 
Despite the different constructions, in \cite{williams-quadratic} Williams predicted, in the context of groups rings over $\ZZ[\tfrac{1}{2}]$, the existence of a fibre sequence
\[
\KK(R,M)_{\hC} \to \KGW(R,M) \to \L^{\cotower{-\infty}}(R,M)
\]
involving the Karoubi-Grothendieck-Witt spectrum defined by Karoubi and Schlichting.
In the terminology of the previous subsection, his conjecture can be reformulated simply as an equivalence $\L^{\cotower{-\infty}}(R,M) \simeq \KL(R,M)$.

Our approach to proving this conjecture in general proceeds by adapting the construction above to produce an endofunctor $\mathcal T$ on $\Fun(\Catp,\Sp)$, equipped with a natural transformation from the identity, and set
\[
\F^{\langle -n \rangle} = \mathcal T^n(\mathcal F(-^\natural)) \quad \text{and} \quad \F^{\langle -\infty \rangle} = \colim_{n \in \mathbb N} \F^{\langle -n \rangle},
\]
for an arbitrary Verdier-localising $\F \colon \Catp \rightarrow \Sp$. When applied to $\L$ and $\K$ one recovers the decoration sequences described above.
But now the construction also applies to Grothendieck-Witt-spectra, where to the best of our knowledge it has not been considered before. As a simultaneous extension of Proposition~\reffour{proposition:intro-D}  and the first part of our \reffour{theorem:pfourmaintheorem} (which concern the cases $n=0$ and $n=\infty$, respectively), we then show:

\begin{theoremintrofour}
\label{theorem:intro-E}%
For any Poincaré $\infty$-category $(\C,\QF)$ and $n\geq 0$ there are exact squares
\[
\begin{tikzcd}
\GW^{\langle -n \rangle}(\C,\QF) \ar[r] \ar[d,"\fgt"] & \KGW(\C,\QF) \ar[d,"\fgt"] && \L^{\langle -n \rangle}(\C,\QF) \ar[r] \ar[d,"\Xi"] & \KL(\C,\QF) \ar[d,"\Xi"] \\
(\tau_{\geq -n} \KK(\C,\QF))^\hC \ar[r] & \KK(\C,\QF)^\hC && (\tau_{\geq -n} \KK(\C,\QF))^\tC \ar[r] & \KK(\C,\QF)^\tC
\end{tikzcd}
\]
and consequently the canonical maps $\GW^{\langle -\infty \rangle}(\C,\QF) \to \KGW(\C, \QF)$ and $\L^{\langle-\infty\rangle}(\C,\QF)\to \KL(\C,\QF)$ are equivalences.
\end{theoremintrofour}

The proof we provide in fact establishes the following generalisation of Williams conjecture: We show that the initial Karoubi-localising functor under $\F$ guaranteed by Theorem~\reffour{theorem:intro-B} 
coincides with
$\F^{\langle -\infty \rangle}$ for \emph{any} Verdier-localising $\F \colon \Catp \rightarrow \Spa$. However, in general the intermediate objects $\F^{\langle -n \rangle}$ do not seem to admit similarly simple descriptions.

\introsubsection{The circle decomposition theorem}
The workhorse result underlying Theorem \reffour{theorem:intro-E} is a unification (and far reaching generalisation) of the Bass-Heller-Swan and Ranicki-Shaneson decompositions, whose proof and discussion forms a substantial part of the present paper. To state it in full, for a functor $\G \colon \Catx \rightarrow \Sp$ and stable $\infty$-category $\D$, we set 
$\mathrm N\mathcal G(\D)$ to be the complementary summand of $\G(\C)$ in $\G(\C_{\mathrm B\mathbb N})$, where $\C_{\mathrm{B}\mathbb{N}}$ denotes the tensor of the stable $\infty$-category $\C$ with the category $\mathrm{B}\mathbb{N}$.
For example, $\Dfree(R)_{\mathrm{B}\mathbb{N}} \simeq \Dfree(R[\NN])$ while $\Dperf(R)_{\mathrm{B}\mathbb{N}}$ and $\Dperf(R[\mathbb{N}])$ still agree up to idempotent completion.
For $\G = \KK$, the nil term $\N\KK$ then coincides with the previously mentioned one when evaluated on perfect derived $\infty$-categories of rings, while for $\G=\K$ we get $\N\K = \tau_{\geq 1} \N\KK$ (which is not quite the complementary summand of $\K(R)$ in $\K(R[\NN])$, but rather the connected cover thereof). %
Replacing $\mathbb N$ by $\mathbb Z$, the operation of tensoring with $\B\ZZ$ can now also be performed on Poincaré $\infty$-categories, and one has for example an equivalence of Poincar\'e $\infty$-categories
\[
\big(\Dfree(R),\QF_M^{\mathrm{g}\lambda}\big)_{\mathrm{B}\mathbb Z}  \simeq \big(\Dfree(R[\ZZ]),\QF_{M[\ZZ]}^{\g\lambda[\ZZ]}\big)
\]
for any form parameter $\lambda=[M_{\Ct} \to Q \to M^{\Ct}]$ on $R$, where $\lambda[\ZZ]$ is the form parameter on $R[\ZZ]$ whose associated module with involution is $M[\ZZ] = \oplus_{n \in \ZZ}M$ as above (i.e.\ the involution exchanges the $n$th and $-n$th summands) and quadratic module given by $Q \oplus \big[\oplus_{n > 0}M\big]$. Note that this transforms the genuine symmetric, even, and quadratic form parameters over $R$ to genuine symmetric, even, and quadratic form parameters over $R[\ZZ]$, respectively. Similarly, one has for any spherical fibration $\xi \colon X \to \Pic(\SS)$ an equivalence of visible Poincar\'e $\infty$-categories
\[
((\Sp_{/X})^\mathrm{f},\QF^\vis_\xi)_{\mathrm B \mathbb Z} \simeq ((\Sp_{/X \times \mathrm S^1})^\mathrm{f},\QF_{\mathrm{pr}^\ast\xi}^{\vis})
\]
where $\mathrm{pr} \colon X \times \mathrm S^1 \to X$ denotes the projection.
We then generally show the following circle decomposition theorem:

\begin{theoremintrofour}
\label{theorem:intro-F}%
Let $\mathcal F \colon \Catp \rightarrow \Sp$ be Verdier-localising and $(\C,\QF)$ a Poincaré $\infty$-category. Then there is a natural equivalence
\[
\F((\C,\QF)_{\mathrm B\mathbb Z}) \simeq \F(\C,\QF) \oplus \Sigma \mathcal F(\C^\natural,\QF) \oplus \mathrm N\F^\hyp(\C).
\]
\end{theoremintrofour}

Under this decomposition the assembly map 
$
\mathrm B\mathbb Z \otimes \F(\C,\QF) \to \F((\C,\QF)_{\mathrm B\mathbb Z})
$
corresponds to the evident map
\[
\mathcal F(\C,\QF) \oplus \Sigma \mathcal F(\C,\QF) \longrightarrow \mathcal F(\C,\QF) \oplus \Sigma \mathcal F(\C^\natural,\QF) \oplus\mathrm N\F^\hyp(\C)
\]
in the first two components. %
In particular, if $\F$ is bordism invariant and $(\C,\QF)$ is idempotent complete then the assembly map
is an equivalence.

Let us now specialise the circle decomposition theorem to a few key examples.
First, for $\K \colon \Catp \rightarrow \Sp^{\mathrm{B}\Ct}$, using $\K^\hyp \simeq \Ct\otimes \K$ as functors to $\Ct$-spectra, we obtain a splitting
\[
\K((\C,\QF)_{\mathrm B\mathbb Z}) \simeq \K(\C,\QF) \oplus \Sigma \K(\C^\natural,\QF) \oplus [\Ct \otimes \mathrm{NK}(\C)],
\]
of $\Ct$-spectra. Disregarding involutions such a splitting was recently established by Saunier \cite{saunier} (starting with an arbitrary stable but idempotent complete $\D$; our proof, which is %
different from Saunier's, in fact works verbatim for all stable $\D$). If we take $\KK$ instead of $\K$ the splitting simplifies to
\[
\KK((\C,\QF)_{\mathrm B\mathbb Z}) \simeq \KK(\C,\QF) \oplus \Sigma \KK(\C,\QF) \oplus [\Ct \otimes \N\KK(\C)]
\]
as one can ignore idempotent completions. %
This recovers the classical Bass-Heller-Swan decomposition, as well as its $\mathrm A$-theoretic analogue from \cites{fundamental-spaces, fundamental-spaces2}, including an identification of the involution on the right hand side.

Second, for $\L \colon \Catp \rightarrow \Sp$, and $\KL$ using $\L^\hyp \simeq 0 \simeq \KL^\hyp$ by bordism invariance, the circle decomposition provides splittings
\[
\L((\C,\QF)_{\mathrm B\mathbb Z}) \simeq \L(\C,\QF) \oplus \Sigma \L(\C^\natural,\QF) \quad \text{ and } \quad \KL((\C,\QF)_{\mathrm B\mathbb Z}) \simeq \KL(\C,\QF) \oplus \Sigma \KL(\C,\QF).
\]
This statement generalises the Ranicki-Shaneson decomposition 
for %
quadratic and genuine symmetric $\L$-theory of ordinary rings from \cite{shaneson,Milgram-Ranicki} to arbitrary (connective in the latter case) ring spectra, and also allows to start with an arbitrary form parameter. In this regard it is perhaps of note that our approach uniformly treats all these cases. %
Let us also mention that Ranicki observed that 
there are two natural such decompositions in $\L$-theory; our general one corresponds to what he termed the algebraic decomposition. In the body of the paper we also perform a detailed analysis of the second one, the geometric decomposition, defined in terms of a certain transfer map using the visible signature of the circle. To this end we also revisit the visible signature from a new perspective in Appendix \reffour{section:appendix-visible-revisited}, which we believe to be of independent interest.

Finally, for $\GW\colon \Catp \to \Sp$, %
we obtain splittings
\[
\GW((\C,\QF)_{\mathrm B\mathbb Z}) \simeq \GW(\C,\QF) \oplus \Sigma \GW(\C^\natural,\QF) \oplus \mathrm{NK}(\C) \quad\text{and}\quad \KGW((\C,\QF)_{\mathrm B\mathbb Z}) \simeq \KGW(\C,\QF) \oplus \Sigma \KGW(\C,\QF) \oplus \N\KK(\C).
\]
that seem to be entirely new. For example, for a ring $R$ with form parameter $\lambda$ we obtain canonical equivalences
\begin{align*}
\GW\bigl(\Dfree(R[\ZZ]),\QF_{M[\ZZ]}^{\g\lambda[\ZZ]}\bigr) & \simeq \GW\bigl(\Dfree(R),\QF^{\g\lambda}_{M}\bigr) \oplus \Sigma \GW\bigl(\Dperf(R),\QF^{\g\lambda}_M\bigr) \oplus \N\K(R) \quad\text{and}\\
\KGW^{\g\lambda[\ZZ]}(R[\ZZ],M[\ZZ]) & \simeq \KGW^{\g\lambda}(R,M) \oplus \Sigma \KGW^{\g\lambda}(R,M) \oplus \N\KK(R).
\end{align*}
Furthermore, through 
the identification of the Grothendieck--Witt spectra of visible Poincar\'e $\infty$-categories with Weiss--Williams' LA-spectra from \papertwo, we in particular obtain equivalences
\begin{align*}
\LA^{\f}(X \times \mathrm S^1, \mathrm{pr}^\ast \xi) &\simeq \LA^{\f}(X, \xi) \oplus \Sigma  \LA(X, \xi) \oplus \mathrm N\mathrm A(X) \quad\text{and}\\
\mathbb{LA}(X \times \mathrm S^1, \mathrm{pr}^\ast \xi) &\simeq \mathbb{LA}(X, \xi) \oplus \Sigma  \mathbb{LA}(X, \xi) \oplus \mathrm N\mathbb A(X),
\end{align*}
where $\LA^{\f}$ is the finitely presented and $\mathbb{LA}$ the Karoubi version of the Weis-Williams $\LA$-spectra, 
which again seems to be new.

Finally, we note that Theorem~\reffour{theorem:intro-F} may also be applied to further invariants such as real (topological) Hochschild and real (topological) cyclic homology, once their domain of definition is extended to Poincar\'e $\infty$-categories in \cite{Harpaz-Nikolaus-Shah}.

\introsubsection{Multiplicative structures on Grothendieck-Witt- and L-spectra}
As the final set of results which we establish in this paper, we turn to multiplicative structures. Our main result in this direction, subsuming in particular the final assertion of the \reffour{theorem:pfourmaintheorem}, is the following: 
\begin{theoremintrofour}
\label{theorem:intro-G}%
The cube
\[
\begin{tikzcd}
[column sep=0.5em, row sep=0.7em]
  & \KGW\ar[rr]\ar[dd] && \KL\ar[dd]\\
\GW\ar[ru]\ar[rr, crossing over]\ar[dd] && \L\ar[ur]\\
  & \KK^\hC\ar[rr] && \KK^\tC\\
\K^\hC\ar[rr] \ar[ru] && \K^\tC\ar[uu,<-, crossing over] \ar[ru]\\
\end{tikzcd}
\]
refines to a cube of lax symmetric monoidal functors $\Catp \rightarrow \Sp$, in a way that makes the functors $\GW$ and $\KGW$ into the initial Verdier- and Karoubi-localising lax symmetric monoidal functors $\Catp \to \Sp$, respectively, and similarly $\L$ and $\KL$ into the initial ones that are in addition bordism invariant. 
\end{theoremintrofour}

In particular, if $(\C,\QF)$ is a (commutative) algebra in $\Catp$, that is, if $\C$ is a stably (symmetric) monoidal $\infty$-category and $\QF \colon \C\op \rightarrow \Sp$ is a lax (symmetric) monoidal Poincaré structure, then for any of the functors $\mathcal F$ in the cube above, the spectrum $\F(\C,\QF)$ acquires the structure of a (commutative) ring spectrum. Likewise, if $(\D,\Phi)$ is a $(\C,\QF)$-module in $\Catp$, then $\F(\D,\Phi)$ acquires the structure of an $\F(\C,\QF)$-module spectrum. 
For example, this applies when $\C = \Dperf(R)$ for a commutative ring $R$ and $\QF$ is one of $\QF^\sym_R$ or $\QF_R^\mathrm{gs}$, recovering and extending results of Laures and McClure from \cite{laures} on multiplications for $\L$-spectra (though we do not 
compare our structures to theirs), and when $\D=\Dperf(A)$ for $A$ an $R$-algebra and $\Phi$ one of $\QF^\qdr_A$ or $\QF^\mathrm{gq}_A$. We spell out precise conditions on form parameters encoding that one is a module over another in Appendix \reffour{section:appendix-form-parameters}. The front face of the cube above was also very recently refined to a lax symmetric monoidal one in \cite{HSV26} (and here our universal property immediately implies compatibility with their structure); our results on the back face, in particular for $\LL = \L^{\langle -\infty \rangle}$ seem to be new.

We go on to show that the association $(\C,\QF) \mapsto [(\C,\QF^\qdr) \to (\C,\QF)]$ refines to a lax symmetric monoidal functor $\Catp \to \Fun([1],\Catp)$, where the target is equipped with Day convolution with respect to the minimum on $[1]$.
It then follows that the functor $\NL \colon \Catp \rightarrow \Sp$, given by the common value (as a consequence of the fundamental fibre sequence and Proposition \reffour{proposition:intro-D}) of 
\[
(\C,\QF) \longmapsto \cof(\F(\C,\QF^\qdr) \rightarrow \F(\C,\QF))
\]
for any of $\F = \GW, \KGW, \L$ or $\KL$ referred to as normal or hyperquadratic $\L$-theory by Ranicki, inherits a lax symmetric monoidal structure. Combining such arguments with the fact that $\NL_0(\Dperf(\mathbb Z),\QF^\mathrm{gs}) \cong \mathbb Z/8$, one obtains for example:

\begin{corollaryintrofour}
\label{corollary:intro-H}%
For any discrete ring $R$ and invertible module with involution $M$ over $R$, the kernel and cokernel of each of the polarisation maps
\[
\GW_i^\qdr(R,M) \longrightarrow \GW_i^\mathrm{gq}(R,M) \longrightarrow \GW_i^\mathrm{gs}(R,M) \longrightarrow \GW_i^\mathrm{s}(R,M)
\]
 and their composites is $8$-torsion for every $i \in \ZZ$, and similarly for $\KGW$, $\L$, or $\KL$ in place of $\GW$.
\end{corollaryintrofour}

This statement really does not extend to arbitrary form parameters or (connective) ring spectra, as Weiss' and Williams' explicit calculations for the sphere spectrum from \cite{WWIII} and our calculation for the genuine Frobenius parameter from \paperthree show.

Our proof of Theorem~\reffour{theorem:intro-G} relies on a general analysis of the Day convolution $\infty$-operads $\Fun(\C,\D)^\Day$ and their localisation properties, to which we devote Appendix \reffour{section:appendix-localisations}. The colours of this operad are the functors $\C \to \D$ and its commutative algebra objects are the lax symmetric monoidal such functors. Our arguments are then based on the following observations:
\begin{enumerate}
\item
The Day convolution $\infty$-operad $\Fun(\Catp,\Spa)^{\Day}$, though not a symmetric monoidal $\infty$-category due to size issues, still admits a unit object, given explicitly by $\SS[\Poinc]$.
\item
Each of the functors $\GW,\KGW,\L,\KL$ is initial under $\SS[\Poinc]$ with respect to a certain property: being additive, Karoubi-localising, or a combination of one of these with bordism invariance. Each of these properties can be phrased as sending a certain collection $\mathcal A$ of cone diagrams to limit diagrams.
\item
When the collection $\mathcal A$ is suitably compatible with the tensor product on $\Catp$, such universal approximations of the unit in $\Fun(\Catp,\Spa)^{\Day}$ become themselves the units of the full suboperad spanned by the functors which send all diagrams in $\mathcal{A}$ to limit diagrams. In particular, these universal approximations refine, in a suitably unique manner, to initial commutative algebra objects, hence initial lax symmetric monoidal functors of the given type. %
\item When $\mathcal A$ consists of finite diagrams and is furthermore stable under forming tensors and cotensors with arbitrary stable $\infty$-categories, these commutative algebra objects are even idempotent and therefore yield smashing localisations onto the category of functors $\Catp \to \Spa$ taking all diagrams in $\mathcal A$ to limit diagrams. 
\end{enumerate}

The last point applies to the collections of diagrams defining additivity and bordism invariance, giving:

\begin{propositionintrofour}
\label{proposition:intro-J}%
Both functors $\GW,\L \colon \Catp \to \Sp$ are idempotent with respect to Day convolution and convolving with them yields a smashing localisation onto the full subcategory of additive functors $\Catp \to \Sp$ in the case of $\GW$ and additive bordism invariant functors $\Catp \rightarrow \Sp$ in the case of $\L$. %
\end{propositionintrofour}

A completely analogous statement holds for $\K \colon \Catx \rightarrow \Sp$
and it follows that all functors from $\Catp$ or $\Catx$ to spectra admit additive (and bordism invariant) approximations. By contrast, the collection defining Karoubi-localising functors is not closed under cotensors and indeed we do not know whether any of the functors $\KK, \KGW$ or $\KL$ is idempotent. Following a suggestion of Ramzi, we do, however, exhibit non-Karoubi-localising modules over all three of these functors in the body of the paper, which rules out that Karoubi-localising functors form a smashing localisation of all functors $\Catp \rightarrow \Spa$. We consequently do not know whether all functors to spectra admit Karoubi-localising approximations (with or without bordism invariance).

There are many variants of the above that can be studied, for example by considering the functors $\K$, $\GW$, $\L$ and their Karoubi-localisations also after postcomposition with the forgetful functors $\Sp \to \Grp_{\Einf} \to \Sps$. For instance, we show that the functor $\Lspace\colon\Catp \to \Grp_{\Einf}$ is idempotent and localises onto bordism invariant and direct sum preserving functors (in particular bordism invariant additive functors to $\Grp_{\Einf}$ also do not form a smashing localisation inside all functors). In upcoming work, Saunier shows that $\Kspace\colon \Catx \to \Sps$ and $\GWspace\colon \Catp \to \Grp_{\Einf},\Sps$ are also idempotent and identifies the associated localised subcategories.

Let us finally mention one application outside the realm of hermitian K-theory, as a further advertisement of our localisation machinery: If $\C$ is a closed stably symmetric monoidal small $\infty$-category and $\D$ is a %
closed and cocomplete stable symmetric monoidal $\infty$-category then our criteria apply to show that the full subcategory $\Exc^n_*(\C,\D)$ of reduced $n$-excisive functors $\C \to \D$ is a smashing localisation of $\Fun(\C,\D)$. The case of $\C=\Spaf$ and $\D=\Spa$ was recently obtained in~\cite{ABHS} and to put matters into perspective, note that even $\Exc^1_*(\C,\Spa)$ is typically not a smashing localisation when $\C$ is not stable. %

\introsubsection{Organisation of the paper}
In \S\reffour{section:localising}, we discuss Karoubi-localising functors and prove %
the second part of Theorem~\reffour{theorem:intro-B} in Theorem ~\reffour{proposition:localising-are-stable}, and Proposition~\reffour{proposition:intro-D} as Theorem~\reffour{theorem:Karoubi-pullback-squares}, respectively. \S\reffour{section:shaneson} is then devoted to the circle decomposition theorem, Theorem~\reffour{theorem:intro-F}, and the subsequent \S\reffour{section:karoubification} etablishes Theorem~\reffour{theorem:intro-E}, and thus Williams' conjecture as Propositions~\reffour{proposition:decorated-GW} and~\reffour{proposition:decorated-L} along with the remaining part of Theorem~\reffour{theorem:intro-B} in Corollary~\reffour{corollary:karoubification}.
Finally, \S\reffour{section:day} discusses multiplicative structures on all Grothendieck-Witt and L-functors constructed in this paper. We prove Theorem~\reffour{theorem:intro-G} as Corollary~\reffour{corollary:kgw-kl-monoidal}, Corollary~\reffour{corollary:intro-H} as Corollary~\reffour{corollary:normalL8tors}, and Proposition~\reffour{proposition:intro-J} in Corollaries~\reffour{corollary:GW-universal-multiplicative} and~\reffour{corollary:lmonoidal}. 

Finally, the paper has four appendices. In \S\reffour{section:appendix-localisations}, we discuss multiplicative and smashing localisations of $\infty$-operads which lay the foundations of our results in \S\reffour{section:day}. %
In \S\reffour{section:presentability} we show that $\Catp$ and $\Cath$ are compactly generated, a technical result we use a few times along the way, in \S\reffour{section:appendix-visible-revisited} we revisit the visible symmetric signature of a Poincar\'e complex to facilitate the discussion of the two possible splittings in the circle decomposition theorem, and in \S\reffour{section:appendix-form-parameters} we unwind in detail our mulitplicitativity results into the more classical language of form parameters. %

\introsubsection{Acknowledgements}
For useful discussions about our project,
we heartily thank
Tobias Barthel,
Clark Barwick,
Patrick Bieker,
Lukas Brantner,
Mauricio Bustamante,
Denis-Charles Cisinski,
Dustin Clausen,
Diarmuid Crowley,
Uriya First,
Thorger Geiß,
Rune Haugseng,
André Henriques,
Lars Hesselholt,
Gijs Heuts,
Manuel Hoff,
Geoffroy Horel,
Marc Hoyois,
Andreas Huber,
Max Karoubi,
Daniel Kasprowski,
Ben Knudsen,
Manuel Krannich,
Achim Krause,
Henning Krause,
Sander Kupers,
Wolfgang Lück,
Ib Madsen,
Cary Malkiewich,
Mike Mandell,
Akhil Mathew,
Lennart Meier,
Vikram Nadig,
Irakli Patchkoria,
Nathan Perlmutter,
Phil Pützstück,
Maxime Ramzi,
Oscar Randal-Williams,
Andrew Ranicki,
George Raptis,
Victor Saunier,
Marco Schlichting,
Peter Scholze,
Stefan Schwede,
Graeme Segal,
Markus Spitzweck,
Jan Steinebrunner,
Georg Tamme,
Ulrike Tillmann,
Michael Weiss,
Christoph Winges,
and
Maria Yakerson.

Furthermore, we owe a tremendous intellectual debt to Jacob Lurie for creating the framework we exploit here, and to Søren Galatius for originally spotting the overlap between various separate projects of ours; his insight ultimately led to the present collaboration.

\medskip

The authors would also like to thank the Hausdorff Center for Mathematics at the University of Bonn, the Newton Institute at the University of Cambridge, the University of Copenhagen and the Mathematical Research Institute Oberwolfach for hospitality and support while parts of this project were undertaken.
\medskip

No AI was used in the production of this paper.
\medskip

BC was supported by the French National Centre for Scientific Research (CNRS) through a ``délégation'' at LAGA, University Paris 13, and by the French National Research Agency (ANR) project ``Motivic homotopy, quadratic invariants and diagonal classes'' (ANR grant no.\ 21-CE40-0015).
BC and YH were supported by the ANR project ``K-theory, Actions \& stable Homotopy'' (ANR grant no.\ 25-CE40-2861) at LAGA, University Paris 13 and LML, University of Artois.
ED was supported by the German Research Foundation (DFG) through the priority program ``Homotopy theory and Algebraic Geometry'' (DFG grant no.\ SPP 1786) at the University of Bonn and WS by the priority program ``Geometry at Infinity'' (DFG grant no.\ SPP 2026) at the University of Augsburg. ED was further supported by the Engineering and Physical Sciences Research Council (EPSRC) through the grant ``Characteristic polynomials for symmetric forms'' (EPSRC grant no.\ EP/W019620/1) at the University of Warwick.
YH was supported by the European Research Council (ERC) through the project ``Foundations of Motivic Real K-Theory'' (ERC grant no.\ 949583) at Paris Cité University.
YH and DN were supported by the ANR project ``Chromatic Homotopy and K-theory'' (ANR grant no.\ 16-CE40-0003) at LAGA, University of Paris 13.
FH was a member of the Hausdorff Center for Mathematics (DFG grant no.\ EXC 2047 390685813) at the University of Bonn and TN of the cluster ``Mathematics Münster: Dynamics-Geometry-Structure'' (DFG grant no.\ EXC 2044 390685587) at the University of Münster. 
FH was further supported by the DFG through the collaborative research center ``Integral structures in Geometry and Representation Theory'' (DFG grant no.\ TRR 358 4913924) at the University of Bielefeld, TN through the centre ``Geometry: Deformations and Rigidity" (DFG grant no.\ SFB 1442 427320536) at the University of Münster, ML and DN through the centre ``Higher Invariants'' (DFG grant no.\ SFB 1085 224262486), and ML and WS through the centre ``Generalised Motivic Methods in Geometry'' (DFG grant no.\ SFB 1785 560334683) both at the University of Regensburg.
FH was also supported by the ERC through the project ``Moduli spaces, Manifolds and Arithmetic'' (ERC grant no.\ 682922) of Søren Galatius and KM by the project ``$\K$-theory, $\L^2$-invariants, manifolds, groups and their interactions'' (ERC grant no.\ 662400) of Wolfgang Lück. 
FH, TN and WS were further supported by the EPSRC through the program ``Homotopy harnessing higher structures'' at the Isaac Newton Institute for Mathematical Sciences (EPSRC grants no.\ EP/K032208/1 and EP/R014604/1). 
ML was supported by the research fellowship ``New methods in algebraic K-theory'' (DFG grant no.\ 424239956) and by the Danish National Research Foundation (DNRF) through the Center for Symmetry and Deformation (DNRF grant no.\ 92) and the Copenhagen Centre for Geometry and Topology (DNRF grant no.\ 151) at the University of Copenhagen. 
KM was also supported by the K\&A Wallenberg Foundation at the University of Stockholm.

\section{Karoubi-Grothendieck-Witt spectra}
\label{section:localising}%

In this section we study Karoubi-localising functors.
We begin in \S\reffour{subsection:cofinality} with a preliminary discussion of the behaviour of $\GW$ and $\L$ under idempotent completions. We show that Schlichting's cofinality theorem for Grothendieck-Witt theory 
holds in our setting 
and deduce the existence of Rothenberg sequences for $\L$-theory. 

Our main result in \S\reffour{subsection:localising-invariants} is Proposition~\reffour{proposition:localising-are-stable}, stating that the $\infty$-category of Karoubi-localising functors valued in an additive $\infty$-category with finite limits is stable. We use this result to construct the Karoubi-Grothendieck-Witt functor $\KGW$ as the unique Karoubi-localising delooping of $\GWspace(-^{\natural})$, and deduce that $\KGW$ is the initial Karoubi-localising functor under $\GW$. We then define the Karoubi $\L$-theory functor $\KL = \KGW^{\bord}$ as the bordification of $\KGW$, and show that it is at the same time also the initial Karoubi-localising functor under $\L$. %
We end this section by a comparison result that describes the fibres of $\GW\to \KGW$ and $\L\to \KL$ in terms of algebraic $\K$-theory.

\subsection{Hermitian cofinality and Rothenberg sequences}
\label{subsection:cofinality}%

We recall from Remark~\reftwo{remark:idem-poincare} that if $(\C,\QF)$ is a Poincar\'e $\infty$-category, then $\QF$ extends (via left Kan extension) uniquely to a Poincar\'e structure on the idempotent completion $\C^\natural$ of $\C$. We therefore, by slight abuse of notation, write $(\C^\natural,\QF)$ for the resulting Poincar\'e $\infty$-category and refer to the association $(\C,\QF) \mapsto (\C^\natural,\QF)$ as idempotent completion of Poincar\'e $\infty$-categories.
We then say that a Poincaré functor $(\C,\QF) \to (\D,\QFD)$ is a Karoubi equivalence if it induces an equivalence on idempotent completions in this sense. Concretely, this means that $\C \to \D$ is a dense inclusion of stable $\infty$-categories, i.e.\ a fully faithful functor such that every object in $\D$ is a retract of an object in $\C$, and $\QFD|_{\C} = \QF$.

Our goal in this subsection is to study the behaviour %
of the Grothendieck-Witt spectrum along Karoubi equivalences $(\C,\QF) \to (\D,\QFD)$.
For this, let us first recall Thomason's cofinality theorem. To ease notation, in the following, for an abelian group $A$ we shall also denote by $A$ its associated Eilenberg-Mac Lane spectrum.

\begin{proposition}
\label{proposition:cofinality-K-theory}%
Let $\C \to \D$ be a Karoubi equivalence of stable $\infty$-categories (e.g.\ $\D = \C^{\natural}$).
Then the induced map 
\[
\K_n(\C)\to \K_n(\D)
\]
is an isomorphism for $n>0$ and a monomorphism for $n=0$. In particular, the square of spectra
\[
\begin{tikzcd}
\K(\C)\ar[r]\ar[d] & \K(\D)\ar[d]\\
\K_0(\C)\ar[r] & \K_0(\D)
\end{tikzcd}
\]
is exact. %
In addition, an object $x\in \D$ is in the essential image of $\C$ %
if and only if its $\K$-class $[x] \in \K_0(\D)$ is in the image of $\K_0(\C)$. 
\end{proposition}
\begin{proof}
The first claim is %
\cite{thomason-trobaugh}*{Theorem 1.10.1} (revisited in the setting of stable $\infty$-categories in \cite{HLS}*{Theorem 7.1}). The last claim follows from Thomason's classification of dense subcategories \cite{thomason-classification}*{Theorem 2.1}. 
\end{proof}

The following result extends~\cite{schlichting-exact}*{Corollary 5.2} to the setting of Poincaré $\infty$-categories.
\begin{proposition}
\label{proposition:weak-cofinality}%
Let $(\C,\QF) \to (\D,\QFD)$ be a Karoubi equivalence of Poincaré $\infty$-categories (e.g.\ $(\D,\QFD) = (\C^\natural,\QF)$).
Then the induced map
\[
\GW_n(\C,\QF)\to \GW_n(\D,\QFD)
\]
is an isomorphism for $n>0$ and a monomorphism for $n=0$, whose image is the subgroup
\[
\GW_0(\D,\QFD) \times_{\K_0(\D)}\K_0(\C)\subseteq \GW_0(\D,\QFD) .
\]
\end{proposition}

\begin{proof}%
We will prove the statement by induction on $n$. Let us first deal with the case $n=0$. By Corollary~\reftwo{corollary:pi-0gw}, the group $\GW_0(\C,\QF)$ is the quotient of the monoid of equivalence classes of Poincaré objects under the relations
\[
[x,q]=[\hyp(y)]
\]
for any Poincar\'e object $(x,q)$ with Lagrangian $y\to x$. In particular, a class $[x,q]$ is 0 if and only if there exist Poincar\'e objects $(x_1,q_1),\,(x_2,q_2)$ with Lagrangians $y_i \to x_i$ for $i=1,2$ such that
\[
(x,q)\oplus (x_1,q_1)\oplus \hyp(y_2)\cong (x_2,q_2)\oplus \hyp(y_1)\,.
\]
To simplify the presentation, let us implicitly identify $\C$ with its essential image in $\D$.
Suppose that $[x,q]$ is in the kernel of the map $\GW_0(\C,\QF)\to \GW_0(\D,\QFD)$. Then we can find $(x_1,q_1)$, $(x_2,q_2)$ and $y_1,y_2$ as above in $\D$. But then there exist $y_1'$ and $y_2'$ in $\D$ such that for $j=1,2$ we have that $y_j'' \coloneq y_j\oplus y_j'$ are in $\C$. Set $(x_j',q_j') = (x_j,q_j) \oplus \hyp(y_j')$. Then $(x_j',q_j')$ is Poincar\'e with Lagrangian $y_j''  \in \C$, which implies that $x_j'$ is in $\C$ as well (since it sits in a fibre sequence between $y_j''$ and its dual). As
\[
(x,q)\oplus (x_1',q_1')\oplus \hyp(y''_2)\cong (x_2',q_2')\oplus \hyp(y_1'')
\]
we have that 
$[x,q]$ vanishes in $\GW_0(\C,\QF)$. We thus conclude that the map $\GW_0(\C,\QF) \to \GW_0(\D,\QFD)$ is injective.
To see that the image is as claimed, suppose $[x,q]\in \GW_0(\D,\QFD)$ is an object such that the $\K_0$-class of $x$ lies in $\K_0(\C)$. Then $x$ is equivalent to an object in $\C$ by Proposition~\reffour{proposition:cofinality-K-theory}. Hence $[x,q]$ is in the image of $\GW_0(\C,\QF) \to \GW_0(\D,\QFD)$.
    
For the induction step, we let $n\geq 1$, and 
contemplate the map of metabolic fibre sequences of Example~\reftwo{example:metabolicfseq}
\[
\begin{tikzcd} 
(\C,\QF\qshift{-1}) \ar[r] \ar[d] & \Met(\C,\QF) \ar[r] \ar[d] & (\C,\QF) \ar[d] \\
(\D,\QFD\qshift{-1}) \ar[r] & \Met(\D,\QFD)\ar[r] & (\D,\QFD) \ .
\end{tikzcd}
\]
By Corollary~\reftwo{corollary:GWHyp},
applying $\GW$ to this diagram produces a diagram of long exact sequences of which a fragment is 
\[
\begin{tikzcd}
[column sep=1em] 
\GW_n(\C,\QF\qshift{-1}) \ar[r]\ar[d] & \K_n(\C)\ar[r]\ar[d,"\cong"] &\GW_n(\C,\QF)\ar[r]\ar[d] & \GW_{n-1}(\C,\QF\qshift{-1}) \ar[r]\ar[d]\arrow[dr, phantom, "\lrcorner", very near start]& \K_{n-1}(\C)\ar[d]\\
\GW_n(\D,\QFD\qshift{-1}) \ar[r] & \K_n(\D)\ar[r] & \GW_n(\D,\QFD)\ar[r] & \GW_{n-1}(\D, \QFD\qshift{-1}) \ar[r]& \K_{n-1}(\D)
\end{tikzcd}
\]
In this diagram, the second vertical arrow is an isomorphism for $n \geq 1$ by Proposition~\reffour{proposition:cofinality-K-theory}. In addition, 
by induction hypothesis, 
the rightmost square is cartesian. A diagram chase then shows that the middle map is surjective. Applying this result to the shifted Poincar\'e structures, we find that also the leftmost map is surjective. Another diagram chase then shows that the middle map is injective.
\end{proof}

We now arrive to the main result of the present subsection, which extends the cofinality theorem of \cite{schlichting-derived}*{Theorem~7.9} 
to the Poincaré categorical setup. In what follows, for an abelian group $A$ with $\Ct$-action, we generally write $A_\hC$, $A^\hC$ and $A^\tC$ for the corresponding construction applied to the Eilenberg-MacLane spectrum of $A$.

\begin{theorem}[Hermitian cofinality]
\label{theorem:strong-cofinality}%
Let $(\C,\QF) \to (\D,\QFD)$ be a Karoubi equivalence of Poincaré $\infty$-categories.
Then all the squares in the following diagrams are cartesian:
\[
\begin{tikzcd}
\GW(\C,\QF)\ar[r]\ar[d] & \GW(\D,\QFD) \ar[d] \\
\K(\C,\QF)^\hC\ar[r]\ar[d] & \K(\D,\QFD)^\hC \ar[d] \\
\K_0(\C,\QF)^\hC\ar[r] & \K_0(\D,\QFD)^\hC
\end{tikzcd}
\qquad
\begin{tikzcd}
\L(\C,\QF)\ar[r]\ar[d] & \L(\D,\QFD)\ar[d] \\
\K(\C,\QF)^\tC\ar[r]\ar[d] & \K(\D,\QFD)^\tC \ar[d] \\
\K_0(\C,\QF)^\tC\ar[r] & \K_0(\D,\QFD)^\tC
\end{tikzcd}
\]
\end{theorem}

\begin{proof}%
The bottom squares are cartesian by Proposition~\reffour{proposition:cofinality-K-theory} and the fact that $(-)^\hC$ and $(-)^\tC$ are exact functors, so it is enough to prove that the top squares are cartesian. Since they are the front and back faces of the cube
\[
\begin{tikzcd}
[column sep=0.5em, row sep=0.7em]
  & \GW(\C,\QF)\ar[rr]\ar[dd] \ar[ld]&& \GW(\D,\QFD)\ar[dd]\ar[ld]\\
\L(\C,\QF)\ar[rr, crossing over]\ar[dd] && \L(\D,\QFD)\\
  & \K(\C,\QF)^\hC\ar[rr]\ar[ld] && \K(\D,\QFD)^\hC\ar[ld]\\
\K(\C,\QF)^\tC\ar[rr] && \K(\D,\QFD)^\tC\ar[uu,<-, crossing over]\\
\end{tikzcd}
\]
whose left and right faces are cartesian, there is a natural equivalence between the total fibres of the two squares under consideration. Let us denote this common total fibre by $T$, viewed as a functor from the $\infty$-category of Karoubi equivalences $(\C,\QF) \to (\D,\QFD)$ to spectra,
and note that the task is to show that $T = 0$. Now since $T$ is the total fibre of the back face, Proposition~\reffour{proposition:weak-cofinality} and Proposition~\reffour{proposition:cofinality-K-theory} imply that $T$ takes values in coconnected spectra. On the other hand, since $T$ is also the total fibre of the front face which consists of bordism invariant functors, we have that $T(-\qshift{n}) = \Sigma^nT$. Together, this shows that $\Sigma^nT$ takes values in coconnective for all integers $n$, and hence that $T=0$ as desired.
\end{proof}

Let us highlight that the claim for L-theory is deduced from the simpler one for Grothendieck-Witt theory via the fundamental exact sequence, contrary to usual flow of information utilised in this paper series.

\begin{remark}
\label{remark:Rothenberg_map}%
On $\pi_n$, the map $\xi \colon \L(\C, \QF)\to \K_0(\C,\QF)^{\tC}$ appearing in Theorem~\reffour{theorem:strong-cofinality}
sends the class of $(x,q)\in \Poinc(\C, \QF\qshift{-n})$ to the class of $x$ in the Tate cohomology group
\[
\widehat \mathrm{H}^{-n}(\Ct; \K_0(\C)) \cong \frac{\{a\in \K_0(\C)\mid \Dual_\QF a = (-1)^n a\}}{\{a+(-1)^n\Dual_\QF a\mid a\in \K_0(\C)\}} ;
\]
indeed, this is true by definition in degree $0$, and in degree $n$ follows by identifying the above map with the map induced on $\pi_0$ by $\L(\C, \QF\qshift{-n})\to \K_0(\C,\QF\qshift{-n})^{\tC}$, where we note that the $\Ct$-action on $\K_0(\C, \QF\qshift{-n})$ differs from that of $\K_0(\C, \QF)$ by multiplication with $(-1)^n$ (e.g., by Theorem~\reftwo{theorem:kr-intro-two}). We denote the resulting $\Ct$-action by $\K_0(\C,\QF)^{(-1)^n}$. Now, the %
periodicity isomorphism
\[
\widehat\mathrm{H}^{-n}(\Ct; \K_0(\C)) = \pi_{n}\K_0(\C,\QF)^{\tC} \cong \pi_0\K_0(\C,\QF\qshift{-n})^{\tC} = \widehat\mathrm{H}^0(\Ct; \K_0(\C)^{(-1)^{n}})
\]
afforded by bordism invariance of $\K_0(-)^{\tC}$ coincides with the standard periodicity isomorphism in Tate cohomology; indeed, by naturality and the classification of dense Poincaré subcategories, it is enough to note that if $\K_0(\C,\QF)$ is generated by an element $[x]$ with dual $(-1)^n[x]$ then $\widehat{\mathrm{H}}^{-n}(\Ct;\K_0(\C,\QF))$ is either $0$ or $\ZZ/2$, and in any case admits no non-identity automorphisms.
\end{remark}

For a Poincaré $\infty$-category $(\C,\QF)$ and a subgroup $c\subseteq \K_0(\C)$ closed under the duality, let us write $\C^c \subseteq \C$ for the dense subcategory spanned by all objects whose K-class lies in $c$, and abbreviate $(\C^c,\QF|_{\C^c})$ to $(\C^c,\QF)$.
We refer to the $\GW$- and $\L$-spectra of $(\C^c,\QF)$ as the $c$-controlled $\GW$- and $\L$-spectra of $(\C,\QF)$. 

\begin{corollary}
\label{cor:rothenbergconc}%
Given a Poincaré $\infty$-category $(\C,\QF)$ and a pair of consecutive subgroups $c \subseteq d \subseteq \K_0(\C)$ closed under the duality, we obtain exact sequences
\[
\GW(\C^c,\QF)\to \GW(\C^d,\QF)\xrightarrow{\xi} (d/c)^\hC \quad\text{and}\quad \L(\C^c,\QF)\to \L(\C^d,\QF)\xrightarrow{\xi} (d/c)^\tC,
\]
\end{corollary}
They translate to long exact sequences relating the corresponding controlled Grothendieck-Witt- and $\L$-groups of $(\C,\QF)$ with group- and Tate-cohomology of $\Ct$ with coefficients in $d/c$, respectively. 
In the case of $\L$-theory, the resulting long exact sequence 
\begin{equation}
\label{equation:rothenberg}%
\dots \to \L_n(\C^c, \QF) \to \L_n(\C^d, \QF) \xrightarrow\xi \widehat{\mathrm{H}}^{-n}(\Ct; d/c) \to \L_{n-1}(\C^c, \QF) \to \dots
\end{equation}
is a generalisation of the \defi{Rothenberg sequence} from~\cite{Ranickiyellowbook}*{Proposition~1.10.1}. Following Remark~\reffour{remark:Rothenberg_map}, the map $\xi$ can be explicitly described as sending $(x,q)\in \Poinc(\C^d, \QF\qshift{-n})$ to the element of $\widehat{\mathrm{H}}^{-n}(\Ct; d/c)$ represented by the $\K$-class of $x$ in $d/c$.

\begin{example}
\label{example:free-dedekind}%
Let $R$ be a Dedekind ring whose field of fractions $K$ is a global field of characteristic not $2$. Recall from Corollary~\refthree{corollary:compute-L-groups-global-field} that its symmetric $\L$-groups are given by
\[
\L_n(\Dperf(R),\QF^\s_R) = \L_n^\s(R) \cong 
\begin{cases}
\W^\s(R) & \text{ for } n \equiv 0 (4) \\
\displaystyle\mathop{\oplus}\limits_{\fp \in \mathcal{I}} \ZZ/2 & \text{ for } n \equiv 1 (4) \\
0 & \text{ for } n \equiv 2 (4) \\
\Pic(R)/2 & \text{ for } n \equiv 3 (4) \ ,
\end{cases}
\]
where $\mathcal{I}$ is the set of prime ideals of $R$ which contain the ideal $(2)$.
To compare this to $\L(\Dfree(R),\QF^\s_R)$, we first recall that 
\[
\tilde \K_0(R) \coloneq \K_0(\Dperf(R))/\K_0(\Dfree(R))) \cong \ker(\K_0(R) \to \K_0(K)) \cong \Pic(R)
\]
where
the induced $\Ct$-action on $\Pic(R)$ is the sign action since every line bundle is inverse to its dual. The maps into Tate cohomology in the Rothenberg sequence are thus maps
\[
\xi_{2k}\colon \L^{\sym}_{2k}(R) \to \Pic(R)[2] \quad\text{and}\quad \xi_{2k+1}\colon \L^{\sym}_{2k+1}(R) \to \Pic(R)/2
\]
into the $2$-torsion part and the quotient modulo $2$ respectively. 

Let us now unwind what these maps are (under the isomorphism of Corollary~\refthree{corollary:compute-L-groups-global-field}):
\begin{enumerate}
\item
For $r = 3\;(\mathrm{mod}\;4)$ the map $\xi_r$ is the identity $\Pic(R)/2 \to \Pic(R)/2$. To see this, note that the $\L$-group $\L^{\sym}_{-1}(R)$ is generated by the images of the unit $[R/\fp,q^{\uni}]$ of $\L^\sym_{0}(R/\fp)$ under the various inflation functors $\Dperf(R/\fp,\QF^{\sym}_{R/\fp}) \to \Dperf(R,(\QF^{\sym}_R)\qshift{1})$ for all primes $\fp$ (after choosing a uniformiser for each $\fp$). 
These images map
to the class mod $2$ of $[\fp]$ in $\L^{\sym}_{-1}(R) = \Pic(R)/2$, while its underlying reduced $\K$-theory class is $[R/\fp] = -[\fp] \in \K_0(\Dperf(R))/\K_0(\Dfree(R))) = \Pic(R)$. Since $[\fp]=-[\fp] \;\mathrm{mod}\;2$ the claim follows from Remark~\reffour{remark:Rothenberg_map}.
\item
For $r = 1\;(\mathrm{mod}\;4)$ the map $\xi_r\colon \oplus_{\fp \in \I}\ZZ/2 \to \Pic(R)/2$ sends the generator of the $\fp$ summand to the class of $[\fp]$ mod $2$.
Here the $\L$-group $\L_1^{\sym}(R)$ is generated by shifted Poincaré objects in the image of the inflation functor $\Dperf(R/\fp,\QF^{\sym}_{R/\fp}) \to \Dperf(R,(\QF^{\sym}_R)\qshift{1})$ for $R/\fp$ of characteristic $2$; the generator of $\L^\sym_{2}(R/\fp) \cong \ZZ/2$ then maps to the generator of the corresponding $\ZZ/2$ term in $\L^{\sym}(R)_{1} = \oplus_{\fp \in \I}\ZZ/2$, while as above its underlying reduced $\K$-class is $[\fp]=-[\fp] \;\mathrm{mod}\;2$, and we use again Remark~\reffour{remark:Rothenberg_map}.
\item
For $r = 2\;(\mathrm{mod}\;4)$ the source of $\xi_r$ is the 0 group.
\item
For $r = 0\;(\mathrm{mod}\;4)$ we do not have an explicit formula, but we claim that $\xi_r\colon \W^{\sym}(R) \to \Pic(R)[2]$ is surjective. Indeed, any element in $\Pic(R)[2]$ can be represented by an ideal $\mathfrak{a} \subseteq R$ such that $\mathfrak{a}\cdot \mathfrak{a} = (a)$ is principal. Then $\mathfrak{a}$ is in particular a projective $R$-module and we can endow it with the unimodular symmetric bilinear form $b_a\colon \mathfrak{a} \times \mathfrak{a} \to R$ defined by
\[
b_a(x,y) = \frac{xy}{a} .
\]
The Poincaré object $(\mathfrak{a},b_a) \in \Poinc(\Dperf(R),\QF^{\sym}_R)$ then has underlying reduced $\K$-class $[\mathfrak{a}] \in \Pic(R)[2]$, so that the map in question is indeed surjective. 
\end{enumerate}

We hence conclude from the Rothenberg sequence above %
that $\L_n^\s(\Dfree(R)) = 0$ for $n=2,3 \;(\mathrm{mod}\; 4)$. Identifying $\L_{4k}^{\s}(\Dfree(R))$ with the Witt group $\W^{\sym}_{\f}(R)$ of free symmetric $R$-modules (using 
Corollary~\refthree{corollary:d-zero}\;\refthreeitem{item:cor-two}; the statement there is for $\Dperf(R)$ but the same argument works for any dense subcategory thereof)
we find that that the Rothenberg sequence in degrees $n=0,1 \;(\mathrm{mod}\; 4)$ takes the following form:
\[
0 \to \Pic(R)[2] \to \L_{4k+1}^\s(\Dfree(R)) \to \displaystyle\mathop{\oplus}\limits_{\fp \in \mathcal{I}} \ZZ/2 \to \Pic(R)/2 \to\W^{\sym}_{\f}(R) \to \W^{\sym}(R) \to \Pic(R)[2] \to 0 .
\]
In particular, if $2$ is invertible in $R$, i.e.\ if $\I = \emptyset$, then $\L_{4k+1}^\s(\Dfree(R)) \cong \Pic(R)[2]$, and the free Witt group sits in an exact sequence
\[
0 \to \Pic(R)/2 \to \W^{\sym}_{\f}(R) \to \W^{\sym}(R) \to \Pic(R)[2] \to 0.
\]
\end{example}

\subsection{Karoubi-localising invariants}
\label{subsection:localising-invariants}%

Recall from Definition~\reftwo{definition:additive} that for an $\infty$-category $\E$ with finite limits, a functor $\F\colon \Catp\to\E$ is called additive, Verdier-localising, or Karoubi-localising if it is reduced and sends split Poincaré-Verdire squares, Poincaré-Verdier squares, or Poincaré-Karoubi squares, to cartesian squares, respectively.
Also recall from Proposition~\reftwo{proposition:localizing-criterion} that the condition that $\F$ is Karoubi-localising is equivalent to the condition that $\F$ is Verdier-localising and Karoubi invariant, that is, inverts Karoubi equivalences, %
a criterion which we will employ several times in the discussion below.
We will denote by 
\[
\Funkloc(\Catp,\E) \subseteq \Funadd(\Catp,\E)
\]
the full subcategory spanned by the Karoubi-localising functors, and note that this subcategory is closed under limits, and also under colimits if $\E$ is stable.

\begin{example}
The spectrum valued functor $(\C,\QF) \mapsto \fib[\GW(\C,\QF) \to \K(\C,\QF)^{\hC}] = \fib[\L(\C,\QF) \to \K(\C,\QF)^{\tC}]$ is Karoubi-localising: Indeed, it is Karoubi-invariant by Theorem~\reffour{theorem:strong-cofinality} and Verdier-localising as it is the fibre of two Verdier-localising functors.
\end{example}

None of the spectrum valued functors $\K, \GW$ and $\L$ are Karoubi-localising. Indeed, by Proposition~\reffour{proposition:cofinality-K-theory}, Proposition~\reffour{proposition:weak-cofinality} and Theorem~\reffour{theorem:strong-cofinality} none of these functors are Karoubi-invariant. %
For the same reason, their space valued counterparts $\Kspace = \Om^{\infty}\K, \GWspace = \Om^{\infty}\GW$ and $\Lspace = \Om^{\infty}\L$ are also not Karoubi-localising.

One may be tempted to remedy this by instead considering the functors $\K(-^\natural), \GW(-^\natural)$, and $\L(-^\natural)$. Clearly, these functors are now Karoubi-invariant. However, they are not Verdier-localising. To see this, let for instance $R$ be the nodal curve over a field $k$ of characteristic different from 2, i.e.\ the pullback of the evaluation at $0,1$ map $k[t] \to k \times k$ along the diagonal, compare Remark~\reftwo{remark:no-verdier} for the computations about $R$ we are about to use. Then the sequence 
\[
\Dperf(R \textrm{ on } t) \to \Dperf(R[t]) \to \Dperf(R[t^{\pm1}])
\]
is a Karoubi sequence which is not a Verdier sequence, see Example~\reftwo{example:app-Verdier-not-karoubi}, and is hence the idempotent completion of the Verdier sequence obtained by replacing $\Dperf(R[t^{\pm1}])$ with its dense subcategory given by the essential image of the localisation functor. 
Therefore, either of the above three functors is Verdier-localising only if it sends this dense inclusion to an equivalence. By the Bass-Heller-Swan decomposition, this dense inclusion induces the direct summand inclusion $\K_0(R) \to \K_0(R) \oplus \K_{-1}(R)$ on $\K_0$; here, we use that $\N\K_0(R) = 0$ as follows from the Bass-Milnor-Murthy excision sequence for the above mentioned pullback description of the nodal curve upon passing to polynomial rings. The (hermitian) cofinality theorems say that this is the case if and only if $\K_{-1}(R)$ vanishes (after applying $\Ct$-fixed points or the $\Ct$-Tate construction, respectively). But $\K_{-1}(R) \cong \ZZ$, so none of these three terms vanish.

Moreover, from the natural equivalence $\L(\C,\QF\qshift{n}) \simeq \Sig^n \L(\C,\QF)$ and the fact that $(-)\qshift{n}$ preserves both Verdier and Karoubi sequences, we deduce that the space-valued functor $\Lspace(-^\natural)$ is also not Verdier (and hence not Karoubi) localising. In contrast, the space-valued functor $\Kspace(-^\natural)$ \emph{is} Karoubi-localising, see \cite{HLS}*{Corollary 7.2}. Similarly, we find:

\begin{proposition}
\label{example:GWspace_Karoubi_localising}%
The space-valued functor $\GWspace(-^{\natural})$ 
is Karoubi-localising.
\end{proposition}
\begin{proof}
The claim is equivalent to saying that $\GWspace$ itself sends all Poincaré-Karoubi squares
\[
 \begin{tikzcd}
  (\C, \QF) \ar[r] \ar[d] & (\C', \QF') \ar[d]\\
  (\D, \QFD) \ar[r] & (\D', \QFD'),
 \end{tikzcd}
\]
consisting of idempotent complete Poincaré $\infty$-categories to exact squares. 
Given such a square we may factor it into a pasting of two pullback squares
\[
\begin{tikzcd}
  (\C, \QF) \ar[r] \ar[d] & (\C', \QF') \ar[d]\\
  (\D_0, \QFD_0) \ar[r] \ar[d] & (\D'_0, \QFD'_0) \ar[d]\\
  (\D, \QFD) \ar[r] & (\D', \QFD')
 \end{tikzcd}
\]
where $\D_0 \subseteq \D$ and $\D_0'\subseteq \D'$ denote the essential images of $\C$ and $\C'$, respectively, with $\QFD_0 \coloneq \QFD|_{\D_0}$ and $\QFD'_0 \coloneq \QFD'|_{\D'_0}$.
By Corollaries~\reftwo{corollary:criterion-poincare-karoubi-projection} and \reftwo{corollary:criterion-poincare-karoubi-is-verdier} 
the vertical arrows in the bottom square are Karoubi equivalences and the vertical arrows the top square are then Poincaré-Verdier projections.
We apply $\GW$ and then consider the map between vertical cofibres in the lower square: According to Proposition~\reffour{proposition:weak-cofinality}, this is a map of
$0$-truncated spectra, which on $\pi_0$ is the map
\[
\GW_0(\D,\QFD) / [(\GW_0(\D,\QFD) \times_{\K_0(\D)} \K_0(\D_0)] \to \GW_0(\D',\QFD') / [(\GW_0(\D',\QFD') \times_{\K_0(\D')} \K_0(\D_0')]
\]
which is itself injective by Thomason's classification of dense subcategories and the fact that the lower square is cartesian. %
We conclude that the total cofibre of the lower square (still after applying $\GW$) is 
0-truncated, and the same holds for the total cofibre of the outer square, as the upper square is exact upon applying $\GW$. The outer square thus has 
(-2)-truncated total fibre after applying $\GW$, and is hence cartesian after applying $\GWspace = \Om^{\infty}\GW$.
\end{proof}

The following result lifts $\GWspace(-^\natural)$ to a spectrum valued Karoubi-localising functor.    
Its analogue for $\Catx$ was obtained previously %
by Cisinski--Khan \cite{CK}*{Theorem 4.4.3} using a different proof. 
Let us also point out that, ignoring Poincaré structures everywhere, 
the proof of Proposition~\reffour{proposition:localising-are-stable} applies mutatis mutandis to prove its analogue for $\Catx$ (compare \cite{Efimov}*{Proposition 1.64} and \cite{stable-envelope} in which the case of localising invariants of exact $\infty$-categories is treated.)

\begin{proposition}
\label{proposition:localising-are-stable}%
Let $\A$ be an additive $\infty$-category with finite limits. The $\infty$-category $\Funkloc(\Catp,\A)$ of $\A$-valued Karoubi-localising functors is stable. In particular, the functor
\[
\Omega^\infty \circ (-)\colon \Funkloc(\Catp,\Spa(\A))\to \Funkloc(\Catp,\A)
\]
is an equivalence. 
\end{proposition}

The proof of Proposition~\reffour{proposition:localising-are-stable} will make use of the following lemma, which will also be used later in \S\reffour{subsection:multiplicative-GW}.

\begin{lemma}
\label{lemma:tensor-preserves-karoubi-squares}%
For every $(\C,\QF) \in \Catp$, the functor $(\C,\QF)\otimes (-)\colon \Catp \to \Catp$ preserves Poincaré-Karoubi squares.
\end{lemma}
\begin{proof}
By Proposition~\reftwo{lemma:Morita-exact-criterion}\;\reftwoitem{item:preserves-K-equivalences}
we have that $(\C,\QF)\otimes (-)$ preserves Karoubi equivalences, and so it will suffice to show that it sends Poincaré-Verdier squares to Poincaré-Karoubi squares. Now given a Poincaré-Verdier square on the left, consider the associated square on the right:
\[
\begin{tikzcd}
(\D,\QFD) \ar[r]\ar[d] & (\E,\QFE) \ar[d] \\
(\D',\QFD') \ar[r] & (\E',\QFE')
\end{tikzcd}
\quad\quad\quad\quad
\begin{tikzcd}
(\C,\QF)\otimes^{\flat}(\D,\QFD) \ar[r]\ar[d] & (\C,\QF)\otimes^{\flat}(\E,\QFE) \ar[d] \\
(\C,\QF)\otimes(\D',\QFD') \ar[r] & (\C,\QF)\otimes(\E',\QFE'),
\end{tikzcd}
\]
where $(\C,\QF)\otimes^{\flat}(\D,\QFD) := [(\C,\QF)\otimes(\D',\QFD')] \times_{[(\C,\QF)\otimes(\D',\QFD')]^{\natural}} [(\C,\QF)\otimes(\D,\QFD)]^{\natural}$, and similarly for $(\C,\QF)\otimes^{\flat}(\E,\QFE)$. The square on the right is then Karoubi equivalent to the image under $(\C,\QF)\otimes (-)$ of the square on the left, and it will hence suffice to show that the former is a Poincaré-Verdier square.
Now the vertical legs in the square on the left are Poincaré-Verdier projections by assumption; they hence remain so after tensoring with $(\C,\QF)$ by Proposition~\reftwo{lemma:Morita-exact-criterion}\;\reftwoitem{item:preserves-PV-proj}, and become Poincaré-Karoubi projections after further taking idempotent completions. Since Poincaré-Karoubi projections are preserved under base change (Corollary~\reftwo{corollary:PK-projpullbackstable}) we conclude that the vertical legs in the right square are Poincaré-Karoubi projections. But they are actually Poincaré-Verdier projections since they are essentially surjective (their sources contain the tensor with $(\C,\QF)$ of the vertical sources in the left square). At the same time, by Proposition~\reftwo{lemma:Morita-exact-criterion}\;\reftwoitem{item:sends-PV-squares-to-PK-squares} the induced map on vertical fibres in the right square is obtained from the corresponding one in the left square by applying $[(\C,\QF)\otimes (-)]^{\natural}$, and is hence an equivalence. By Proposition~\reftwo{proposition:poincare-verdier-square-char} it will now suffice to show that the right square is adjointable, or equivalently, that the square
\[
\begin{tikzcd}
\Ind\C\otimes\Ind\D \ar[r]\ar[d] & \Ind\C\otimes\Ind\E \ar[d] \\
\Ind\C\otimes \Ind\D' \ar[r] & \Ind\C \otimes \Ind\E'
\end{tikzcd}
\]
obtained by taking inductive completions is both horizontally and vertically right adjontable, where we have used that $\Ind$ inverts Karoubi equivalences and is symmetric monoidal when viewed as a functor from $\Catx$ to $\Pr^\mathrm L$. Now recall that the right adjoints arising after applying $\Ind$ are colimit preserving, and so the resulting adjunctions in $\Pr^\mathrm L$ are internal ones. It will hence suffice to show that tensoring with $\Ind(\C)$ preserves vertically/horizontally internally right adjointable squares in $\Pr^\mathrm L$. Indeed, this is because $\Ind(\C)\otimes (-)$ is a $\Pr^\mathrm L$-linear endomorphism of $\Pr^\mathrm L$, and hence also $\Cat$-linear.
\end{proof}

\begin{proof}[Proof of Proposition~\reffour{proposition:localising-are-stable}]
The full subcategory $\Funkloc(\Catp,\A) \subseteq \Fun(\Catp,\A)$ is closed under limits, and so the loop functor on this category is given by postcomposition with the loop functor on $\A$. The second part of the proposition hence follows from the first. We now proceed to prove that $\Funkloc(\Catp,\A)$ is stable. Since it has finite limits and is pointed, it suffices to show that $\Omega$ is an equivalence, see \cite{HA}*{Proposition~1.4.2.11(3)}. We will produce an explicit inverse to it.

To this end, let $(\Spa^{\omega_1},\hat{\QF}^{\uni})$ be the hermitian $\infty$-category of $\omega_1$-compact spectra equipped with the Poincaré structure $\hat{\QF}^{\uni}$ right Kan extended along $(\Spaf)\op \subseteq (\Spa^{\omega_1})\op$ from the universal Poincaré structure. Then 
\[
\Om^{\infty}\hat{\QF}^{\uni}\colon (\Spa^{\omega_1})\op \to \Sp
\]
preserves countable cofiltered limits,
and so by~\cite{HTT}*{Lemma 6.3.3.6} the $\infty$-category $\catforms(\Spa^{\omega_1},\hat{\QF}^{\uni})$ admits countable filtered colimits, and by a similar reasoning $\Funh((\Spa^{\omega_1},\hat{\QF}^{\uni}),(\Spa^{\omega_1},\hat{\QF}^{\uni}))$ admits countable filtered colimits. Taking the sequential colimit of the functors $\oplus_{\{1,...n\}}\id$ we hence obtain an endofunctor $S\colon (\Spa^{\omega_1},\hat{\QF}^{\uni}) \to (\Spa^{\omega_1},\hat{\QF}^{\uni})$ which satisfies the
Eilenberg swindle (in the semi-additive $\infty$-category $\Cath$), 
that is, $S\oplus \id \simeq S$. Applying the right adjoint $\Pairings\colon \Cath \to \Catp$, we find that the Poincar\'e category $\Latt^{\omega_1}(\Spaf, \QF^{\uni}) \coloneq \Pairings(\Spa^{\omega_1},\hat{\QF}^{\uni})$ supports an Eilenberg swindle in the semi-additive $\infty$-category $\Catp$. For any Karoubi-localising functor $\FF\colon \Catp\to \A$ we therefore have $\FF(\Latt^{\omega_1}(\Spaf, \QF^{\uni}))\simeq 0$ (since $\A$ is assumed additive). Now the composite
\[
(\Spaf,\QF^\uni) \to\Pairings(\Spaf,\QF^\uni)\to \Pairings(\Spa^{\omega_1},\hat{\QF}^{\uni}) = \Latt^{\omega_1}(\Spaf, \QF^{\uni}),
\]
is a Poincaré-Verdier inclusion (see Remark~\refone{remark:unit-fully-faithful}) and the corresponding Poincaré-Verdier quotient %
\[
\Tate^{\omega_1}(\Spaf, (\QF^{\uni})\qshift{1}) \coloneq \Latt^{\omega_1}(\Spaf, \QF^{\uni})/(\Spaf,\QF^\uni)
\]
is an $\omega_1$-bounded variant of the Tate category, see \S\reftwo{subsection:poincare-verdier-sequences}. %
Now for any Poincaré $\infty$-category $(\C,\QF)$ the sequence
\[
(\C,\QF) = (\C,\QF) \otimes (\Spaf,\QF^{\uni}) \to (\C,\QF) \otimes \Latt^{\omega_1}(\Spaf, \QF^{\uni}) \to (\C,\QF) \otimes \Tate^{\omega_1}(\Spaf, (\QF^{\uni})\qshift{1})
\]
is a Poincaré-Karoubi sequence by Proposition~\reftwo{lemma:Morita-exact-criterion}\;\reftwoitem{item:sends-PV-squares-to-PK-squares}.
As $(\C,\QF) \otimes \Latt^{\omega_1}(\Spaf, \QF^{\uni})$ inherits from $\Latt^{\omega_1}(\Spaf, \QF^{\uni})$ an Eilenberg swindle endomorphism
we obtain an equivalence
\[
\FF\simeq \Omega\FF((-) \otimes \Tate^{\omega_1}(\Spaf,(\QF^{\uni})\qshift{1}))
\]
of functors $\Catp\rightarrow \A$. %
By Lemma~\reffour{lemma:tensor-preserves-karoubi-squares} tensoring with $\Tate^{\omega_1}(\Spaf,(\QF^{\uni})\qshift{1})$ preserves Poincar\'e-Karoubi squares, 
and so precomposing with it defines an endofunctor of $\Funkloc(\Catp,\A)$, inverse to $\Omega$.

\end{proof}

We now arrive to the main definition of this paper, which is afforded by Proposition~\reffour{proposition:localising-are-stable} and Example \reffour{example:GWspace_Karoubi_localising}.

\begin{definition}
\label{definition:KGW}%
We define the \defi{Karoubi-Grothendieck-Witt} spectrum $\KGW\colon \Catp \to \Spa$ to be the essentially unique Karoubi-localising functor whose infinite loop space is $\GWspace(-^\natural)$. 
\end{definition}

Inspecting the proof of Proposition~\reffour{proposition:localising-are-stable} we see that the $\KGW(\C,\QF)$ can be represented by the $\Omega$-spectrum object whose $n$'th space is $\GWspace((\Tate^{\omega_1}(\Spaf,(\QF^{\uni})\qshift{1})^{\otimes n}\otimes (\C,\QF))^\natural)$. This Lattice-Tate construction of $\KGW$ can be viewed as a Poincaré categorical analogue of the cone-suspension ring approach of Karoubi and Schlichting. In fact, the latter is recovered by the present construction in the case of rings, as shown by the following remark.

\begin{remark}
Recall from~\cite{comparison}*{Theorem B} that for a ring $R$ equipped with an invertible module with involution $M$, the genuine symmetric Grothendieck-Witt space $\GWspace^{\gs}(R,M)$ recovers the classical symmetric Grothendieck-Witt space $\GWspace_{\cl}^{\sym}(R,M)$. 
We claim that this comparison extends to Karoubi-Grothendieck-Witt spectra. To see this, recall the classical construction of the Karoubi-Grothendieck spectrum by Karoubi and Schlichting using cone and suspension rings, see e.g.\ \cite{schlichting-mv}*{\S 10.1}.
Here, the cone ring $C(R)$ of $R$ is the infinite matrix ring $C(R)$ consisting of the $R$-valued $\NN\times \NN$-matrices each of whose rows and each of whose columns has finitely many non-zero entries, and the suspension ring $S(R)$ is its quotient by the idempotent two-sided ideal $I(R) \subseteq C(R)$ consisting of the matrices which globally have finitely many non-zero entries. Then $S(R)$ is an Ore-localisation of $C(R)$ and consequently $S(R)$ and $I(R)$ are flat $C(R)$-modules. In particular, $I(R)$ is also idempotent in the derived sense and $C(R) \to S(R)$ is a derived localisation.
The rings $C(R)$ and $S(R)$ inherit invertible modules with involution $C(M)$ and $S(M)$, respectively, where $C(M)$ is given by $M$-valued $\NN\times\NN$-matrices with the same vanishing properties as $C(R)$, and $S(M) = C(M)/I(M)$ where $I(M)$ is spanned by the $\NN\times \NN$-matrices with the same vanishing properties as $I(R)$, the involution being induced in both cases by transposition followed by the pointwise involution of $M$.
Writing $P = \oplus_{i \in \NN} R$ for the $C(R)$-projective $(C(R),R)$-bimodule of column vectors one then constructs a sequence of duality preserving functors
\[
(\Proj(R),\Dual_M) \xrightarrow{P\otimes_R (-)} (\Proj(C(R)),\Dual_{C(M)}) \longrightarrow (\Proj(S(R)),\Dual_{S(M)}) ,
\]
which one can verify to be a fibre sequence of additive categories with duality, using the fact that $P^{\vee}\otimes_{C(R)} P \cong R$ and $P\otimes_R P^{\vee} \cong I(R)$, where $P^{\vee} := \Hom_{C(R)}(P,C(R))$ can be identified with the $(R,C(R))$-bimodule of row vectors. By~\cite{schlichting-mv}*{Theorem 10.2, Lemma 10.7} this sequence induces a fibre sequence with contractible middle term on classical symmetric $\GW$-spaces, yielding an equivalence 
\[
\GWspace_{\cl}^{\sym}(R,M) \xrightarrow{\sim} \Om\GWspace_{\cl}^{\sym}(S(R),S(M)).
\]
Iterating this construction produces the connecting maps in a spectrum, which we will refer to as the classical symmetric Karoubi-Grothendieck-Witt spectrum $\KGW_{\cl}^{\sym}(R,M)$ of $(R,M)$. To obtain an equivalence with the Karoubi-Grothendieck-Witt spectrum of Definition~\reffour{definition:KGW} we simply add that, in light of Proposition~\refone{proposition:restriction-of-2-excisive-functors}, the construction of the above fibre sequence of additive categories works the same to produce a fibre sequence
\[
\bigl(\Dperf(R),\QF^{\gs}_M\bigr) \xrightarrow{P\otimes_R (-)} \bigl(\Dperf(C(R)),\QF^{\gs}_{C(M)}\bigr) \longrightarrow \bigl(\Dperf(S(R)),\QF^{\gs}_{S(M)}\bigr)
\]
of Poincaré $\infty$-categories. It then follows from Proposition~\reftwo{proposition:left-Kan-extension-specific-structures} that it is in fact a Poincar\'e-Karoubi sequence as one may apply Example~\reftwo{example:involution-induction-compatible}~\reftwoitem{item:ore-compatible} to see that $C(M)$ is compatible with the derived localisation $C(R) \to S(R)$.
Similar to \cite{schlichting-mv}*{Lemma 10.7}, by Eilenberg swindle one can show that $\KGW^\gs(C(R),C(M))$ vanishes. Hence, the above sequence determines an equivalence
\[
\Om^{\infty-1}\KGW^{\gs}(R,M) \simeq \Om^{\infty}\KGW^{\gs}(S(R),S(M)) = \GWspace^{\gs}(S(R),S(M)).
\]
In particular, we obtain an equivalence $\KGW^{\gs}(R,M) \simeq \KGW_{\cl}^{\sym}(R,M)$
as claimed. A similar argument applies for the quadratic and even form parameters.
\end{remark}

Recall the universal properties of $\GWspace = \Om^{\infty}\GW$ and $\GW$ established in \papertwo, namely, $\GWspace$ is the initial group-like additive functor under $\Poinc$ (Observation~\reftwo{observation:univ-GW-space}) and $\GW$ is the initial additive functor to spectra equipped with a map $\Poinc \Rightarrow \Omega^\infty\GW$ (Corollary~\reftwo{corollary:univ-GW-spectrum}). %
Since $\KGW$ is in particular additive the composite transformation $\Poinc \Rightarrow \GWspace\Rightarrow \GWspace(-^\natural) = \Omega^\infty \KGW$ determines a map
\[
\GW\Rightarrow \KGW
\]
which agrees with the transformation $\GWspace\Rightarrow\GWspace(-^\natural)$ on infinite loop spaces.
In particular, by Proposition~\reffour{proposition:weak-cofinality} the map
\[
\pi_n \GW \to \pi_n \KGW
\]
is an isomorphism for $n>0$ and a monomorphim for $n=0$. 

\begin{proposition}
\label{proposition:universal-KGW}%
The transformations $\SS[\Poinc] \Rightarrow \GW\Rightarrow \KGW$ exhibit $\KGW$ as the initial Karoubi-localising functor under $\SS[\Poinc]$ and $\GW$.
\end{proposition}

The proof of Proposition~\reffour{proposition:universal-KGW} will make use of the following lemma:
\begin{lemma}
\label{lemma:idem-invariant}%
Let $\F\colon \Catp \to \E$ be a functor. Then $\F(-^\natural)$ is the initial Karoubi invariant functor under $\F$.
\end{lemma}
\begin{proof}
The inclusion of Karoubi invariant functors inside all functors can be realised as the restriction along the localisation $\Catp \to\Catpi$ of $\Catp$ by the Karoubi equivalences. This localisation is left Bousfield with fully-faithful right adjoint given by the inclusion $\Catpi \subseteq \Catp$ as idempotent complete $\infty$-categories. Pre-composition along this right adjoint hence gives a left adjoint to the inclusion of Karoubi invariant functors, as desired.
\end{proof}

\begin{proof}[Proof of Proposition~\reffour{proposition:universal-KGW}]
The two statements are equivalent since the map $\SS[\Poinc] \Rightarrow \GW$ induces an equivalence on mapping spaces into any additive functor (and in particular into any Karoubi-localising functor) by the universal property of $\GW$. It will hence suffice to show that $\KGW$ is the initial Karoubi-localising functor under $\SS[\Poinc]$. By adjunction, this is equivalent to showing that for every Karoubi-localising functor $\GG \colon \Catp \to \Sp$ the composite functor
\[
\Nat(\KGW,\GG) \to \Nat(\Om^{\infty}\KGW,\Om^{\infty}\GG) \to \Nat(\GWspace,\Om^{\infty}\GG) \to \Nat(\Poinc,\Om^{\infty}\GG)
\]
is an equivalence. Here, we note that for the second and third term, it does not matter if we consider the two functors as additive functors to $\Grp_{\Einf}$ or as group-like additive functors to $\Sps$, since the forgetful functor between these two types of functors is an equivalence. The first map is then an equivalence by Proposition~\reffour{proposition:localising-are-stable} applied with $\A = \Grp_{\Einf}$, the second is an equivalence by Lemma~\reffour{lemma:idem-invariant} since $\Om^{\infty}\KGW = \GWspace(-^{\natural})$, and the third is an equivalence by the universal property of $\GWspace$.
\end{proof}

We write 
\[
\KK \coloneq \KGW^{\hyp} \colon \Catp \to \Fun(\B\Ct, \Sp)
\]
for the functor obtained by pointwise composing $\KGW$ with $\Hyp\colon \Catp\to \Fun(\B\Ct, \Catp)$. By abuse of notation we will also use this symbol for the underlying functor $\Catx \to \Sp$, which agrees with non-connective, or Bass $\K$-theory: By the $\Catx$-analogue of Proposition~\reffour{proposition:localising-are-stable} established by Cisinski--Khan, it suffices to show that their infinite loop spaces agree, and both are equivalent to $\Kspace(-^\natural)$.

\begin{proposition}
\label{proposition:universal-KK}%
The transformations $\SS[\core] \Rightarrow \K \Rightarrow \KK$ of functors $\Catp\to \Spa$ exhibit $\KK$ as the initial Karoubi-localising functor under $\SS[\core]$ and $\K$. The same statement holds with $\Catp$ replaced by $\Catx$.
\end{proposition}
\begin{proof}
As in the proof of Proposition~\reffour{proposition:universal-KGW}, the claim for $\Catp$ follows from Proposition~\reffour{proposition:localising-are-stable}
and Lemma~\reffour{lemma:idem-invariant} using the universal properties of $\K$ and $\Kspace$ (which in this setting follow from the universal properties of $\GW$ and $\GWspace$; indeed, since $\Hyp\colon \Catx \to \Catp$ is both left and right adjoint to the forgetful functor $U\colon \Catp \to \Catx$ by Corollary~\refone{corollary:hyp-is-adjoint} we have that $\Hyp\circ U$ is both left and right adjoint to itself, and since pre-composition with $\Hyp\circ U$ preserves additive functors it also preserves universal additive approximations). The analogous claim %
for $\Catx$ can be obtained by a similar argument using the corresponding statements for $\Catx$ (and in particular \cite{CK}*{Theorem 4.4.3}).
\end{proof}

\medskip 

We now turn to the Karoubi-invariant version of $\L$-theory. Recall from Lemma~\reftwo{lemma:criterion}
that a spectrum-valued additive functor is bordism invariant if and only if it vanishes on hyperbolic (or equivalently metabolic) Poincaré $\infty$-categories, and that the inclusion of spectrum-valued bordism invariant additive functors inside all spectrum-valued additive functors admits both a left and a right adjoint, called bordification and cobordification respectively, see Theorem~\reftwo{theorem:localisation}. The map $\GW \Rightarrow \L$ then exhibits $\L$ as the bordification of $\GW$.

\begin{definition}
We define the \emph{Karoubi L-theory functor} $\KL\colon \Catp\to\Spa$ to be the bordification of $\KGW$.
\end{definition}

\begin{lemma}
\label{lemma:bordification-again-KL}%
Let $\A$ be a (co)complete stable $\infty$-category. Then the (co)bordification %
of a Karoubi-localising functor $\FF\colon \Catp \to \A$ is again Karoubi-localising. In particular, $\KL$ is Karoubi-localising. %
\end{lemma}
\begin{proof}
According to Proposition~\reftwo{proposition:phi-bord} (whose proof carries over verbatim from $\A=\Sp$ to an arbitrary stable $\infty$-category $\A$ equipped with the relevant (co)limits), an explicit model for the bordification of $\FF$ is the cofibre of the canonical map $\FF^\hyp(-)_\hC\to \FF$, and an explicit model for the cobordification of $\FF$ is the fibre of the canonical map $\FF\to \FF^\hyp(-)^\hC$. Since $\Hyp\colon \Catp\to \Catp$ preserves Karoubi sequences, both of them have the required property. %
\end{proof}

By the universal property of $\L$, the composite functor $\SS[\Poinc] \Rightarrow \KGW \Rightarrow \KL$ factors through $\L$.
\begin{proposition}
\label{proposition:universal-KL}%
The square
\[
\begin{tikzcd}
\SS[\Poinc]\ar[r]\ar[d] & \KGW \ar[d] \\
\L \ar[r] & \KL
\end{tikzcd}
\]
exhibits $\KL$ as initial in three manners, namely, as a bordism invariant functor under $\KGW$, as a Karoubi-localising functor under $\L$, and as a Karoubi-localising and bordism invariant functor under $\SS[\Poinc]$.
\end{proposition}

\begin{proof}
The first of these manners is by definition, and the last follows from the first by
the universal property of $\KGW$ from Proposition~\reffour{proposition:universal-KGW}. %
Combining this with the universal property of $\L$ we conclude that the map $\L \Rightarrow \KL$ exhibits $\KL$ as the initial bordism invariant and Karoubi-localising functor under $\L$. Since cobordification preserves Karoubi-localising functors by Lemma~\reffour{lemma:bordification-again-KL} we have that $\KL$ is also initial as a Karoubi-localising functor under $\L$.
\end{proof}

Applying Corollary~\reftwo{corollary:tate-square} we now find that $\KGW$ and $\KL$ fit into the Karoubi-localising avatar of the fundamental fibre square from Corollary~\reftwo{corollary:tate-square-L}.
\begin{corollary}
\label{corollary:Karoubi-fund-square}%
There is an exact square 
\begin{equation}
\label{equation:fundamental-karoubi}%
 \begin{tikzcd}
  \KGW \ar[r,"\bord"] \ar[d,"\fgt"'] & \KL \ar[d,"\Xi"] \\
  \KK^{\hC} \ar[r,"\mathrm{can}"] & \KK^{\tC},
 \end{tikzcd}
\end{equation}
of Karoubi-localising functors $\Catp \rightarrow \Spa$. 
\end{corollary}

\begin{theorem}[Karoubi cofinality]
\label{theorem:Karoubi-pullback-squares}%
The following two squares are cartesian.
\[
\begin{tikzcd}
\GW(\C,\QF)\ar[r]\ar[d] & \KGW(\C,\QF)\ar[d] & \L(\C) \ar[r] \ar[d] & \KL(\C,\QF) \ar[d] \\
\K(\C,\QF)^\hC\ar[r] & \KK(\C,\QF)^\hC & \K(\C,\QF)^\tC \ar[r] & \KK(\C,\QF)^{\tC}.
\end{tikzcd}
\]
In particular, we have 
\begin{enumerate}
\item $\pi_*\cof(\GW(\C,\QF)[\frac{1}{2}] \to \KGW(\C,\QF)[\frac{1}{2}]) \cong 
\begin{cases}
0 & \text{if } *  >0 \\
\big(\KK_0(\C,\QF) / \K_0(\C,\QF)\big) [\frac{1}{2}]^{\Ct} & \text{if } * =0 \\
\KK_*(\C,\QF)[\frac{1}{2}]^{\Ct} & \text{if } *  <0
\end{cases} $
and
\item the map $\L(\C,\QF)[\frac{1}{2}] \to \KL(\C,\QF)[\frac{1}{2}]$ is an equivalence.
\end{enumerate}
\end{theorem}
\begin{proof}
Arguing as in the proof of Theorem~\reffour{theorem:strong-cofinality}, consider the commutative cube
\[
\begin{tikzcd}
[column sep=0.5em, row sep=0.7em]
  & \GW(\C,\QF)\ar[rr]\ar[dd] \ar[ld]&& \KGW(\C,\QF)\ar[dd]\ar[ld]\\
\L(\C,\QF)\ar[rr, crossing over]\ar[dd] && \KL(\C,\QF)\\
  & \K(\C,\QF)^\hC\ar[rr]\ar[ld] && \KK(\C,\QF)^\hC\ar[ld]\\
\K(\C,\QF)^\tC\ar[rr] && \KK(\C,\QF)^\tC\ar[uu,<-, crossing over]\\
\end{tikzcd}
\]
whose back and front faces are the two squares of the theorem and whose left and right faces are the fundamental fibre squares for $\GW$ and $\KGW$ respectively. %
The total fibres of the back and front faces are hence equivalent, and we may write $T(\C,\QF)$ for the common value. Then $T(-)$ is a bordism invariant functor since the front face consists of such functors and $T(-)$ takes values in $(-1)$-truncated spectra since the square of infinite loop spaces underlying the back face is the square
\[
\begin{tikzcd}
\GWspace(\C,\QF)\ar[r]\ar[d] & \GWspace(\C^{\natural},\QF)\ar[d]  \\
\Kspace(\C,\QF)^\hC\ar[r] & \Kspace(\C^{\natural},\QF)^\hC
\end{tikzcd}
\]
whose horizontal fibres are $(-1)$-truncated by Propositions~\reffour{proposition:weak-cofinality} and~\reffour{proposition:cofinality-K-theory}. We conclude that $T=0$ as desired.

To obtain the ``in particular'', we argue as follows. Put $D=\cof(\K(\C,\QF)\to\KK(\C,\QF))$.
By Proposition~\reffour{proposition:cofinality-K-theory},
$D$ is coconnective, with
$\pi_0 D\cong\KK_0(\C,\QF)/\K_0(\C,\QF)$ and
$\pi_i D\cong\KK_i(\C,\QF)$ for $i<0$.
The two cartesian squares identify the cofibres of
$\GW(\C,\QF)\to\KGW(\C,\QF)$ and
$\L(\C,\QF)\to\KL(\C,\QF)$ with
$D^{\hC}$ and $D^{\tC}$, respectively.
Now we recall that $(-)^{\hC}$ commutes with filtered colimits of uniformly bounded above spectra with $\Ct$-action. For spectra with a single homotopy group, this boils down to the commutation of $\Ct$-group cohomology with filtered colimits, which can for instance be deduced from its identification with shifted $\Ct$-group homology. Inductively over the steps in the Whitehead tower, one obtains the desired result for bounded spectra and by passing to the colimit over the Whitehead tower using \cite{NS}*{Lemma I.2.6}, one obtains the desired result for bounded above spectra.
Since $(-)_{\hC}$ commutes with all colimits, we deduce that the canonical comparison maps 
\[
D^{\hC}[\tfrac{1}{2}]
\longrightarrow (D[\tfrac{1}{2}])^{\hC}
\quad\text{and}\quad
D^{\tC}[\tfrac{1}{2}]
\longrightarrow (D[\tfrac{1}{2}])^{\tC}
\]
are equivalences, so the desired result follows from the homotopy fixed point and Tate spectral sequences, respectively.
\end{proof}

\begin{corollary}
\label{corollary:regular-noetherian}%
Let $A$ be a connective ring spectrum such that $\pi_0(A)$ is regular noetherian. Then the map
$\K(A) \to \KK(A)$ is an equivalence. As a consequence, for any Poincaré structure $\QF$ on $\Modp{A}$ the maps
\[
\GW(\Modp{A},\QF) \to \KGW(\Modp{A},\QF) \quad\text{and}\quad \L(\Modp{A},\QF) \to \KL(\Modp{A},\QF)
\]
are equivalences as well.
\end{corollary}
\begin{proof}
By Theorem \reffour{theorem:Karoubi-pullback-squares} the claim for Grothendieck-Witt- and $\L$-theory follows from the claim for $\K$-theory.
The latter should be considered as well-known, and follows for example from combining the following two claims:
\begin{enumerate}
\item
\label{item:a}%
For every connective ring spectrum $A$ and every $n\geq 0$, the map $\K_{-n}(A) \to \K_{-n}(\pi_0(A))$ is an isomorphism. %
\item
\label{item:b}%
For every regular noetherian ring $R$ and all $n> 0$, one has $\K_{-n}(R) = 0$. %
\end{enumerate}
The first statement above appears for example in \cite{BGT}*{Theorem 9.53}. The second is due to Bass. We include a proof for the convenience of the reader.

To see \reffouritem{item:a}, we note that the Bass-Heller-Swan decomposition (see~\cite{saunier}*{Theorem 1.3} for a categorical version that is applicable in particular to ring spectra, or the next section for a proof in the Poincare setting which equally applies to stable categories) %
implies that the map $\K_{-n}(A) \to \K_{-n}(\pi_0(A))$ is a retract of the map $\K_0(A[\ZZ^n]) \to \K_0(\pi_0(A)[\ZZ^n]) \cong \K_0(\pi_0(A[\ZZ^n]))$, so it suffices to prove the case $n=0$, which follows immediately from $A$ and $\pi_0(A)$ having the same projective modules up to isomorphism, see e.g.\ \cite{Lurie-K-theory}*{Corollary 3}. For \reffouritem{item:b}, again using the Bass-Heller-Swan decomposition
we find that $\K_{-n}(R)$ is a retract of $\K_{-1}(R[\ZZ^{n-1}])$. Since with $R$ being so, also $R[\ZZ^{n-1}]$ is regular noetherian, 
and so it suffices to prove the claim for $n=1$. But another application of Bass-Heller-Swan %
shows that 
\[
\K_{-1}(R) \cong \mathrm{coker}\big[\K_0(R[t]) \oplus \K_0(R[t]) \xrightarrow{t \mapsto (t,t^{-1})} \K_0(R[t^{\pm 1}])\big]
\]
(this is of course Bass' original definition of negative $\K$-groups) and since $R[t]$ is also regular noetherian, 
it suffices to show that $\K_0(R[t]) \to \K_0(R[t^{\pm1}])$ is surjective. This is true more generally for any localisation $S \rightarrow S[f^{-1}]$ of a regular coherent ring $S$ at a central element $f \in S$: it clearly suffices to lift all finitely generated projective modules over $S[f^{-1}]$ to perfect modules over $S$. But for a regular coherent ring all finitely presented module are perfect, so such a lift can be given simply by lifting a presentation and clearing denominators, see \papertwo Example 4.5.2.
\end{proof}

\begin{corollary}
\label{corollary:filtcolimpreskgwspectrum}%
The functors $\KGW$ and $\KL$ preserve filtered colimits.
\end{corollary}
\begin{proof}
This follows from Karoubi cofinality (Theorem~\reffour{theorem:Karoubi-pullback-squares}), using that the functors $[\KK(-)/\K(-)]^{\hC}$ and $[\KK(-)/\K(-)]^{\tC}$ preserve filtered colimits (see the proof of loc.\ cit.) and the earlier established facts that the functors $\GW$ and $\L$ preserve filtered colimits, see Corollaries~\reftwo{corollary:filtcolimpresgwspectrum} and~\reftwo{corollary:Ladd}.
\end{proof}

We close with a simple application of Karoubi cofinality regarding the commutation of Grothendieck-Witt and $\L$-spectra with infinite products. Recall that both $\K, \KK \colon \Catx \rightarrow \Spa$ do commute with arbitrary products; this was originally established by Carlsson in the context of exact categories in \cite{carlsson-infinite} and for stable $\infty$-categories in \cite{kaswing} and with a more conceptual proof again in \cite{Efimov}.

\begin{corollary}
The functors $\GW,\KGW,\L \colon \Catp \rightarrow \Sp$ commute with arbitrary products, and $\KL \colon \Catp \rightarrow \Sp$ commutes with products over families $\C_i$ with $\KK(\C_i)$ uniformly bounded below.
\end{corollary}

\begin{proof}
That $\L_i \colon \Catp \rightarrow \Ab$ commutes with products follows immediately from the cycle description of $\L_i(\C,\QF)$ as Poincar\'e objects in $(\C,\QF\qshift{i})$ modulo metabolics. The claim for $\L$ thus follows from homotopy groups commuting with products. 

From the fundamental fibre sequence $\K(\C)_\hC \rightarrow \GW(\C,\QF) \rightarrow \L(\C,\QF)$ it then follows that $\GW$ commutes with infinite products if $\K(-)_\hC$ does. Dually to the earlier argument on the commutation of $(-)^{\hC}$ with filtered colimits, the homotopy orbits functor $(-)_\hC$ commutes with cofiltered limits, and consequently arbitrary products, of uniformly bounded below collections of spectra. The claim for $\GW$ thus follows from the one for $\K$. 

For $\KGW$ the statement is now a direct consequence of the left cartesian square in Theorem~\reffour{theorem:Karoubi-pullback-squares}. The claim about $\KL$ follows similarly from the right square in Theorem~\reffour{theorem:Karoubi-pullback-squares} since with $(-)_{\hC}$ also $(-)^{\tC}$ commutes with infinite products of uniformly bounded below spectra. 
\end{proof}

\begin{remark}
The functor $\KL \colon \Catp \rightarrow \Sp$ does not commute with arbitrary infinite products as we show next. Using the identification $\KL = \L^{\cotower{-\infty}}$ that we will establish in Proposition~\reffour{proposition:decorated-L} below, a counterexample was first constructed by Winges in \cite{winges-Ltheory}. 

In fact, the Karoubi-localising fundamental fibre square~\eqreffour{equation:fundamental-karoubi}
implies that $\KL(\prod_{i \in I} (\C_i,\QF_i)) \to \prod_{i \in I} \KL(\C_i,\QF_i)$ is an equivalence if and only if the map 
\[
\left[ \prod_{i \in I} \KK(\C_i,\QF_i)\right]^{\tC} \to \prod_{i \in I} \KK(\C_i,\QF_i)^{\tC}
\]
is one. A different (possibly simpler) example than the one provided by Winges where this is not the case is the following: Let $\mathscr{O}_\infty$ denote the infinite Cuntz algebra, a complex $\mathrm C^*$-algebra (and thus in particular a ring with involution) with $\KK(\mathscr{O}_\infty) \simeq \K^\mathrm{top}(\mathscr{O}_\infty) \simeq \KU$ with trivial $\Ct$-action; this equivalence has first been observed by Cortinas--Phillips \cite{CP}, see e.g.\ \cite{BLP}*{\S 5.3} for some general explanations. Considering then a countable product of $(\Dperf(\mathscr{O}_\infty),\QF^\sym_{\mathscr O_\infty})$ we need to show that the canonical map 
\[
\left[ \prod_{\mathbb{N}} \KU\right]^{\tC} \to \prod_{\mathbb{N}} \KU^{\tC}
\]
is not an equivalence. To that end, we recall that $\pi_*(\KU^\tC) \cong \mathbb{Q}_2[t^{\pm 1}]$, see e.g.\ \cite{NS}*{Example I.2.3 (iii)} with $|t|=-2$ and $\beta \cdot t = 2$ where $\beta$ denotes the image of the Bott element in $\pi_2(\KU)$ under the canonical map $\KU \to \KU^{\tC}$. Moreover, we have $\pi_*(\KU^{\hC}) = \KU_*[[t]]/(t^2-2t)$ and for any $\KU$-module $M$, the map $M^{\hC} \to M^{\tC}$ is a localisation at $t$. 
Hence, the target of the above comparison map evaluates to $\prod_{\mathbb{N}} \mathbb{Q}_2[t^{\pm 1}]$ while the source is the localisation of $\prod_{\mathbb{N}} \KU_*[[t]]/(t^2-2t)$ at $t$ and the sequence $\{\frac{1}{2^n} = \frac{1}{\beta^nt^n}\}_{n \in \mathbb{N}} \in \prod_\mathbb{N} \mathbb{Q}_2[t^{\pm 1}]$ does not lie in the image of the comparison map. 
\end{remark}

\section{The circle decomposition theorem}
\label{section:shaneson}%
In this section, we state and prove the circle decomposition theorem, Theorem~\reffour{theorem:intro-F} from the introduction, which describes the behaviour of Verdier-localising functors under the operation $(\C,\QF) \mapsto (\C,\QF)_{\B\ZZ}$ of tensoring with the circle. Our framework reproduces both the Bass-Heller-Swan and the Ranicki-Shaneson decompositions in a unified manner. 
Our main result is stated and put into perspective in \S\reffour{subsection:shaneson} and proven in \S\reffour{subsection:proof}. In \S\reffour{subsection:splittings}, we revisit Ranicki's observation about the difference between the algebraic and geometric Ranicki-Shaneson decompositions in the setting of Poincaré $\infty$-categories.

\subsection{Statement of the circle decomposition theorem}
\label{subsection:shaneson}%

To explain the precise statement of the circle decomposition theorem, we recall that for a monoid $M$, we may form the pointed category $\B M$ with a single object whose endomorphisms are $M$, naturally in the monoid $M$. When $M$ is a group we will also view $\B M$ as a pointed space. The examples we have in mind are the monoid $\NN$ of natural numbers and the group $\ZZ$ of the integers.
In particular, 
for a stable $\infty$-category we may form $\C_{\B\NN}$ and $\C_{\B\ZZ}$, the tensor of $\C$ with $\B\NN$ and $\B\ZZ$, respectively. The monoid inclusion $\NN \hrar \ZZ$ then determines a map $q\colon\B\NN \to \B\ZZ$ and hence an exact functor $\C_{\B\NN} \to \C_{\B\ZZ}$ for every stable $\infty$-category $\C$. 
Similarly, we may form these tensors with $(\C,\QF)$ a hermitian $\infty$-category and if $(\C,\QF)$ is Poincar\'e, then so is $(\C,\QF)_{\B\ZZ}$, while $(\C,\QF)_{\B\NN}$ typically is not. Finally, we note that the functors $\ast \to \B\NN \to \B\ZZ \to \ast$ make $\C$ and $(\C,\QF)$ a retract of $\C_{\B\NN}$ and $(\C,\QF)_{\B\ZZ}$, respectively. To keep in mind the effect of tensoring with $\B\ZZ$, we recall the following.

\begin{example}
\label{example:tensor-BZ}%
Suppose $A$ is an $\Eone$-ring, $(M,N,\alpha \colon N \to M^{\tC})$ a module with genuine involution over $A$, and $(\Mod^\free(A),\QF^\alpha_M)$ the associated Poincar\'e $\infty$-category. Then we have 
\[
(\Mod^\free(A),\QF^\alpha_M)_{\B\ZZ} \simeq (\Mod^\free(A[\ZZ]),\QF^{\alpha[\ZZ]}_{M[\ZZ]})
\]
where $\alpha[\ZZ]$ is the map $N \simeq N \otimes \SS \to M^{\tC} \otimes \SS[\ZZ]^{\tC} \to M[\ZZ]^{\tC}$, see Example \refone{example:tensor-of-rings}. In particular, we have:
\begin{enumerate}
\item
\label{item:ex-tensor-BZ-genuine-family}%
If $(R,\sigma)$ is a ring with involution, then $(\Dfree(R),\QF^{\geq m}_\sigma)_{\B\ZZ} \simeq (\Dfree(R[\ZZ]),\QF^{\geq m}_{\bar{\sigma}})$, where $\bar{\sigma}$ is the unique extension of the involution $\sigma$ on $R$ which is the inversion map on $\ZZ$; under the canonical isomorphism $R[\ZZ] \cong R[t^{\pm 1}]$, this means $\bar{\sigma}(t) = t^{-1}$.
\item 
For a ring $R$ with form parameter $\lambda=[M_{\Ct} \to Q \to M^{\Ct}]$ and associated Poincar\'e structure $\QF^{\g\lambda}_M$, we have an equivalence
\[
(\Dfree(R),\QF_M^{\g\lambda})_{\mathrm{B}\mathbb Z}  \simeq (\Dfree(R[\ZZ]),\QF_{M[\ZZ]}^{\g\lambda[\ZZ]})
\]
where $\lambda[\ZZ]$ is the form parameter on $R[\ZZ]$ whose associated module with involution is $M[\ZZ] = \oplus_{n \in \ZZ}M$ as above and quadratic module given by $Q \oplus \big[\oplus_{n > 0}M\big]$. While the operation $\lambda \mapsto \lambda[\ZZ]$ preserves quadratic, even, and symmetric form parameters %
by \reffouritem{item:ex-tensor-BZ-genuine-family} above, there are other classes of form parameters which do not obey such a compatibility. For instance, if $(R,\sig)$ is a commutative ring with involution and we consider the Frobenius form parameter $\lambda = R^\fr$ as in Example~\reffour{example:form_parameter_over_Tate_Poincare}, then $\lambda[\ZZ]$ is not the Frobenius form parameter on $R[\ZZ]$, unless $R$ is zero. Otherwise put, on the level of associated Poincaré structures we have that 
$\QF^{\g\fr[\ZZ]}_{\sig[\ZZ]} \not{\simeq} \QF^{\g\fr}_{\sigma[\ZZ]}$ 
as Poincaré structures on $\Dfree(R[\ZZ])$ %
(indeed, the $\pi_0$ of the $R[\ZZ]$-complex representing the linear part is $R$ in the former case and $R[\ZZ]$ in the latter, see Example~\reffour{example:modules-over-genuine} below).
\item For a spherical fibration $\xi \colon X \rightarrow \Pic(\SS)$ on a space $X$ and the associated visible Poincar\'e category $(\Sp_X^\free,\QF^\vis_\xi)$ we have 
$(\Sp_X^\free,\QF^\vis_\xi)_{\B\ZZ} \simeq (\Sp_{X\times S^1}^\free,\QF^\vis_{p^*(\xi)})$
where $p\colon X \times S^1 \to X$ is the projection.
\end{enumerate}
\end{example}

\begin{construction}
\label{construction:nil}%
\label{remark:splitting}%
Let $\G\colon \Catx \to \A$ be a functor from $\Catx$ to a stable $\infty$-category $\A$. We define
$\N\G\colon \Catx \to \A$ to be the functor given by the cofibre of the map
$\G(\C) \to \G(\C_{\B \NN})$ or equivalently the fibre of the map $\G(\C_{\B\NN}) \to \G(\C)$ and obtain a natural splitting 
\[
\G(\C_{\B\NN}) \simeq \G(\C) \oplus \N\G(\C)
\]
and refer to it as the \defi{nil term} of $\G$. The terminology comes from the well-known fact that for $\G$ Karoubi-localising, $\N\G(\C)$ is also (a shift of) the fibre $\G(\mathrm{Nil}(\C)) \to \G(\C)$ where $\mathrm{Nil}(\C)$ is the category of nilpotent endomorphisms on $\C$, see e.g.\ \cite{LT2}*{Corollary 4.5 \& Prop.\ 4.7}.
In particular, if $\F\colon \Catp \to \A$ is a functor defined on Poincaré $\infty$-categories, we may form $\N\F^{\hyp}$, where $\F^\hyp = \F \circ \Hyp\colon \Catx \to \A$.
\end{construction}

\begin{definition}
\label{definition:assembly}%
The \defi{assembly map} for $\B\ZZ$ is the canonical colimit-exchange map associated to the constant diagram on $(\C,\QF)$ indexed on $\B\ZZ$:
\[
{\B\ZZ} \otimes \F(\C,\QF) := \colim\limits_{\B\ZZ} \F(\C,\QF) \xrightarrow{\ass} \F(\colim\limits_{\B\ZZ} (\C,\QF)) = \F((\C,\QF)_{\B\ZZ}).
\]
\end{definition}

\begin{remark}
\label{remark:circle-decompose-1}%
We always view $\B\ZZ$ as a pointed space.
The generator $1 \in \ZZ$ determines a pointed self homotopy $\tau$ of the composite $\ast_+ \to \ast \to \B\ZZ$, yielding a commutative square
\[
\begin{tikzcd}
*_+ \ar[r]\ar[d] & \ast \ar[d] \\
* \ar[r] & \B\ZZ 
\end{tikzcd}
\]
in pointed spaces, where we view $\tau$ as the commuting homotopy from the left-bottom composite to the top-right composite. This square is cocartesian, exhibiting $\B\ZZ$ as the suspension of $*_+$ in pointed spaces. 
For any stable $\infty$-category $\A$ and any object $x \in \A$, the above construction 
determines a decomposition
\[
\B\ZZ \otimes x = x \oplus \Sig x .
\]
In particular, if $(\C,\QF)$ is any Poincaré $\infty$-category then we can identify the assembly map with the map $\F(\C,\QF) \oplus \Sig\F(\C,\QF) \to \F((\C,\QF)_{\B\ZZ})$ induced by the resulting commutative square
\[
\begin{tikzcd}
\F(\C,\QF)\oplus \F(\C,\QF) \ar[r]\ar[d] & \F(\C,\QF) \ar[d] \\
\F(\C,\QF) \ar[r] & \F((\C,\QF)_{\B\ZZ}) 
\end{tikzcd}
\]
in $\A_{\F(\C,\QF)/}$, where we view the top left corner as living under $\F(\C,\QF)$ via the left summand.
\end{remark}

Let us now describe our circle decomposition theorem for Poincaré $\infty$-categories.
Given a Verdier-localising functor $\F\colon \Catp \to \A$ valued in a stable $\infty$-category $\A$ and a Poincaré $\infty$-category $(\C,\QF)$, 
consider the commutative diagram
\[
\begin{tikzcd}
\F(\Hyp(\C)) \ar[rr]\ar[dd]\ar[dr] && \F(\Hyp(\C_{\B\NN}))\ar[dd, dashed]\ar[dr] &\\
& \B\ZZ \otimes \F(\C,\QF) \ar[rr, "\ass", pos=0.4]\ar[dd] && \F((\C,\QF)_{\B\ZZ})\ar[dd]\ar[r] & \F(\C,\QF) \ar[dd] \\ 
\F(\Hyp(\C^{\natural})) \ar[rr, dashed]\ar[dr] && \F(\Hyp((\C^{\natural})_{\B\NN}))\ar[dr] &\\
& \B\ZZ \otimes \F(\C^{\natural},\QF) \ar[rr, "\ass",pos=0.4] && \F((\C^{\natural},\QF)_{\B\ZZ}) \ar[r] & \F(\C^{\natural},\QF)
\end{tikzcd}
\]
where the right diagonal arrows are instances of the map induced by the composite %
\[
\Hyp(\D_{\B\NN}) \to \Hyp(\D_{\B\ZZ}) \xrightarrow{\hyp} (\D,\Phi)_{\B\ZZ}
\]
and the left diagonal arrows are instances of the composite 
\[
\F(\Hyp(\D)) \xrightarrow{\hyp} \F(\D,\Phi) \xrightarrow{i} {\B\ZZ} \otimes \F(\D,\Phi).
\]
The main result of this section is then:

\begin{theorem}[The circle decomposition theorem]
\label{theorem:circle-decomposition}%
Let $\F\colon \Catp \to \A$ and $\G \colon \Catx \to \A$ be Verdier-localising functors valued in a stable $\infty$-category $\A$. %
Then the following hold:
\begin{enumerate}
\item
\label{item:main-face}%
The bottom face in the above diagram is cartesian. 
\item
\label{item:relative-circle}%
The functor $(\C,\QF) \mapsto \fib[\F((\C,\QF)_{\B\ZZ}) \to \F(\C,\QF)]$ on $\Catp$ is invariant under idempotent completion. In particular, 
the front right square in the above diagram is cartesian. 
\item
\label{item:nil}%
The functor $\D \mapsto \N\G(\D)$ on $\Catx$ is invariant under idempotent completion. In particular, the back face in the above diagram is cartesian. 
\end{enumerate}
In particular, for a Poincar\'e $\infty$-category $(\C,\QF)$ we have a natural decomposition
\[
\F((\C,\QF)_{\B\ZZ}) \simeq \F(\C,\QF) \oplus \Sig\F(\C^{\natural},\QF) \oplus \N\F^{\hyp}(\C).
\]
\end{theorem}

\begin{remark}
Unwinding definitions we find the map from right to left in the above natural decomposition is given by the sum of the following maps:
\begin{enumerate}
\item the tautological map $\F(\C,\QF) \rightarrow \F((\C,\QF)_{\B\ZZ})$,
\item the composite 
\begin{align*}
\Sig \F(\C^\natural,\QF) & \simeq \fib(\B\ZZ \otimes \F(\C^\natural,\QF) \rightarrow \F(\C^\natural,\QF)) \xrightarrow{\ass} \fib[\F((\C^\natural,\QF)_{\B\ZZ}) \rightarrow \F(\C^\natural,\QF)] \\
		& \simeq \fib[\F((\C,\QF)_{\B\ZZ}) \rightarrow \F(\C,\QF)] \longrightarrow \F((\C,\QF)_{\B\ZZ}),
\end{align*}
\item the map $\N\F^\hyp(\C) \rightarrow \F\Hyp(\C_{\B\NN}) \rightarrow \F((\C,\QF)_{\B\ZZ})$.
\end{enumerate}
Moreover, by construction the assembly map $\B\ZZ \otimes \F(\C,\QF) \to \F((\C,\QF)_{\B\ZZ})$ corresponds, under the splitting of its source as in 
Remark~\reffour{remark:circle-decompose-1},
to the composite of
\[
\F(\C,\QF) \oplus \Sig\F(\C,\QF) \to \F(\C,\QF) \oplus \Sig\F(\C^\natural,\QF) \to \F(\C,\QF) \oplus \Sig\F(\C^\natural,\QF) \oplus \N\F^{\hyp}(\C)
\]
where the first map is induced by the inclusion $(\C,\QF) \to (\C^\natural,\QF)$ and the second map is the summand inclusion. 
\end{remark}

\begin{corollary}[Ranicki-Shaneson decomposition]
\label{corollary:algebraic-shaneson-splitting}%
Let $(\C,\QF)$ be a Poincaré $\infty$-category and $\F\colon \Catp \to \A$ a bordism invariant Verdier-localising functor valued in a stable $\infty$-category $\A$. Then we have a natural decomposition
\[
\F((\C,\QF)_{\B\ZZ}) \simeq \F(\C,\QF) \oplus \Sig\F(\C^\natural,\QF)
\]
under which the assembly map from $\B\ZZ \otimes \F(\C,\QF) = \F(\C,\QF) \oplus \Sig\F(\C,\QF)$
is given by the identity on the first term and by the map induced by the inclusion $(\C,\QF) \to (\C^\natural,\QF)$ on the second term.
\end{corollary}

Arguing by induction we then conclude:
\begin{corollary}
\label{corollary:iterated-shaneson-splitting}%
Let $(\C,\QF)$ be a Poincaré $\infty$-category, $\F\colon \Catp \to \A$ a bordism invariant Verdier-localising functor valued in a stable $\infty$-category $\A$ and $n \geq 1$ an integer. Assume either that $\C_{\B\ZZ^{n-1}}$ is idempotent complete or that $\F$ is Karoubi-localising. Then the torus
assembly map
\[
\B\ZZ^n \otimes \F(\C,\QF) \to \F((\C,\QF)_{\B\ZZ^n})
\]
is an equivalence.
\end{corollary}

\begin{example}
\label{example:shaneson-splitting-Ghyp}%
If $\F\colon \Catp \to \A$ is a Verdier-localising functor we may form its genuine hyperbolisation $\F^\ghyp\colon \Catp \to \A^\gCt$, the genuine $\Ct$-$\A$-object valued functor whose underlying functor is $\F^\hyp = \F \circ \Hyp \colon \Catp \to \A^{\B\Ct}$, whose genuine fixed points is $\F$ and whose geometric fixed points is $\F^\bord$; see \S\reftwo{subsection:Mackey} for details. From loc.\ cit.\ we find $(\F^\ghyp)^\hyp \simeq \Ct \otimes \F^\hyp$ as functors $\Catx \to \A^{\gCt}$ where $C_2 \otimes - $ denotes the genuine induction functor, that is, the left adjoint of the underlying object functor $\A^{\gCt} \to \A$. Consequently,
Theorem~\reffour{theorem:circle-decomposition} gives an equivalence
\[
\F^{\ghyp}((\C,\QF)_{\B\ZZ}) \simeq \F^{\ghyp}(\C,\QF) \oplus \Sigma \F^{\ghyp}(\C^\natural,\QF) \oplus (\Ct \otimes \mathrm \N\F^\hyp(\C)).
\]
For $\F= \GW$, we have $\F^{\ghyp}= \KR$ and thus
\[
\KR((\C,\QF)_{\B\ZZ}) \simeq \KR(\C,\QF) \oplus \Sigma \KR(\C^{\natural},\QF) \oplus (\Ct \otimes \N\K(\C)),
\]
a Bass-Heller-Swan type decomposition for the real K-spectra of Poincaré categories. On underlying spectra with $\Ct$-action, this gives
\[
\K((\C,\QF)_{\B\ZZ}) \simeq \K(\C,\QF) \oplus \Sig \K(\C^\natural,\QF) \oplus [\Ct \otimes \N\K(\C)]
\]
with $\Ct$-action summand-wise. In particular, the $\Ct$-action swaps the two nil terms appearing in the classical Bass-Heller-Swan decomposition. On fixed points, we obtain
\[
\GW((\C,\QF)_{\B\ZZ}) = \GW(\C,\QF) \oplus \Sigma \GW(\C^{\natural},\QF) \oplus \N\K(\C)
\]
and on geometric fixed points, we obtain
\[
\L((\C,\QF)_{\B\ZZ}) = \L(\C,\QF) \oplus \Sigma \L(\C^{\natural},\QF).
\]
The same holds true for $\F=\KGW$, in which case the terms $\KK$ and $\KL$ appear instead of $\K$ and $\L$.
One may then insert the Poincar\'e $\infty$-categories appearing in Example~\reffour{example:tensor-BZ} and obtain formulas for the Grothendieck-Witt and L-theory of Laurent polynomial rings $R[t^{\pm 1}]$ equipped with involution sending $t$ to $t^{-1}$ or form parameter $\lambda[\ZZ]$ induced from a form parameter on $R$, as well as for the visible Poincar\'e structure on parametrised spectra over $X \times S^1$ with spherical fibration factoring through the projection $X \times S^1 \to X \to \Pic(\SS)$, see the circle decomposition part of the introduction.
\end{example}

\begin{remark}
\label{remark:all-involutions}%
Given a commutative ring $R$ with the identity involution, Example~\reffour{example:shaneson-splitting-Ghyp} gives in particular formulas for the Grothendieck--Witt- and L-spectra associated to the ring with $R[t^{\pm 1}]$ with involution given by $t\mapsto t^{-1}$. This is only one of the four canonical involutions $t \mapsto \pm t^{\pm 1}$. A number of things are known in the other cases as well:
\begin{enumerate}
\item
\label{item:twisted}%
The case $t\mapsto -t^{-1}$ was studied in~\cite{Milgram-Ranicki}, 
and fits in the more general setting of the twisted circle decomposition obtained in~\cite{twisted-shaneson}*{Theorem 1.5}, a special case of which yields the formulas
\[
\KGW^{\sym}(R[t^{\pm 1}];t\mapsto -t^{-1}) = \coeq[\KGW^{\sym}(R) \rightrightarrows \KGW^{\sym}(R)] \oplus \N\KK(R)
\]
and 
\[
\KL^{\sym}(R[t^{\pm 1}],t\mapsto -t^{-1}) = \coeq[\KL^{\sym}(R) \rightrightarrows \KL^{\sym}(R)] ,
\]
where one may, similarly to our circle decomposition theorem, replace the symmetric Poincaré structure
on $\Dperf(R)$ with the one determined by a form parameter on $R$, and the Poincaré structure on $\Dperf(R[t^{\pm 1}])$ with a corresponding induced form parameter on $R[t^{\pm1}]$.
The coequalisers in the above formulas are taken along the automorphisms $\id,\tau$ of $\KGW^{\sym}(R)$ and $\KL^{\sym}(R)$, respectively, where $\tau$ is the automorphism induced by the sign involution of $\QF^{\sym}_R$; in the case of $\KL^{\sym}(R)$ this is just the sign involution of $\KL^{\sym}(R)$, so that the coequaliser in question becomes the cofibre of multiplication by $2$ (compare~\cite{Milgram-Ranicki}). This coequaliser term can again be described as the image of a the assembly map, this time for a non-constant family of Poincaré $\infty$-categories indexed by the circle $\B\ZZ$.
\item
\label{item:motivic}%
The case where the involution is the identity can be deduced from the projective line formula of~\cite{motives}*{\S 6.1}. One obtains the formulas
\[
\KGW^{\sym}(R[t^{\pm 1}]) = \KGW^{\sym}(R) \oplus \Sig\KGW^{\sym,\scalebox{0.5}{[-1]}}(R) \oplus \N\KGW^{\sym}(R)\oplus \N\KGW^{\sym}(R)
\]
and 
\[
 \KL^{\sym}(R[t^{\pm 1}]) = \KL^{\sym}(R) \oplus \KL^{\sym}(R) \oplus \N\KL^{\sym}(R)\oplus \N\KL^{\sym}(R),
\]
where $\KGW^{\sym,\scalebox{0.5}{[-1]}}(R):=\KGW\big(\Dperf(R),(\QF^{\sym}_R)\qshift{-1}\big)$ and $\N\KGW^{\sym}(R)$ and $\N\KL^{\sym}(R)$ are the cofibres of the map $\KGW^{\sym}(R) \to \KGW^{\sym}(R[t])$ and $\KL^{\sym}(R) \to \KL^{\sym}(R[t])$, respectively. The situation becomes more transparent when passing to real $\K$-theory, in which case the above formulas can be combined into a single formula %
\[
\KKR^{\sym}(R[t^{\pm1}]) = [S^{\sig} \otimes \KKR^{\sym}(R)] \oplus \N\KKR^{\sym}(R) \oplus \N\KKR^{\sym}(R),
\]
where the assembly map term is replaced with tensoring with the sign representation sphere $S^{\sig}$, viewed as a genuine $\Ct$-space.
By \cite{motives}*{Theorem 6.3.1} the above nil groups vanish 
for $R$ regular noetherian of finite Krull dimension.
The analogous nil groups for quadratic $\L$-theory  %
have been shown to agree with Cappell's unitary nil terms, and are non-zero even for $R=\ZZ$, see \cite{Connolly-Ranicki}. 
\item One can expect to obtain a similar formula for the final involution $t\mapsto -t$, which can most likely be extracted from a suitable projective line formula in the setting of schemes with involution. This is the topic of work in progress~\cite{schemes-with-involution}.
\end{enumerate}
We remark that in the upcoming PhD thesis of Andreas Huber, a general formula for the value of Karoubi-localising invariants on certain pushouts of Poincaré $\infty$-categories is studied. 
One can then show that at least the involutions from~\reffouritem{item:twisted} and~\reffouritem{item:motivic}
give rise to Poincar\'e $\infty$-categories which can be described as pushouts of easier pieces referring only to $R$. In these cases, the results in Huber's thesis
apply to reproduce the above formulas in a unified manner. 
\end{remark}

\subsection{Proof of the circle decomposition theorem}
\label{subsection:proof}%
This subsection is dedicated to the proof of Theorem~\reffour{theorem:circle-decomposition}.
We shall use the following auxiliary construction.
Let \(\I\) be the category with objects \(\{-1,1\}\) and which has exactly two non-identity arrows \(\alp^+,\alp^-\), both going from \(-1\) to \(1\) (in particular, \(1\) and \(-1\) have no non-identity endomorphisms and there are no arrows from \(1\) to \(-1\)). 
The functor
\[
p\colon \I \to \B \ZZ
\]
sending $\alp^+$ to $1 \in \ZZ = \Hom_{\B\ZZ}(\ast,\ast)$ and $\alp^-$ to $0$ then exhibits $\B\ZZ$ as a localisation of $\I$.
\begin{remark}
\label{remark:unique}%
The functor $p$ above is completely determined by $p(1) = \ast$ and $p(\alp^+) \circ p(\alp^{-})^{-1} = 1 \in \ZZ$: More precisely the subgroupoid of $\Fun_*((\I,1),(\B\ZZ,\ast))$ spanned by those pointed functors $(\I,1) \to (\B\ZZ,\ast)$ such that $p(\alp^+) \circ p(\alp^-)^{-1} = 1$ is contractible. Indeed, this follows from the universal property of localisations and the fact that the group of pointed orientation preserving self-homotopy equivalences of the circle is contractible.

\end{remark}

Given a Poincaré $\infty$-category $(\C,\QF)$ we may tensor it with $\I$ in $\Cath$ to obtain a hermitian $\infty$-category $(\C,\QF)_\I$. Since the category \(\I\) is strongly finite,
Lemma~\refone{lemma:finite} and Proposition~\refone{proposition:tensor-strongly-finite} give that we may identify \(\C_\I\) with the stable \(\infty\)-category of diagrams \(\I\op \to \C\), such that the hermitian functor \(\QF_\I\) is given by
\[
\QF_{\I}(\vphi) = \colim_{i \in \I}\QF(\vphi(i)) 
\]
and its bilinear part of $\QF_{\I}$ is represented by the %
duality
\[
\Dual_{\I}(\vphi)(j) = \colim_{i \in \I}\Dual(\vphi(i))^{\Hom_{\I}(i,j)}.
\]
Calculating this explicitly gives that $\QF_\I$ takes the coequaliser of the two maps obtained by applying $\QF$, and 
\begin{equation}
\label{equation:explicit-dual-I}%
\Dual_{\I}
\big[
\begin{tikzcd}
\x  & \ar[l,shift right]\ar[l,shift left] \y 
\end{tikzcd}
\Big] = 
\Big[
\begin{tikzcd}
\Sig\Dual\x  & \ar[l,shift right]\ar[l,shift left] \Sig[\Dual\x \times_{\Dual \y}\Dual\x] 
\end{tikzcd}
\big],
\end{equation}
where the pullback $\Dual \x\times_{\Dual \y} \Dual \x$ is to be taken along the two evident maps, and the two arrows in the second diagram are given by the two projections from the fibre product. Direct inspection shows that this duality is perfect, so that $(\C,\QF)_\I$ is Poincaré.
We will also deduce this fact from abstract principles in Proposition~\reffour{proposition:I-tensor-poincare} below.

Let $\NN^{\pm} \subseteq \ZZ$ be the submonoid spanned by $\pm 1$, respectively. Let $p^{\pm}\colon \I \to \B\NN^{\pm}$ be the functor sending $\alp^{\pm}$ to $\pm 1 \in \NN^{\pm}$ and $\alp^{\mp}$ to $0$. Then the composites 
\[
\I \xrightarrow{p^{\pm}} \B\NN^{\pm} \to \B\ZZ
\]
both send the formal difference $\alp^+ \circ (\alp^-)^{-1}$ to $1 \in \ZZ$, and are hence both equivalent to $p\colon \I \to \B\ZZ$ by its uniqueness, see Remark~\reffour{remark:unique}. %
We thus obtain a commutative square of $\infty$-categories
\[
\begin{tikzcd}
\I \ar[r, "p^+"]\ar[d, "p^-"'] & \B \NN^+ \ar[d] \\
\B\NN^- \ar[r] & \B\ZZ
\end{tikzcd}
\]
in which all arrows are localisation functors: the horizontal functors exhibit their targets as localisations by $\alp^-$ and $-1$, respectively, while the vertical functors exhibit their targets as localisations by $\alp^+$ and $1$, respectively. We note that unlike $(\C,\QF)_\I$, the tensored hermitian $\infty$-categories $(\C,\QF)_{\B\NN^{\pm}}$ are generally not Poincaré. However, by Proposition~\reftwo{proposition:localizationforp4}\;\reftwoitem{item:induced-functor-verdier}, the functor $(\C,\QF)_\I \to (\C,\QF)_{\B\ZZ}$ is a Poincar\'e-Verdier projection, so we may use this projection, together with an analysis of its kernel, to infer information from $\F((\C,\QF)_\I)$ to $\F((\C,\QF)_{\B\ZZ})$.

The study of $(\C,\QF)_\I$ itself will in turn rest on a portion of %
the dual $\Q$-construction which we now recall. We consider the hermitian tensor $\dualQ_1(\C,\QF) = (\C,\QF)_{\Twar(\Del^1)}$ of $(\C,\QF)$ by the twisted arrow category of $\Del^1$. Since $\Twar(\Del^1)$ is the category classifying pushouts, for ease of notation we write $\mathcal{P} = \Twar(\Del^1)$ and denote its initial object by $-1$ and the two other objects by $1$ and $1'$.
Then $\dualQ_1(\C,\QF) = (\C,\QF)_{\mathcal{P}}$ is Poincaré by Lemma~\reftwo{lemma:dqexists} and similarly as with the tensor with $\I$, its objects are given by $\cP\op$-indexed diagrams $[x \rightarrow z \leftarrow y]$ with the Poincaré structure given by
\[
\QF_{\mathcal{P}}([x \rightarrow z \leftarrow y]) = \QF(x) \oplus_{\QF(z)} \QF(y)
\]
and duality given by 
\[
\Dual_{\mathcal{P}}([x \xrightarrow{f} z \xleftarrow{g} y]) = [\cof(\Dual(f)) \rightarrow \Sig \Dual(z) \leftarrow \cof(\Dual(g))].
\]
In fact, $\dualQ_1(\C,\QF)$ is part of a cosimplicial object $[n] \mapsto \dualQ_n(\C,\QF)$ in $\Catp$ with $\dualQ_0(\C,\QF) = (\C,\QF)$, see Lemma~\reftwo{lemma:dqexists} and what follows. The two face maps $[0] \to [1]$ (with images $\{0\}$ and $\{1\}$ respectively) and the degeneracy map $s\colon [1] \to [0]$ then induce Poincaré functors
\[
l,r\colon (\C,\QF) \to (\C,\QF)_{\mathcal{P}} \quad\quad s\colon (\C,\QF)_{\mathcal{P}} \to (\C,\QF)
\]
which are explicitly given by $l(x) = [x \rightarrow 0 \leftarrow 0]$ and $r(y) = [0 \rightarrow 0 \leftarrow y]$ and $s[x \rightarrow z \leftarrow y] = x \times_z y$. We note that $s$ admits a fully-faithful left adjoint given by $x \mapsto [x = x = x]$, and is hence in particular a split Verdier projection, in fact a split Poincar\'e-Verdier projection as one checks readily. On the other hand, the face functors $l,r$ are both split Poincaré-Verdier inclusions, with left adjoints given by evaluation at $1,1' \in \mathcal{P}$ respectively. %

\begin{construction}
\label{construction:twar}%
Consider the commutative diagram of categories
\begin{equation}
\label{equation:pre-square-I-Q1}%
\begin{tikzcd}
\{\ast, \ast\} \ar[r,"{1,1'}"] \ar[d]  & \mathcal{P} \ar[d] & \ar[l,"-1"'] \ast \ar[d,equal] \\
\ast \ar[r,"1"] & \I & \ar[l,"-1"'] \ast 
\end{tikzcd}
\end{equation}
where the left square is the pushout square exhibiting $\I$ as obtained from $\mathcal{P}$ by collapsing the objects $1,1'$ in $\mathcal{P}$ to the object $1 \in \I$, and the right square is a pullback square. By direct inspection the map $\mathcal{P} \to \I$ is a cocartesian fibration and hence the opposite of the right square satisfies the Beck-Chevalley property for right Kan extensions along the vertical maps. For a given stable $\infty$-category $\D$ we hence obtain a commutative diagram
\begin{equation}
\label{equation:stable-square-I-Q1}%
\begin{tikzcd}
\D \oplus \D \ar[r,"l \oplus r"]\ar[d] & \D_{\mathcal{P}} \ar[d] \ar[r,"\ev^{\cP}_{-1}"] & \D \ar[d,equal] \\
\D \ar[r,"i_1"] & \D_\I \ar[r,"\ev^{\I}_{-1}"] & \D
\end{tikzcd}
\end{equation}
where the left square is induced on tensors by the left square in the previous diagram and the right square is obtained by doing the same to the right square in the previous diagram and then replacing the horizontal arrows by their left adjoints. 
\end{construction}

When identifying tensors with diagram categories as above, all functors in the left square of~\eqreffour{equation:stable-square-I-Q1} are induced by right Kan extension and the horizontal maps in the right square are induced by evaluation 
at $-1$ in $\mathcal{P}$ and $\I$, respectively. Concretely, for the left square, $(l\oplus r)(x,y) = [x \rightarrow 0 \leftarrow y]$ and $i_1(x) = [0 \leftleftarrows x]$. 
Since tensoring with $\D$ is a colimit preserving operation the left square is cocartesian in $\Cat$.

\begin{proposition}
\label{proposition:I-tensor-poincare}%
The rows in diagram~\eqreffour{equation:stable-square-I-Q1} are split Verdier sequences. In addition, if $(\C,\QF)$ is a Poincar\'e $\infty$-category, the hermitian $\infty$-category $(\C,\QF)_\I$ is Poincar\'e and diagram~~\eqreffour{equation:stable-square-I-Q1} for $\D=\C$ refines to a commutative diagram in $\Catp$ whose rows are split Poincar\'e-Verdier sequences
\begin{equation}
\label{equation:square-I-Q1}%
\begin{tikzcd}
	(\C,\QF) \oplus (\C,\QF) \ar[r,"l\oplus r"]\ar[d] & (\C,\QF)_{\mathcal{P}} \ar[d] \ar[r,"\ev^{\cP}_{-1}"] & (\C,\QF\qshift{1}) \ar[d,equal]\\
	(\C,\QF) \ar[r,"i_1"] & (\C,\QF)_\I \ar[r,"\ev^{\I}_{-1}"] & (\C,\QF\qshift{1})
\end{tikzcd}
\end{equation}
where on underlying hermitian $\infty$-categories the left square is obtained from the left square of~\eqreffour{equation:pre-square-I-Q1} by tensoring with $(\C,\QF)$. In addition, the left square is a pushout in $\Catp$.
In particular, any additive functor $\F\colon \Catp \to \A$ to a stable $\infty$-category $\A$ sends the left square %
to an exact square. 
\end{proposition}
\begin{proof}
By direct inspection the rows of~\eqreffour{equation:stable-square-I-Q1} are fibre sequences. In addition, the evaluation functors $\ev^{\cP}_{-1}$ and $\ev^{\I}_{-1}$ admit fully-faithful left and right adjoints given by left and right Kan extension along the inclusion of $\{-1\}$. We hence conclude that the rows of~\eqreffour{equation:stable-square-I-Q1} are split Verdier sequences.

Given a Poincaré $\infty$-category $(\C,\QF)$, the left square of~\eqreffour{equation:stable-square-I-Q1} for $\D=\C$ refines to the left square of~\eqreffour{equation:square-I-Q1} on the level of hermitian $\infty$-categories by performing the tensor there, yielding again a cocartesian square since tensoring with $(\C,\QF)$ is colimit preserving. Since the top horizontal and left vertical arrows in the left square of~\eqreffour{equation:square-I-Q1} are Poincaré functors between Poincaré $\infty$-categories and $\Catp$ is closed in $\Cath$ under colimits by Proposition~\refone{proposition:Catp-cocomplete}, the left square in~\eqreffour{equation:square-I-Q1} is contained in $\Catp$ and is a pushout square there. In particular, $(\C,\QF)_\I$ is Poincaré.  
Now the horizontal arrows in the left square are Poincaré-Verdier inclusions by Proposition~\reftwo{proposition:example-cosieve}.
Combining this with the fact that the left square is cocartesian we conclude that the induced map $(\ev^{\cP}_{-1})_!\QF_{\cP} \Rightarrow (\ev^{\I}_{-1})_!\QF_\I$ is an equivalence. 
Concretely, using the left adjoints $x \mapsto [0 \rightarrow x \leftarrow 0]$ and $x \mapsto [0 \leftleftarrows x]$ of $\ev^{\cP}_{-1}$ and $\ev^{\I}_{-1}$, respectively, this map can be identified with the map
\[
\colim\big[\begin{tikzcd}
0 & \QF(x) \ar[r]\ar[l] & 0 
\end{tikzcd}\big] \xrightarrow{\simeq}
\colim\big[\begin{tikzcd}
\QF(x) \ar[r,shift right]\ar[r,shift left] & 0 
\end{tikzcd}\big] 
\]
induced on colimits by $\cP \to \I$. We now finish the proof by identifying the left hand side with $\QF\qshift{1}$.
\end{proof}

\begin{remark}
\label{rem:variant}%
\label{remark:dmet}%
The split Poincaré-Verdier projection $\ev^{\cP}_{-1}\colon (\C,\QF)_{\mathcal{P}} \to (\C,\QF\qshift{1})$ of Proposition~\reffour{proposition:I-tensor-poincare} factors through the Poincaré-Verdier quotient $\dMet_1(\C,\QF)=(\C,\QF)_{\cP}/l(\C,\QF)$,
see Definition~\reftwo{definition:dualmet}. Concretely, we may identify objects of $\dMet_1(\C,\QF)$ with arrows $z \leftarrow y$, such that its Poincaré structure is given by $[z \leftarrow y]\mapsto \cof[\QF(z) \to \QF(y)]$. Under this identification the projection $(\C,\QF)_{\cP} \to \dMet_1(\C,\QF)$ is given by
\[
[x \rightarrow z \leftarrow y] \mapsto [z \leftarrow y] \quad\quad \colim\big[\QF(x) \leftarrow \QF(z) \rightarrow \QF(y)\big]  \to \cof[\QF(z) \to \QF(y)] .
\]
The resulting Poincaré functor $\dMet_1(\C,\QF) \to (\C,\QF\qshift{1})$ is then the canonical projection 
\[
[z \leftarrow y] \mapsto z \quad\quad \cof[\QF(z) \to \QF(y)] \mapsto \QF\qshift{1}(z)
\]
from the dual metabolic category. Alternatively, identifying $\cof[\QF(x) \to \QF(z)]$ with $\fib[\QF\qshift{1}(x) \to \QF\qshift{1}(z)]$ we may also identify $\dMet_1(\C,\QF)$ with $\Met(\C,\QF\qshift{1})$ and the above projection with $\met\colon \Met(\C,\QF\qshift{1}) \to (\C,\QF\qshift{1})$. Summarising the situation, Diagram~\eqreffour{equation:square-I-Q1} can be extended to a diagram
\[
\begin{tikzcd}
	(\C,\QF) \ar[r,"l"] \ar[d,"i_l"] & (\C,\QF)_{\mathcal{P}} \ar[r,"\ev"] \ar[d,equal] & \Met(\C,\QF\qshift{1}) \ar[d,"\met"] \\
	(\C,\QF) \oplus (\C,\QF) \ar[r,"l\oplus r"]\ar[d] & (\C,\QF)_{\mathcal{P}} \ar[d] \ar[r,"\ev^{\cP}_{-1}"] & (\C,\QF\qshift{1}) \ar[d,equal]\\
	(\C,\QF) \ar[r,"i_1"] & (\C,\QF)_\I \ar[r,"\ev^{\I}_{-1}"] & (\C,\QF\qshift{1})
\end{tikzcd}
\]
in which also the top horizontal sequence is a split Poincar\'e-Verdier sequence, where $i_l$ is the inclusion of the left summand. 
\end{remark}

Having analysed $\D_\I$ and $(\C,\QF)_\I$, we now address the kernel of the induced functor $\D_\I \to \D_{\B\ZZ}$, which is a Verdier projection by Proposition~\reftwo{proposition:localizationforp4}.

\begin{construction}
\label{construction:E}%
Let $\D$ be a stable $\infty$-category.
For \(i \in \I\) and \(\x \in \D\) we let \(R_{i,\x}\colon \I\op \to \D\) denote the right Kan extension of \(\x\) along the inclusion \(\{i\} \subseteq \I\op\). Explicitly, the diagrams \(R_{1,\x}\) and \(R_{-1,\x}\) can be depicted as
\[
\begin{tikzcd}[row sep=2pt]
R_{1,x} = \big[0 & \ar[l,shift left]\ar[l,shift right] \x \big] && \text{and} && R_{-1,x} = \big[\x & \ar[l,shift left]\ar[l,shift right] \x \oplus \x \big], \
\end{tikzcd}
\]
respectively. 
The maps \(\alp^{\pm}\colon -1 \to 1\) in $\I$ now induce maps \(\alp^{\pm}_{\ast}\colon R_{-1,\x} \to R_{1,\x}\), whose fibres \(z^+_{\x} := \fib[\alp^+_{\x}]\) and \(z^-_x := \fib[\alp^-_{\x}]\) are depicted by the diagrams
\[
\begin{tikzcd}[row sep = 2ex]
z^+_x = \big[\x & \ar[l,shift left,"0"]\ar[l,shift right,"\id"'] \x\big] && \text{and} && z^-_x = \big[ \x & \ar[l,shift left,"\id",pos=0.45]\ar[l,shift right,"0"',pos=0.45] \x\big] , \
\end{tikzcd}
\]
respectively.
We then write \(\cN^{\pm}(\D) \subseteq \D_\I\) for the smallest stable and retract-closed subcategory of $\D_\I$ containing the objects \(z^{\pm}_{\x}\) for every \(\x \in \D\), and \(\cN(\D) \subseteq \D_{\I}\) for the smallest stable retract-closed subcategory of \(\D_{\I}\) containing $\cN^+(\D)$ and $\cN^-(\D)$. The categories $\cN^{\pm}(\D)$ are equivalent to the categories of nilpotent endomorphisms in $\D$ which we have briefly mentioned in Construction~\reffour{construction:nil}. If $(\C,\QF)$ is a Poincar\'e $\infty$-category,
the explicit formula above for \(\Dual_{\I}\) shows that \(\Dual_{\I}(z^+_{\x}) \simeq \Sigma z^-_{\Dual\x}\) and \(\Dual_{\I}(z^-_{\x}) \simeq \Sigma z^+_{\Dual \x}\), and so:
\begin{enumerate}
\item
\(\cN(\C)\) is %
invariant under the duality of $(\C,\QF)_\I$.
\item
The induced duality of $\cN(\C)$ switches between $\cN^+(\C)$ and $\cN^-(\C)$, and in particular determines an equivalence $\cN^-(\C) \simeq (\cN^+(\C))\op$.
\end{enumerate}
We will then denote by $\cN(\C,\QF) = (\cN(\C),(\QF_\I)_{|\cN(\C)})$ the Poincaré $\infty$-category obtained by restricting the Poincaré structure $\QF_\I$ to $\cN(\C)$.
\end{construction}

\begin{remarks}\
\label{remark:cofree}%
\begin{enumerate}
	\item The functor $\D_{\mathcal{P}} \to \D_\I$ in \eqreffour{equation:stable-square-I-Q1} sends $[x = x \leftarrow 0]$ to $z^+_x$ and $[0 \rightarrow x = x]$ to $z^-_x$. Indeed, this can be seen by direct inspection, or by noting that in both cases the object in question is obtained as a fibre of a map between two right Kan extended diagrams, and using the fact that the above induced functor is realised via right Kan extension.
	\item Since $\alp^{\mp}$ maps to an equivalence by $p^{\pm}\colon \I \to \B\NN^{\pm}$ %
and the induced map on tensors $p^{\pm}_*\colon \D_\I \to \D_{\B\NN^{\pm}}$ is realised via right Kan extension, we find for all $x \in \D$ that 
\[
p^+_*(z_x^-) = 0 \in \D_{\B\NN^+} \quad\text{and}\quad p^-_*(z^+_x) = 0 \in \D_{\B\NN^-} \, .
\]
\end{enumerate}
\end{remarks}

\begin{lemma}
\label{lemma:nil-term-hyperbolic}%
For a stable $\infty$-category $\D$, the full subcategories $\cN^+(\D),\cN^-(\D) \subseteq \cN(\D)$ induce a fully faithful and dense inclusion $\cN^-(\D) \oplus \cN^+(\D) \to \cN(\D)$. Moreover, for a Poincar\'e $\infty$-category $(\C,\QF)$, under the above identification $\cN^+(\C) = (\cN^-(\C))\op$, the full subcategory $\cN^-(\C) \oplus \cN^+(\C)$ of $\cN(\C,\QF)$ with the restricted Poincaré structure is naturally equivalent to $\Hyp(\cN^-(\C))$. 
In particular, if $\C$ is idempotent complete then the induced map 
\[
\Hyp(\cN^-(\C)) \to \cN(\C,\QF)
\]
is an equivalence.
\end{lemma}
\begin{proof}
We first show that for all $x^\pm \in \cN^\pm(\C)$, we have $\hom(x^+,x^-) = 0 = \map(x^-,x^+)$ from which it follows that the functor $\cN^-(\D) \oplus \cN^+(\D) \to \cN(\D)$ is fully faithful. Since its image is dense by definition of $\cN(\D)$, the first part of the lemma follows. To see the claim about the mapping spectra, it suffices to show $\map(z_x^-,z_y^+) = 0 = \map(z_x^+,z_y^-)$, since the mapping functor preserves finite colimits in each variable and retracts of $0$ are $0$. From the pushout formula $\I = \{1\} \amalg_{\{1,1'\}} \mathcal{P}$, recall diagram~\eqreffour{equation:pre-square-I-Q1}, we obtain the following pullback square
\[
\begin{tikzcd}
	\map_{\D_\I}(z_x^+,z_y^-) \ar[r] \ar[d] & \map_{\D_{\mathcal{P}}}([x=x\xleftarrow{0} x],[y\xrightarrow{0} y = y]) \ar[d]\\
	\map_\D(x,y) \ar[r,"\Delta"] & \map_\D(x,y) \times \map_\D(x,y)
\end{tikzcd}
\]
Unravelling the definitions, we find that the right vertical map is an equivalence in the first component of the target and the zero map in the second component of the target. As a consequence, the pullback vanishes as claimed. The same argument also shows $\map(z_x^-,z_y^+) = 0$.

For $(\C,\QF)$ a Poincar\'e $\infty$-category, using that the duality of $\cN(\C,\QF)$ interchanges $\cN^-(\C)$ and $\cN^+(\C)$, we deduce from the previous result on mapping spectra that $\B_\I$ vanishes upon restriction to $\cN^{\pm}(\C)$. As a result, $\QF_\I$ restricts to a linear functor on $\cN^\pm(\C)$, so we deduce from the computation
\[
\QF_\I(z^{+}_{\x}) = \colim \Big[
\begin{tikzcd}
\QF(\x) \ar[r,shift right,"0"']\ar[r,shift left,"\id"] & \QF(\x) 
\end{tikzcd}
\Big] = 0,
\]
and similarly for \(z^-_{\x}\), that $\QF_\I$ vanishes upon restriction to $\cN^{\pm}(\C)$. 
Hence, the restriction of $\QF_\I$ on $\cN^-(\C) \oplus \cN^+(\C)$ vanishes on each of the components and has a duality which swaps them, so we obtain that the induced Poincaré functor 
\[
\Hyp(\cN^-(\C)) \to (\cN^-(\C) \oplus \cN^+(\C),(\QF_\I)_{|\cN^+(\C)\oplus \cN^-(\C)})
\]
is an equivalence as desired.
\end{proof}

\begin{proposition}
\label{proposition:PV-sequences-Q1-hyp}%
For a stable $\infty$-category $\D$ the induced functor $p^+_*\colon \D_\I \to \D_{\B\NN^+}$ fits into a commutative diagram
\begin{equation}
\label{equation:tw-BN-C}%
\begin{tikzcd}
\D \ar[rr,"{x \mapsto [x = x \leftarrow 0]}"]\ar[d,"x \mapsto z_x^-"'] && \D_{\mathcal{P}} \ar[r]\ar[d] & \D_{\Del^1} \ar[d] \\
\cN^-(\D) \ar[rr] && \D_\I \ar[r, "p^+_*"] & \D_{\B\NN^+}
\end{tikzcd}
\end{equation}
in which the rows are Verdier sequences (with the top a split one) and the right square is induced on tensors by the square
\begin{equation}
\label{equation:tw-BN}%
\begin{tikzcd}
\mathcal{P} \ar[r,"c"]\ar[d] & \Del^1\ar[d] \\
\I \ar[r] & \B\NN^+ 
\end{tikzcd}
\end{equation}
where $c(-1) = c(1') =-1$ and $c(1) = 1$ (which admits a unique splitting).
\end{proposition}
\begin{proof}
Since the horizontal arrows in~\eqreffour{equation:tw-BN} are localisation functors, Proposition~\reftwo{proposition:localizationforp4} gives that the horizontal arrows in the right square of~\eqreffour{equation:tw-BN-C} are Verdier projections whose kernels are the smallest stable subcategories closed under retracts containing $[x = x \leftarrow 0] \in \D_{\mathcal{P}}$ and $z^-_x \in \D_\I$ for every $x \in \D$, respectively. In the former case the collection of cospans in $\D$ of the form $[x = x \leftarrow 0]$ already form a reflective and coreflective stable subcategory of $\D_{\mathcal{P}}$ which is equivalent to $\D$ itself via the displayed exact functor $x \mapsto [x = x \leftarrow 0]$. In the latter case the relevant smallest full stable subcategory closed under retracts is $\cN^-(\D)$. The left most vertical leg in the commutative diagram is necessarily the induced map on fibres, which by Remark~\reffour{remark:cofree} is given by the exact functor $\D \to \cN^-(\D)$ sending $x$ to $z_x^-$.
\end{proof}

\begin{proposition}
\label{proposition:PV-sequences-Q1}%
For a Poincaré $\infty$-category $(\C,\QF)$ the induced hermitian functor $p_*\colon (\C,\QF)_\I \to (\C,\QF)_{\B\ZZ}$ is Poincaré and furthermore fits in a commutative diagram
\begin{equation}
\label{equation:tw-S1}%
\begin{tikzcd}
\Hyp(\C) \ar[r]\ar[d] & (\C,\QF)_{\mathcal{P}} \ar[r,"s"]\ar[d] & (\C,\QF) \ar[d] \\
\cN(\C,\QF) \ar[r] & (\C,\QF)_\I \ar[r,"p_*"] & (\C,\QF)_{\B\ZZ} 
\end{tikzcd}
\end{equation}
whose rows are Poincaré-Verdier sequences (with the top a split one), whose right square is induced on tensors by the square
\[
\begin{tikzcd}
\mathcal{P} \ar[r]\ar[d] & \ast \ar[d] \\
\I \ar[r] & \B\ZZ \ ,
\end{tikzcd}
\]
and such that the two Poincaré functors out of $\Hyp(\C)$ are induced by the exact functors $x \mapsto [x \rightarrow x \leftarrow 0]$ and $x \mapsto z^-_x$, respectively.
\end{proposition}
\begin{proof}
Since $\cN(\C) \subseteq \C_\I$ is closed under the duality of $(\C,\QF)_\I$ we have by Proposition~\reftwo{proposition:localizationforp4}\;\reftwoitem{item:induced-functor-verdier} that $(\C,\QF)_\I \to (\C,\QF)_{\B\ZZ}$ is a Poincaré-Verdier projection whose kernel coincides with $\cN(\C,\QF)$. This gives the bottom row of~\eqreffour{equation:tw-S1}.
For the top row, note that the hermitian functor induced on tensors by the terminal map $\mathcal{P} \to \ast$ is the degeneracy functor 
$s$, which we observed earlier is a split Poincaré-Verdier projection. 
Its kernel consists of cospans of the form
$[x \rightarrow x \oplus y \leftarrow y] = [x = x \leftarrow 0] \oplus [0 \rightarrow y = y]$.
Since the Poincaré structure of $(\C,\QF)_{\mathcal{P}}$ vanishes on cospans of the form $[x = x \leftarrow 0]$ and on cospans of the form $[0 \rightarrow y = y]$, and at the same time the duality $\Dual_{\mathcal{P}}$ interchanges these two types of cospans, we conclude that the kernel of the Verdier projection $s\colon \C_{\mathcal{P}} \to \C$ is closed under the duality and identifies with $\Hyp(\C)$ when equipped with the restricted Poincaré structure. Example~\reftwo{example:kan-poincare} then yields the top Poincaré-Verdier sequence in \eqreffour{equation:tw-S1} and left most vertical leg in \eqreffour{equation:tw-S1} is necessarily the induced map on fibres, and by Remark~\reffour{remark:cofree} 
it is induced by the exact functor $\C \to \cN(\C)$ which sends $x$ to $z_x^-$.
\end{proof}

\begin{remark}
\label{remark:equ}%
Composing the Poincaré-Verdier projection $(\C,\QF)_\I \to (\C,\QF)_{\B\ZZ}$ of Proposition~\reffour{proposition:PV-sequences-Q1} with the induced Poincaré functor $(\C,\QF)_{\B\ZZ} \to (\C,\QF)$ yields a Poincaré functor
\[
\equ\colon (\C,\QF)_\I \to (\C,\QF)
\]
whose underlying hermitian functor is induced on tensors by the terminal map $\I \to \ast$. The Poincaré functor $\equ$ then provides a retraction to the split Poincaré-Verdier inclusion $i_1\colon (\C,\QF) \to (\C,\QF)_\I$ of Proposition~\reffour{proposition:I-tensor-poincare}, whose cofibre is identified with $(\C,\QF\qshift{1})$ as part of that proposition. In particular, 
for any additive functor $\F\colon \Catp \to \A$ we obtain a decomposition 
\[
\F((\C,\QF)_\I) = \F(\C,\QF) \oplus \F(\C,\QF\qshift{1}) .
\]
\end{remark}

\begin{remark}
\label{remark:I-is-assembly}%
Given a Poincaré $\infty$-category $(\C,\QF)$, the left square of Proposition~\reffour{proposition:I-tensor-poincare} and the right square of Proposition~\reffour{proposition:PV-sequences-Q1} fit into the following commutative diagram.
\[
\begin{tikzcd}
(\C,\QF) \oplus (\C,\QF) \ar[d]\ar[r] & (\C,\QF)_{\mathcal{P}} \ar[d]\ar[r] & (\C,\QF) \ar[d] \\
(\C,\QF) \ar[r] & (\C,\QF)_\I \ar[r] & (\C,\QF)_{\B\ZZ}
\end{tikzcd}
\]
Applying a Verdier-localising functor $\F\colon \Catp \to \cA$ with a stable target to this diagram sends the left square to a cartesian square by Proposition~\reffour{proposition:I-tensor-poincare}, and hence induces an equivalence
\[
\begin{tikzcd}
\F(\C,\QF) \displaystyle\mathop{\oplus}_{\F(\C,\QF) \oplus \F(\C,\QF)} \F(\C,\QF)\ar[dr,"\ass"]\ar[rr,"\simeq"] && \F(\C,\QF) \displaystyle\mathop{\oplus}_{\F((\C,\QF)_{\mathcal{P}})} \F((\C,\QF)_\I)  \ar[dl] \\
& \F((\C,\QF)_{\B\ZZ}) &
\end{tikzcd}
\]
from the assembly map on the left to the diagonal map on the right (see Remark~\reffour{remark:circle-decompose-1}). If furthermore $\F$ is bordism invariant then the map $\F((\C,\QF)_{\mathcal{P}}) \to \F(\C,\QF)$ is an equivalence by Proposition~\reffour{proposition:PV-sequences-Q1} so that the assembly map reduces to the map induced by the Poincaré functor $(\C,\QF)_\I \to (\C,\QF)_{\B\ZZ}$.

We then note that the inclusion $i_l\colon (\C,\QF) \to (\C,\QF) \oplus (\C,\QF)$ of the left summand induces a split Poincar\'e-Verdier inclusion of $(\C,\QF)$ into all four terms of the left square in the first displayed diagram above. Taking the quotients of each of these split Poincar\'e-Verdier inclusions, using Remark~\reffour{rem:variant}, we obtain the metabolic Poincar\'e-Verdier sequence
\[
\begin{tikzcd}
	(\C,\QF) \ar[r] \ar[d] & \Met(\C,\QF\qshift{1}) \ar[d] \\
	0 \ar[r] & (\C,\QF\qshift{1})
\end{tikzcd}
\]
It follows that both terms of the horizontal equivalence in the above commutative triangle receive a map from $\F(\C,\QF)$ which is split (see Remark~\reffour{remark:equ}) and whose induced map between cofibres is the equivalence
\[
\begin{tikzcd} \Sigma \F(\C,\QF) \ar[r,"\simeq"] & \cof[\Met(\C,\QF\qshift{1}) \to (\C,\QF\qshift{1})] \end{tikzcd}
\]
induced by the metabolic fibre sequence.
\end{remark}

We are now ready to establish the first part of the circle decomposition theorem (Theorem~\reffour{theorem:circle-decomposition}\;\reffouritem{item:main-face}):
	
\begin{theorem}
\label{theorem:circle-decomposition-idem}%
Let $\F\colon \Catp \to \A$ a Verdier-localising functor valued in a stable $\infty$-category $\A$ and $(\C,\QF)$ an idempotent complete Poincaré $\infty$-category. Then the square
\[
\begin{tikzcd}
\F(\Hyp(\C)) \ar[r]\ar[d] & \F(\Hyp(\C_{\B\NN^+}))\ar[d] \\
{\B\ZZ} \otimes \F(\C,\QF) \ar[r, "\ass"] & \F((\C,\QF)_{\B\ZZ}),
\end{tikzcd}
\]
is cartesian.
\end{theorem}
\begin{proof}
Consider the following vertical transformation between two diagrams of Poincaré $\infty$-categories:
\[
\begin{tikzcd}
0 \ar[rr]\ar[dd]\ar[dr] && \Hyp(\C\oplus\C)\ar[dd,dashed] \ar[rr,equal]\ar[dr] && \Hyp(\C \oplus \C)\ar[dd,dashed] \ar[dr] & \\
& \Hyp(\C)\ar[rr]\ar[dd] && \Hyp(\C_{\mathcal{P}}) \ar[rr]\ar[dd] && \Hyp(\C_{\Del^1})\ar[dd] \\
0 \ar[rr, dashed]\ar[dr] && \Hyp(\C)  \ar[rr,equal,dashed]\ar[dr] && \Hyp(\C)\ar[dr]  & \\
& \Hyp(\cN^-(\C))\ar[rr] && \Hyp(\C_\I) \ar[rr] && \Hyp(\C_{\B\NN^+}) \\
&& \color{white}{\;} \ar[r,phantom,"\rotatebox{270}{$\Rightarrow$}"] & \color{white}{\;}  && \\
0 \ar[rr]\ar[dd]\ar[dr] && (\C,\QF) \oplus (\C,\QF)\ar[dd,dashed] \ar[rr,equal]\ar[dr] && (\C,\QF) \oplus (\C,\QF)\ar[dd,dashed] \ar[dr] & \\
& \Hyp(\C)\ar[rr]\ar[dd] && (\C,\QF)_{\mathcal{P}} \ar[rr]\ar[dd] && (\C,\QF)\ar[dd] \\
0 \ar[rr, dashed]\ar[dr] && (\C,\QF)  \ar[rr,equal,dashed]\ar[dr] && (\C,\QF)\ar[dr]  & \\
& \cN(\C,\QF)\ar[rr] && (\C,\QF)_\I \ar[rr] && (\C,\QF)_{\B\ZZ}.
\end{tikzcd}
\]
The front horizontal rows are the Poincaré-Verdier sequences of Propositions \reffour{proposition:PV-sequences-Q1-hyp} (after applying $\Hyp$) and \reffour{proposition:PV-sequences-Q1}, and the middle squares in the top and bottom diagrams are those of Construction~\reffour{construction:twar} (after applying $\Hyp$) and Proposition~\reffour{proposition:I-tensor-poincare}, respectively.
The functor $\F$ hence sends each of the rows to a fibre sequence. In addition, the transformation between the top and bottom diagrams %
is an equivalence on left most faces by Lemma~\reffour{lemma:nil-term-hyperbolic} (which is where we use the assumption that $\C$ is idempotent complete), while the middle faces are both sent to cartesian squares via $\F$ by Proposition~\reffour{proposition:I-tensor-poincare}. %
After applying $\F$, the transformation from the top diagram to the bottom one thus induces an equivalence on total cofibres of the right most faces. At the same time, the bottom right most face encodes the assembly map, see Remark~\reffour{remark:circle-decompose-1}.
We conclude that the resulting square
\[
\begin{tikzcd}
\F(\Hyp(\C)) \displaystyle\mathop{\oplus}_{\F(\Hyp(\C \oplus \C))} \F(\Hyp(\C_{\Del^1})) \ar[r]\ar[d] & \F(\Hyp(\C_{\B\NN^+}))\ar[d] \\
{\B\ZZ} \otimes \F(\C,\QF) \ar[r] & \F((\C,\QF)_{\B\ZZ}) 
\end{tikzcd}
\]
is cartesian.
To finish the proof it will now suffice to show that the map
\[
\F((\Hyp(\C \oplus \C)) \to \F(\Hyp(\C_{\Del^1})),
\]
induced by the boundary inclusion $\{\ast,\ast\} \to \Del^1$, is an equivalence.
Indeed, this map 
is obtained by applying the Verdier-localising functor $\F^{\hyp}$ to the middle vertical map in the commutative rectangle
\[
\begin{tikzcd}
\C \ar[rr, "{x\mapsto (x,0)} "]\ar[d,equal] && \C \oplus \C\ar[rr, "{(x,y) \mapsto y}"]\ar[d] && \C \ar[d,equal] \\
\C \ar[rr, "{x \mapsto [x \to 0]}"] && \C_{\Del^1} \ar[rr, "{[x\to y] \mapsto y}"] && \C
\end{tikzcd}
\]
whose rows are (split) Verdier sequences. 
\end{proof}

\begin{remark}
\label{remark:proof-for-bordism-invariant}%
The proof of Theorem~\reffour{theorem:circle-decomposition-idem} described above simplifies when $\F$ is bordism invariant: in this case, $\F$ vanishes on the top diagram
as well as on the left most face in the bottom diagram. Since the rows are Poincaré-Verdier sequences we obtain an equivalence between the value of $\F$ on the remaining faces:
\[
\begin{tikzcd}
\F(\C,\QF) \oplus \F(\C,\QF)\ar[dd] \ar[dr] && \F(\C,\QF) \oplus \F(\C,\QF)\ar[dd] \ar[dr] & \\
& \F((\C,\QF)_{\mathcal{P}}) \ar[dd]\ar[dr,phantom,"\simeq"] && \F(\C,\QF)\ar[dd] \\
\F(\C,\QF)  \ar[dr] && \F(\C,\QF)\ar[dr]  & \\
& \F((\C,\QF)_\I)  && \F((\C,\QF)_{\B\ZZ}).
\end{tikzcd}
\]
Since the left one is exact so is the right, and we obtain that the assembly map is an equivalence. 
\end{remark}

Next, we address the parts of the circle decomposition theorem pertaining to the behaviour under idempotent completion, that is Theorem~\reffour{theorem:circle-decomposition}\;\reffouritem{item:relative-circle} and~\reffouritem{item:nil}.

\begin{lemma}
\label{lemma:S1-natural}%
For a stable $\infty$-category $\D$ and a Poincar\'e $\infty$-category $(\C,\QF)$, the squares
\[
\begin{tikzcd}
\D_{\B\NN} \ar[r] \ar[d] & \D \ar[d] & (\C,\QF)_{\B\ZZ} \ar[r] \ar[d] & (\C,\QF) \ar[d] \\
(\D^{\natural})_{\B\NN} \ar[r] & \D^{\natural} & (\C^\natural,\QF)_{\B\ZZ} \ar[r] & (\C^\natural,\QF)
\end{tikzcd}
\]
are cartesian in $\Catx$ and $\Catp$, respectively. 
\end{lemma}
\begin{proof}
The vertical maps in the squares are dense Karoubi and Poincaré-Karoubi inclusions, respectively.
At the same time, the horizontal arrows induce bijections on $\K_0$: indeed, they admit sections induced on tensors by the base point inclusions $\ast \to \B\NN$ and $\ast \to \B\ZZ$, and these section are surjective on $\K_0$ since their essential image generates their respective targets under finite colimits and desuspensions.
The desired result hence follows from Thomason's classification of dense subcategories.
\end{proof}

\begin{construction}
\label{construction:PV-square}%
Let $\D$ be a stable $\infty$-category. We let $\D^{\flat}_{\I} \subseteq \D_\I^\natural = (\D^{\natural})_\I$ be the dense subcategory spanned by those objects whose equaliser lies in $\D$. Similarly, given a Poincaré $\infty$-category 
$(\C,\QF)$ we define $(\C,\QF)^{\flat}_\I$ to be 
the Poincaré subcategory of $(\C^\natural,\QF)_\I$ spanned by $\C^{\flat}_\I$. We note that this subcategory is invariant under the duality since taking equalisers is duality invariant, see Remark~\reffour{remark:equ}.
By Lemma~\reffour{lemma:S1-natural}, the pasting lemma, and the fact that Poincaré-Verdier projections are closed under base change, see Remark~\reftwo{remark:cocartesian}~\reftwoitem{item:PV-stable-by-pullback}, 
the Poincaré-Verdier projection $(\C,\QF)^{\flat}_{\I} \to (\C^\natural,\QF)_{\B\ZZ}$ of Proposition~\reffour{proposition:PV-sequences-Q1} and the Verdier projection $\D^{\flat}_\I \to (\D^{\natural})_{\B\NN}$ from Proposition~\reftwo{proposition:localizationforp4}\;\reftwoitem{item:induced-functor-verdier} extend to the following Verdier and Poincaré-Verdier squares, respectively.
\[
\begin{tikzcd}
\D^{\flat}_\I \ar[r, dashed]  \ar[d, hook] &  \D_{\B\NN}\ar[d,hook] & (\C,\QF)^{\flat}_\I \ar[r, dashed]  \ar[d, hook] & (\C,\QF)_{\B\ZZ} \ar[d,hook] \\
\D^\natural_\I \ar[r] & (\D^\natural)_{\B\NN} &  (\C^\natural,\QF)_\I \ar[r] & (\C^\natural,\QF)_{\B\ZZ}. 
\end{tikzcd}
\]
\end{construction}

\begin{lemma}
\label{lemma:PV-sequence-3}%
The Poincaré $\infty$-category $(\C,\QF)^{\flat}_\I$ sits in a commutative diagram
\[
\begin{tikzcd}
	(\C,\QF) \ar[r,"i_1"] \ar[d, hook] & (\C,\QF)^{\flat}_\I \ar[r,"\ev_{-1}"] \ar[l,bend right,"\equ"']  \ar[d,hook] & (\C^\natural,\QF\qshift{1}) \ar[d, equal] \\
	(\C^\natural,\QF) \ar[r,"i_1"]& (\C^\natural,\QF)_\I \ar[r,"\ev_{-1}"] \ar[l,bend right,"\equ"'] & (\C^\natural,\QF\qshift{1})
\end{tikzcd}
\]
whose top row is a Poincaré-Verdier sequence, and whose bottom row is the split Poincaré-Verdier sequence of
Proposition~\reffour{proposition:I-tensor-poincare}.
Similarly, if $\D$ is a stable $\infty$-category then we have a diagram in $\Catx$ with the analogous properties involving $\D^{\flat}_\I$ and $\D^{\natural}_\I$.
\end{lemma}
\begin{proof}
Since the vertical maps are Poincar\'e-Karoubi equivalences and the bottom row is a split Poincaré-Verdier sequence, we find that the top row is a Poincar\'e-Karoubi sequence. Since the top right horizontal functor is also essentially surjective and its kernel is precisely the image of $(\C,\QF)$ under $i_1$,
we conclude that the top row is a Poincaré-Verdier sequence. The same argument also applies on the level of stable $\infty$-categories.
\end{proof}

The following two corollaries finish the proof of Theorem~\reffour{theorem:circle-decomposition}.

\begin{corollary}
\label{corollary:S1-natural}%
Let $\F\colon \Catp \to \A$ be a Verdier-localising invariant with stable target $\A$ and $(\C,\QF)$ a Poincar\'e $\infty$-category. Then both squares in the diagram  
\[
\begin{tikzcd}
	(\C,\QF) \ar[r]\ar[d]  & (\C,\QF)_{\B\ZZ} \ar[d]\ar[r] & (\C,\QF) \ar[d] \\
	(\C^\natural,\QF) \ar[r] & (\C^\natural,\QF)_{\B\ZZ} \ar[r] & (\C^\natural,\QF)
\end{tikzcd}
\]
are sent to cartesian squares by $\F$. In particular, Theorem~\reffour{theorem:circle-decomposition}\;\reffouritem{item:relative-circle} holds true.
Similarly, if $\D$ is a stable $\infty$-category then the analogous diagram in $\Catx$ is sent to a diagram with cartesian squares by any Verdier-localising functor $\G\colon \Catx \to \A$.
\end{corollary}
\begin{proof}%

Since the horizontal composites in the diagram under investigation are the identity, it suffices to consider the left square, which in turn is obtained from pasting the following two squares:
\[
\begin{tikzcd}
	(\C,\QF) \ar[r,"i_1"] \ar[d] & (\C,\QF)^{\flat}_\I \ar[r] \ar[d] & (\C,\QF)_{\B\ZZ} \ar[d] \\
	(\C^\natural,\QF) \ar[r,"i_1"] & (\C^\natural,\QF)_\I \ar[r] & (\C^\natural,\QF)_{\B\ZZ}  
\end{tikzcd}
\]
Here, the left square 
is sent to a cartesian square by $\F$ as a consequence of Lemma~\reffour{lemma:PV-sequence-3}. The right square is a Poincaré-Verdier square, see Construction~\reffour{construction:PV-square}, and is hence sent to a cartesian square since $\F$ is Verdier-localising, giving the desired result. The same argument applies when $\D$ is a stable $\infty$-category and $\F$ is a Verdier-localising functor on $\Catx$.
\end{proof}

\begin{corollary}
\label{corollary:nil-is-invariant}%
Let $\G\colon \Catx \to \cA$ be a Verdier-localising functor valued in a stable $\infty$-category $\cA$ and $\D$ a stable $\infty$-category. Then both squares in the diagram
\[
\begin{tikzcd}
	\D \ar[r]\ar[d] & \D_{\B\NN} \ar[r] \ar[d] & \D \ar[d] \\
	\D^\natural \ar[r] & (\D^{\natural})_{\B\NN} \ar[r] & \D^\natural
\end{tikzcd}
\]
are sent to cartesian squares by $\G$. In particular, Theorem~\reffour{theorem:circle-decomposition}\;\reffouritem{item:nil} holds true.
\end{corollary}
\begin{proof}
As before, it suffices to consider the left square, which is again obtained from pasting the following two squares:
\[
\begin{tikzcd}
	\D \ar[r] \ar[d] & \D^\flat_\I \ar[r] \ar[d] & \D_{\B\NN} \ar[d] \\
	\D^\natural \ar[r] & \D^\natural_\I \ar[r] & (\D^\natural)_{\B\NN}
\end{tikzcd}
\]
The left square is sent to a cartesian square by $\G$ as a consequence of (the stable $\infty$-category version of) Lemma~\reffour{lemma:PV-sequence-3}. The right square is a Verdier square, see Construction~\reffour{construction:PV-square}, and is hence sent to a cartesian square since $\G$ is Verdier-localising, giving the desired result.
\end{proof}

\subsection{The extended assembly map}
\label{subsection:extended}%

When $\C$ is idempotent complete or $\F$ Karoubi-localising Theorem~\reffour{theorem:circle-decomposition} can be expressed as saying that the assembly map $\B\ZZ \otimes \F(\C,\QF) \to \F((\C,\QF)_{\B\ZZ})$ is a split inclusion with complementary summand given by the nil term $\N\F^{\hyp}(\C)$ and therefore an equivalence if $\F$ is in addition bordism invariant. 

When $\C$ is not idempotent complete one can still obtain such a formulation by replacing the assembly map with the base change  
\[
[\B\ZZ \otimes \F(\C^{\natural},\QF)] \times_{(\C^{\natural},\QF)_{\B\ZZ}} (\C,\QF)_{\B\ZZ} \to (\C,\QF)_{\B\ZZ}
\]
of the assembly map at $(\C^{\natural},\QF)$.
Our goal in this subsection is to give a $\Catp$-level description for the domain of this map, which we shall henceforth call the \defi{extended assembly map}.
We will make use of this description when %
analysing the difference between the algebraic and geometric decompositions in \S\reffour{subsection:splittings}.

\begin{construction}
\label{construction:ext-ass}%
Let $\F\colon \Catp \to \A$ be a direct sum preserving functor to a stable $\infty$-category $\A$. 
Denote by 
\[
(\C,\QF) \oplus^{\flat} (\C,\QF) \subseteq (\C^{\natural},\QF) \oplus (\C^{\natural},\QF)
\]
the full Poincaré subcategory spanned by those pairs $(x,y) \in \C^{\natural} \oplus \C^{\natural}$ such that $x \oplus y \in \C$, so that the sum map $(\C^\natural,\QF) \oplus (\C^{\natural},\QF) \to (\C^{\natural},\QF)$ restricts to a map $\oplus\colon(\C,\QF) \oplus^{\flat} (\C,\QF) \to (\C,\QF)$. We then define a new functor ${\B\ZZ} \otimes^{\flat}\F \colon \Catp \rightarrow \A$ via the pushout in the square
\[
\begin{tikzcd}
\F((\C,\QF) \oplus^{\flat} (\C,\QF)) \ar[r,"\F(\oplus)"]\ar[d,"\F(\oplus)"'] & \F(\C,\QF) \ar[d] \\
\F(\C,\QF) \ar[r] & {\B\ZZ} \otimes^{\flat} \F(\C,\QF) \ .
\end{tikzcd}
\]
Let us make explicit again that ${\B\ZZ} \otimes^{\flat} \F(\C,\QF)$ is not given by some operation applied to $\F(\C,\QF)$. We note that the inclusion $(\C,\QF) \oplus (\C,\QF) \to (\C,\QF) \oplus^\flat (\C,\QF)$ induces a map $\B\ZZ \otimes \F(\C,\QF) \to \B\ZZ \otimes^\flat \F(\C,\QF)$ which is an equivalence if $\C$ is idempotent complete.
The square
\[
\begin{tikzcd}
(\C,\QF) \oplus^{\flat} (\C,\QF) \ar[r]\ar[d] & (\C,\QF) \ar[d] \\
(\C,\QF) \ar[r] & (\C,\QF)_{\B\ZZ} \ ,
\end{tikzcd}
\]
obtained via restriction from the pushout square of $(\C^{\natural},\QF)_{\B\ZZ}$ obtained from writing $\B\ZZ$ as the pushout $\ast \amalg_{\{\ast,\ast\}} \ast$, then determines a map
\[
\ass^{\flat}\colon {\B\ZZ} \otimes^{\flat} \F(\C,\QF) \to \F((\C,\QF)_{\B\ZZ}) ,
\]
whose precomposition with the above mentioned tautological map $\B\ZZ \otimes \F(\C,\QF) \to \B\ZZ \otimes^\flat \F(\C,\QF)$ is the assembly map. 
\end{construction}

We %
call $\ass^{\flat}$ %
the extended assembly map. It coincides with base change of the assembly map of $(\C^{\natural},\QF)$ by the next remark.
If $\C$ is idempotent complete, the extended assembly map recovers the assembly map.

\begin{remark}
\label{remark:circle-decompose-2}%
We have a Poincaré-Verdier sequence
\[
(\C,\QF) \xrightarrow{x \mapsto (x,0)} (\C,\QF) \oplus^{\flat} (\C,\QF) \xrightarrow{(x,y) \mapsto y} (\C^{\natural},\QF)
\]
whose inclusion is retracted by the Poincaré functor $\oplus \colon (\C,\QF) \oplus^{\flat} (\C,\QF) \to (\C,\QF)$. We thus obtain natural decompositions
\[
\B\ZZ \otimes^{\flat} \F(\C,\QF) = \F(\C,\QF) \oplus \Sig\F(\C^{\natural},\QF) \quad\text{and}\quad \F((\C,\QF) \oplus^{\flat} (\C,\QF)) = \F(\C,\QF) \oplus \F(\C^{\natural},\QF)
\]
and conclude that the external rectable in %
\[
\begin{tikzcd}
\B\ZZ \otimes^{\flat} \F(\C,\QF) \ar[r]\ar[d]& \F((\C,\QF)_{\B\ZZ})\ar[r]\ar[d] & \F(\C,\QF) \ar[d] \\
\B\ZZ \otimes \F(\C^\natural,\QF) \ar[r] & \F((\C^{\natural},\QF)_{\B\ZZ})\ar[r]& \F(\C^\natural,\QF)
\end{tikzcd}
\]
is cartesian with horizontal fibres given by $\Sig\F(\C^{\natural},\QF)$. It then follows that both squares are cartesian by Corollary~\reffour{corollary:S1-natural}.

\end{remark}

Unravelling the definitions, we find that the circle decomposition theorem now implies: 
\begin{corollary}
For any Poincar\'e $\infty$-category $(\C,\QF)$ and any Verdier-localising functor $\F$, the extended assembly map 
\[
\B\ZZ \otimes^\flat \F(\C,\QF) \to \F((\C,\QF)_{\B\ZZ})
\]
is a split inclusion with complementary term the nil term $\N\F^\hyp(\C)$. In particular, the extended assembly map is an equivalence if $\F$ is bordism invariant.
\end{corollary}

\begin{remark}
\label{remark:algebraic-splitting-2}%
\label{remark:algebraic-splitting}%
Recall from 
Remark~\reffour{remark:I-is-assembly} that for $\F$ a Verdier-localising functor and $(\C,\QF)$ a Poincaré $\infty$-category 
the assembly map can be identified with the induced map
\[
\F(\C,\QF) \displaystyle\mathop{\oplus}_{\F((\C,\QF)_{\mathcal{P}})} \F((\C,\QF)_\I)  \to \F((\C,\QF)_{\B\ZZ}).
\]
Let us complement this by explaining how to do the same for the extended assembly map.
For this, let $(\C,\QF)^{\flat}_{\mathcal{P}}$ be the dense Poincaré subcategory of $(\C^\natural,\QF)_{\mathcal{P}}$ spanned by those objects $[x \rightarrow z \leftarrow y] \in \C^{\natural}_{\mathcal{P}}$ whose image $x \times_z y \in \C^{\natural}$ under the degeneracy functor $s$ belongs to $\C$. By construction, the diagram of Proposition~\reffour{proposition:I-tensor-poincare} and Remark~\reffour{rem:variant} then restricts to a commutative diagram
\[
\begin{tikzcd}
(\C,\QF) \ar[r] \ar[d,"i_r"] & (\C,\QF)^\flat_{\mathcal{P}} \ar[r,"\ev"] \ar[d,equal] & \Met(\C^\natural,\QF\qshift{1}) \ar[d] \\
(\C,\QF)\oplus^{\flat} (\C,\QF) \ar[d]\ar[r] & (\C,\QF)^{\flat}_{\mathcal{P}} \ar[d]\ar[r,"\ev_{-1}"] & (\C^{\natural},\QF\qshift{1}) \ar[d,equal] \\
(\C,\QF) \ar[r] & (\C,\QF)^{\flat}_\I \ar[r,"\ev_{-1}"] & (\C^{\natural},\QF\qshift{1})
\end{tikzcd}
\]
in which the rows are Poincaré-Verdier sequences by Lemma~\reffour{lemma:PV-sequence-3} (and the analogous argument for the top rows).
As in Remark~\reffour{remark:I-is-assembly}, we obtain a commutative triangle 
\[
\begin{tikzcd}
\F(\C,\QF) \displaystyle\mathop{\oplus}_{\F((\C,\QF) \oplus^{\flat} (\C,\QF))} \F(\C,\QF)\ar[dr,"\ass^{\flat}"]\ar[rr,"\simeq"] && \F(\C,\QF) \displaystyle\mathop{\oplus}_{\F((\C,\QF)^{\flat}_{\mathcal{P}})} \F((\C,\QF)^{\flat}_\I)  \ar[dl] \\
& \F((\C,\QF)_{\B\ZZ}) &
\end{tikzcd}
\]
whose horizontal functor induces an equivalence 
from the extended assembly map on the left to the diagonal map on the right. In particular, if $\F$ is bordism invariant then the extended assembly map is equivalent to the map induced by the Poincaré functor $(\C,\QF)^{\flat}_\I \to (\C,\QF)_{\B\ZZ}$. Again, just as in Remark~\reffour{remark:I-is-assembly}, it follows that splitting off the canonical summand $\F(\C,\QF)$ in the above horizontal equivalence, we obtain the canonical equivalence
\[
\Sigma\F(\C^\natural,\QF) \xrightarrow{\simeq} \cof[\F(\Met(\C^\natural,\QF\qshift{1})) \to \F(\C^\natural,\QF\qshift{1})]
\]
induced by the metabolic fibre sequence.
\end{remark}

\subsection{The framed transfer for the circle}
\label{subsection:splittings}%

In this section, we discuss an observation due to Ranicki that for L-spectra, in fact more generally for any bordism invariant Verdier-localising invariant $\F$, in addition to the Ranicki-Shaneson decomposition obtained in \reffour{corollary:algebraic-shaneson-splitting}, there is a second such decomposition which Ranicki refers to as the geometrically relevant decomposition. It is based on studying a suitable product with the visible signature of the circle $S^1 \simeq \B\ZZ$, thought of as a Poincar\'e element of $(\Spa^\free,(\QF^\uni)\qshift{-1})_{\B\ZZ}$ using a stable framing of $S^1$. Following Ranicki, we make explicit the difference between the two decompositions. It follows from this analysis that the decompositions agree when one uses the null-bordant stable framing and restricts attention to idempotent complete Poincar\'e $\infty$-category, but not in general.

Let us now explain the situation in more detail.
\begin{construction}
\label{const:transfer}%
To begin, recall that a finite connected Poincar\'e complex $X \in \Sps$, with dualising complex $\nu_X\colon X\to \Pic(\SS) \subset \Spa$, has a visible signature $\vsign(X)$, a Poincaré object in $(\Spa_X^\free,\QF^\vis_{\nu_X})$, see \S\refone{subsection:parametrised-spectra} or Appendix~\reffour{section:appendix-visible-revisited} below for a second perspective.
Given a Poincar\'e $\infty$-category $(\C,\QF)$, let us write 
\[
(\C, \QF)_{(X, \nu_X)} \coloneq \colim_{x\in X} (\C, \QF\otimes \nu_X(x)) = (\C,\QF) \otimes (\Spaf, \QF^\uni)_{(X, \nu_X)} = (\C,\QF) \otimes (\Spa_X^\free,\QF^\vis_{\nu_X})
\]
for the colimit of the canonical diagram induced by $\nu_X$. Recall that the tensor unit $(\Spaf,\QF^\uni)$ of $\Catp$ corepresents the Poincar\'e objects functor; In particular, $\vsign(X)$ is equivalently described by a Poincar\'e functor $(\Spaf,\QF^\uni) \to (\Spa_X^\free,\QF^\vis_{\nu_X})$. 
Tensoring the visible signature $\vsign(X)$ of $X$ with 
$(\C,\QF)$ then gives the Poincar\'e functor 
\[
(\C, \QF) \to (\C, \QF)\otimes (\Spaf, \QF^\uni)_{(X, \nu_X)} \simeq (\C, \QF)_{(X, \nu_X)},
\]
which we refer to as the \defi{parametrised transfer functor} associated to the Poincar\'e complex $X$. 
\end{construction}

\begin{remark}
\label{rem:framings}%
Suppose that the dualising complex $\nu_X$ is equipped with a framing, that is, an equivalence  %
\[
\tau\colon \nu_X \simeq \SS_X\qshift{-n}
\]
where $n=\dim(X)$ and $\SS_X\qshift{-n}$ is the constant local system on $X$ with value $\SS\qshift{-n} = \Omega^n\SS$.
Then we may form the composite, often referred to as \emph{spherical fundamental class} %
\[
c\colon \SS \longrightarrow \Th(\nu_X) \xrightarrow{\Th(\tau)} \Omega^n\Th(\SS_X) = \SS\qshift{-n}[X], 
\]
where for a stable spherical fibration $\xi\colon X\to \Pic(\SS)$, we write $\Th(\xi) = \colim_X \xi$ for its Thom spectrum and as above $\SS\qshift{-n}[X] = \Omega^n\SS[X]$. The first map in the above composite, often call the collapse map of $X$, is the composite of the constant inclusion $\SS\to \map(X, \SS)$ and the equivalence $\map(X, \SS)=\lim_X \SS\simeq \Th(\nu_X)$ issued from the defining property of the dualising complex.

Conversely, any map $c\colon \SS \to \SS\qshift{-n}[X]$ which induces an equivalence
\[
\lim_X E \simeq \hom_X(\SS^{-n}_X,\SS^{-n}_X \otimes E) \xrightarrow{\colim} \hom(\SS^{-n}[X], \colim_X \SS_X^{-n} \otimes E) \xrightarrow{c} \colim_X \mathbb S^{-n}_X \otimes E
\]
for all $E \in \Sp_X$ determines a framing of $\nu_X$ by the uniqueness of the dualising complex; its associated spherical fundamental class then coincides with $c$. Passing between $\tau$ and $c$ in this manner, one obtains an equivalence between the $\infty$-groupoid framings of $\nu_X$ and the full subgroupoid of $\Omega^{\infty+n} \SS[X] = \Hom(\SS,\SS\qshift{-n}[X])$ spanned by the spherical fundamental classes. %
We also note that any spherical fundamental class $c$ induces, %
via the projection $q\colon X\to *$, an element 
\[
q_*[c]\in \pi_n(\SS)
\]
which we call the Pontryagin-Thom class of $c$.
If $X$ is a stably framed smooth manifold, then $q_*[c]$ is the element determined by the usual Pontryagin-Thom isomorphism between stable homotopy and framed bordism groups.  
\end{remark}

\begin{construction}
\label{construction:visible-sig-c}%
Given a spherical fundamental class $c\colon \SS \to \SS\qshift{-n}[X]$, we may view $X$ as a framed Poincar\'e complex as described in Remark~\reffour{rem:framings} and obtain an induced identification %
$(\C, \QF)_{(X, \nu_X)}\simeq (\C, \QF\qshift{-n})_X$ for any Poincaré $\infty$-category $(\C,\QF)$. We may hence view the visible signature of $X$ as a Poincar\'e object $\vsign_c(X)$ in $(\Sp^\free, (\QF^\uni)\qshift{-n})_X$ where we have added the subscript $c$ to keep track of the dependence of this construction on the spherical fundamental class $c$. 
By Proposition~\reffour{proposition:visible-signature-revisited},
under the embedding $\Spaf_X \subseteq \Fun(X\op,\Pro(\Spaf))$ (see Construction~\refone{construction:tensor} in the case of tensoring by an $\infty$-groupoid) the Poincaré object $\vsign_c(X)$ corresponds to the 
constant diagram $\SS_X$ equipped with the Poincar\'e form
\[
\beta_c\colon \SS \xrightarrow{\bet_{\uni}} \QF^{\uni}(\SS) \xrightarrow{c\otimes \QF^{\uni}(\SS)} \SS\qshift{-n}[X] \otimes \QF^{\uni}(\SS) = \colim_{X} \QF^{\uni}(\SS)\qshift{-n} = (\QF^\uni)_X\qshift{-n}(\SS_X) ,
\]
where $\bet_{\uni}$ is the universal Poincaré form on $\SS$. 
The unique map $q\colon X \to \ast$ induces a Poincar\'e functor
\[
q_*\colon (\Spaf,(\QF^{\uni})\qshift{-n})_{X} \longrightarrow (\Spaf,(\QF^{\uni})\qshift{-n})
\]
which we may apply to the visible signature to obtain a 
Poincaré object 
$q_*\vsign_c(X) = (q_*\SS_X,q_*\beta_c) = (\SS^{X},q_*\beta_c)$
in $(\Sp^\free,(\QF^\uni)\qshift{-n})$. 
Tensoring with $\vsign_c(X)$ and $q_*\vsign_c(X)$ then yields 
Poincar\'e functors
\[
\tr_{(X,c)}\colon (\C, \QF\qshift{n})\to (\C, \QF)_X \quad \text{and}\quad \tr_{(X,c)}^q \colon \quad (\C,\QF\qshift{n}) \to (\C,\QF)
\]
which we refer to as the \defi{(parametrised) framed transfer} associated to $X$ and its spherical fundamental class $c$. We note that the framed transfer is obtained from the parametrised framed transfer by post-composition with the map induced by $q\colon X \to \ast$.%
\end{construction}

We now apply these observations to $X=\B\ZZ$. Using the splitting 
\[
\SS\qshift{-1}[\B\ZZ] = \SS \oplus \SS\qshift{-1}
\]
of Remark~\reffour{remark:circle-decompose-1}
the $\infty$-groupoid of spherical fundamental classes identifies with the 
subgroupoid of $\Hom(\SS,\SS \oplus \SS\qshift{-1})$
spanned by those maps $c\colon \SS \to \SS \oplus \SS\qshift{-1}$ whose post-composition with the projection to $\SS$ is an equivalence. This $\infty$-groupoid therefore splits as a product $\gl_1(\SS) \times \Om^{\infty+1}\SS$, canonically pointed by the pair $(\id,0)$. The set of connected components of this $\infty$-groupoid is
\[
\pi_0\gl_1(\SS) \times \pi_1(\SS) = \{\pm1\} \times \ZZ/2 ,
\]
so that the equivalence type of any fundamental class is encoded by a sign 
$\eps(c) \in \{\pm1\}$,
which encodes the gap between the orientation determined by $c$ and our fixed orientation of $\B\ZZ$, and the associated Pontryagin-Thom class $q_*[c] \in \pi_1(\SS) = \ZZ/2$.

\begin{remark}
\label{remark:pontryagin}%
Geometrically, the Pontryagin-Thom class determines whether or not 
the corresponding stable framing of the 1-circle extends to a stable framing of the 2-disc. 
From the present perspective, this is equivalent to the question whether or not the the map
$\SS \to \SS^{\B\ZZ}$ induced by $q\colon \B\ZZ \to \ast$ refines to a Lagrangian of the Poincaré object $q_*\vsign_c(\B\ZZ)$ in $(\Spaf,(\QF^{\uni})\qshift{-1})$.
Indeed, unwinding the definitions we get that the restriction of the form $q_*\beta_c \in (\QF^{\uni})\qshift{-1}(\SS^{\B\ZZ})$ %
along the constant inclusion $\SS\to \SS^{\B\ZZ}$ %
is given by the composite
\[
\SS \xrightarrow{q_*[c]} \SS^{-1} \to (\QF^{\uni})\qshift{-1}(\SS)
\]
where the second map is given by the universal Poincar\'e form. This composite is then null-homotopic if and only if the first map is, since $\SS \to \QF^{\uni}(\SS)$ admits a retraction $\QF^{\uni}(\SS) \to  \Lam^{\uni}(\SS) = \SS$ by projecting to the linear part of the universal Poincaré structure. In addition, when this composite is null-homotopic, any null-homotopy exhibits $\SS$ as a Lagrangian in $q_*\vsign_c(\B\ZZ)$, just like any extension of a stable framing of the 1-circle to the 2-disc constitutes a stably framed null-cobordism. In particular, if $q_*[c] = 0$, then $q_*\vsign_c(\B\ZZ)$ is null-cobordant. Conversely, it can also be shown that if $q_*[c] \neq 0$, then $q_*\vsign_c(\B\ZZ)$ is not null-cobordant.
\end{remark}

A specific feature of $\B\ZZ$ is that the $\K$-class underlying the visible signature $\vsign_c(\B\ZZ)$ is trivial. As a consequence, the (parametrised) framed transfers on $(\C,\QF)$ and $(\C^\natural,\QF)$, respectively, factor as follows:
\[
\begin{tikzcd}[column sep=large]
	(\C,\QF\qshift{1}) \ar[r] & (\C^\natural, \QF\qshift{1}) \ar[r,"{\tr^\flat_{(\B\ZZ,c)}}", dashed] & (\C, \QF)_{\B\ZZ} \ar[r,hook]  & (\C^\natural, \QF)_{\B\ZZ}
\end{tikzcd}
\]
We call $\tr^\flat_{(\B\ZZ,c)}$ the \defi{extended parametrised framed transfer} associated to $(\B\ZZ,c)$ and its composite $\tr^{\flat,q}_{(\B\ZZ,c)}$ with the tautological map $(\C,\QF)_{\B\ZZ} \to (\C,\QF)$ induced by $q \colon \B\ZZ \to \ast$ the \defi{extended framed transfer}.

\begin{proposition}
\label{proposition:geometric}%
Let $\F\colon \Catp \to \A$ be a Verdier-localising functor valued in a stable $\infty$-category $\A$. Under the equivalence of Theorem~\reffour{theorem:circle-decomposition}, the extended parametrised framed transfer
\[
\F(\C^\natural,\QF\qshift{1}) \xrightarrow{\tr^{\flat}_{(\B\ZZ,c)}} \F((\C, \QF)_{\B\ZZ})\simeq \F(\C,\QF) \oplus \Sig\F(\C^\natural,\QF) \oplus \N\F^{\hyp}(\C)
\]
has the following components:
\begin{enumerate}
\item
\label{item:base-component}%
Its component at $\F(\C,\QF)$ is induced by the extended framed transfer $\tr_{(\B\ZZ,c)}^{\flat,q}$ associated to $(\B\ZZ,c)$.
\item
\label{item:nil-component}%
Its component at $\N\F^{\hyp}(\C)$ is zero.
\item
\label{item:sig-component}%
Its component at $\Sig\F(\C^\natural,\QF)$ is given by $\eps(c)\cdot \partial$, where $\partial$ is the boundary map of the metabolic sequence. 
\end{enumerate}
\end{proposition}

\begin{corollary}%
\label{corollary:geometric-shaneson-splitting}%
Let $\F\colon \Catp \to \A$ be a bordism invariant and Verdier-localising functor valued in a stable $\infty$-category $\A$. 
Then the Poincaré functor $(\C,\QF) \to (\C,\QF)_{\B\ZZ}$ induced by 
the base point $\ast \to \B\ZZ$ and the extended framed transfer $\tr_{\B\ZZ}^{\flat,q}$ of $(\B\ZZ,c)$ induce an equivalence
\[
\F(\C,\QF) \oplus \Sig\F(\C^\natural,\QF) \simeq \F(\C,\QF) \oplus \F(\C^\natural,\QF\qshift{1}) \xrightarrow{\simeq} \F(\C,\QF)_{\B\ZZ} .
\]
\end{corollary}

We warn the reader again that this \defi{geometric Ranicki-Shaneson decomposition} is generally different from the one appearing in Corollary~\reffour{corollary:algebraic-shaneson-splitting}, as will discuss below.

\begin{proof}[Proof of Proposition~\reffour{proposition:geometric}]
Statement~\reffouritem{item:base-component} holds by definition. Next, we recall 
the Poincaré-Verdier projection
\[
p_*\colon (\Spaf,(\QF^{\uni})\qshift{-1})_{\I}^\flat \to (\Spaf,(\QF^{\uni})\qshift{-1})_{\B\ZZ}.
\]
We will show momentarily that the visible signature $\vsign_c(\B\ZZ)$ lifts 
to a Poincar\'e object $\bar{\sigma}^\vis_c(\I)$ in $(\Spaf,(\QF^{\uni})\qshift{-1})_{\I}^\flat$. As a consequence, the refined parametrised transfer factors as a composite
\[
(\C^\natural,\QF\qshift{1}) \to (\C,\QF)^\flat_\I \to (\C,\QF)_{\B\ZZ}.
\]
By Remark~\reffour{remark:algebraic-splitting-2}, we find that second displayed functor, after applying $\F$, factors as
\[
\F((\C,\QF)^\flat_\I) \to \B\ZZ\otimes^{\flat}\F(\C,\QF) \to \F((\C,\QF)_{\B\ZZ})
\]
where the second map is the extended assembly map of \S\reffour{subsection:extended}.
Since $\N\F^\hyp(\C)$ is the cofibre of the extended assembly map, we deduce statement~\reffouritem{item:nil-component}. 
For~\reffouritem{item:sig-component}, note that, by the last part of Remark~\reffour{remark:algebraic-splitting-2}, the left map in the last factorisation induces, after splitting off the summand $\F(\C,\QF)$ on both sides, the boundary map
\[
\F(\C,\QF\qshift{1}) \to \Sig\F(\C^{\natural},\QF)
\]
arising from the metabolic sequence.
To show statement~\reffouritem{item:sig-component} we show momentarily that the composite
\[
(\C^\natural,\QF\qshift{1}) \to (\C,\QF)_\I^\flat \xrightarrow{\ev_{-1}} (\C^\natural,\QF\qshift{1})
\]
is given by the Poincar\'e functor $(\id,\eps(c)\cdot \id)$. Taking this for granted for the moment, the result follows immediately in case $\eps(c) =1$. If instead $\eps(c) = -1$, we recall from Corollary~\reftwo{corollary:inversion} that $\F(\id,-\id) = -\id_{\F(\C^\natural,\QF\qshift{1})} - \F(\mathrm{N}\Omega)$, where $\mathrm{N}\Omega$ is the norm, in the sense of Construction~\reftwo{construction:norm}, of the exact functor $\Omega$ on $\C^\natural$. Concretely, $\mathrm{N}\Omega$ is given by the composite
\[
(\C^\natural,\QF\qshift{1}) \xrightarrow{\Omega^\hyp} \Hyp(\C^\natural) \to (\C^\natural,\QF\qshift{1}).
\]
where the latter functor is the canonical functor using that $\Hyp(-)$ is left adjoint to the forgetful functor. This functor in turn factors as 
\[
\Hyp(\C^\natural) \xrightarrow{\mathrm{can}} \Met(\C^\natural,\QF\qshift{1}) \to (\C^\natural,\QF\qshift{1})
\]
whose second map is the Poincar\'e-Verdier projection of the metabolic sequence.
As a result, the composite 
\[
\F(\C^\natural,\QF\qshift{1}) \xrightarrow{\mathrm{N}\Omega} \F(\C^\natural,\QF\qshift{1}) \xrightarrow{\partial} \Sigma \F(\C^\natural,\QF)
\]
vanishes canonically, and the result follows.

Now, the remaining claim is implied by the statement that the
composite
\[
(\Spaf,\QF^\uni) \to (\Spaf,(\QF^\uni)\qshift{-1})_\I^\flat \xrightarrow{\ev_{-1}} (\Spaf,\QF^\uni)
\]
classifies $(\SS,\eps(c)\cdot \beta_u)$. To proceed we now have to make explicit the construction of $\bar{\sigma}^\vis_c(\I)$.
To that end, we note that the constant diagram $[\begin{tikzcd} \SS & \SS \ar[l,shift left,"\id"] \ar[l, shift right,"\id"'] \end{tikzcd}]$
maps to 
the constant diagram $\SS_{\B\ZZ}$ under the functor $\Spaf_\I \to \Spaf_{\B\ZZ}$, since $\I \to \B\ZZ$ is a localisation. %
Moreover, the induced map
\[
(\QF^{\uni})\qshift{-1}_\I(\SS \leftleftarrows \SS) =  \colim_\I (\QF^{\uni})\qshift{-1}(\SS) \to  \colim_{\B\ZZ}(\QF^{\uni})\qshift{-1}(\SS) = (\QF^{\uni}_{\B\ZZ})\qshift{-1}(\SS_X)
\]
is then an equivalence since any localisation is cofinal. %
We may consequently (essentially uniquely) lift the visible signature form $\bet_{c}$ on $\SS_{\B\ZZ}$ to a hermitian form $\bet^\I_c \in (\QF^{\uni})\qshift{-1}_\I([\SS \leftleftarrows \SS])$. Concretely, $\beta^\I_c$ is given by the composite
\begin{equation}
\label{equation:beta-I-c}%
\SS \xrightarrow{\bet_{\uni}} \QF^{\uni}(\SS) \xrightarrow{c\otimes \QF^{\uni}(\SS)} \SS\qshift{-1}[\B\ZZ] \otimes \QF^{\uni}(\SS) = \colim_{\B\ZZ} (\QF^{\uni})\qshift{-1}(\SS) \simeq \colim_{\I} (\QF^{\uni})\qshift{-1}(\SS),
\end{equation}
where $\bet_{\uni}$ is the universal Poincaré form on $\SS$. Given this, it is readily computed that the Poincar\'e functor $\ev_{-1}$ takes $\beta_c^\I$ to $\eps(c)\cdot \beta_u$.
It remains to note that the hermitian object $(\SS \leftleftarrows \SS,\bet^\I_c)$ in $(\Spaf,(\QF^{\uni})\qshift{-1})_\I$ is in fact Poincar\'e. For this, one can employ the embedding $\Spaf_\I \subseteq \Fun(\I\op,\Pro(\Spaf))$, which takes $\SS \leftleftarrows \SS$ to the constant diagram on $\SS$ and observe that also its dual is constant by the explicit formula \eqreffour{equation:explicit-dual-I}. But under the analogous embedding for $\Spaf_{\B\ZZ}$ the functor $p_* \colon \Spaf_\I \rightarrow \Spaf_{\B\ZZ}$ is induced by right Kan extension along $p\colon \I\op \to \B\ZZ\op$, which has restriction as a fully faithful adjoint. Since $\SS \leftleftarrows \SS$ and its dual both lie in the image of this restriction functor and the map connecting them is taken to the equivalence $\beta_c$ by $p_*$, it is itself an equivalence as claimed.
\end{proof}

We finish this subsection by discussing the difference between the geometric and algebraic decompositions of $\F(\C,\QF)_{\B\ZZ}$ into a direct sum of $\F(\C,\QF)$ and $\Sig\F(\C^{\natural},\QF)$, for a bordism invariant and Verdier-localising functor $\F$. The above shows that under the canonical equivalence $\Sigma \F(\C,\QF) \simeq \F(\C,\QF\qshift{1})$ the composite
\[
\begin{tikzcd}[column sep=large]
	\F(\C,\QF) \oplus \Sigma \F(\C^\natural,\QF) \ar[r,"\simeq"] & \F((\C,\QF)_{\B\ZZ}) \ar[r,"\simeq","{\reffour{corollary:algebraic-shaneson-splitting}}"'] & \F(\C,\QF) \oplus \Sigma \F(\C^\natural,\QF) 
\end{tikzcd}
\]
where the first and second equivalence are the ones from Corollaries~\reffour{corollary:geometric-shaneson-splitting} and~\reffour{corollary:algebraic-shaneson-splitting}, respectively,
is given by the following matrix:
\[
\begin{pmatrix}
	\id & 0 \\ \F(\tr^{\flat,q}_{\B\ZZ,c}) & \eps(c)\cdot \id
\end{pmatrix}
\]
Now, by construction, $\tr^{\flat,q}_{\B\ZZ,c}$ is given by multiplication by $q_*\vsign_c(\B\ZZ)$ and if the spherical fundamental class $c$ has Pontryagin-Thom class $q_*[c] = 0$, then $q_*(\vsign_c(\B\ZZ))$ is cobordant to $0$ in $(\Spaf,(\QF^\uni)\qshift{-1})$ by Remark~\reffour{remark:pontryagin}. Note that it is not evidently also cobordant to $0$ in the minimalisation of $\Spaf$, in fact, we will show later that it is \emph{not} null-cobordant therein. As a consequence of the above, we nevertheless obtain:

\begin{corollary} 
\label{cor:geometric=algebraic-decomposition}%
If the spherical fundamental class $c$ satisfies $q_*[c]=0$, the map 
\[
\F(\tr^{\flat,q}_{\B\ZZ,c})\colon \F(\C^\natural,\QF\qshift{1}) \to \F(\C,\QF)
\]
is null-homotopic after pre-composition with to $\F(\C, \QF\qshift{1})\to \F(\C^\natural,\QF\qshift{1})$ and after post-composition with $\F(\C,\QF) \to \F(\C^\natural,\QF)$. In particular, for idempotent complete Poincar\'e $\infty$-categories and for a spherical fundamental class $c$ with $(\eps(c),q_*[c]) = (1,0)$ in $\{\pm1\} \times \ZZ/2$, the algebraic and geometric decompositions of $\F(\C,\QF)_{\B\ZZ}$ coincide.
\end{corollary}

Let us now investigate further the case where $(\eps(c),q_*[c]) = (1,0)$ but $\C$ is not necessarily idempotent complete. For this, let $(\C,\QF)\subseteq (\D,\QFD)$ be a Poincar\'e-Karoubi equivalence, so that $\QF = \QFD|_{\C}$ and $\C$ contains the minimal dense subcategory $\D^0\subseteq \D$ spanned by the objects with trivial K-class. Our goal is to give a more explicit description of the associated extended framed transfer
\[
-\otimes q_*\vsign_c(\B\ZZ)\colon \F(\D, \QFD\qshift 1) \to \F(\C, \QF)
\]
As a consequence of Corollary~\reffour{cor:geometric=algebraic-decomposition}, this map is canonically nullhomotopic if we restrict the source to $\F(\C, \QF\qshift 1)$ and if we map the target into $\F(\D, \QFD)$. This simple observation admits the following strengthening:

\begin{proposition}
\label{proposition:tau-is-equiv}%
In the above situation, the extended framed transfer factors canonically as 
\[
\F(\D,\QFD\qshift{1}) \to \cof[\F(\C,\QF\qshift{1}) \to \F(\D,\QFD\qshift{1})]
\xrightarrow\simeq  \fib[\F(\C,\QF) \to \F(\D,\QFD)] \to \F(\C,\QF).
\]
where the  outer maps are the canonical ones, and the middle one is an equivalence. 
\end{proposition}

\begin{proof}
Write $\Met(\D,\C,\QFD) \subseteq \Met(\D,\QFD)$ for the dense Poincaré subcategory spanned by those $\cob \to x$ such that $x \in \C$, and $\Ar(\D,\C,\QFD) \subseteq \Ar(\D,\QFD)$ for the dense Poincaré subcategory spanned by those $f\colon y \to z$ such that $\cof(f) \in \C$. We note that $\Met(\D,\C,\QFD) \simeq \Ar(\D,\C,\Phi\qshift{-1})$ via the 
equivalence of Lemma~\refone{lemma:arrow-is-met}. In addition, these Poincaré $\infty$-categories then fit into the solid commutative squares 
\[
\begin{tikzcd}
(\C,\QF\qshift{1}) \ar[r]\ar[d,"\id_{(-)}"'] 
  & (\D,\QFD\qshift{1}) \ar[d,"\id_{(-)}"]  
  && \Met(\D,\C,\QFD) \ar[d,"\met"']\ar[r] 
  & \Met(\D,\QFD) \ar[d,"\met"] 
  \\
  \Ar(\C,\QF\qshift{1}) \ar[r] 
  & \Ar(\D,\C,\QFD\qshift{1})\ar[urr,dotted,"{-\otimes z}"] && (\C,\QF) \ar[r] & (\D,\QFD)
\end{tikzcd}
\]
in which the solid horizontal arrows are dense Poincaré-Karoubi inclusions. We now observe that the right solid square in the diagram is cartesian by construction and its right vertical arrow is a split Poincaré-Verdier projection, so that the whole square is a split Poincaré-Verdier square by Proposition~\reftwo{corollary:PV-proj-pullback}.
On the left hand side we have the dual phenomenon: the vertical arrows are both split Poincaré-Verdier inclusions (with left and right adjoint sending an arrow to its target and source, respecitvely), and the square is vertically left and right adjointable by direct inspection. The induced functor on vertical cofibres is an equivalence (where both sides can be identified with $(\C,\QF)$ via $\cof$), and so by Proposition~\reftwo{proposition:poincare-verdier-square-char} we conclude that the left solid square is cocartesian. 
The functor $\F$ therefore takes both solid squares to cartesian squares by additivity, and since $\F$ vanishes on $\Ar(\C,\QF)$ and $\Met(\C,\QF)$ by bordism invariance, we obtain identifications
\[
\F(\Ar(\D,\C,\QFD\qshift{1})) = \cof[\F(\C,\QF\qshift{1}) \to \F(\D,\QFD\qshift{1})] \quad\text{and}\quad \F(\Met(\D,\C,\QFD)) = \fib[\F(\C,\QF) \to \F(\D,\QFD)] .
\]

We next construct the diagonal Poincar\'e functor marked by a dotted arrow in the diagram: For this, we identify its target with $\Ar(\D, \C, \QF\qshift{-1})$ as above and will use that $\Ar\colon \Catp\to \Catp$ is lax symmetric monoidal for the tensor product of Poincar\'e categories, as a composite of the strong symmetric monoidal inclusion $\Catp\to \Cath$ and the lax symmetric monoidal right adjoint $\Pairings\colon \Cath\to \Catp$. On underlying categories, this lax symmetric monoidal structure is pointwise, i.e.\ given by 
\[
\Ar(\C)\otimes \Ar(\D)\to \Ar(\C\otimes \D), \quad (c\to c')\otimes (d\to d') \mapsto (c\otimes d\to c'\otimes d').
\]
By passing to dense subcategories, this symmetric monoidal refinement of $\Ar(-)$ induces a Poincar\'e functor 
\[
\Ar(\D, \C, \QFD\qshift n)\otimes \Ar(\Spaf, \Spa^0, (\QF^{\uni})\qshift m)\to \Ar(\D, \C, \QFD\qshift{n+m}),
\]
where we write $\Spa^0$ for the minimalisation $(\Spaf)^{\min}$ of $\Spaf$ consisting of all objects with trivial $\K$-class.

Let us now fix a null-homotopy of Pontryagin-Thom map %
$q_*[c]\colon \SS \to \SS\qshift{-1}$,
which by Remark~\reffour{remark:pontryagin}
determines %
a Lagrangian $(\SS \to \SS^{\B\ZZ},\bet_c,\eta)$, and let us write 
\[
z \in \Poinc(\Met(\Spaf, \Spa^0, (\QF^{\uni})\qshift{-1}))\simeq \Poinc(\Ar(\Spaf, \Spa^0, (\QF^{\uni})\qshift{-2}))
\]
for this corresponding Poincaré object.
In particular, the underlying object in $\Met(\Spaf)$ is given by the constant inclusion $\SS\to \SS^{\B\ZZ}$, and the corresponding object in $\Ar(\Spaf)$ is obtained by passing to fibres, giving a map $\SS\qshift{-2}\to \SS$. 
The diagonal Poincar\'e functor $\Ar(\D, \C, \QF\qshift{1}) \to \Ar(\D, \C, \QF\qshift{-1})$ in the above diagram is defined by tensoring with this Poincar\'e object $z = (\SS\qshift{-2} \to \SS,\gam)$. The composite functor
\[
(\D,\QFD\qshift{1})  \to \Ar(\D,\C,\QFD\qshift{1}) \xrightarrow{(-)\otimes z} \Ar(\D,\C,\QFD\qshift{1}) \xrightarrow{\cof} (\C,\QF)
\]
is then given by tensoring with $q_*\vsign_c(\B\ZZ)$.

We can therefore finish the proof by showing that the diagonal functor is a bordism equivalence. 
For this, it will suffice to show that the Poincaré object $z$
is graded bordism invertible, that is, it admits a tensor inverse in $\Poinc(\Ar(\Spaf,\Spa^0,(\QF^{\uni})\qshift{2}))$ up to cobordism. Let $(\SS^2,0) \in \spsforms(\Spaf,(\QF^{\uni})\qshift{2})$ be the 2-shifted hermitian object whose underlying object is the $2$-shifted sphere spectrum and whose hermitian form is $0$. Under proposition Proposition~\refone{proposition:thom} this corresponds to a hyperbolic Poincaré object of the form $([\SS^2 \xrightarrow{0} \SS],\hyp) \in \Poinc(\Ar(\Spaf,(\QF^{\uni})\qshift{2}))$, and since $\SS^2$ and $\SS$ have the same $\K$-class this Poincaré object actually lies in $\Ar(\Spaf,\Spa^0,(\QF^{\uni})\qshift{2})$. Its tensor product with $z=(\SS\qshift{-2}\to \SS, \gam)$ is
\[
([\SS^2 \xrightarrow{0} \SS],\hyp) \otimes ([\SS^{-2} \to \SS],\gam) = ([\SS \xrightarrow{0} \SS],\hyp),
\]
the hyperbolic Poincaré object corresponding to the trivial (unshifted) form on $\SS$. We now claim
that $([\SS \xrightarrow{0} \SS],\hyp)$ is cobordant to the unit Poincaré object $([\SS \xrightarrow{\id} \SS],u)$ of $\Ar(\Spaf,\Spa^0,\QF^{\uni})$. Indeed, this follows from the cofinality theorem \reffour{theorem:strong-cofinality}, whose $\L$-theory version gives 
\[
\L_0(\Ar(\Spaf, \Spa^0, \QF^{\uni}) \simeq \pi_0 (\L(\Spaf, \QF^{\uni})/\L(\Spa^0, \QF^{\uni})) \simeq \pi_0 (\K_0(\Spaf)/\K_0(\Spa^0))^\tC \simeq \ZZ/2
\]
where the composite map explicitly sends $c\to d$ to the $\K_0$-class of $d$. 
\end{proof}

Let us finally examine the special case where $\F$ is the $\L$-spectrum functor. 
By Theorem~\reffour{theorem:strong-cofinality}
the equivalence of Proposition~\reffour{proposition:tau-is-equiv} then takes the form
\[
\Sig[\K_0(\D)/\K_0(\C)]^{\tC} \to \Om[\K_0(\D)/\K_0(\C)]^{\tC} .
\]
which on $\pi_{-n}$ is the canonical isomorphism in Tate cohomology
\[
\hat{\mathrm{H}}^{n+1}(\Ct,\K_0(\D)/\K_0(\C)) \xrightarrow{\cong} \hat{\mathrm{H}}^{n-1}(\Ct,\K_0(\D)/\K_0(\C)) .
\]
Indeed, by naturality and the classification of dense Poincaré subcategories, it is enough to note that if $\K_0(\C,\QF)$ is cyclic then all its Tate cohomologies are
either $0$ or $\ZZ/2$, 
and in any case admits no non-identity automorphisms (compare Remark~\reffour{remark:Rothenberg_map}).
We may hence view Proposition~\reffour{proposition:tau-is-equiv} as a generalisation of Ranicki's \cite{Ranicki-splitting}*{Proposition 6.2}, where this statement is proven specifically for quadratic L-theory of group rings.

\begin{example}
\label{example:concrete-difference}%
Suppose that $R$ is a number ring in which 2 is invertible and such that $\Pic(R)/2$ is non-trivial. For example, we may take $R = \ZZ[\sqrt{-pq},\tfrac{1}{2}]$ where $p$ and $q$ are two distinct positive odd primes; the 2-rank of $\Pic(\ZZ[\sqrt{-pq}])$ is then $2$ by Gauss' genus formula (see, e.g., \cite{frohlich-taylor}*{\S V.I, Theorem 39}), and so the 2-rank of $\Pic(R)$ is at least 1. %
As explained in Example~\reffour{example:free-dedekind}, in this case the map
\[
\L_{2k+1}(\Dperf(R),\QF^{\sym}_R) \to \pi_{2k+1}[\K_0(\Dperf(R))/\K_0(\Dfree(R))]^{\tC} = \Pic(R)/2
\]
is an isomorphism when $k$ is odd and zero when $k$ is even (since the assumption that $2$ is invertible in $R$ implies that $R$ has no dyadic primes).
It then follows that the composite
\[
\begin{tikzcd}
\L_{-1}(\Dperf(R),\QF^{\sym}_R) \ar[r, ->>] & \Pic(R)/2 \ar[rrr,"\tau_{(\Dperf(R),\Dfree(R),\QF^{\sym})}","\simeq"'] &&& \Pic(R)/2 \ar[r,"\partial", hook] & \L_0(\Dfree(R), \QF^{\sym}_R)
\end{tikzcd}
\]
is non-zero as soon as $\Pic(R)/2$ is non-zero. In particular, in this case the geometric and algebraic decompositions of $\L((\Dfree(R),\QF^{\sym}_R)_{\B\ZZ})$ into a direct sum $\L(\Dfree(R),\QF^{\sym})$ and $\Sig\L(\Dperf(R),\QF^{\sym})$ do not coincide.
\end{example}

\section{Karoubi-localising approximations and Williams' conjecture}
\label{section:karoubification}%
In the present section we employ the circle decomposition theorem to study the decoration sequence, most importantly for $\L$-spectra. We will to this end extend the construction to arbitrary functors and analyse the resulting sequence for Grothendieck-Witt spectra using hermitian cofinality. This will result in one proof of Williams' conjecture from the introduction. The second half of the section then analyses the decoration construction more broadly and show that it, alongside the Lattice-Tate-construction already employed in the construction of $\KGW$ in Proposition \reffour{proposition:localising-are-stable}, in the limit gives rise to Karoubi-localising approximations for arbitrary Verdier-localising functors (and even a few more), without a detailed analysis of the entire decoration sequence, in particular resulting in a second proof of Williams' conjecture.

\subsection{Ultimate decorations and Williams' conjecture}
\label{subsection:universal-dec}%

The purpose of this subsection is to introduce the decoration construction and then use the circle decomposition formula established in \S\reffour{subsection:shaneson} together with hermitian cofinality in order to give a formula for $\KGW$ and $\KL$. In the case of $\KL$, this reproduces what is classically known as $\L$-theory with ultimate decoration, thus showing that the latter coincides with Karoubi $\L$-theory, in particular establishing a Williams' conjecture.

\begin{definition}
\label{definition:decoration}%
Let $\F \colon \Catp \to \A$ be a functor valued in a stable $\infty$-category $\A$. We define $\T\F\colon \Catp \to \A$ by the formula
\[
\T\F(\C,\QF) = \Sig\TotFib\left[
\begin{tikzcd}
\F(\Hyp(\C)) \ar[r]\ar[d] & \F(\Hyp(\C_{\B\NN})) \ar[d] \\
\F(\C,\QF) \ar[r] & \F((\C,\QF)_{\B\ZZ})
\end{tikzcd}\right] \ ,
\]
where $\TotFib$ denotes the total fibre of the given square. $\T\F$ is therefore the cofibre of the induced map between the horizontal fibres. There is a natural map $\F \Rightarrow \T\F$ given as a composite
\[
\F(\C,\QF) = \fib[\F(\C,\QF) \to \B\ZZ \otimes \F(\C,\QF)] \Rightarrow \fib[\F(\C,\QF) \Rightarrow \F((\C,\QF)_{\B\ZZ})] \Rightarrow \T\F(\C,\QF).
\]
For $n \geq 0$ we then define $\F^{\cotower{-n}}\colon \Catp \to \A$ inductively as follows. We set $\F^{\cotower{0}} = \T\F$, and for $n \geq 1$ we set $\F^{\cotower{-n}} = \T\F^{\cotower{-n+1}}$. We then set
\[
\F^{\cotower{-\infty}} = \colim[\F \Rightarrow \F^{\cotower{0}} \Rightarrow \F^{\cotower{-1}} \Rightarrow ... ].
\]
We will refer to this as the \defi{decoration sequence} of $\F$, and to $\F^{\cotower{-\infty}}$ as \defi{ultimately decorated} approximation of $\F$. 
\end{definition}

Applied to $\F = \K$ this reproduces the Bass delooping of the $\K$-spectrum functor, so that $\K^{\cotower{-\infty}} = \KK$, see Lemma \reffour{proposition:connective-covers} below.

\begin{lemma}
\label{lemma:preserves-limits}%
Let $\A$ be a stable $\infty$-category. Then the following hold:
\begin{enumerate}
\item
\label{item:T-preserves-limits-and-colimits}%
The functor $\T \colon \Fun(\Catp,\A) \to \Fun(\Catp,\A)$ preserves additive functors as well as all limits and all colimits that exist in $\A$.
\item
\label{item:T-commutes-translation}%
There is an equivalence
\[
m_{(\D,\QFD)}^*\T(\F) \simeq \T(m_{(\D,\QFD)}^*\F)
\]
natural in $(\D,\QFD) \in \Catp$ and $\F\in \Fun(\Catp,\A)$, where $m_{(\D,\QFD)}^*$ is the endo-functor of $\Fun(\Catp,\A)$ given by $m_{(\D,\QFD)}^*\F \simeq \F((\D,\QFD) \otimes -)$. \item
\label{item:TFhyp-equivalence}%
There is a natural equivalence $(\T\F)^{\hyp} \simeq \T(\F^{\hyp})$ in $\Fun(\Catp,\A^{\B\Ct})$. 
\end{enumerate}
\end{lemma}
\begin{proof}
The first half of~\reffouritem{item:T-preserves-limits-and-colimits} follows directly from Proposition~\reftwo{lemma:Morita-exact-criterion}. The other statement follow from the following analysis: 
Using the projection formula from Corollary~\refone{corollary:hyp-projection} and the compatibility of tensors and the monoidal product from Remark~\refone{remark:comparing-2}, we may equivalently write the functor $\T\colon \Fun(\Catp,\A) \to \Fun(\Catp,\A)$ as
\[
\T\F(-) = \Sig\TotFib\left[
\begin{tikzcd}
\F((-) \otimes \Hyp(\Spaf)) \ar[r]\ar[d] & \F((-) \otimes \Hyp(\Spaf_{\B\NN})) \ar[d] \\
\F((-) \otimes (\Spaf,\QF^{\uni})) \ar[r] & \F((-) \otimes (\Spaf,\QF^{\uni})_{\B\ZZ})
\end{tikzcd}\right] \ .
\]
More precisely, $\T$ is the suspension of the total fibre of the square $\tau\colon \Del^1 \times \Del^1 \to \Fun(\Catp,\A)$ which is the image of the square 
\[
\begin{tikzcd}
\Hyp(\Spaf) \ar[r]\ar[d] & \Hyp(\Spaf_{\B\NN}) \ar[d] \\
(\Spaf,\QF^{\uni}) \ar[r] & (\Spaf,\QF^{\uni})_{\B\ZZ})
\end{tikzcd}
\]
under the map
\[
\Catp \to \Fun(\Fun(\Catp,\A),\Fun(\Catp,\A)) \quad\quad (\D,\QFD) \mapsto m_{(\D,\QFD)}^*
\]
Since each of the functors $m_{(\D,\QFD)}^*$ preserves limits and colimits and $\A$ is stable we conclude that $\T$ preserves limits and colimits, showing the second half of~\reffouritem{item:T-preserves-limits-and-colimits}. The equivalence $m_{(\Dtwo,\QFDtwo)} \circ m_{(\D,\QFD)} = m_{(\D,\QFD)}\circ m_{(\Dtwo,\QFDtwo)}$ determined by the symmetric monoidal structure of $\Catp$ then induces the equivalence of~\reffouritem{item:T-commutes-translation}. Finally, \reffouritem{item:TFhyp-equivalence} is essentially a special case of \reffouritem{item:T-commutes-translation} for $(\D,\QFD) = \Hyp(\Spaf)$, using again the projection formula from Corollary~\refone{corollary:hyp-projection}), only that the result is interpreted as taking values in $\A^{\B\Ct}=\Fun(\B\Ct,\A)$ via the $\Ct$-action on $\Hyp(\Spaf)$ as a Poincaré $\infty$-category. 
\end{proof}

\begin{proposition}
\label{proposition:T-of-localising}%
Let $\F\colon \Catp \to \A$ be a Verdier-localising functor valued in a stable $\infty$-category $\A$. Then $\T\F$ is naturally equivalent, as a functor under $\F$, to the functor $(\C,\QF) \mapsto \F(\C^\natural,\QF)$. In particular, if $\F$ is Karoubi-localising then the map $\F \Rightarrow \T\F$ is an equivalence.
\end{proposition}
\begin{proof}
Consider the commutative diagram
\[
\begin{tikzcd}
\F\Hyp(\C)\ar[r,equal]\ar[d]  & \F\Hyp(\C) \ar[r]\ar[d] & \F(\Hyp(\C_{\B\NN}))\ar[d] \\
\F(\C,\QF) \ar[r] & {\B\ZZ} \otimes^{\flat} \F(\C,\QF) \ar[r, "\ass^{\flat}"] & \F((\C,\QF)_{\B\ZZ}) \ .
\end{tikzcd}
\]
By Theorem~\reffour{theorem:circle-decomposition} and Remark~\reffour{remark:circle-decompose-2} the right square is exact, and hence the total cofibre of the external rectangle is naturally equivalent to the total cofibre of the left square, which in turn is given by $\cof[\F(\C,\QF) \to {\B\ZZ} \otimes^{\flat} \F(\C,\QF)] = \Sig\F(\C^\natural,\QF)$. %
We conclude that $\T\F(\C,\QF) = \F(\C^{\natural},\QF)$, as desired.
\end{proof}

Suppose that $\F\colon \Catp \to \A$ is a functor and $\F \Rightarrow \FF$ a natural transformation such that $\FF$ is Karoubi-localising. Then in the commutative square
\[
\begin{tikzcd}
\F \ar[r]\ar[d] & \FF \ar[d, "\simeq"] \\
\F^{\cotower{-n}} \ar[r]\ar[ur,dotted] & \FF^{\cotower{-n}}
\end{tikzcd}
\]
the right vertical map is an equivalence by Proposition~\reffour{proposition:T-of-localising}, and hence the dotted lift exists and is essentially unique. In particular, any map $\F \Rightarrow \FF$ with a Karoubi-localising target yields an augmentation 
\[
\F \Rightarrow \F^{\cotower{0}} \Rightarrow \F^{\cotower{-1}} \Rightarrow \cdots \Rightarrow \F^{\cotower{-n}} \Rightarrow \cdots \Rightarrow \FF
\]
of the decoration sequence of $\F$, and so induces a map
\[
\F^{\cotower{-\infty}} \Rightarrow \FF .
\]

\begin{proposition}
\label{proposition:connective-covers}%
Let $\FF\colon \Catp \to \Spa$ be a Karoubi-localising functor and let $\F \Rightarrow \FF$ be a map which exhibits $\F$ as the pointwise $m$-connective cover of $\FF$. Then the induced map $\T\F \Rightarrow \T\FF$ exhibits $\T\F$ as the pointwise $(m-1)$-connective cover of $\T\FF \simeq \FF$. 
As a consequence, for $\F = \tau_{\geq m}\FF$ we have $\F^{\cotower{-n}} = \tau_{\geq m-n-1}\F$ and the associated map
$\F^{\cotower{-\infty}} \Rightarrow \FF$ is an equivalence.
\end{proposition}
\begin{proof}
For $(\C,\QF) \in \Catp$ consider the commutative diagram
\[
\begin{tikzcd}
\FF\Hyp(\C)\ar[r,equal]\ar[d]  & \FF\Hyp(\C) \ar[r]\ar[d] & \FF(\Hyp(\C_{\B\NN}))\ar[d] \\
\FF(\C,\QF) \ar[r] & {\B\ZZ} \otimes \FF(\C,\QF) \ar[r, "\ass"] & \FF((\C,\QF)_{\B\ZZ}) \ .
\end{tikzcd}
\]
Since $\FF$ is Karoubi-localising, Theorem~\reffour{theorem:circle-decomposition} implies that the right square is exact. Furthermore, the horizontal arrows in that square admit retractions, and are hence injective on homotopy groups. We conclude that the right square remains exact after applying $\tau_{\geq m}$. The total fibre of the external rectangle then coincides with that of the left square, so the map $\T\F = \T(\tau_{\geq m}\FF) \Rightarrow \T\FF$ can be identified with the map induced on horizontal fibres in the square
\[
\begin{tikzcd}[column sep=5pt]
\tau_{\geq m}\FF(\C,\QF) \ar[rr]\ar[d] && \tau_{\geq m}[{\B\ZZ} \otimes \FF(\C,\QF)] \ar[r,phantom, "="]\ar[d] & \tau_{\geq m}\FF(\C,\QF) \oplus \tau_{\geq m}(\Sig\FF(\C,\QF))\ar[d] \\
\FF(\C,\QF) \ar[rr] && \B\ZZ \otimes \FF(\C,\QF) \ar[r,phantom, "="] & \FF(\C,\QF) \oplus \Sig\FF(\C,\QF) \\
\end{tikzcd}
\]
which is the map $\Om\tau_{\geq m}(\Sig\FF) \Rightarrow \Om\Sig\FF=\FF$. This map is indeed an $(m-1)$-connective cover, as desired.
\end{proof}

\begin{corollary}
\label{corollary:decorated-K}%
For $n \geq 0$ the map $\K^{\cotower{-n}} \Rightarrow \KK$ exhibits $\K^{\cotower{-n}}$ as the $(-n)$-connective cover $\tau_{\geq -n}\KK$ of $\KK$. 
In particular, the induced map 
$\K^{\cotower{-\infty}} \Rightarrow \KK$ is an equivalence.
\end{corollary}
\begin{proof}
We have $\K^{\cotower{0}} = \K((-)^{\natural}) = \tau_{\geq 0}\KK$ by Proposition~\reffour{proposition:T-of-localising} and the definition of $\KK$. The desired result hence follows from Proposition~\reffour{proposition:connective-covers}.
\end{proof}

Let us now bootstrap the above result from $\K$- to Grothendieck-Witt-spectra. 

\begin{proposition}
\label{proposition:decorated-GW}%
For all $n\geq 0$, we have that the square 
\[
\begin{tikzcd} \GW^{\cotower{-n}} \ar[r] \ar[d] & \KGW \ar[d] \\
(\tau_{\geq -n} \KK)^\hC \ar[r] & \KK^\hC
\end{tikzcd}
\]
is exact. In particular, the map $\GW^{\cotower{-\infty}} \Rightarrow \KGW$ is an equivalence.
\end{proposition}
\begin{proof}
Consider the square
\[
\begin{tikzcd}
\GW \ar[r]\ar[d] & \KGW\ar[d] \\
\K^{\hC} \ar[r] & \KK^{\hC} 
\end{tikzcd} \ ,
\]
which is exact by Theorem~\reffour{theorem:Karoubi-pullback-squares}. Applying the functor $(-)^{\cotower{-n}}$ using Proposition~\reffour{proposition:T-of-localising} and Lemma~\reffour{lemma:preserves-limits}\;\reffouritem{item:T-preserves-limits-and-colimits} we obtain the exact square
\[
\begin{tikzcd}
\GW^{\cotower{-n}} \ar[r]\ar[d] & \KGW\ar[d] \\
(\K^{\hC})^{\cotower{-n}} \ar[r] & \KK^{\hC} \ .
\end{tikzcd} 
\]
Corollary~\reffour{corollary:decorated-K} together with the fact that the operation $\T(-)$ preserves all limits by Lemma~\reffour{lemma:preserves-limits}~\reffouritem{item:T-preserves-limits-and-colimits} identifies this with the square from the statement. Finally, taking the colimit over $n \in \mathbb N$, the bottom horizontal map becomes the map
\[
\colim_{n \geq 0}(\KK^{\hC})^{\cotower{-n}} = \colim_{n \geq 0}(\tau_{\geq -n} \KK)^{\hC} \longrightarrow \KK^{\hC},
\]
which is an equivalence by \cite{NS}*{Lemma I.2.6 (ii)} and hence the top horizontal arrow is an equivalence as well.
\end{proof}

\begin{proposition}
\label{proposition:decorated-L}%
For all $n\geq 0$, we have that the square 
\[
\begin{tikzcd} \L^{\cotower{-n}} \ar[r] \ar[d] & \KL \ar[d] \\
(\tau_{\geq -n} \KK)^\tC \ar[r] & \KK^\tC
\end{tikzcd}
\]
is exact. In particular, the induced map $\L^{\cotower{-\infty}} \Rightarrow \LL$ is an equivalence.
\end{proposition}
\begin{proof}
The proof proceeds as that of Proposition~\reffour{proposition:decorated-GW} by applying the functor $(-)^{\cotower{-n}}$, this time to the exact square
\[
\begin{tikzcd}
\L \ar[r]\ar[d] & \KL\ar[d] \\
\K^{\tC} \ar[r] & \KK^{\tC} 
\end{tikzcd} \ ,
\]
equally established in Theorem~\reffour{theorem:Karoubi-pullback-squares}. The resulting square can then identified with the square in the proposition by using Proposition~\reffour{proposition:T-of-localising}, Proposition~\reffour{corollary:decorated-K} and the fact that 
$(-)^{\cotower{-n}}$ 
commutes with the Tate construction by 
Lemma~\reffour{lemma:preserves-limits}\;\reffouritem{item:T-preserves-limits-and-colimits}. The last part follows again from \cite{NS}*{Lemma I.2.6 (ii)}.
\end{proof}

\begin{remark}
\label{remark:all-verdier}%
Corollary~\reffour{corollary:decorated-K} and Propositions \reffour{proposition:decorated-GW} and \reffour{proposition:decorated-L} imply that the natural transformation $\F \Rightarrow \F^{\cotower{-\infty}}$ exhibits $\F^{\cotower{-\infty}}$ as the initial Karoubi-localising functor under $\F$ for $\F=\K,\GW$ and $\L$, respectively. In Theorem~\reffour{theorem:karoubification}
below 
we will show that this in fact holds for any Verdier-localising functor $\F\colon \Catp \to \Spa$.
\end{remark}

\begin{remark}
\label{rem:spectral_rothenberg_lower_decoration}%
Propositions \reffour{proposition:decorated-GW} and \reffour{proposition:decorated-L}
give rise to exact squares
\[
\begin{tikzcd}
	\GW^{\cotower{-n+1}} \ar[r] \ar[d] & \GW^{\cotower{-n}} \ar[d] && \L^{\cotower{-n+1}} \ar[r] \ar[d] & \L^{\cotower{-n}} \ar[d] \\
	(\tau_{\geq -n+1}\KK)^{\hC} \ar[r] & (\tau_{\geq -n}\KK)^{\hC} && (\tau_{\geq -n+1}\KK)^{\tC} \ar[r] & (\tau_{\geq -n}\KK)^{\tC}
\end{tikzcd}
\]
and consequently, one obtains fibre sequences
\[
\begin{tikzcd}[row sep=tiny]
	\GW^{\cotower{-n+1}}(-) \ar[r] & \GW^{\cotower{-n}}(-) \ar[r] & \mathbb{S}^{-n} \otimes \KK_{-n}(-)^{\hC}\,, \\
	\L^{\cotower{-n+1}}(-) \ar[r] & \L^{\cotower{-n}}(-) \ar[r] & \mathbb{S}^{-n} \otimes  \KK_{-n}(-)^{\tC}
\end{tikzcd}
\]
the second of which on homotopy groups induces Ranicki's Rothenberg sequences for lower L-theory \cite{Ranicki-lower}*{Theorem 17.2}. 
\end{remark}

\subsection{Karoubi approximation via Day cotensoring}
\label{subsection:karoubi-approximation-verdier}%

In this subsection we turn to the behaviour of the ultimate decoration construction for general Verdier-localising functors, and aim to show that it always provides an initial Karoubi-localising approximation. The same in fact goes for the Lattice-Tate construction used in the proof of Proposition~\reffour{proposition:localising-are-stable}, and our approach is to consider both of them as instances of Day cotensoring (the right adjoint to Day convolution), and then formulate sufficient conditions under which iterative Day cotensoring leads to Karoubi approximation. Let us begin by setting up the stage.

In what follows, recall that a functor $\C \to \D$ among locally small $\infty$-categories is said to be small if it is left Kan extended from a small full subcategory $\C_0 \subseteq \C$. Given a locally small $\infty$-category $\C$, its Ind-completion $\Ind(\C)$ is by definition the smallest full subcategory of $\Fun(\C\op,\Sps)$ containing the representable presheaves and closed under (small) filtered colimits. In particular, $\Ind(\C)$ is again locally small. If $\C$ admits finite colimits then we can also identify $\Ind(\C)$ with the $\infty$-category of finite-limit-preserving small presheaves $\C\op \to \Sps$.

\begin{construction}
\label{const:day-conv}%
Let $\cB$ be a locally small and idempotent complete symmetric monoidal $\infty$-category acting on an idempotent complete locally small $\infty$-category $\E$, and assume that $\cB$ and $\E$ admit finite colimits compatibly with the tensor product and action. The $\infty$-categories $\Ind(\cB)$ and $\Ind(\E)$ are then cocomplete and locally small, and $\Ind(\cB)$ inherits a symmetric monoidal structure and an action on $\Ind(\E)$, both compatible with small colimits. 

By Corollary~\reffour{corollary:persist-locally-small}, the $\infty$-category $\Fun^{\sm}(\Catp,\Ind(\cB))$ of small functors $\Catp \to \Ind(\cB)$ admits a symmetric monoidal structure by Day convolution, and acts via Day convolution on $\Fun(\Catp,\Ind(\E))$, where in the latter we do not have to restrict attention to small functors. Here we recall that the Day convolution of a small functor $\T\colon \Catp \to \Ind(\cB)$ and a functor $\F$ from $\Catp$ to either $\Ind(\cB)$ or $\Ind(\E)$ is given by the left Kan extension along $\otimes\colon \Catp \times \Catp \to \Catp$ of the functor $((\C,\QF),(\D,\QFD)) \mapsto \T(\C,\QF)\otimes \F(\C,\QF)$, and the condition that $\T$ is small and $\Catp$ is closed ensures that this left Kan extension is well-defined.

Let us write $\Fun^{\omega}(\Catp,\Ind(\cB)) \subseteq \Fun^{\sm}(\Catp,\Ind(\cB))$ for the full subcategory spanned by the compact objects in $\Fun^{\sm}(\Catp,\Ind(\cB))$. We claim that a small $\T\colon \Catp \to \Ind(\cB)$ is compact if and only if it is a retract of a functor of the form
\[
\colim_{i \in \I}\Gamma(i)_!\tau(i) 
\]
for some finite $\infty$-category $\I$, some functor $\Gamma\colon \I\op \to \Catp$, and some functor $\tau\colon \I \to \cB$, where $\Gamma(i)_!$ means left Kan extension along the functor $\ast \to \Catp$ picking out $\Gamma(i)$. 
Indeed, on the one hand, any such functor is compact, being a retract of a finite colimit of compacts. On the other hand, if $\T$ is a compact small functor then pick a %
small full subcategory $\mathcal{U} \subseteq \Catp$ from which $\T$ left Kan extends, so that we can write $\T$ as a small colimit of the form
\[
\T = \colim_{\stackrel{[(\C,\QF) \to (\D,\QFD)]}{\in \Twar(\mathcal{U})}}(\D,\QFD)_!\T(\C,\QF) .
\]
By~\cite{HTT}*{Variant 4.2.3.15} there exists a poset $\J$ and a cofinal map $\J \to \Twar(\mathcal{U})$, so that the above colimit can be computed after restricting to $\J$. %
Writing $\J$ as a colimit of finite subposets $\I \subseteq \J$ and using the compactness of $\T$ we obtain that
$\T$ is a retract of a functor of the form 
\[
\colim_{i \in \I}\Gamma(i)_!\bar{\tau}(i)
\]
for some finite poset $\I$, some functor $\Gamma\colon \I\op \to \Catp$, and some functor $\bar{\tau}\colon \I \to \Ind(\cB)$. Finally, note that the above expression is colimit preserving in $\bar{\tau}$; writing the latter as a filtered colimit of diagrams which take values in $\cB$ using~\cite{HTT}*{Proposition 5.3.5.15} %
and invoking again the compactness of $\T$ we may assume without loss of generality that $\bar{\tau}$ takes values in $\cB$, and so we get the claim.

Now since $(\Spaf,\QF^{\uni})_!\one_{\cB}$ is compact and $(\C,\QF)_!b \otimes_{\Day} (\C',\QF')_!b' = ((\C,\QF)\otimes(\C',\QF'))_!(b\otimes b')$ we conclude that the collection of compact functors contains the Day unit and is closed under Day convolution, therefore inheriting from $\Fun^{\sm}(\Catp,\Ind(\cB))$ a symmetric monoidal structure and an action on $\Fun(\Catp,\Ind(\E))$. In addition, since $(\C,\QF)_!b \otimes_{\Day} \F = b\otimes \F(\Funx((\C,\QF),-))$ and $\E$ is idempotent complete we see that the Day action of compact functors preserves $\Fun(\Catp,\E) \subseteq \Fun(\Catp,\Ind(\E))$. All in all, this construction results in a symmetric monoidal $\infty$-category $\Fun^{\omega}(\Catp,\Ind(\cB))$ acting on the $\infty$-category $\Fun(\Catp,\E)$.
\end{construction}

In what follows, we would like to construct Karoubi approximations to Verdier-localising functors $\Catp \to \E$, where $\E$ is as in Construction~\reffour{const:day-conv} for some $\cB$. We will hence fix the assumption that $\E$ has finite limits. We recall that Verdier-localising functors are by definition reduced, that is, send the zero Poincaré $\infty$-category to the terminal object $\ast \in \E$, and hence uniquely lift to functors $\Catp \to \E_{\ast/}$. Replacing $\E$ with $\E_{\ast/}$, we may as well assume that $\E$ is pointed, in which case we may as well assume that $\cB$ is pointed as well. 
We then write $\Fun_{\ast}(\Catp,\E) \subseteq \Fun(\Catp,\E)$ for the full subcategory spanned by the reduced functors.
This full subcategory is reflective, with left adjoint given by $\F^{\red} = \cof[(0,0)_!\F(0) \to \F(-)]$. Similarly, let $\Fun_{\ast}^{\omega}(\Catp,\Ind(\cB)) \subseteq \Fun^{\omega}(\Catp,\Ind(\cB))$ be the full subcategory spanned by the compact reduced functors. 
Given a functor of the form $(\C,\QF)_!b$ as above, its reduction is a finite colimit of such, and in particular compact. We conclude that the reduction preserves compact functors, and so $\Fun_{\ast}^{\omega}(\Catp,\Ind(\cB))$ is again a reflective subcategory of $\Fun^{\omega}(\Catp,\Ind(\cB))$.
Using that $(0,0) \otimes (\D,\QFD) = (0,0)$ for every $(\D,\QF)$ we then get that this reflective localisation is multiplicative with respect to Day convolution (see Proposition~\reffour{proposition:multiplicative-day}), and so we have an induced symmetric monoidal structure on $\Fun_{\ast}^{\omega}(\Catp,\Ind(\cB))$ and an induced action of $\Fun_{\ast}^{\omega}(\Catp,\Ind(\cB))$ on $\Fun_{\ast}^{\omega}(\Catp,\E)$. 
The unit $\one^{\red}$ of $\Fun_{\ast}^{\omega}(\Catp,\Ind(\cB))$ is given by $((\Spaf,\QF^{\uni})_!\one_{\cB})^{\red} = \Poinc \tilde{\otimes} \one_{\cB}$, where $\Poinc$ is viewed as a functor to pointed spaces and $\tilde{\otimes}$ stands for the action of pointed spaces on $\Ind(\cB)$ afforded by the fact that it is pointed with finite colimits. More generally, $((\C,\QF)_!b)^{\red} = \Map((\C,\QF),-)\tilde{\otimes} b$, where $\Map((\C,\QF),-)$ is viewed as pointed by the zero Poincaré functor.

Given the above, we adopt the following assumption:

\begin{assumption}
\label{assumption:cotensor}%
Let $\cB$ and $\E$ be as in Construction~\reffour{const:day-conv} and assume in addition
\begin{enumerate}
\item
$\E$ and $\cB$ are pointed.
\item
$\E$ admits finite limits.
\item
For every $b \in \cB$ the action $b\otimes(-)\colon \E \to \E$ admits a right adjoint $(-)^b\colon \E \to \E$.
\end{enumerate}
\end{assumption}

Assumption~\reffour{assumption:cotensor} ensure that the Day convolution action of $\Fun^{\omega}_{\ast}(\Catp,\Ind(\cB))$ on $\Fun_{\ast}(\Catp,\E)$ is closed, that is, for any reduced compact functor $\T\colon \Catp \to \Ind(\cB)$ the reduced Day convolution action $\T\otimes_{\Day}(-)\colon \Fun_{\ast}(\Catp,\E) \to \Fun_{\ast}(\Catp,\E)$ admits a right adjoint. 
For example if $\T = \colim_{i \in \I}\Gamma(i)_!\tau(i)$ for some finite $\infty$-category $\I$ and functors $\Gamma\colon \I\op \to \Catp$ and $\tau\colon \I \to \cB$ then 
\[
\G^{\T}(\C,\QF) = \lim_{i \in \I\op}\G(\Gamma(i))^{\bet(i)} \in \E,
\]
and for a general $\T$ we will get a retract of such a formula.

\begin{definition}
Let $\T\colon \Catp \to \Ind(\cB)$ be a compact reduced functor equipped with a reduced augmentation, that is, a map $\eps\colon \T \Rightarrow \one^{\red} =  \Poinc\tilde{\otimes} \one_{\cB}$ to the reduced unit of Day convolution. We will say that $(\T,\eps)$ is \defi{Karoubi preserving} if for any Karoubi-localising functor $\G\colon \Catp \to \E$ the induced map $\eps^*\colon \G \Rightarrow \G^{\T}$ is an equivalence. 
\end{definition}

\begin{example}
\label{example:lattice-tate}%
Take $\cB = \Sps^{\fin}_*$ be the $\infty$-category of finite pointed spaces and $\E$ any pointed idempotent complete $\infty$-category with finite limits and colimits, so that $\E$ is canonically a $\cB$-module and the action of $\cB$ on $\E$ is closed.
Consider the reduced compact functor
\[
\T_{\Tate} := \cof[\Map(\Tate^{\omega_1}(\Spaf, (\QF^{\uni})\qshift{1}),-) \to \Map(\Latt^{\omega_1}(\Spaf, \QF^{\uni}),-)],
\]
where the mapping spaces in $\Catp$ are viewed as pointed via the zero Poincaré functor.
Reduced augmentations $\eps\colon\T_{\Tate} \Rightarrow \Poinc$ then correspond to Poincaré functors 
\[
(\Spaf,\QF^{\uni}) \to \fib\bigl[\Latt^{\omega_1}(\Spaf, \QF^{\uni}) \to \Tate^{\omega_1}(\Spaf, (\QF^{\uni})\qshift{1})\bigr] ,
\]
and we pick the one determined by the Lattice-Tate sequence (as in the proof of Proposition~\reffour{proposition:localising-are-stable}):
\[
(\Spaf,\QF^{\uni}) \to \Latt^{\omega_1}(\Spaf, \QF^{\uni}) \to \Tate^{\omega_1}(\Spaf, (\QF^{\uni})\qshift{1}).
\]
Then for any reduced $\F\colon \Catp \to \E$ we have 
\[
\G^{\T_{\Tate}}(\C,\QF) = \fib\big[\G((\C,\QF)\otimes \Latt^{\omega_1}(\Spaf, \QF^{\uni})) \to \G((\C,\QF) \otimes \Tate^{\omega_1}(\Spaf, (\QF^{\uni})\qshift{1}))\big],
\]
and if in addition $\E$ is additive and $\G$ preserves finite direct sums this simplifies to $\G^{\T_{\Tate}}(\C,\QF) = \Om\G((\C,\QF) \otimes \Tate^{\omega_1}(\Spaf, (\QF^{\uni})\qshift{1}))$. 
Since the above Poincaré-Verdier sequence remains at least a Karoubi-Verdier sequence after tensoring with any Poincaré $\infty$-category $(\C,\QF)$ we conclude that $\T_{\Tate}$ is Karoubi preserving.
\end{example}

\begin{example}
\label{example:decoration}%
Take $\cB = \Spa^{\fin}$ be the $\infty$-category of finite spectra and $\E$ any idempotent complete stable $\infty$-category, so that $\E$ is canonically a $\cB$-module and the action of $\cB$ on $\E$ is closed.

Let $\T_{\dec}\colon \Catp \to \Spa$ be given by the desuspension of the total cofibre of 
\[
\begin{tikzcd}
\Sig^{\infty}\Map((\Spaf,\QF^{\uni})_{\B\ZZ},-) \ar[r]\ar[d] & \Sig^{\infty}\Map(\Hyp(\Spaf_{\B\NN}),-) \ar[d] \\
\Sigma^{\infty}\Map((\Spaf,\QF^{\uni}),-)\ar[r] & \Sigma^{\infty}\Map(\Hyp(\Spaf),-),
\end{tikzcd}
\]
where again all mapping spaces in $\Catp$ are canonically viewed as pointed spaces and $\Sigma^{\infty}$ is the reduced suspension spectrum. Equivalently, $\T_{\dec}$ is the fibre of the induced map on vertical cofibres in the above square.
We then define a reduced augmentation $\eps\colon \T_{\dec} \Rightarrow \Sig^{\infty}\Poinc$ by the composite
\begin{align*}
\T &\Rightarrow \cof\big[\Sig^{\infty}\Map((\Spaf,\QF^{\uni})_{\B\ZZ},-) \Rightarrow \Sig^{\infty}\Map((\Spaf,\QF^{\uni}),-)\big] \\
& \Rightarrow \cof\big[(\Sig^{\infty} \Map_*((\Spaf,\QF^{\uni}),-))^{\B\ZZ} \Rightarrow \Sig^{\infty}\Map((\Spaf,\QF^{\uni}),-)\big] \\
& = \Sig^{\infty}\Map((\Spaf,\QF^{\uni}),-)\\
& = \Sig^{\infty}\Poinc.
\end{align*}
For any reduced $\F\colon \Catp \to \E$ the functor
\[
\G^{\T}(\C,\QF) = \Sig\TotFib\left[
\begin{tikzcd}
\G(\Hyp(\C)) \ar[r]\ar[d] & \G(\Hyp(\C_{\B\NN})) \ar[d] \\
\G(\C,\QF) \ar[r] & \G((\C,\QF)_{\B\ZZ})
\end{tikzcd}\right] = \T\G(\C,\QF)
\]
coincides with the decoration construction of Definition~\reffour{definition:decoration} such that the map $\G \Rightarrow \G^{\T}$ identifies with the map $\G\Rightarrow \T\G$ from loc.\ cit., and is therefore an equivalence whenever $\G$ is Karoubi-localising by Proposition~\reffour{proposition:T-of-localising}. In particular, $\T_{\dec}$ is Karoubi preserving.
\end{example}

\begin{lemma}
\label{lemma:T-Karoubi-examples}%
For $\cB$ and $\E$ be as in Example~\reffour{example:lattice-tate} we have $\G^{\T_{\Tate}} = \G(-^{\natural})$ for every Verdier-localising $\G\colon \Catp \to \E$. Similarly, 
For $\cB$ and $\E$ be as in Example~\reffour{example:decoration} we have $\G^{\T_{\dec}} = \G(-^{\natural})$ for every Verdier-localising $\G\colon \Catp \to \E$.
\end{lemma}
\begin{proof}
The claim for $\T_{\mathrm{dec}}$ holds by Proposition~\reffour{proposition:T-of-localising}.
For $\T_{\Tate}$, it will suffice to show that for any Poincaré $\infty$-category $(\C,\QF)$ the kernel of the Verdier projection 
\[
(\C,\QF) \otimes \Latt^{\omega_1}(\Spaf, \QF^{\uni}) \to (\C,\QF) \otimes \Tate^{\omega_1}(\Spaf, (\QF^{\uni})\qshift{1})
\]
is given by $(\C^{\natural},\QF)$. To that end, consider the induced map of Karoubi sequences
\[
\begin{tikzcd}
(\C,\QF)\ar[r]\ar[d] & (\C,\QF) \otimes \Latt^{\omega_1}(\Spaf, \QF^{\uni}) \ar[r]\ar[d] & (\C,\QF) \otimes \Tate^{\omega_1}(\Spaf, (\QF^{\uni})\qshift{1})\ar[d] \\
(\C^{\natural},\QF)\ar[r] & (\C^{\natural},\QF) \otimes \Latt^{\omega_1}(\Spaf, \QF^{\uni}) \ar[r] & (\C^{\natural},\QF) \otimes \Tate^{\omega_1}(\Spaf, (\QF^{\uni})\qshift{1}) 
\end{tikzcd}
\]
in which the bottom one is a fibre sequence in $\Catp$ since its left most term is idempotent complete. Since tensoring with a fixed Poincaré $\infty$-category preserves Karoubi equivalences (see Proposition~\reftwo{lemma:Morita-exact-criterion}\;\reftwoitem{item:preserves-K-equivalences}) 
the vertical maps are Karoubi equivalences, and in particular fully-faithful. The induced map on horizontal kernels in the right square is then equivalent to the base changed Karoubi equivalence
\[
(\C^{\natural},\QF) \times_{(\C^{\natural},\QF) \otimes \Latt^{\omega_1}(\Spaf, \QF^{\uni})}[(\C,\QF) \otimes \Latt^{\omega_1}(\Spaf, \QF^{\uni})] \hookrightarrow (\C^{\natural},\QF).
\]
What we need to show is that this map is actually an equivalence in $\Catp$. Now recall that $\Latt^{\omega_1}(\Spaf,\QF^{\uni})$ is defined as the pairings $\infty$-category 
\[
\Latt^{\omega_1}(\Spaf,\QF^{\uni}) = \Pairings(\Ind^{\omega_1}(\Spaf),\hat{\QF}^{\uni}) ,
\]
and so by Lemma~\refone{lemma:pairing-projection} we have a natural equivalence
\[
(\C,\QF) \otimes \Latt^{\omega_1}(\Spaf,\QF^{\uni}) \simeq \Pairings((\C,\QF) \otimes (\Ind^{\omega_1}(\Spaf),\hat{\QF}^{\uni})) .
\]
Being in particular a right adjoint, $\Pairings$ preserves fibre products. Consequently, we obtain a cartesian square
\[
\begin{tikzcd}
(\C^{\natural},\QF)\times_{\Pairings((\C^{\natural},\QF) \otimes (\Ind^{\omega_1}(\Spaf),\hat{\QF}^{\uni}))} [\Pairings((\C,\QF) \otimes (\Ind^{\omega_1}(\Spaf),\hat{\QF}^{\uni}))] \ar[d]\ar[r,hookrightarrow] & (\C^{\natural},\QF) \ar[d] \\
\Pairings\big((\C^{\natural},\QF) \times_{(\C^{\natural},\QF) \otimes (\Ind^{\omega_1}(\Spaf),\hat{\QF}^\uni)} [(\C,\QF) \otimes (\Ind^{\omega_1}(\Spaf),\hat{\QF}^\uni)] \big) \ar[r,hookrightarrow] & \Pairings(\C^{\natural},\QF) 
\end{tikzcd}
\]
in which the horizontal maps are Karoubi equivalences. To show that the top horizontal map is an equivalence it will hence suffice to show that the Karoubi equivalence of hermitian $\infty$-categories 
\[
(\C^{\natural},\QF) \times_{(\C^{\natural},\QF) \otimes (\Ind^{\omega_1}(\Spaf),\hat{\QF}^\uni)} [(\C,\QF) \otimes (\Ind^{\omega_1}(\Spaf),\hat{\QF}^\uni)] \hookrightarrow (\C^{\natural},\QF)
\]
is an equivalence. But this is detected on the level of underlying stable $\infty$-categories since in a Karoubi equivalence the hermitian structure on the domain is restricted from that of the target.  Hence, we are reduced to showing that the dense full subcategory inclusion
\[
\C^{\natural}  \times_{\C^{\natural} \otimes \Ind^{\omega_1}(\Spaf)} [\C \otimes \Ind^{\omega_1}(\Spaf)] \hookrightarrow \C^{\natural}
\]
is essentially surjective.

Let $\NN$ be the set of non-negative integers, $\SS[\NN] \in \Ind^{\omega_1}(\Spaf)$ the $\NN$-indexed direct sum of copies of $\SS$, and $T\colon \SS[\NN] \to \SS[\NN]$ the map induced by the map $(-)+1\colon \NN \to \NN$.
Suppose given a retract diagram $d \to c \to d$ in $\C^{\natural}$ with $c \in \C$ and let $\alp\colon c \to d \to c$ be the corresponding idempotent endomorphism of $c$. We need to show that $d \otimes \SS \in \C^{\natural} \otimes \Ind^{\omega_1}(\Spaf)$ is contained in $\C \otimes \Ind^{\omega_1}(\Spaf)$.
We claim that $d \otimes \SS$ is in fact equivalent to the coequaliser 
\[
e := \coeq\big[\begin{tikzcd}c \otimes \SS[\NN] \ar[r,yshift=-2pt,"\id"']\ar[r,yshift=2pt,"\alp\otimes T"] & c \otimes \SS[\NN] \end{tikzcd}\big] \in \C \otimes \Ind^{\omega_1}(\Spaf) .
\]
To see this, let us decompose $c$ as a direct sum $c = d \oplus d'$ for some $d'\in \C^{\natural}$. Then $\alp\colon c \to c$ decomposes as a direct sum of the identity on $d$ and $0$ on $d'$, yielding a decomposition 
\begin{align*}
e &= \coeq\big[\begin{tikzcd}[ampersand replacement=\&]d \otimes \SS[\NN] \ar[r,yshift=-2pt,"\id"']\ar[r,yshift=2pt,"\id\otimes T"] \& d \otimes \SS[\NN] \end{tikzcd}\big] \oplus \coeq\big[\begin{tikzcd}[ampersand replacement=\&]d' \otimes \SS[\NN] \ar[r,yshift=-2pt,"\id"']\ar[r,yshift=2pt,"0\otimes T"] \& d' \otimes \SS[\NN] \end{tikzcd}\big]\\
&= d\otimes \coeq\big[\begin{tikzcd}[ampersand replacement=\&, column sep=18pt]\SS[\NN] \ar[r,yshift=-2pt,"\id"']\ar[r,yshift=2pt,"T"] \& \SS[\NN] \end{tikzcd}\big]  \\
&= d \otimes \SS
\end{align*}
as claimed.
\end{proof}

\begin{construction}
\label{construction:karoubification}%
Let $\cB,\E$ be as in Assumption~\reffour{assumption:cotensor} and assume in addition that $\E$ admits sequential colimits. Let $\T\colon \Catp \to \Ind(\cB)$ be a compact reduced functor equipped with a reduced augmentation $\eps\colon\T \Rightarrow \one^{\red}$. Iterating the natural transformation $\G \Rightarrow \G^{\T}$ we obtain a sequence
\[
\G \Rightarrow \G^{\T} \Rightarrow \G^{\T\otimes \T} \Rightarrow ... \Rightarrow \G^{\T^{\otimes n}} \Rightarrow ...
\]
and we write $\G^{\T^{\infty}} = \colim_n\G^{\T^{\otimes n}}$ for its colimit.
We note that, as usual with these constructions, there are a-priori $n+1$ different maps $\G^{\T^{\otimes n}} \Rightarrow \G^{\T^{\otimes n+1}}$ one could consider here, indexed by $i \in \{0,...,n\}$, where the $i$'th one is induced by identifying the source with $\G^{\T^{\otimes i}\otimes \one^{\red} \otimes \T^{\otimes n-i}}$ and using the natural transformation $\eps$ on the component $\one^{\red}$. 
Nonetheless, the automorphisms of $\G^{\T^{\otimes n+1}}$ induced by the permutation automorphisms of $\T^{\otimes n+1}$ act transitively on the set of choices, and so the colimit $\G^{\T^{\infty}}$ (as an object of $\Fun(\Catp,\E)_{\G/}$) is in fact independent on the specific choices we make.  
\end{construction}

For instance, we have $\G^{\T_{\mathrm{dec}}^{\infty}} = \G^{\cotower{-\infty}}$ as functors under $\G$. Let us also write $\G^{\Tate} := \G^{\T_{\Tate}^{\infty}}$. We now come to the main result of this subsection:

\begin{theorem}
\label{theorem:karoubification}%
Let $\cB$ and $\E$ be as in Assumption~\reffour{assumption:cotensor} and assume in addition that $\E$ is differentiable (that is, admits finite limits and sequential colimits, and these two commute) and that the cotensoring $(-)^b\colon \E \to \E$ preserves sequential colimits for any $b \in \cB$.
Let $\T\colon \Catp \to \Ind(\cB)$ be a compact reduced functor and $\F\colon \Catp \to \E$ a Verdier-localising functor. If $\T$ is Karoubi preserving and $\F^{\T}$ is Karoubi invariant then $\F^{\T^{\infty}}$ is Karoubi-localising and the map 
\[
\F \Rightarrow \F^{\T^{\infty}}
\]
exhibits its target as initial in Karoubi-localising functors under $\F$.
\end{theorem}

Combining the above theorem with Lemma~\reffour{lemma:T-Karoubi-examples} we conclude:

\begin{corollary}
\label{corollary:tate-dec}%
If $\E$ is a pointed differentiable $\infty$-category with finite colimits and $\F\colon \Catp \to \E$ is a Verdier-localising functor then $\F \Rightarrow \F^{\Tate}$ is the initial Karoubi-localising functor under $\F$. If $\E$ is furthermore stable then the same holds for $\F \Rightarrow \F^{\cotower{-\infty}}$.
\end{corollary}

\begin{corollary}
\label{corollary:karoubification}%
If $\E$ is a pointed differentiable $\infty$-category with finite colimits then the inclusion 
\[
\Funkloc(\Catp,\E) \to \Funvloc(\Catp,\E)
\]
admits a left adjoint $(-)^{\mathrm{kar}}$. In addition, if $R\colon \E \to \E'$ is any functor between pointed differentiable $\infty$-categories with finite colimits such that $R$ commutes with finite limits and sequential colimits then the square
\[
\begin{tikzcd}
\Funkloc(\Catp,\E)\ar[r]\ar[d,"R\circ (-)"] & \Funvloc(\Catp,\E) \ar[d,"R\circ (-)"] \\
\Funkloc(\Catp,\E')\ar[r] & \Funvloc(\Catp,\E')
\end{tikzcd}
\]
is horizontally left adjointable, that is, the canonical map $(R\F)^{\mathrm{kar}} \Rightarrow R(\F^{\mathrm{kar}})$ is an equivalence for all Verdier-localising functors $\F \colon \Catp \to \E$.
\end{corollary}
\begin{proof}
We set $\cB = \Sps^{\fin}_*$ and form Karoubi-localising approximations via $(-)^{\Tate}$ (using Corollary~\reffour{corollary:tate-dec}) to obtain the horizontal left adjoints.
To show the adjointability property we need to show that if $\F\colon \Catp \to \E$ is a Verdier-localising functor then the map $R\circ \F \Rightarrow R\circ \F^{\Tate}$ induced by $\F \Rightarrow \F^{\Tate}$ is initial in Karoubi-localising functors $\Catp \to \E'$ under $R\circ \F$. Note that we do not assume $R$ preserves finite colimits. In particular, $R$ is not a $\cB$-linear, but only lax $\cB$-linear, and similarly the post-composition functor $R\circ (-)\colon \Fun_{\ast}(\Catp,\E) \to \Fun_{\ast}(\Catp,\E)$ is lax $\Fun^{\omega}_{\ast}(\Catp,\Sps_*)$-linear with respect to Day convolution. This means that on the level of Day cotensoring we have an oplax comparison map
\[
R \circ \F^{\T_{\Tate}} \Rightarrow (R\circ \F)^{\T_{\Tate}}
\]
for $\F \in \Fun_{\ast}(\Catp,\E)$.
The condition that $R$ commutes with finite limits implies that this comparison map is an equivalence. The additional assumption that $R$ commutes with sequential colimits now implies that the left facing map in the induced zig-zag
\[
R \circ \F^{\Tate} = R \circ \bigl[\colim_n\F^{\T_{\Tate}^{n}}\bigr] \Leftarrow \colim_n \bigl[R \circ \F^{\T_{\Tate}^{n}}\bigr] \stackrel{\sim}{\Rightarrow} \colim_n[R\circ \F]^{\T_{\Tate}^{n}} = (R \circ \F)^{\Tate}
\]
is an equivalence as well. By Corollary~\reffour{corollary:tate-dec} we conclude that the induced map
\[
R \circ \F \Rightarrow R \circ \F^{\Tate}
\]
exhibits its target as the initial Karoubi-localising functor under the source, as desired.
\end{proof}

\begin{remark}
The proof of Theorem~\reffour{theorem:karoubification} given below works verbatim if we replace $\Catp$ by $\Catx$, and is applicable to the non-hermitian analogues of Examples~\reffour{example:lattice-tate} and~\reffour{example:decoration}. In the former case, one can either work with the underlying stable $\infty$-categories $\Latt^{\omega_1}(\Spaf)$ and $\Tate^{\omega_1}(\Spaf)$, or with their one-sided counterparts $\Ind^{\omega_1}(\Spaf)$ and $\mathrm{Calc}^{\omega_1}(\Spaf) = \Ind^{\omega_1}(\Spaf)/\Spaf$. In all cases the assumptions of Theorem~\reffour{theorem:karoubification} hold for any Verdier-localising functor, as follows from the argument in the proof of Lemma~\reffour{lemma:T-Karoubi-examples}.
In particular, we get that every Verdier-localising functor on $\Catx$ with a finitely cocomplete pointed differentiable target admits a universal Karoubi-localising approximation, and that the latter can be obtained by the Lattice-Tate (or Ind-Calc) construction, or, alternatively, by the decoration construction if the target is stable. 
Related results can also be found in Blumberg-Gepner-Tabuada~\cite{BGT}*{\S 9}.
\end{remark}

The proof of Theorem~\reffour{theorem:karoubification} is obtained by combining Lemma~\reffour{lemma:T-Karoubi-examples} with the following two, essentially independent results:

\begin{proposition}
\label{proposition:is-karoubi}%
Let $\cB$ and $\E$ be as in Theorem~\reffour{theorem:karoubification}, $\T\colon \Catp \to \Ind(\cB)$ a compact reduced functor which is Karoubi preserving and $\F\colon \Catp \to \E$ a Verdier-localising functor.
If $\F^{\T}$ is Karoubi invariant then $\F^{\T^{\infty}}$ is Karoubi-localising.
\end{proposition}

\begin{proposition}
\label{proposition:is-universal}%
Let $\cB$ and $\E$ be as in Theorem~\reffour{theorem:karoubification}, $\T\colon \Catp \to \Ind(\cB)$ a compact reduced functor which is Karoubi preserving and $\F\colon \Catp \to \E$ a functor.
If $\F^{\T^{\infty}}$ is Karoubi-localising then the map $\F \Rightarrow \F^{\T^{\infty}}$ exhibits its target as the Karoubi approximation of its source.
\end{proposition}

\begin{proof}[Proof of Proposition~\reffour{proposition:is-universal}]
This is a well-known argument that goes back to the work of Schwede on symmetric spectra. A proof in the generality we need is recorded in~\cite{motives}*{Proposition 7.1.2}.

Let $\mathrm{Inj}^{\otimes} = \mathrm{Env}(\mathbb{E}_0^{\otimes})$ denote the symmetric monoidal enveloping $\infty$-category of the $\infty$-operad $\mathbb{E}_0^{\otimes}$. Concretely, the objects of $\mathrm{Inj}$ are finite sets, the morphisms are injective maps of such, and the symmetric monoidal structure is given by disjoint union. 
The augmented object $\eps \colon \T \Rightarrow \one^{\red}$ then extends in an essentially unique manner to a symmetric monoidal functor 
\[
\mathrm{Inj}^{\otimes} \to (\Fun^{\omega}_{\ast}(\Catp,\Ind(\cB))\op)^{\Day} \quad\quad I \mapsto \T^{\otimes I} ,
\]
where on the right hand side we have the symmetric monoidal structure induced on opposites by (localised) Day convolution on $\Fun^{\omega}_{\ast}(\Catp,\Ind(\cB))$.
The action of $\Fun^{\omega}_{\ast}(\Catp,\Ind(\cB))\op$ on $\Fun_{\ast}(\Catp,\E)$ via Day cotensoring then restricts to give an action $\beta$ of $\mathrm{Inj}$ on $\Fun_{\ast}(\Catp,\E)$, where $I \in \mathrm{Inj}$ acts on $\Fun_{\ast}(\Catp,\E)$ by $\G \mapsto \G^{\T^{\otimes I}}$. 
Our assumption that cotensors of the form $(-)^b$ for $b \in \cB$ and finite limits commute with sequential colimits in $\E$ then implies $\mathrm{Inj}$ acts on $\Fun_{\ast}(\Catp,\E)$ via sequential-colimit-preserving functors.

A $\G \in \Fun(\Catp,\E)$ is called $\beta$-stable if $\G \to \G^{\T}$ is an equivalence. In particular, every Karoubi-localising $\G$ is $\bet$-stable by our assumptions on $\T$.
Let now $\F\colon \Catp \to \E$ be such that $\F^{\T^{\infty}}$ is Karoubi-localising. Then $\F^{\T^{\infty}}$ is $\beta$-stable and hence $\F$ is $\beta$-semistable in the sense of \cite{motives}*{Proposition 7.1.2}. Loc.\ cit.\ then implies that the map $\F \to \F^{\T^{\infty}}$ is the initial map to a $\beta$-stable object, and in particular to a Karoubi-localising functor.
\end{proof}

\begin{proof}[Proof of Proposition~\reffour{proposition:is-karoubi}]
Since $\T$ is compact there exists a finite $\infty$-category $\I$ and diagrams $\tau\colon \I \to \Ind(\cB)$ and $\Gamma\colon \I\op \to \Catp$ such that $\T$ is a retract of
\[
\T_{\fin} := \colim_{i \in \I} \Gamma(i)_!\tau(i) ,
\]
and so $\G^{\T}$ is a retract of 
\[
\G^{\T_{\fin}} =  \lim_{i \in \I} \G(-\otimes \Gamma(i))^{\tau(i)}
\]
for every reduced $\G\colon \Catp \to \E$. In particular, if $\G$ is Karoubi-invariant then $\G^{\T}$ is again Karoubi-invariant, since tensoring with any Poincaré $\infty$-category preserves Karoubi equivalences by Proposition~\reftwo{lemma:Morita-exact-criterion}\;\reftwoitem{item:preserves-K-equivalences}. 

Now in the case at hand we start with a Verdier-localising functor $\F\colon \Catp \to \E$, in particular reduced, and $\F^\T$ is Karoubi-invariant by assumption. We hence conclude that $\F^{\T^{\otimes n}}$ is reduced and Karoubi-invariant for every $n \geq 1$ and so $\F^{\T^{\infty}}$ is reduced and Karoubi-invariant. To show that $\F^{\T^{\infty}}$ is Karoubi-localising it will hence suffice to show that it sends Poincar\'e-Verdier sequences to fibre sequences. To that end, let $\alp\colon \Del^1\times \Del^1 \to \Catp$ be a Poincaré-Verdier square
\[
\begin{tikzcd}
(\C,\QF)\ar[r]\ar[d] & (\D,\QFD) \ar[d] \\
(\cB,\Pi)\ar[r] & (\E,\QFE)
\end{tikzcd}
\]
and for $n \in \NN \cup \{\infty\}$ write
\[
\Xi(\alp,n) = \bigl[\F^{\T^{\otimes n}}(\C,\QF)  \to \F^{\T^{\otimes n}}(\cB,\Pi)\times_{\F^{\T^{\otimes n}}(\E,\QFE)} \F^{\T^{\otimes n}}(\D,\QFD)\bigr] \in \Ar(\E)
\]
for the arrow whose invertiblity is equivalent to $\F^{\T^{\otimes n}}(\alp)$ being cartesian.
Write $\F^{\T^{\infty}}$ as the colimit of the sequence
\[
\F \Rightarrow \F^{\T} \Rightarrow \F^{\T\otimes \T} \Rightarrow ... \Rightarrow \F^{\T^{\otimes n}} \Rightarrow ...
\]
where in the $n$'th step the map $\F^{\T^{\otimes n}} \Rightarrow \F^{\T^{\otimes n+1}}$ is induced by the map $\F \to \F^{\T}$ via cotensoring with $\T^{\otimes n}$.
Since sequential colimits in $\E$ commute with finite limits by the assumption that $\E$ is differentiable we have $\Xi(\alp,\infty) = \colim_n \Xi(\alp,n)$.
Our goal is hence to show that 
$\colim_n \Xi(\alp,n)$ is an invertible arrow in $\E$, and we will do so by showing that $\Xi(\alp,n) \to \Xi(\alp,n+1)$ factors through an invertible arrow for every $n \geq 0$. 

Now since $\F^{\T}$ is Karoubi invariant the map $\F \Rightarrow \F^{\T}$ factors as a composite $\F \Rightarrow \F(-^{\natural}) \Rightarrow \F^{\T}$, and so the map 
$\F^{\T^{\otimes n}} \Rightarrow \F^{\T^{\otimes n+1}}$ factors through 
$[\F(-^{\natural})]^{\T^{\otimes n}}$. At the same time, the induced map 
$\F^{\T^{\otimes n}} \to [\F(-^{\natural})]^{\T^{\otimes n}}$ is a retract of the induced 
map $\F^{\T^{\otimes n}_{\fin}} \to [\F(-^{\natural})]^{\T^{\otimes n}_{\fin}}$.
To finish the proof it will hence suffice to show that the map of squares
\[
\begin{tikzcd}
\lim\limits_{I \in \I^n} \F((\C,\QF)\otimes \Gamma(I))^{\beta(I)} \ar[r]\ar[d] & \lim\limits_{I \in \I^n} \F((\D,\QFD)\otimes \Gamma(I))^{\beta(I)} \ar[d] \\
\lim\limits_{I \in \I^n} \F((\cB,\Pi)\otimes \Gamma(I))^{\beta(I)}\ar[r] & \lim\limits_{I \in \I^n} \F((\E,\QFE)\otimes \Gamma(I))^{\beta(I)}
\end{tikzcd}
\]
\[
\Downarrow
\]
\[
\begin{tikzcd}
\lim\limits_{I \in \I^n} \F([(\C,\QF)\otimes \Gamma(I)]^{\natural})^{\beta(I)} \ar[r]\ar[d] & \lim\limits_{I \in \I^n} \F([(\D,\QFD)\otimes \Gamma(I)]^{\natural})^{\beta(I)} \ar[d] \\
\lim\limits_{I \in \I^n} \F([(\cB,\Pi)\otimes \Gamma(I)]^{\natural})^{\beta(I)}\ar[r] & \lim\limits_{I \in \I^n} \F([(\E,\QFE)\otimes \Gamma(I)]^{\natural})^{\beta(I)}
\end{tikzcd}
\]
factors through a cartesian square. Now by Lemma~\reffour{lemma:tensor-preserves-karoubi-squares}  each of the squares
\[
\begin{tikzcd}
{[(\C,\QF)\otimes \Gamma(I)]^{\natural}} \ar[r]\ar[d] & {[(\D,\QFD)\otimes \Gamma(I)]^{\natural}} \ar[d] \\
{[(\cB,\Pi)\otimes \Gamma(I)]^{\natural}}\ar[r] & {[(\E,\QFE)\otimes \Gamma(I)]^{\natural}}
\end{tikzcd}
\]
is a Poincaré-Karoubi square among idempotent complete Poincaré $\infty$-categories, hence in particular cartesian. Note that the collection of cartesian squares forms a reflective subcategory inside all squares, with reflection given by replacing the top left corner by the fibre product of the remaining terms. We conclude that the above map of squares factors through the square
\[
\begin{tikzcd}
\lim\limits_{I \in \I^n} \F((\C,\QF)(I))^{\beta(I)} \ar[r]\ar[d] & \lim\limits_{I \in \I^n} \F((\D,\QFD)\otimes \Gamma(I))^{\beta(I)} \ar[d] \\
\lim\limits_{I \in \I^n} \F((\cB,\Pi)\otimes \Gamma(I))^{\beta(I)}\ar[r] & \lim\limits_{I \in \I^n} \F((\E,\QFE)\otimes \Gamma(I))^{\beta(I)},
\end{tikzcd}
\]
where
\[
(\C,\QF)(I) := [(\cB,\Pi) \otimes \Gamma(I)] \times_{(\E,\QFE) \otimes \Gamma(I)} [(\D,\QFD) \otimes \Gamma(I)] .
\]
We now observe that the last square is cartesian since each of the squares
\[
\begin{tikzcd}
(\C,\QF)(I) \ar[r]\ar[d] & (\D,\QFD)\otimes \Gamma(I) \ar[d] \\
(\cB,\Pi)\otimes \Gamma(I)\ar[r] & (\E,\QFE)\otimes \Gamma(I),
\end{tikzcd}
\]
is a Poincaré-Verdier square by Proposition~\reftwo{lemma:Morita-exact-criterion}\;\reftwoitem{item:preserves-PV-proj},
and $\F$ is assumed Verdier-localising.
\end{proof}

It is natural to ask if Theorem~\reffour{theorem:karoubification} can be extended beyond the realm of Verdier-localising functors.
To this end, for $\cB,\E$ and $\T$ as in Theorem~\reffour{theorem:karoubification}, let us denote by $\cK_\T$ the class of functors $\F\colon \Catp \to \E$ for which the map $\F \Rightarrow \F^{\T^{\infty}}$ is a Karoubi-localising approximation. 
By Proposition~\reffour{proposition:is-universal}, this is exactly the class of functors for which $\F^{\T^{\infty}}$ is Karoubi-localising, which means that the class $\cK_\T$ is closed under all operations which commute with $(-)^{\T^{\infty}}$ and preserves Karoubi-localising functors, notably: 
\begin{enumerate}
\item
Pre-composition with $(\C,\QF)\otimes (-)$ for some $(\C,\QF) \in \Catp$.
\item
Finite limits.
\item
All colimits that commute with finite limits in $\E$. 
\end{enumerate}
Now in the case of $\T=\T_{\Tate},\T_{\dec}$, the class $\cK_\T$ contains all Verdier-localising functors by Theorem~\reffour{theorem:karoubification}
and Lemma~\reffour{lemma:T-Karoubi-examples}, and hence also contains all functors of the form $\F((\C,\QF)\otimes (-))$ for Verdier-localising $\F$. We note that $\F((\C,\QF)\otimes (-))$ are not Verdier-localising in general (see Remark~\reftwo{remark:no-verdier}). 

\begin{corollary}
\label{corollary:translated}%
For any Poincaré $\infty$-category $(\C,\QF)$ the maps
\[
\GW((\C,\QF)\otimes(-)) \to \KGW((\C,\QF)\otimes (-)) \quad\text{and}\quad \L((\C,\QF)\otimes (-)) \to \KL((\C,\QF)\otimes(-))
\]
exhibit their targets as initial Karoubi-localising functors under the sources. 
\end{corollary}

We note that the statement of Corollary~\reffour{corollary:translated} is not a formal consequence of the analogous statement for $\KGW$ and $\KL$:
While pre-composition with $(\C,\QF)\otimes (-)$ preserves Karoubi-localising functors, there is a-priori no reason for it to preserve Karoubi-localising approximations, and we do not know if this is true in general. 

Corollary~\reffour{corollary:translated} provides a non-trivial class of non-Verdier-localising additive functors for which a Karoubi-localising approximation can be obtained via the Tate or decoration construction. This does not hold for all additive functors, as the following example shows:

\begin{example}
\label{example:not-all-additive}%
Let $\cB,\E$ be as in Assumption~\reffour{assumption:cotensor} with $\E = \Spa$, and let $\T\colon \Catp \to \Ind(\cB)$ be a Karoubi preserving compact reduced functor (e.g., we can take $\cB$ and $\T$ to be as in Example~\reffour{example:lattice-tate}, or as in Example~\reffour{example:decoration}). In particular, $\T$ is a retract of $\colim_{i \in \I}\Gamma(i)_!\tau(i)$ for some finite $\infty$-category $\I$ and functors $\tau\colon \I \to \cB$ and $\Gamma\colon \I\op \to \Catp$.
Then we can find a large enough uncountable cardinal $\kappa$ such that $\Gamma$ takes values in Poincaré $\infty$-categories with $\kappa$-small mapping spectra.
Let $R = \QQ[x_i]_{i \in \I}$ be the polynomial ring over the rational numbers with set of indeterminates $\{x_i\}_{i \in I}$ of cardinality $|I| = \kappa$, $K$ be the fraction field of $R$, and set 
\[
\F := \KGW(\Funx((\Dperf(K),\QF^{\sym}),-))\colon \Catp \to \Spa,
\]
so that $\F$ is an additive functor. 
We claim that $\F^{\T^{\infty}}$ is not Karoubi-localising. More precisely, we claim that $\F^{\T^{\infty}}$ sends the Karoubi sequence
\[
(\Dperf_K(R),\QF^{\sym}) \to (\Dperf(R),\QF^{\sym}) \to (\Dperf(K),\QF^{\sym})
\]
of Corollary~\reftwo{corollary:ore-localisation} to a sequence of the form $0 \to 0 \to E$ for some non-zero spectrum $E$. To see this, we first note that if $\A$ is a stable $\infty$-category whose mapping spectra are all $\kappa$-small then every exact functor $\Dperf(K) \to \A$ is zero. Indeed, if we let $x \in \A$ be the image of $K \in \Dperf(K)$ and $x$ is not zero then the induced map $K \to \pi_0\map_{\A}(x,x)$ is injective and this is not possible since $K$ is not $\kappa$-small. Secondly, if $\A$ is a stable $\infty$-category and $M \in \Dperf(R)\otimes \A$ is an object whose image in $\Dperf(R')\otimes \A$ is zero for every ring homomorphism $R \to R'$ to a countable ring $R'$, then $M=0$. Indeed, $R$ is a filtered colimit of the countable polynomial rings $\QQ[x_{i_1},...,x_{i_n}]$ ranging over the finite subsets $\{i_1,...,i_n\} \subseteq I$, and so $\Dperf(R)$ is a filtered colimit of the various $\Dperf(\QQ[x_{i_1},...,x_{i_n}])$, and similarly $\Dperf(R)\otimes \A$ is a filtered colimit of the various 
$\Dperf(\QQ[x_{i_1},...,x_{i_n}]) \otimes \A$; the claim then follows from the fact that each of the inclusions $\QQ[x_{i_1},...,x_{i_n}] \subseteq R$ admits a retraction (by sending the other variables to zero, for example).
Combining the above two points with the fact that the mapping spectra in $\Dperf(R')\otimes \A$ are $\kappa$-small whenever $R'$ is countable and $\A$ has $\kappa$-small mapping spectra (indeed, it is enought to consider pure tensors, in which case the mapping spectra are tensor products of a $\kappa$-small spectrum and a countable spectrum, hence again $\kappa$-small) we conclude that if $\A$ has $\kappa$-small mapping spectra then any exact functor $\Dperf(K) \to \Dperf(R)\otimes \A$ is zero, and so $\F^{\T^{\infty}}$ vanishes on $(\Dperf(R),\QF^{\sym})$, and similarly on any Poincaré subcategory of it. We now show that $\F^{\T^{\infty}}(\Dperf(K),\QF^{\sym}) \neq 0$. For this, note that $(\Dperf(K),\QF^{\sym})$ is a symmetric monoidal Poincaré $\infty$-category, so that restriction along the unit Poincaré functor $u\colon (\Spaf,\QF^{\uni}) \to (\Dperf(K),\QF^{\sym})$ gives a natural transformation $u^*\colon \F \Rightarrow \KGW$ fitting in a commutative square
\[
\begin{tikzcd}
\F\ar[r,"u^*"]\ar[d] & \KGW \ar[d,"\simeq"] \\
\F^{\T^{\infty}}\ar[r] & \KGW^{\T^{\infty}} 
\end{tikzcd}
\]
where the right vertical map an equivalence (since we assumed that $\T$ is Karoubi preserving).
It will hence suffice to show that the map
\[
\KGW_0(\Funx((\Dperf(K),\QF^{\sym}),(\Dperf(K),\QF^{\sym}))) \to \KGW_0(\Dperf(K),\QF^{\sym}) ,
\]
induced by the component $u^*_{(\Dperf(K),\QF^{\sym})}$, has a non-zero image. Indeed, the action of $(\Dperf(K),\QF^{\sym})$ on itself via the monoidal product determines a Poincaré section 
\[
(\Dperf(K),\QF^{\sym}) \to \Funx((\Dperf(K),\QF^{\sym}),(\Dperf(K),\QF^{\sym}))
\]
of $u^*_{(\Dperf(K),\QF^{\sym})}$, and $\KGW_0(\Dperf(K),\QF^{\sym})$ is a non-zero ring (it supports a ring homomorphism to $\KK_0(K) = \ZZ$, for example). 
\end{example}

\begin{remark}
\label{rem:accessible-functors-have-Karoubi-approximations}%
We do not know whether every functor $\Catp \to \Spa$ admits a Karoubi-localising approximation. However, in \paperfive 
we will show that if $\kappa$ is an uncountable regular cardinal then every Poincaré-Karoubi square is a $\kappa$-filtered colimits of Poincaré-Karoubi squares among $\kappa$-compact Poincaré $\infty$-categories.
Since $\Catp$ is compactly generated (see Proposition~\reffive{proposition:catp-cath-are-compact}) it then follows from the adjoint functor theorem that any $\kappa$-accessible functor $\Catp \to \Spa$ (hence, eventually, any accessible functor by choosing a large enough $\kappa$) admits a Karoubi-localising approximation. 
\end{remark}

\subsection{The two towers}
\label{subsection:the-two-towers}%

In this final subsection we show that the Tate towers of $\K,\GW$ and $\L$, that is, the towers (or sequences) obtained by iteratively cotensoring with $\T_{\Tate}$, in fact coincide with their decoration towers described in \S\reffour{subsection:universal-dec}. We then revisit the resulting tower from a different perspective, showing that it coincides with the decalage of the stabilisation tower of Definition~\reftwo{definition:stab-construction}.

We begin by proving an analogue of Proposition~\reffour{proposition:connective-covers}:

\begin{proposition}
\label{proposition:connective-covers-tate}%
Let $\FF\colon \Catp \to \Spa$ be a Karoubi-localising functor and let $\F \Rightarrow \FF$ be a map which exhibits $\F$ as the pointwise $m$-connective cover of $\FF$. Then the induced map $\F^{\T_{\Tate}}\Rightarrow \FF^{\T_{\Tate}}$ exhibits $\F^{\T_{\Tate}}$ as the pointwise $(m-1)$-connective cover of $\FF^{\T_{\Tate}} \simeq \FF$. 
As a consequence, for $\F = \tau_{\geq m}\FF$ we have $\F^{\T^{\otimes n}_{\Tate}} = \tau_{\geq m-n-1}\F$ and the associated map
$\F^{\Tate} \Rightarrow \FF$ is an equivalence.
\end{proposition}
\begin{proof}
By definition, the map $\F^{\T_{\Tate}}(\C,\QF)\to \FF^{\T_{\Tate}}(\C,\QF)$ is the map induced on horizontal fibres in the square
\[
\begin{tikzcd}
\F((\C,\QF) \otimes \Latt^{\omega_1}(\Spaf, \QF^{\uni}))\ar[r]\ar[d] & \F((\C,\QF) \otimes \Tate^{\omega_1}(\Spaf, (\QF^{\uni})\qshift{1}) \ar[d] \\
\F((\C,\QF) \otimes \Latt^{\omega_1}(\Spaf, \QF^{\uni}))\ar[r] & \FF((\C,\QF) \otimes \Tate^{\omega_1}(\Spaf, (\QF^{\uni})\qshift{1}) .
\end{tikzcd}
\]
Since $\FF$ is Karoubi-localising it is in particular additive and hence $\F$ is additive as well, so that the terms on the left column vanish by the Eilenberg swindle (as in the proof of Proposition~\reffour{proposition:localising-are-stable}). We may hence identify the map in question with the induced map $\Om\F((\C,\QF) \otimes \Tate^{\omega_1}(\Spaf, (\QF^{\uni})\qshift{1}) \to \FF((\C,\QF) \otimes \Tate^{\omega_1}(\Spaf, (\QF^{\uni})\qshift{1})$ which is indeed an $(m-1)$-connective cover, being the loop of an $m$-connective cover. 
\end{proof}

\begin{corollary}
\label{corollary:two-towers}%
The Tate tower $\F^{\T_{\Tate}^{\otimes n}}$ of $\F$ coincides with the decoration tower $\F^{\<n\>}$ for $\F$ any of the functors $\K,\GW,\L,\K^{\hC},\K^{\tC}$.
\end{corollary}
\begin{proof}
By Lemma~\reffour{lemma:T-Karoubi-examples}
we have $\F^{\T_{\Tate}} \simeq \F^{\T_{\dec}} \simeq \F(-^{\natural})$ for $\F$ any of these three. Now $\K(-^{\natural}) = \tau_{\geq 0}\KK$ and so by Theorem~\reffour{theorem:karoubification} and Proposition~\reffour{proposition:connective-covers-tate} the tower
\[
\K^{\T_{\Tate}} \to \K^{\T_{\Tate}^{\otimes 2}} \to ... \to \K^{\T_{\Tate}^{\otimes n}} \to ... \to \K^{\Tate}
\]
is equivalent to the tower
\[
\tau_{\geq 0}\KK \to \tau_{\geq -1}\KK \to ... \to \tau_{\geq -n+1}\KK \to ... \to \KK ,
\]
and hence coincides with the decoration tower of $\K$ by Corollary~\reffour{corollary:decorated-K}. This holds whether we consider $\K$ as a functor $\Catp \to \Sp$ or $\Catp \to \Sp^{\B\Ct}$. 
Now the Day cotensor operation $\F \mapsto \F^{\T_{\Tate}}$ preserves limits (being right adjoint to the Day convolution with $\T_{\Tate}$), but also colimits by stability, and similarly for $\T_{\dec}$ (as we already observed in Lemma~\reffour{lemma:preserves-limits}). We conclude that the Tate towers of $\K^{\hC}$ and $\K^{\tC}$ are given by
\[
(\tau_{\geq 0}\KK)^{\hC} \to (\tau_{\geq -1}\KK)^{\hC} \to ... \to (\tau_{\geq -n+1}\KK)^{\hC} \to ... \to \KK^{\hC}
\]
and
\[
(\tau_{\geq 0}\KK)^{\tC} \to (\tau_{\geq -1}\KK)^{\tC} \to ... \to (\tau_{\geq -n+1}\KK)^{\tC} \to ... \to \KK^{\tC} ,
\]
and similarly for the corresponding decoration towers (as we already used in the proof of Proposition~\reffour{proposition:decorated-GW}). Finally, for $\GW$ and $\L$ the result follows from the fibre squares
\[
\begin{tikzcd}
\GW(\C,\QF)\ar[r]\ar[d] & \KGW(\C,\QF)\ar[d] & \L(\C) \ar[r] \ar[d] & \KL(\C,\QF) \ar[d] \\
\K(\C,\QF)^\hC\ar[r] & \KK(\C,\QF)^\hC & \K(\C,\QF)^\tC \ar[r] & \KK(\C,\QF)^{\tC}.
\end{tikzcd}
\]
of Theorem~\reffour{theorem:Karoubi-pullback-squares}. In particular, the Tate towers of $\GW$ and $\L$ fit in the same form of fibre squares as we recorded from the decoration tower in Theorem~\reffour{proposition:decorated-GW} and Proposition~\reffour{proposition:decorated-L}.
\end{proof}

We finally offer a third perspective on the equivalent towers of Corollary~\reffour{corollary:two-towers}. Namely, recall from Definition~\reftwo{definition:stab-construction} the two-sided sequence
\[
\dots \longrightarrow \SS^{-1} \otimes \KGW(\C,\QF\qshift{1}) \longrightarrow \KGW(\C,\QF) \longrightarrow \SS^1 \otimes \KGW(\C,\QF\qshift{-1}) \longrightarrow \dots
\]
whose maps are the boundary maps of the Bott-Genauer sequences associated to the various $(\C,\QF\qshift{i})$. By Proposition~\reftwo{proposition:Q-bordism}, the colimit of this sequence is canonically identified with $\KGW^\bord(\C,\QF) \simeq \KL(\C,\QF)$, while the limit is the cobordification $\KGW^\mathrm{cbord}(\C,\QF) \simeq \fib[\LL(\C,\QF) \rightarrow \KK(\C,\QF)^\tC]$. We can therefore regard this sequence as a (non-separated) filtration $\mathrm{stab}_\KGW$ of $\KL$, the stabilisation filtration, and consider its décalage $\dec(\mathrm{stab}_\KGW)$; we will recall the definition in detail below. Roughly, the décalage of a filtered spectrum is a new filtered spectrum, whose associated graded arises from that of the original filtration by taking homology with respect to the canonical differential on it (that also gives rise to the first differential of the associated spectral sequence). In particular, it has the same limit and colimit as the original filtration. 

We claim that the filtration formed by $\L^{\cotower {-n}}(\C, \QF)$ is the décalage of the stabilisation filtration for $n\geq 0$. More precisely, we have:

\begin{proposition}
\label{prop:dec-stab}%
There are cartesian squares
\[
\begin{tikzcd} 
\dec_n(\mathrm{stab}_\KGW) \ar[d] \ar[r] & \LL \ar[d] \\
(\tau_{\geq -n} \KK)^\tC \ar[r] & \KK^\tC
\end{tikzcd}
\]
for all $n \in \mathbb Z$, so that in particular we have $\dec_{n}(\mathrm{stab}_\KGW) \simeq \L^{\cotower{-n}}$ for all $n \geq 0$.
\end{proposition}

To give the definition of the décalage filtration, recall from \cite{ariotta}*{Section 5} the Beilinson $t$-structure on the category $\Fun(\mathbb Z,\Sp)$ of filtered spectra: Its coconnective part is spanned by all $F \colon \mathbb \ZZ \rightarrow \Sp$ for which $F_i$ is $i$-coconnective for all $i \in \mathbb Z$, and its connective part by those $F$ for which $\mathrm{gr}_iF$ is $i$-connective for all $i \in \mathbb Z$, where $\mathrm{gr}_iF = \cof(F_{i-1} \rightarrow F_i)$ is the $i$-th piece of the associated graded spectrum. We will write $\tau_{\geq n}^\mathrm B F$ for the $n$-th connective cover of a filtered spectrum $F$ in this $t$-structure, and denote by $|\tau_{\geq n}^\mathrm B F|$ its colimit.
The décalage $\dec(F)$ of a filtered spectrum $F$ is then given by the sequence
\[
\dots \longrightarrow |\tau_{\geq -n+1}^\mathrm B F| \longrightarrow |\tau_{\geq -n}^\mathrm B F| \longrightarrow \dots
\]
It derives its relevance from the fact that by the main results of \cite{antieau-dec}, the spectral sequence associated to $\dec(F)$ agrees with that associated to $F$ up to a page shift by $1$, that is the $r$-th page for $\dec(F)$ agrees with the $(r+1)$-st page for $F$. This fundamental property allows one to analyse the spectral sequence associated to $F$ entirely in terms of properties of the décalage construction. 

For our purposes the fundamental examples are the following, originally due to Hedenlund, Krause and the eighth author; they are exploited in the forthcoming \cite{KrauseHedenlundNikolaus} to analyse the Tate spectral sequence in detail.

\begin{example}\
\begin{enumerate}
\item Let $X \in \Sp^\hC$ be a $\Ct$-spectrum. Then we have the filtered spectrum, often known as the Tate filtration of $X$, 
\[
\dots \longrightarrow (\SS^{n\sigma} \otimes X)^\hC \longrightarrow (\SS^{(n+1)\sigma} \otimes X)^\hC \longrightarrow \dots,
\]
with maps induced by the cofibre sequence $(\Ct)_+ \rightarrow \ast_+ \rightarrow \mathrm S^\sigma$ of pointed $\Ct$-spaces by tensoring with $\SS^{n\sigma}$. Its limit vanishes and its colimit indeed identifies with $X^\tC$ by noting that the analogous sequence for orbits has vanishing colimit and is constant for the Tate construction since $(\SS^\sigma \otimes X)^\tC \simeq X^\tC$ as a consequence of the cofibre sequence above (and the Tate construction vanishing on induced spectra).

By inspection, the $k$-th connective cover and truncation of this filtration in the Beilinson $t$-structure are then given by 
\[
\dots \longrightarrow (\SS^{n\sigma} \otimes \tau_{\geq k} X)^\hC \longrightarrow (\SS^{(n+1)\sigma} \otimes \tau_{\geq k} X)^\hC \longrightarrow \dots,
\]
and
\[
\dots \longrightarrow (\SS^{n\sigma} \otimes \tau_{\leq k} X)^\hC \longrightarrow (\SS^{(n+1)\sigma} \otimes \tau_{\leq k} X)^\hC \longrightarrow \dots,
\]
respectively: The lower sequence is evidently $k$-coconnective and on account of the cofibre sequence of $\Ct$-spaces above, the $i$-th associated graded of the upper sequence is simply $\mathbb S^i \otimes \tau_{\geq k} X$, which evidently makes it $k$-connective. 

Thus we learn that the décalage of the Tate filtration is given by the Tate construction of the Whitehead filtration:
\[
\dots \longrightarrow (\tau_{\geq -n+1} X)^\tC \longrightarrow (\tau_{\geq -n} X)^\tC \longrightarrow \dots
\]
\item There is a similar filtration 
\[
\dots \longrightarrow (\SS^{n\sigma} \otimes X)^\gC \longrightarrow (\SS^{(n+1)\sigma} \otimes X)^\gC \longrightarrow \dots
\]
for any genuine $\Ct$-spectrum $X \in \Sp^\gC$, with colimit the geometric fixed points $X^\geofix$ and limit $\fib(X^\gC \rightarrow X^\hC) \simeq \fib(X^\geofix \rightarrow X^\tC)$; the argument is the same as above. Again by direct inspection, the $k$-th connective cover and truncation are given by 
\[
\dots \longrightarrow (\SS^{n\sigma} \otimes \tau_{\geq k} X)^\hC \times_{(\SS^{n\sigma} \otimes X)^\hC} (\SS^{n\sigma} \otimes X)^\gC \longrightarrow (\SS^{(n+1)\sigma} \otimes \tau_{\geq k} X)^\hC \times_{(\SS^{(n+1)\sigma} \otimes X)^\hC} (\SS^{(n+1)\sigma} \otimes X)^\gC \longrightarrow \dots,
\]
and
\[
\dots \longrightarrow (\SS^{n\sigma} \otimes \tau_{\leq k} X)^\hC \longrightarrow (\SS^{(n+1)\sigma} \otimes \tau_{\leq k} X)^\hC \longrightarrow \dots,
\]
respectively, since the upper sequence has the same associated graded as in the previous example. It thus follows that the décalage of this filtration of the geometric fixed points has $n$-th term the pullback
\[
\begin{tikzcd}
\dec_n((\SS^{\bullet\sigma} \otimes X)^\gC) \ar[r] \ar[d] & X^\geofix \ar[d] \\
(\tau_{\geq -n} X)^\tC \ar[r] & X^\tC.
\end{tikzcd}
\]
\end{enumerate}
\end{example}
 
Now, for the proof of Proposition \reffour{prop:dec-stab} it only remains to observe that by Karoubi periodicity the stabilisation filtration of $\KL(\C,\QF)$ is precisely the filtration of the geometric fixpoints described in the previous example for the genuine $\Ct$-spectrum $\KKR(\C,\QF)$, see the discussion after \reftwo{remark:compare-bordification}. In fact, this works for an arbitrary additive functor $\F \colon \Catp \rightarrow \Sp$: Recall that for any such $\F$ and $(\C,\QF)$ in $\Catp$ we produced a genuine $\Ct$-spectrum $\F^\ghyp(\C,\QF)$, such that the stabilisation filtration $\mathrm{stab}_\mathcal F$ given by
\[
\dots \longrightarrow \SS^{n} \otimes \F(\C,\QF\qshift{-n}) \longrightarrow \SS^{n+1} \otimes \F(\C,\QF\qshift{-n-1}) \longrightarrow \dots
\]
is precisely the filtration
\[
\dots \longrightarrow (\SS^{n\sigma} \otimes \F^\ghyp(\C,\QF))^\gC \longrightarrow (\SS^{(n+1)\sigma} \otimes \F^\ghyp(\C,\QF))^\gC \longrightarrow \dots
\]
of $\F^\ghyp(\C,\QF)^\geofix \simeq \F^\bord(\C,\QF)$ considered in the previous example by \reftwo{theorem:periodicity}. We thus find:

\begin{proposition}
Let $\mathcal F \colon \Catp \rightarrow \Sp$ be additive and $(\C,\QF) \in \Catp$. Then we have exact squares
\[
\begin{tikzcd}
\dec_n(\mathrm{stab}_\mathcal F) \ar[r] \ar[d] & \F^\bord \ar[d] \\
(\tau_{\geq -n} \F^\hyp)^\tC \ar[r] & (\F^\hyp)^\tC
\end{tikzcd}
\]
for all $n \in \mathbb Z$.
\end{proposition}

Inserting $\F = \KGW$ and using Proposition~\reffour{proposition:decorated-L} we thus obtain Proposition \reffour{prop:dec-stab}. %
For $\F = \GW$ the decalaged stabilisation filtration
agrees with that of $\KGW$ %
in negative degrees, %
then terminates at $\L$ in degree $0$, and is constant in positive degrees.
For $\F$ bordism invariant both the stabilisation filtration and its décalage are constant, and finally for $\F = \K$, we find that the stabilisation filtration has vanishing connecting maps and the décalage vanishes.

\section{Multiplicative structures on Grothendieck-Witt- and L-spectra}
\label{section:day}%

In this section we establish universal multiplicative properties for the Grothendieck-Witt space and spectrum, the $\L$-space and spectrum, and the Karoubi-localising variants $\KGW$ and $\KL$ constructed in \S\reffour{subsection:localising-invariants}, compatibly with the respective fundamental fibre squares. The main results, which rely to a large extent on Appendix~\reffour{section:appendix-localisations}, are collected in \S\reffour{subsection:multiplicative-GW}. We then dedicate \S\reffour{subsection:examples} to a comprehensive review of examples, mostly from the world of commutative rings and ring spectra. In \S\reffour{subsection:PeriodicityviaMultiplication} we discuss the compatibility of the above multiplicative structures with Ranciki periodicity and Karoubi's fundamental theorem, while in \S\reffour{subsection:mult-shaneson} we discuss their compatibility with the circle decomposition theorem of \S\reffour{subsection:shaneson}. In \S\reffour{subsection:normal-L} we discuss the induced multiplicative structure on normal L-theory. Finally, in \S\reffour{subsection:automorphism-GW-L} we use the above multiplicative structures to describe the ring of homotopy classes of self-equivalences of $\GW$ and $\L$.

\subsection{Day convolution of additive and localising functors}
\label{subsection:multiplicative-GW}%
\label{subsection:multiplicative-L}%
\label{subsection:mult-tate}%

Our construction of the desired multiplicative structures starts by interpreting a lax symmetric monoidal structure on a functor $F\colon \C\to \D$ between symmetric monoidal categories $\C^\otimes$ and $\D^\otimes$ as a commutative algebra structure on $F$ in the Day convolution operad $\Fun(\C, \D)^\Day$, whose underlying $\infty$-category is $\Fun(\C, \D)$ and whose $n$-ary morphisms $(F_1, \dots, F_n)\to G$ are given by the multi-natural transformations $\bigotimes_{i=1}^n F_i(c_i) \to G(\bigotimes_{i=1}^n c_i)$. In the case at hand we have $\C=\Catp$ equipped with its symmetric monoidal structure from~\S \refone{subsection:tensor-product} and \S \refone{subsection:monoidal-structure}, and $\D = \Sps, \Sp$ or $\Grp_\Einf$. Upon restriction to finitary (or more generally accessible) functors, the resulting Day convolution operad turns out to be a symmetric monoidal category. In the presence of such a symmetric monoidal structure, the arguments that we are about to give are well-established, and build on the notion of multiplicative Bousfield localisations. However, the restriction to accessible functors turns out to be unnecessary, because the arguments can be completely performed on an operadic level. We have recorded in Appendix \reffour{section:appendix-localisations} %
the requisite extension of the theory of multiplicative and smashing localisations. In what follows, we shall freely use constructions and results from that appendix.

Taking the target $\D=\Spa$ for illustration, our argument is based on the following key observations:
\begin{enumerate}
\item
The Day convolution $\infty$-operad $\Fun(\Catp,\Spa)^{\Day}$, though not a symmetric monoidal $\infty$-category, still admits a unit object, given explicitly by $\SS[\Poinc]$.
\item
Each of the functors $\GW,\KGW,\L,\KL$ is initial under $\SS[\Poinc]$ with respect to a certain property: being additive, Karoubi-localising, or a combination of one of these with bordism invariance. Each of these properties can be phrased as sending a certain collection of cone diagrams to limit diagrams.
\item
When the collection of cone diagrams is suitably compatible with the tensor product on $\Catp$, such universal approximations of the unit in $\Fun(\Catp,\Spa)^{\Day}$ become themselves the units of the full suboperad spanned by the functors which send the relevant collection of cone diagrams to limit diagrams. In particular, these universal approximations refine to initial commutative algebra objects, hence initial lax symmetric monoidal functors of the given type. %
\end{enumerate}
Replacing the target with $\Grp_{\Einf}$ or $\Sps$ one can also consider in this manner the functors $\GWspace$ and $\Lspace$. On the other hand, one can also replace $\Catp$ with $\Catx$ and consider functors such as $\K,\KK$ or $\Kspace$, thus recovering classical multiplicative properties of $\K$-theory.

To perform our arguments simultaneously for all our examples, we start with symmetric monoidal $\infty$-categories $\C$ and $\D$, %
with $\D$ cocomplete and such that its tensor product preserves colimits in each argument. 
In this situation, the Day convolution operad admits a unit object $\one_\Day = [c \mapsto (\one_\D)_{\Hom_\C(\one_\C, c)}]$, which consequently refines to form the initial lax symmetric monoidal functor $\C \to \D$. 
Let $\cR$ be a collection of cone diagrams $K^{\triangleleft} \rightarrow \C$ (with varying $K \in \Cat$) and denote by $\Fun^\cR(\C,\D) \subseteq \Fun(\C,\D)$ the full subcategory spanned by those functor that take each diagram in $\cR$ to a limit cone. We will refer to such functors as $\cR$-local.
Here are the examples that we have in mind:

\begin{definition}
\label{def:collections-R}%
We define the following collection of cone diagrams in $\Catp$:
\begin{enumerate}
\item We write $\cR_{\add}$ for the collection consisting of $0 \colon \ast \rightarrow \Catp$ and all split Poincaré-Verdier squares.
\item We write $\cR_{\kloc}$ for the collection consisting of $0 \colon \ast \rightarrow \Catp$ and all Poincaré-Karoubi squares.
\item We write $\cR_{\binv}$ for the collection consisting of all cylinder inclusions $\cyl\colon (\C,\QF) \to \Q_1(\C,\QF)$.
\item We write $\cR_{\bord}$ for $\cR_\add \cup \cR_{\binv}$.
\end{enumerate}
\end{definition}

Then, by definition, for an $\infty$-category $\E$, %
a functor $\F\colon \Catp \to \E$ is
\begin{enumerate}
\item $\cR_{\add}$-local if and only if it is additive,
\item $\cR_{\kloc}$-local if and only if it is Karoubi-localising,
\item $\cR_{\binv}$-local if and only if it is bordism invariant, and hence
\item $\cR_\bord$-local if and only if it is additive and bordism invariant, and similarly
\item $(\cR_{\kloc} \cup \cR_{\bord})$-local if and only if it is Karoubi-localising and bordism invariant.
\end{enumerate}

\begin{remark}
\label{remark:binv-is-not-add}%
Let us point out that while in \papertwo we only considered the notion of bordism invariance for additive functors, here we separate additivity and bordism invariance. This leads to a better behaved theory for functors with unstable targets, such as spaces. For example,  
the second part of Corollary~\reffour{corollary:universal-multiplicative-L-space} fails for $\cR_{\bord}$ in place of $\cR_{\binv}$.
\end{remark}

The following result is obtained by applying 
Corollary~\reffour{corollary:multiplicative-day} 
of Appendix \reffour{section:appendix-localisations} 
to the $\infty$-operads $\cO=\Com$ and $\cO=\MCom$. %

\begin{proposition}
\label{proposition:universal-unit}%
Let $\C$ and $\D$ be symmetric monoidal $\infty$-categories %
such that $\D$ is cocomplete and closed, %
and let $\cR$ be a collection of cone diagrams in $\C$ as described above. Assume that $\cR$ is stable under tensor products with arbitrary objects of $\C$, and suppose given a natural transformation $\eta\colon \one_{\Day} \Rightarrow F$ which exhibits $F$ as the initial $\cR$-local functor under $\one_{\Day}$. Then the following holds:
\begin{enumerate}
\item
$F$ carries an essentially unique lax symmetric monoidal structure compatible with $\eta$.
\item
This structure renders $F$ initial in lax symmetric monoidal $\cR$-local functors $\C \to \D$. 
\item
If $\E$ is a cocomplete $\infty$-category equipped with a closed action of $\D$ then any $\cR$-local functor $\C \rightarrow \E$ is a module over $F$ in an essentially unique fashion.
\end{enumerate}
\end{proposition}

It turns out that for some of the functors we consider, notably $\GW$, more is true: not only is $\GW$ a commutative algebra objects with respect to Day convolution, it is actually an idempotent algebra object, and Day convolving with it yields a smashing localisations of $\Fun(\Catp,\Spa)^{\Day}$ into its full suboperad spanned by the additive 
functors. This means in particular that the notion of an additive functor $\G\colon \Catp \to \Spa$ is equivalent to that of $\GW$-module, and every $\G$ admits a universal additive approximation given by Day convolution with $\GW$. The analogous statement also holds for the $\L$-space functor $\Lspace\colon \Catp \to \Sps$ and the $\K$-spectrum and $\K$-space functors on $\Catx$,
for example, but fail for $\KGW,\KK$ and $\KL$, as well as for $\GWspace$ (see Remark~\reffour{remark:not-convolving}). The required sufficient conditions for a smashing localisation to arise are summarised in the following result, %
which is a special case of Proposition~\reffour{proposition:smashing-day} of Appendix \reffour{section:appendix-localisations}. Here, recall that a functor $F\colon \C \to \D$ is said to be \defi{small} if it is left Kan extended from a small subcategory of $\C$. As we establish in Corollary~\reffour{corollary:persist-locally-small}, under the conditions of the following proposition (for instance) small functors are symmetric monoidal objects in $\Fun(\C,\D)^{\Day}$,
that is, Day convolution of a small functor against any functor is well-defined.

\begin{proposition}
\label{proposition:sample-day-smashing}%
Let $\C$ and $\D$ be closed symmetric monoidal $\infty$-categories %
such that $\D$ is cocomplete, %
and let $\cR$ be a collection of cone diagrams in $\C$ as described above. 
Suppose given a natural transformation $\eta\colon \one_{\Day} \Rightarrow F$ 
such that $F$ is small and $\cR$-local,
and such that the following hypotheses hold:
\begin{enumerate}
\item
The collection $\cR$ is stable under tensoring and cotensoring with objects in $\C$. 
\item
The class of $\cR$-local functors is closed under colimits and tensoring with objects of $\D$.
\end{enumerate}
Then $\eta$ exhibits $F$ as the initial $\cR$-local functor under $\one_{\Day}$ if and only if the induced map $G \mapsto F \otimes_{\Day} G$ is an equivalence for every $\cR$-local $G$, in which case
$F$ is idempotent with respect to Day convolution and convolving with it yields a smashing localisation of $\infty$-operads
\[
\Fun(\C,\D)^{\Day} \longrightarrow \Fun^\cR(\C,\D)^{\Day}, \quad G \longmapsto F \otimes_\Day G .
\]
In particular, in this case $\Fun^\cR(\C,\D)^{\Day}$ identifies with the $\infty$-operad of $F$-modules in $\Fun(\C,\D)^{\Day}$.
\end{proposition}

To apply these general results we now verify the relevant hypotheses for the cone diagrams of Definition~\reffour{def:collections-R}. Note that the tensoring and cotensoring mentioned in these hypotheses correspond to the operations $(\C,\QF) \otimes -$ and $\Funx((\C,\QF),-)$ of taking tensor products with, and internal mapping objects from, a given $(\C,\QF) \in \Catp$.

\begin{proposition}
\label{prop:closure-properties-of-Rs}%
The above defined collections satisfy the following properties:
\begin{enumerate}
\item
\label{item:list1}%
The collections $\cR_{\add}$ and $\cR_{\binv}$ are closed under (co)tensoring with arbitrary $(\C,\QF) \in \Catp$.
\item
\label{item:list2}%
The collection $\cR_{\kloc}$ is closed under tensoring with arbitrary $(\C,\QF) \in \Catp$.
\end{enumerate}
\end{proposition}
\begin{proof}
First we note that $0 \otimes (\C,\QF) = 0 = \Funx((\C,\QF),0)$. 
Combined with Proposition~\reffour{lemma:tensor-preserves-karoubi-squares}, this establishes Part \reffouritem{item:list2}. We now address Part \reffouritem{item:list1}.

To obtain the claim for $\cR_{\add}$ we show that the functors $(\C,\QF)\otimes -$ and $\Funx((\C,\QF),-)$ preserve split Poincar\'e-Verdier squares.
For the latter, since it preserves limits and adjunctions on the level of underlying $\infty$-categories, it suffice to check that it preserves Poincaré-Verdier projections, which follows from Lemma~\reftwo{lemma:and-smashing}.
For the former, we have established in Proposition~\reftwo{lemma:Morita-exact-criterion}\;\reftwoitem{item:preserves-split-PV-sequences} that $(\C,\QF) \otimes -$ preserves split Poincaré-Verdier sequences.
By Proposition~\reftwo{proposition:poincare-verdier-square-char} it hence suffices to show that it also preserves vertically left and right adjointable squares whose vertical legs admit left and right adjoints. This follows again since $\C\otimes -$ preserves adjunctions, so a left or right adjoint to a given vertical leg induces the same kind of adjoint after tensoring with $\C$, with unit and counit induced by the original ones. The resulting Beck-Chevalley transformations are then also induced from the original ones after tensoring with $\C$, and are hence in particular invertible as soon as the original ones were so. 

To see the claim about $\cR_\binv$, we recall that $\Q_1(\D,\QFD)$ is the cotensor $(\D,\QFD)^{\mathcal{P}}$ where $\mathcal{P}$ is the finite poset classifying pushouts as considered earlier.
Under this identification, the Poincaré functor $\cyl$ is the one induced by terminal map $\mathcal{P} \to \ast$.
From the natural equivalence $\Funx((\C,\QF),(\D,\QFD)^{\mathcal{P}}) = \Funx((\C,\QF),(\D,\QFD))^{\mathcal{P}}$ we find that $\cR_{\binv}$ is closed under cotensors. Finally, it follows from Corollary~\refone{corollary:internally-corepresented} that since $\mathcal{P}$ is a finite poset, $(\Spaf,\QF^\uni)^{\mathcal{P}}$ is dualisable in $\Catp$ with dual $(\Spaf,\QF^\uni)_{\mathcal{P}}$. In particular, the natural map
$(\D,\QFD) \otimes (\Spaf,\QF^{\uni})^{\mathcal{P}} \to (\D,\QFD)^{\mathcal{P}} = \Funx((\Spaf,\QF^\uni)_{\mathcal{P}},(\D,\QFD))$ is an equivalence for all $(\D,\Phi)$, implying the remaining claim. 
\end{proof}

\begin{remark}
\label{rem:contensor-does-not-preserve-KS}%
The collection $\cR_{\kloc}$ is not closed under cotensoring with arbitrary $(\C,\QF) \in \Catp$. For instance, $\Funx((\Dperf(\QQ),\QF^\sym),-)$ does not preserve Poincar\'e-Karoubi sequences as it sends the Poincar\'e-Verdier sequence 
\[
(\Dperf(\ZZ)^\mathrm{tors},\QF^\s) \to (\Dperf(\ZZ),\QF^\s) \to (\Dperf(\QQ),\QF^\s)
\]
to the sequence $0 \to 0 \to (\Dperf(\QQ),\QF^\s)$ which is far from a Poincar\'e-Karoubi sequence. 
\end{remark}

Before moving on to our applications of Propositions~\reffour{proposition:universal-unit} and \reffour{proposition:sample-day-smashing},
we record the following:
\begin{remark}
\label{remark:small-functors}%
All of the functors $\Kspace,\GWspace,\Lspace,\K,\GW,\L,\KK,\KGW,\KL$ preserve filtered colimits. Consequently, since $\Catx$ and $\Catp$ are compactly generated (see Proposition~\reffive{proposition:catp-cath-are-compact}), all of these functors are small. Indeed, that $\Kspace$, $\K$ and $\KK$ preserve filtered colimits is well-known. That $\GWspace,\GW,\Lspace$ and $\L$ preserve filtered colimits was noted in Corollaries \reftwo{corollary:filtcolimpresgwspace}, \reftwo{corollary:filtcolimpresgwspectrum}, Theorem~\reftwo{theorem:L-is-localising} and Corollary~\reftwo{corollary:Ladd}, respectively. Finally, for $\KGW$ and $\KL$ this is established in Corollary \reffour{corollary:filtcolimpreskgwspectrum}.
\end{remark}

\begin{corollary}
\label{corollary:GWspace-universal-multiplicative}%
The natural transformation $\Poinc \Rightarrow \GWspace$
refines in an essentially unique manner to a lax symmetric monoidal natural transformation among lax symmetric monoidal functors, making $\GWspace$ into the initial lax symmetric monoidal group-like additive functor $\Catp \to \Sps$.
Furthermore, any group-like additive functor $\Catp \to \E$ valued in a %
cocomplete and finitely complete $\infty$-category $\E$ is a module over $\GWspace$ in an essentially unique manner.
\end{corollary}
\begin{proof}
Recall that group-like additive functors $\Catp \to \Sps$ can be canonically identified with additive functors $\Catp \to \Grp_{\Einf}$ (or, more generally, group-like additive functors $\Catp \to \E$ with additive functors $\Catp \to \Grp_{\Einf}(\E)$). Viewing $\GWspace$ as a $\Grp_{\Einf}$-valued functor, 
we first deduce from Observation~\reftwo{observation:univ-GW-space} that $\GWspace$ is an $\cR_\add$-localisation of the Day unit 
\[
\Fr_{\Einf-\grp}[\Poinc] \in \Fun(\Catp,\Grp_{\Einf})^\Day ,
\]
where $\Fr_{\Einf-\grp}$ denotes the free $\Einf$-group functor, and then apply Proposition~\reffour{proposition:universal-unit}
with $\cR = \R_\add$ using Proposition~\reffour{prop:closure-properties-of-Rs}, obtaining the stated result with $\Fr_{\Einf-\grp}[\Poinc]$ in place of $\Poinc$. We then finish by noting that natural transformations of space-valued functors $\Poinc \Rightarrow \GWspace$ correspond via adjunction to natural transformations of $\Grp_{\Einf}$-valued functors $\Fr_{\Einf-\grp}[\Poinc] \Rightarrow \GWspace$ (and the same holds for lax symmetric monoidal natural transformations; the latter are contractible in both cases since $\Poinc$ and $\Fr_{\Einf-\grp}[\Poinc]$ are the initial lax symmetric monoidal functors to $\Sps$ and $\Grp_{\Einf}$, respectively).
\end{proof}

\begin{corollary}\
\label{corollary:GW-universal-multiplicative}%
\begin{enumerate}
\item
\label{item:GW-universal}%
The natural transformation $\SS[\Poinc] \Rightarrow \GW$ refines in an essentially unique manner to a lax symmetric monoidal natural transformation among lax symmetric monoidal functors, making $\GW$ into the initial 
lax symmetric monoidal additive functor $\Catp \rightarrow \Spa$. Furthermore, any additive functor $\Catp \to \E$ valued in a stable cocomplete $\infty$-category $\E$ is a module over $\GW$ in an essentially unique manner.
\item
\label{item:GW-idempotent}%
The functor $\GW\colon \Catp \to \Spa$ is Day idempotent and Day convolution with it yields a smashing localisation of $\infty$-operads
\[
\Fun(\Catp,\Spa)^{\Day} \locadjL{\GW \otimes_{\Day} (-)} \Funadd(\Catp,\Spa)^{\Day} .
\]
\end{enumerate}
\end{corollary}
\begin{proof}
We deduce from Corollary~\reftwo{corollary:univ-GW-spectrum} that $\GW$ is an $\cR_\add$-localisation of $\SS[\Poinc]$, so the first part follows similarly as above from Proposition~\reffour{proposition:universal-unit}.
To see the second part, note that the collection of spectrum-valued $\cR_\add$-local functors is closed under colimits and tensors with spectra since $\Spa$ is stable. The claim then follows from
Proposition~\reffour{proposition:sample-day-smashing}
since $\GW$ is small by Remark~\reffour{remark:small-functors}.
\end{proof}

\begin{remark}
It turns out that $\GWspace$ is also Day idempotent, both as a space and an $\Einf$-group valued functor, see~\cite{saunier-motives}, but the category it localises onto is not that of additive functors. Instead, a functor is $\G\colon \Catp \to \Grp_{\Einf}$ is invariant under Day convolving with $\GWspace\colon \Catp \to \Grp_{\Einf}$ if and only if for every Poincaré $\infty$-category $\C$ containing an isotropic subcategory $\Lag \hrar \C$ the induced map $\G(\C) \to \G(\Hyp(\Lag)) \times \G(\Hlgy(\Lag))$ is an equivalence, while a functor $\G\colon \Catp \to \Sps$ is invariant under Day convolving with $\GWspace\colon \Catp \to \Sps$ if and only if for every $\C$ and $\Lag \hrar \C$ as before 
the induced map $\G(\C) \to \G(\Hyp(\Lag) \times \Hlgy(\Lag))$ is an equivalence. %
The $\infty$-category of additive functors $\Catp \to \Grp_{\Einf}$ is, in turn, simply not a smashing localisation (otherwise its idempotent would have to be $\GWspace$), and we do not know whether every functor $\Catp \to \Grp_{\Einf}$ admits an additive approximation.
\end{remark}

\begin{corollary}
\label{corollary:kgw-multiplicative}%
The natural transformation $\SS[\Poinc] \Rightarrow \KGW$ refines in an essentially unique manner to a lax symmetric monoidal natural transformation among lax symmetric monoidal functors, making $\KGW$ into the initial 
lax symmetric monoidal Karoubi-localising functor $\Catp \rightarrow \Spa$. Furthermore, any Karoubi-localising functor $\Catp \to \E$ valued in a stable cocomplete $\infty$-category $\E$ is a module over $\GWspace$ in an essentially unique manner.
\end{corollary}
\begin{proof}
By Proposition~\reffour{proposition:universal-KGW} we find that $\KGW$ is a $\cR_\kloc$-localisation of $\SS[\Poinc]$ so that we may apply Proposition~\reffour{proposition:universal-unit}
with $\cR = \cR_\kloc$.
\end{proof}

\begin{remark}\
\label{remark:not-convolving}%
We do not now whether $\KGW$ is idempotent, nor whether every functor $G\colon \Catp \to \Spa$ admits a universal Karoubi-localising approximation $G \to \mathbb{G}$, but we do know that such an approximation cannot be obtained by Day convolving with $\KGW$. Indeed, if that were true then the Day convolution 
\[
\KGW \otimes_\Day \SS[\Hom_{\Catp}((\C,\QF),-)] \simeq \KGW(\Funx((\C,\QF),-))
\]
would be Karoubi-localising, and this fails, e.g., for $(\C,\QF) = (\Dperf(\QQ),\QF^\sym)$ %
by Remark~\reffour{rem:contensor-does-not-preserve-KS}. 
\end{remark}

Before moving on, we will establish the following strengthening of the universal property of the L-space functor established in \papertwo. To avoid confusion, recall from Remark~\reffour{remark:binv-is-not-add} that, unlike in loc.\ cit., here we do not require bordism invariant functors to be additive. %
 
\begin{proposition}
\label{proposition:universal-L-space}%
The functor $\Lspace \colon \Catp \rightarrow \Sps $ is the initial bordism invariant functor under $\Poinc$.
\end{proposition}
\begin{proof}
The functor $\Lspace$ is bordism invariant by Theorem~\reftwo{theorem:L-is-localising}. 
Now suppose that $\G\colon \Catp \to \Sps$ is a bordism invariant functor. From the equivalence 
$\Lspace= \colim_{n \in \Delta\op} y((\Spaf,\QF^\uni)_{\T_n})$, where $\T_n$ denotes the face poset of $\Delta^n$, we deduce that the map
$\Nat(\Lspace,\G) \to \Nat(\Poinc,\G)$
identifies with the projection
\[
\lim_{n \in \Del}\G((\Spaf,\QF^{\uni})_{\T_n}) \to \G(\Spaf,\QF^{\uni}).
\]
To show that this map is an equivalence, by the bordism invariance of $\G$, it suffices to show that each of the Poincaré functors $(\Spaf,\QF^{\uni})_{\T_n} \to (\Spaf,\QF^{\uni})$ is a bordism equivalence. By Corollary~\refone{corollary:internally-corepresented}, this functor is obtained by applying the functor $\Funx(-,(\Spaf,\QF^{\uni}))$ to $(\Spaf,\QF^{\uni}) \to (\Spaf,\QF^{\uni})^{\T_n}$, and it is enough to prove that the latter is a bordism equivalence.
We then observe that this last functor is the inclusion of the homology $\infty$-category $\Hlgy(\Lag_n)$ for an isotropic subcategory $\Lag_n \subseteq (\Spaf)^{\T_n}$ and then use that such an inclusion is a bordism equivalence, see Example~\reftwo{example:isotropic}. To define $\Lag_n$, let $\dM_n \subseteq (\Spaf)^{\T_n}$ be the full subcategory spanned by those diagrams $\vphi\colon  \T_{n} \rightarrow \Spaf$ such that $\vphi(S \cup \{0\}) \rightarrow \vphi(S)$ is an equivalence for every $S$ (equivalently, which are left Kan extended from their restriction to the subposet $\T^0_n \subseteq \T_n$ spanned by those $S$ which contain $0$). We then define $\Lag_n \subseteq \dM_n$ 
to be the full subcategory spanned by those diagrams which additionally satisfy $\varphi(0) \simeq 0$. 
The inclusions $\Lag_n \subseteq \dM_n \subseteq (\Spaf)^{\T_n}$ then both admit right adjoints and one readily checks that for $\varphi \in \dM_n$ there is an equivalence
$
\QF^{\T_n}(\varphi) \simeq \QF(\varphi(0)),
$
so that $\Lag_n$ is indeed isotropic. Furthermore, inspecting the explicit formula for the duality we find that $\Dual^{\T_n}(\dM_n)$ consists of those diagrams whose restriction to $\T^0_n$ is constant, so that
$(\Lag_n)^\perp = \dM_n$ and 
$\Hlgy(\Lag_n) \simeq \dM_n \cap \Dual^{\T_n}(\dM_n)$ contains precisely the constant diagrams.
\end{proof}

\begin{corollary}
\label{corollary:universal-multiplicative-L-space}%
\
\begin{enumerate}
\item
The natural transformation $\Poinc \Rightarrow \Lspace$ refines in an essentially unique manner to a lax symmetric monoidal natural transformation among lax symmetric monoidal functors, making $\Lspace$ into the initial 
lax symmetric monoidal bordism invariant functor $\Catp \to \Sps$. Furthermore, any bordism invariant functor $\Catp \to \E$ valued in a cocomplete $\infty$-category $\E$ is a module over $\Lspace$ in an essentially unique manner.
\item
The functor $\Lspace\colon \Catp \to \Sps$ is Day idempotent and Day convolution with it yields a smashing localisation of $\infty$-operads
\[
\Fun(\Catp,\Sps)^{\otimes} \locadjL{\Lspace \otimes_{\Day} (-)} \Fun^\binv(\Catp,\Sps)^{\otimes} .
\]
\end{enumerate}
\end{corollary}
\begin{proof}
The first part follows from Proposition~\reffour{proposition:universal-unit}
for $\cR=\cR_\binv$ using Propositions~\reffour{prop:closure-properties-of-Rs} and~\reffour{proposition:universal-L-space}. Since arbitrary colimits of equivalences are equivalences, the second part then follows from Proposition~\reffour{proposition:sample-day-smashing}
since $\Lspace$ is small by Remark~\reffour{remark:small-functors}.
\end{proof}

\begin{variant}
Proposition~\reffour{proposition:sample-day-smashing} can also be applied to $\Lspace$ as an $\Einf$-monoid valued functor. For this, consider the collection $\cR^{\oplus-\binv}$ consisting of the diagrams in $\cR^{\binv}$ together with the diagram $\Del^0 = \emptyset^{\triangleleft}\to \Catp$ picking the zero Poincaré category and the diagrams $(\Del^0 \coprod \Del^0)^{\triangleleft} \to \Catp$ witnessing all cartesian products in $\Catp$. Then a functor $\F\colon \Catp \to \Mon_{\Einf}$ is $\cR^{\oplus-\binv}$-local if and only if it is bordism invariant and direct sum preserving. Since $\Catp$ and $\Mon_{\Einf}$ are semi-additivly symmetric monoidal the first two assumptions of Proposition~\reffour{proposition:sample-day-smashing} hold for $\cR^{\oplus-\binv}$ as well. Finally, since $\Lspace$ and $\Poinc$ are both direct sum preserving as $\Mon_{\Einf}$-valued functors and since direct sum preserving functors $\Catp \to \Mon_{\Einf}$ embed in space-valued functors (as the product preserving ones), Proposition~\reffour{proposition:universal-L-space} implies that $\Lspace$ is also initial as a bordism invariant direct sum preserving functor $\Catp \to \Mon_{\Einf}$ under $\Poinc$.
We conclude from Proposition~\reffour{proposition:sample-day-smashing} that $\Lspace$ is idempotent as a $\Mon_{\Einf}$-valued functor and Day convolving with it localises into the subcategory of bordism invariant direct sum preserving functors. A-posteriori, this implies that all bordism invariant direct sum preserving functors to $\Mon_{\Einf}$ are group-like, namely, taking values in $\Grp_{\Einf}$, and that also as an $\Einf$-group valued functor $\Lspace$ is idempotent and localises to bordism invariant direct sum preserving functors.
\end{variant}

\begin{corollary}
\label{corollary:lmonoidal}%
\
\begin{enumerate}
\item
\label{item:L-universal}%
The natural transformation $\SS[\Poinc] \Rightarrow \L$ refines in an essentially unique manner to a lax symmetric monoidal natural transformation among lax symmetric monoidal functors, making $\L$ into the initial 
lax symmetric monoidal bordism invariant and additive functor to $\Sp$. Furthermore, any additive bordism invariant functor $\Catp \to \E$ valued in a stable cocomplete $\infty$-category $\E$ is a module over $\L$ in an essentially unique manner.
\item
\label{item:L-smashing}%
The functor $\L\colon \Catp \to \Spa$ is Day idempotent and Day convolution with it yields a smashing localisation of $\infty$-operads
\[
\Fun(\Catp,\Spa)^{\otimes} \locadjL{\L \otimes_{\Day} (-)} \Funbord(\Catp,\Spa)^{\otimes}.
\]
\end{enumerate}
\end{corollary}
In particular, the bordification $(-)^{\bord}$ of additive functors studied in \S\reftwo{subsection:bordification} satisfies $(-)^\bord \simeq \L \otimes_\Day -$.
\begin{proof}
The first part follows from Proposition~\reffour{proposition:universal-unit}
for $\cR=\cR_\bord$ using Proposition~\reffour{prop:closure-properties-of-Rs} and Theorem~\reftwo{theorem:lisbordgw}. Since arbitrary colimits of equivalences are equivalences and pushouts in spectra are closed under colimits and tensor products, the second part then follows from 
Proposition~\reffour{proposition:sample-day-smashing}
using in addition Remark~\reffour{remark:small-functors}.
\end{proof}

\begin{corollary}
\label{corollary:kl-monoidal}%
The natural transformation $\SS[\Poinc] \Rightarrow \KL$ refines in an essentially unique manner to a lax symmetric monoidal natural transformation among lax symmetric monoidal functors, making $\KL$ into the initial 
lax symmetric monoidal bordism invariant and Karoubi-localising functor to $\Sp$. Furthermore, any bordism invariant Karoubi-localising functor $\Catp \to \E$ valued in a stable cocomplete $\infty$-category $\E$ is a module over $\KL$ in an essentially unique manner.
\end{corollary}
\begin{proof}
We have established in Proposition~\reffour{proposition:universal-KL} above that $\KL$ is the initial bordism invariant and Karoubi-localising functor under $\SS[\Poinc]$, so the result follows from  Proposition~\reffour{proposition:universal-unit}
with $\cR_{\bord} \cup \cR_{\kloc}$.
\end{proof}

\begin{remark}
\label{remark:KL-not-convolving}%
As with $\KGW$ (see Remark~\reffour{remark:not-convolving}), we do not now whether $\KL$ is idempotent, nor whether every functor $G\colon \Catp \to \Spa$ admits a universal Karoubi-localising bordism invariant 
approximation $G \to \mathbb{G}$, but again as for $\KGW$ we do know that such an approximation cannot be obtained by Day convolving with $\KL$, since by Remark~\reffour{rem:contensor-does-not-preserve-KS} functors of the form $\KL(\Funx((\C,\QF),-))$ are generally not Karoubi-localising. 
\end{remark}

\begin{remark}
One can also apply Proposition~\reffour{proposition:sample-day-smashing}
with $\cR = \cR_\binv$ and $\D=\Spa$. We then find that $\SS[\Lspace] \otimes_\Day -$ implements bordification for functors $\Catp \rightarrow \Sp$. Using Corollary \reffour{corollary:GW-universal-multiplicative} above we thus find that $\SS[\Lspace] \otimes_{\Day} \GW$ is the bordification of $\GW$ inside all functors, but retains its additivity. By Theorem~\reftwo{theorem:lisbordgw} we thus find $\L \simeq \SS[\Lspace] \otimes_\Day \GW$.
\end{remark}

For completeness sake, we also record the multiplicative properties of K-theory.
\begin{corollary}
\label{corollary:K-universal}%
\
\begin{enumerate}
\item 
The natural transformation $\core \Rightarrow \Kspace$
refines in an essentially unique manner to a lax symmetric monoidal natural transformation among lax symmetric monoidal functors, making $\Kspace$ into the initial lax symmetric monoidal group-like additive functor $\Catx \to \Sps$. Furthermore, any group-like additive functor $\Catx \to \E$ valued in a cocomplete $\infty$-category with finite products $\E$ is a module over $\Kspace$ in an essentially unique manner.
\item
\label{item:k-idempotent}%
The functor $\Kspace\colon \Catp \to \Grp_{\Einf}$ is Day idempotent and Day convolution with it yields a smashing localisation of $\infty$-operads
\[
\Fun(\Catx,\Grp_{\Einf})^{\Day} \locadjL{\Kspace \otimes_{\Day} (-)} \Funadd(\Catx,\Grp_{\Einf})^{\Day} .
\]
\item 
The natural transformation $\SS[\core] \Rightarrow \K$ refines in an essentially unique manner to a lax symmetric monoidal natural transformation among lax symmetric monoidal functors, making $\K$ into the initial lax symmetric monoidal additive functor $\Catx \rightarrow \Spa$. Furthermore, any additive functor $\Catx \to \E$ valued in a stable cocomplete $\infty$-category $\E$ is a module over $\K$ in an essentially unique manner. 
\item
\label{item:K-idempotent}%
The functor $\K\colon \Catp \to \Spa$ is Day idempotent and Day convolution with it yields a smashing localisation of $\infty$-operads
\[
\Fun(\Catx,\Spa)^{\Day} \locadjL{\K \otimes_{\Day} (-)} \Funadd(\Catx,\Spa)^{\Day} .
\]
\item 
The natural transformation $\SS[\core] \Rightarrow \KK$ refines in an essentially unique manner to a lax symmetric monoidal natural transformation among lax symmetric monoidal functors, making $\KK$ into the initial lax symmetric monoidal Karoubi-localising functor $\Catx \rightarrow \Spa$. Furthermore, any Karoubi-localising functor $\Catx \to \E$ valued in a stable cocomplete $\infty$-category $\E$ is a module over $\KK$ in an essentially unique manner.
\end{enumerate}
\end{corollary}
\begin{proof}
We apply Propositions~\reffour{proposition:universal-unit} and \reffour{proposition:sample-day-smashing}
with $\C=\Catx$, target either $\Grp_{\Einf}$ or $\Sp$, and the variants of the collections $\cR_{\add}$ and $\cR_\kloc$ from Definition~\reffour{def:collections-R} obtained by discarding Poincar\'e structures. All the arguments here are the same as their hermitian analogues above, with the exception of Statement~\reffouritem{item:k-idempotent} (whose analogue for $\GWspace$ is false), which requires verifying that the collection of additive functors $\F\colon \Catx \to \Grp_{\Einf}$ is closed under colimits and tensoring with objects in $\Grp_{\Einf}$; indeed, since split Verdier projections admit sections in $\Catx$ we have that $\F$ sends Verdier squares to pullback squares if and only if it sends them to pushout squares.
\end{proof}

\begin{remark}
\label{remark:KK-not-convolving}%
As with $\KGW$ and $\KL$ (see Remarks \reffour{remark:not-convolving} and \reffour{remark:KL-not-convolving}), we do not now whether $\KK$ is idempotent, nor whether every functor $G\colon \Catx \to \Spa$ admits a universal Karoubi-localising approximation $G \to \mathbb{G}$, but as with the previous examples we do know that such an approximation cannot be obtained by Day convolving with $\KK$, since by Remark~\reffour{rem:contensor-does-not-preserve-KS} functors of the form $\KK(\Funx(\C,-))$ are generally not Karoubi-localising, e.g., for $\C = \Dperf(\QQ)$.
As pointed to us by Maxime Ramzi, the same example also shows that Day convolving with $\KK$ cannot localise onto the full subcategory of motivic functors $\Catp \to \Spa$, that is, functors which invert motivic equivalences, or equivalently, the functors which factor through $\Mot$, where we recall from~\cite{ramzi-motivic} that the functor $U\colon \Catx \to \Mot$ exhibits $\Mot$ as the localisation of $\Catx$ by the motivic equivalences. Indeed, we claim that if a functor of the form $\KK(\Funx(\C,-))$ is motivic then $\Om^{\infty}\KK(\Funx(\C,-))$ is Karoubi-localising (which already fails for $\C = \Dperf(\QQ)$). 
To see this, note that if $\KK(\Funx(\C,-))$ is motivic then the space valued functor $\Om^{\infty}\KK(\Funx(\C,-))$ is also motivic. Now since $U$ is symmetric monoidal the operation $\G \mapsto \G(\C\otimes-)$ preserves motivic functors and so its left adjoint $\F \mapsto \F(\Funx(\C,-))$ preserves maps which are local equivalences with respect to motivic functors. Since $\core\Rightarrow \Om^{\infty}\KK = \Map_{\Mot}(\one,U(-))$ is one such local equivalence we conclude that $\core(\Funx(\C,-)) \Rightarrow \Om^{\infty}\KK(\Funx(\C,-))$ is also one, and hence by 2-out-of-3 so is the map
\[
\Om^{\infty}\KK(\Funx(\C,-)) \Rightarrow \Map_{\Mot}(U(\C),U(-)) .
\]
In particular, if the left hand side is motivic then this map would have to be an equivalence, hence the left hand side would be Karoubi-localising, as claimed.
\end{remark}

Using the results of the appendix for more general base $\infty$-operads 
one can deduce further results, for example:

\begin{corollary}
\label{corollary:multgw}%
Let $\E$ be a presentable $\infty$-category and consider the full inclusions
\[
\Funadd(\Catp,\E) \subseteq \Fun(\Catp,\E) \quad ,\quad \Funbord(\Catp,\E) \subseteq \Fun(\Catp,\E)
\]
\[
\Funadd(\Catx,\E) \subseteq \Fun(\Catx,\E) \quad \text{and} \quad \Fun^{\mathrm{binv}}(\Catp,\E) \subseteq \Fun(\Catp,\E).
\]
The the following holds:
\begin{enumerate}
\item
The bottom right inclusion is always reflective with reflection given by $\F \mapsto \Lspace \otimes_{\Day} \F$.
\item
If $\E$ is additive then the bottom left inclusion is reflective with reflection given by $\F \mapsto \Kspace \otimes_{\Day} \F$.
\item 
If $\E$ is stable then the top two inclusions are reflective with reflections given by $\F \mapsto \GW \otimes_{\Day} \F$ and $\F \mapsto \L \otimes_\Day \F$, respectively.
\end{enumerate}
In addition, if we replace $\E$ with a presentably $\cO$-monoidal $\infty$-category for some $\infty$-operad $\cO$ (keeping the additivity or stability assumptions in the relevant cases) then the above reflection functors yield multiplicative localisations on the level of Day convolution $\infty$-operads. When $\cO$ admits a unit object these multiplicative localisations are furthermore smashing.
\end{corollary}

\begin{proof}
Since $\E$ is presentable it carries a canonical colimit preserving action of $\Sps$, and if $\E$ is additive (resp. stable) this action factors through a colimit preserving action of $\Grp_{\Einf}$ (resp. $\Sp$). We may then encode this action via presentably $\MCom$-monoidal $\infty$-category $\M^{\otimes} \to \MCom$, %
to which we can apply Proposition~\reffour{proposition:smashing-day}.
For the last part apply Proposition~\reffour{proposition:multiplicative-day} for the multiplicative statement, given that under the assumptions all functors in question admit local approximations by the first part of the corollary.
As for the smashing statement, let us consider for simplicity the case of $\GW$. Let $\cO^{\otimes}$ be an $\infty$-operad with a unit object $\one_{\cO}$ and $\D^{\otimes} \to \cO^{\otimes}$ a presentably $\cO$-monoidal $\infty$-category, so that $\D_{\one_{\cO}}$ acquires the structure of a presentably symmetric monoidal $\infty$-category. Then post-composition with both the unit functor $u\colon \Sp \to \D_{\one_{\cO}}$ and its right adjoint preserves additive functors and hence $u\GW$ as the initial additive functor under the Day unit $u\Sig^{\infty}_+\Poinc\colon \Catp \to \D_{\one_{\cO}}$. In addition, if $a \in \cO$ is any object then $\D_a$ is a stable presentable $\infty$-category and the action of $\Sp$ on $\D_a$ factors through an action of $\D_{\one_{\cO}}$, 
so that the operations $u\GW\otimes_{\Day}(-)$ and $\GW\otimes_{\Day} (-)$ on $\D_a$ coincide. The desired result now follows from Proposition~\reffour{proposition:smashing-day}. The argument in the case of $\L,\Lspace$ and $\Kspace$ are the same.
\end{proof}

Having established the multiplicative properties of the K, GW, and L-functors, let us finally spell out their interaction. To this end, recall from Corollary~\refone{corollary:hyp-is-adjoint} the adjunctions  
\begin{equation}
\label{equation:two-sided}%
\begin{tikzcd}
	\Catx \ar[r, bend right,"\Hyp"'] \ar[r, bend left,"\Hyp"] & \Catp \ar[l, "\mathrm{fgt}"']
\end{tikzcd}
\end{equation}
exhibiting $\Hyp$ as both left and right adjoint to the forgetful functor $\fgt$. Since the functors $\Hyp$ and $\fgt$ preserves zero objects and split (Poincaré-)Verdier squares, restriction along $\Hyp$ and $\fgt$ is an operation which preserves additive functors. We then obtain an induced two sided adjunction 
\[
\begin{tikzcd}
	\Funadd(\Catp,\Spa) \ar[r, "\Hyp^*"]  & \Funadd(\Catx,\Spa) \ar[l, bend right,"\fgt^*"'] \ar[l, bend left,"\fgt^*"]
\end{tikzcd}
\]
on the level of additive functors.
In Proposition~\reftwo{proposition:bordrecoll}, it was shown that this two-sided adjunction yields a stable recollement
\[
\begin{tikzcd}
[column sep=7ex]
\Funbord(\Catp,\Spa) \ar[r] & 
\Funadd(\Catp,\Spa) 
\ar[r] 
\ar[l,bend left=25,shift left=1.5ex,start anchor=west,end anchor=east, "(-)^{\rbord}"] 
\ar[l,bend right=25,shift right=1.5ex,start anchor=west,end anchor=east,"(-)^{\bord}"'] 
\ar[l,phantom,shift left=1.2ex,start anchor=west,end anchor=east,"\perp"] 
\ar[l,phantom,shift right=1.3ex,start anchor=west,end anchor=east,"\perp"]
& \Funadd(\Catx,\Spa)^\hC \ .
\ar[l,bend left=25,shift left=1.5ex,start anchor=west,end anchor=east,"{(- \circ \fgt)^\hC}"] 
\ar[l,bend right=25,shift right=1.5ex,start anchor=west,end anchor=east,"{(- \circ \fgt)_\hC}"']
\ar[l,phantom,shift left=1.2ex,start anchor=west,end anchor=east,"\perp"] 
\ar[l,phantom,shift right=1.3ex,start anchor=west,end anchor=east,"\perp"]
\end{tikzcd}
\]
In~\cite{shah2021recollements}*{Definition 2.20} Shah defines the notion of a multiplicative stable recollement with respect to a symmetric monoidal structure on the ambient $\infty$-category, by requiring that the left Bousfield localisations to the two constituents of the recollement are compatible with said structure. 
His definition readily extends to the case where the ambient $\infty$-category is not a symmetric monoidal $\infty$-category but only an $\infty$-operad, by requiring these left Bousfield localisations to be multiplicative in the sense of Definition~\reffour{definition:multiplicative-localisation}. 
Returning to the above recollement, we may view the ambient term as the $\infty$-category of colours underlying the Day convolution full suboperad $\Funadd(\Catp,\Spa)^{\Day} \subseteq \Fun(\Catp,\Spa)^{\Day}$. The two constituents of the recollement can then be viewed as full suboperads of $\Funadd(\Catp,\Spa)^{\Day}$. For the left constituent this just gives the Day convolution full suboperad $\Funbord(\Catp,\Spa)^{\Day} \subseteq \Fun(\Catp,\Sp)^{\Day}$. As for the right constituent, we have the following:

\begin{lemma}
\label{lemma:sanity-check}%
The full suboperad of $\Funadd(\Catp,\Spa)^{\Day}$ spanned by the image of the fully-faithful right adjoint 
\[
(- \circ \fgt)^\hC \colon \Funadd(\Catx,\Spa)^{\hC} \to \Funadd(\Catp,\Spa)
\]
naturally coincides with the homotopy fixed points $(\Funadd(\Catx,\Spa)^{\Day})^{\hC}$ of the Day convolution full suboperad $\Funadd(\Catx,\Spa)^{\Day} \subseteq \Fun(\Catx,\Spa)^{\Day}$. 
\end{lemma}
\begin{proof}
Recall that the forgetful functor $\fgt\colon \Catp \to \Catx$ is symmetric monoidal and $\Ct$-equivariant, where $\Ct$ acts trivially on $\Catp$. Pre-composition with $\fgt$ hence determines a $\Ct$-equivariant map of $\infty$-operads
\[
(-\circ \fgt)\colon \Fun(\Catx,\Spa)^{\Day} \to \Fun(\Catp,\Spa)^{\Day} ,
\]
with $\Ct$ acting trivially on the target.
Passing to homotopy fixed points and using that $\Fun(\Catp,-)^{\Day}$ is a right adjoint and hence limit preserving we obtain a composite map of $\infty$-operads
\[
(\Fun(\Catx,\Spa)^{\Day})^{\hC} \to (\Fun(\Catp,\Spa)^{\Day})^{\hC}  = \Fun(\Catp,\Spa^{\B\Ct})^{\Day} \xrightarrow{(-)^{\hC}\circ (-)} \Fun(\Catp,\Spa)^{\Day} ,
\]
where in the last step we composed with the $\infty$-operad map induced on Day convolution $\infty$-operads by the lax symmetric monoidal functor $(-)^{\hC}\colon \Spa^{\B\Ct} \to \Spa$ (viewed as an $\infty$-operad map). Now the full suboperad inclusion $\Funadd(\Catx,\Spa)^{\Day} \subseteq \Fun(\Catx,\Spa)^{\Day}$ induces a full suboperad inclusion after taking homotopy fixed points, and the $\infty$-operad map constructed above sends this full suboperad to the full suboperad $\Funadd(\Catp,\Spa)$, thus giving an $\infty$-operad map
\[
(-\circ \fgt)^{\hC}\colon (\Funadd(\Catx,\Spa)^{\Day})^{\hC} \to \Funadd(\Catp,\Spa),
\]
refining the fully-faithful functor $(-\circ \fgt)^{\hC}\colon \Funadd(\Catx,\Spa)^{\hC} \hrar \Funadd(\Catp,\Spa)$ on the level of colours. To finish the proof it will hence suffice to show that $(-\circ \fgt)^{\hC}$ is also fully-faithful as a map of $\infty$-operads, that is, for every sequence of $\Ct$-equivariant additive functors $\F_1,...,\F_n,\G\colon \Catx \to \Spa$ (with respect to the $\mop$ action on $\Catx$ and the trivial action on $\Spa$) the induced map
\[
\Nat(\otimes_n \circ (\F_1 \times ... \times \F_n),\G \circ \otimes_n)^{\hC} \to \Nat([\otimes_n \circ ([\F_1 \circ \fgt]^{\hC} \times ... \times [\F_n \circ \fgt]^{\hC}),[\G \circ \fgt]^{\hC} \circ \otimes_n)
\]
is an equivalence, where we write $\otimes_n$ for the $n$-fold tensor product of either $\Catx,\Catp$ or $\Spa$. 
Now since $\fgt$ is symmetric monoidal we can rewrite this map as
\[
\Nat(\otimes_n \circ (\F_1 \times ... \times \F_n),\G \circ \otimes_n)^{\hC} \to \Nat([\otimes_n \circ ([\F_1 \circ \fgt]^{\hC} \times ... \times [\F_n \circ \fgt]^{\hC}),[\G \circ \otimes_n\circ (\fgt \times ... \times \fgt)]^{\hC}) .
\]
Next, since $\Hyp$ is in particular right adjoint to $\fgt$, the functor 
\[
[(-) \circ (\fgt \times ... \times \fgt)]^{\hC}\colon \Fun(\Catx \times ... \times \Catx,\Spa)^{\hC} \to \Fun(\Catp \times ... \times \Catp,\Spa)
\]
has a left adjoint given by
\[
(-) \circ (\Hyp \times ... \times \Hyp)\colon \Fun(\Catp \times ... \times \Catp,\Spa) \to \Fun(\Catx \times ... \times \Catx,\Spa)^{\hC}.
\]
The map that we want to show is an equivalence can hence be identified via adjunction with the map
\[
\Nat(\otimes_n \circ (\F_1 \times ... \times \F_n),\G \circ \otimes_n)^{\hC} \to \Nat([\otimes_n \circ ([\F_1 \circ \fgt \circ \Hyp]^{\hC} \times ... \times [\F_n \circ \fgt \circ \Hyp]^{\hC}),\G \circ \otimes_n)^{\hC}
\]
induced by pre-composition with the composite 
\[
\otimes_n\circ ([\F_1 \circ \fgt\circ \Hyp]^{\hC},...,[\F_n \circ \fgt \circ \Hyp]^{\hC}) \Rightarrow [\otimes_n\circ (\F_1 \circ \fgt\circ \Hyp,...,\F_n \circ \fgt \circ \Hyp)]^{\hC} \Rightarrow   \otimes_n\circ (\F_1,...,\F_n) ,
\]
where the second map is the counit of the above adjunction. But this composite is simply the map induced after $\otimes_n \circ$ by the tuple of counits $[\F_i \circ \fgt \circ \Hyp]^{\hC} \Rightarrow \F_i$
of the adjunction $(-)\circ \Hyp\dashv [(-)\circ \fgt]^{\hC}$ for $i=1,...,n$. The desired result now follows from the fact that the counit $[\F \circ \fgt \circ \Hyp]^{\hC} \Rightarrow \F$ is an equivalence whenever $\F \colon \Catx \to \Spa$ is an additive $\Ct$-equivariant functor.
\end{proof}

\begin{proposition}
The above stable recollement is multiplicative with respect to the Day convolution $\infty$-operad structure on $\Funadd(\Catp,\Spa)$, yielding two left Bousfield localisations of $\infty$-operads
\[
\begin{tikzcd}
[column sep=7ex]
\Funbord(\Catp,\Spa)^{\Day} \ar[r] & 
\Funadd(\Catp,\Spa)^{\Day} 
\ar[r]  
\ar[l,bend right=25,shift right=1.5ex,start anchor=west,end anchor=east,"(-)^{\bord}"'] 
\ar[l,phantom,shift right=1.3ex,start anchor=west,end anchor=east,"\perp"]
& (\Funadd(\Catx,\Spa)^{\Day})^\hC \ .
\ar[l,bend left=25,shift left=1.5ex,start anchor=west,end anchor=east,"{(- \circ \fgt)^\hC}"] 
\ar[l,phantom,shift left=1.2ex,start anchor=west,end anchor=east,"\perp"] 
\end{tikzcd}
\]
\end{proposition}

\begin{proof}
We need to verify that the left localisations 
\begin{enumerate}
\item[-] $(-)^{\brd}\colon \Funadd(\Catp,\Spa) \to \Funbord(\Catp,\Spa)$ and,
\item[-] $\Hyp^*\colon \Funadd(\Catp,\Spa) \to \Funadd(\Catx,\Spa)^{\hC}$
\end{enumerate}
are multiplicative with respect to Day convolution. 
For $(-)^{\brd}$ this follows from Corollary~\reffour{corollary:lmonoidal}\;\reffouritem{item:L-smashing} and Remark~\reffour{remark:full-suboperad}.
For $\Hyp^*$, it will suffice by Lemma~\reffour{lemma:canonical-colimit} to show that $\Hyp^*$-equivalences are closed under Day convolution with (small) functors of the form $E \otimes \SS[\Map_{\Catp}((\C,\QF),-)]\colon \Catp \to \Spa$. Concretely, we need to for $\F\Rightarrow \F'$ a natural transformation of functors $\Catp \to \Spa$ such that $\F(\Hyp(\D)) \to \F'(\Hyp(\D))$ is an equivalence for every $\D \in \Catx$ also the map $E \otimes \F(\Funx((\C,\QF),\Hyp(\D))) \to E \otimes \F'(\Funx((\C,\QF),\Hyp(\D)))$ is an equivalence, which follows from the natural equivalence 
$\Funx((\C,\QF),\Hyp(-)) \simeq \Hyp(\Funx(\C,-))$.
\end{proof}

Applying~\cite{shah2021recollements}*{Corollary 2.29}, which applies verbatim to $\infty$-operads rather than symmetric monoidal $\infty$-categories, we now obtain the following symmetric monoidal refinement of the structure square for additive functors:

\begin{corollary}
\label{corollary:generic-fracture-monoidal}%
Associating to a functor $\F \colon \Catp \to \Spa$ 
its fundamental fibre square
\[
\begin{tikzcd}
	\F \ar[r] \ar[d] & \F^\brd \ar[d]   \\
	(\F^\hyp)^{\hC} \ar[r] & (\F^\hyp)^{\tC} 
\end{tikzcd}
\]
naturally refines to a map of $\infty$-operads $\Fun(\Catp,\Spa)^\Day \to \Fun(\Delta^1 \times \Delta^1,\Fun(\Catp,\Spa)^\Day)$. In particular, the fundamental fibre square of $\F$ naturally refines to a square of lax symmetric monoidal functors once $\F$ refines to a lax symmetric monoidal functor.
\end{corollary}

\begin{corollary}
\label{corollary:kgw-kl-monoidal}%
\label{corollary:GW-symmetric-monoidal}%
The entire cube
\[
\begin{tikzcd}
[column sep=0.5em, row sep=0.7em]
  & \KGW\ar[rr]\ar[dd] && \KL\ar[dd]\\
\GW\ar[ru]\ar[rr, crossing over]\ar[dd] && \L\ar[ur]\\
  & \KK^\hC\ar[rr] && \KK^\tC\\
\K^\hC\ar[rr] \ar[ru] && \K^\tC\ar[uu,<-, crossing over] \ar[ru]\\
\end{tikzcd}
\]
canonically refines to a diagram of lax symmetric monoidal functors $\Catp \rightarrow \Sp$. Here, the lax symmetric monoidals on the functors in the top face are those of Corollaries~\reffour{corollary:kl-monoidal}, 
\reffour{corollary:lmonoidal}\;\reffouritem{item:L-universal},
\reffour{corollary:kgw-multiplicative} and 
\reffour{corollary:GW-symmetric-monoidal}\;\reffouritem{item:GW-universal}, 
while those on the functors in the bottom face are induced by the lax symmetric monoidal structures of Corollary~\reffour{corollary:K-universal}. %
\end{corollary}
\begin{proof}
Note first that the lax symmetric monoidal structures induced on $\L = \GW^\mathrm{bord}$ and $\KL = \KGW^\mathrm{bord}$ by Corollary~\reffour{corollary:generic-fracture-monoidal} 
necessarily agree with those of Corollaries \reffour{corollary:lmonoidal}
and \reffour{corollary:kl-monoidal} respectively, on account of the uniqueness property of the latter. %
Secondly, the map $\GW \Rightarrow \KGW$ canonically refines to one of lax symmetric monoidal functors, since $\GW$ is initial both as an additive functor under $\SS[\Poinc]$ (Corollary~\reftwo{corollary:univ-GW-spectrum}) and as an additive lax symmetric monoidal functor (Corollary~\reffour{corollary:GW-universal-multiplicative}), while at the same time $\SS[\Poinc]$ is the initial lax symmetric monoidal functor $\Catp \to \Sp$ and the composite $\SS[\Poinc] \Rightarrow \GW \Rightarrow \KGW$ is a map of lax symmetric monoidal functors (Corollary~\reffour{corollary:kgw-multiplicative}). We may then apply Corollary~\reffour{corollary:generic-fracture-monoidal} to the lax monoidal transformation $\GW \to \KGW$.
\end{proof}
Using the fact that the pullback diagram
\[
\begin{tikzcd}
\Sp^\gCt \ar[rr,"{(-)^{\varphi\Ct} \Rightarrow (-)^\tC}"] \ar[d,"\mathrm{fgt}"] && \Ar(\Sp) \ar[d,"t"]\\
\Sp^\hC \ar[rr,"{(-)^\tC}"] && \Sp,
\end{tikzcd}
\]
describing genuine $\Ct$-spectra as Tate-squares as above is one of $\infty$-operads, Corollary~\reffour{corollary:kgw-kl-monoidal} can equivalently be stated as:

\begin{corollary}
\label{corollary:real-is-lax-symmetric-monoidal}%
The natural transformation $\KR \Rightarrow \KKR \colon \Catp \rightarrow \Sp^\gCt$ canonically refines to one of lax symmetric monoidal functors.
\end{corollary}

This result in particular provides an $\Einf$-structure to the genuine $\Ct$-spectra $\KR(\C,\QF)$ and $\KKR(\C,\QF)$ whenever $(\C,\QF)$ is a symmetric monoidal Poincar\'e $\infty$-category. In this paper, we will not discuss the more subtle structure of a normed $\Einf$-algebra structure (sometimes also called a $\Ct$-$\Einf$-structure) on $\KR(\C,\QF)$.

We finish this subsection by explaining how to deduce from the above results multiplicative statements about $\GW$ (or any of the other functors considered above) of rings with form parameters. More precisely, let $(\Ringp)^{\otimes}$ be the symmetric monoidal $\infty$-category sitting in the pullback square
\[
\begin{tikzcd}[column sep=40pt]
(\Ringp)^{\otimes}\ar[r]\ar[d] & (\Catpi)^{\otimes}\ar[d] \\
\Alg_{\Eone}(\Sp)^{\otimes} \ar[r,"\Modp{-}"] & (\Catxi)^{\otimes}
\end{tikzcd}
\]
of symmetric monoidal functors, where $\Catxi$ is the symmetric monoidal localisation of $\Catx$ by the Karoubi equivalences, and similarly for $\Catpi$. A similar construction using $\Alg_{\Einf}$ in place of $\Alg_{\Eone}$ was considered in~\cite{poincare-rings}, who coined the term ``Poincaré rings''.
By the classification of Poincaré structures on module $\infty$-categories (see \S\refone{subsection:classification-in-modules}) we may identify objects in $\Ringp$ with tuples $(A,M,N,\alp\colon N \to M^{\tC})$ where $A$ is an $\Eone$-ring spectrum and $\M := (M,N,\alp)$ is a module with genuine involution over $A$, or, equivalently, an $\N A$-module in the $\infty$-category $\Sp^{\gCt}$ of genuine $\Ct$-spectra, where $\N A \in \CAlg(\Sp^{\g\Ct})$ is the Hill-Hopkins-Ravanel norm of $A$ (whose underlying spectrum is $A \otimes A$ and whose geometric fixed points is $A$), see Remark~\refone{remark:NormModules}. %
The symmetric monoidal product on $(\Ringp)^{\otimes}$ is then given by $(A,\M) \otimes (A',\M') = (A \otimes A',\M \otimes \M')$, where $\M \otimes \M'$ is the tensor product in genuine $\Ct$-spectra, equipped with its induced $A \otimes A'$-module structure, see Example~\refone{example:tensor-of-rings}. Now recall that the $\infty$-category $\Sp^{\gCt}$ admits a t-structure, which we will refer to as the Mackey t-structure, where an $E \in \Sp^{\gCt}$ is Mackey connective (resp. Mackey coconnective) if and only if the underlying spectrum $\bar{E} \in \Sp$ and the genuine fixed points $E^{\Ct} \in \Sp$ are connective (resp. coconnective). In other words, this is the pointwise t-structure after identifying $\Sp^{\gCt}$ with spectral Mackey functors, and its heart consists exactly of the $\Ab$-valued Mackey functors. It is also convenient to note that connectivity (though not coconnectivity) can be tested on underlying objects and geometric fixed points. This t-structure is multiplicative in the sense that the monoidal unit is Mackey connective and the collection of Mackey connective objects is closed under tensor product. In addition, if $A$ is a connective ring spectrum then $\N A$ is Mackey connective.
Finally, let us say that a module with involution over $A$ is projective if it is projective 
as an $A$-modules (that is, is a retract of $A^n$ for some $n$) with respect to either of its $A$-module structures, and that an $\N A$-module is projective if the underlying module with involution is projective in the above sense. We recall from Remark~\refone{remark:form-parameters-and-modules-with-genuine-involution} that when $A$ is a discrete ring the data of a Mackey discrete $\N A$-module corresponds exactly to that of a form parameter over $A$. We say that a form parameter $\lambda$ on $A$ is Poincaré if the associated module with involution $M$ over $A$ is projective and invertible, which is equivalent to saying that 
the hermitian structure $\QF^{\g\lambda}_M$ on $\Dperf(A)$ is Poincaré and its duality preserves $\Proj(R)$ (see~\cite{comparison}*{Example 3.2.4(4)} for a discussion, though note the notation clash with loc.\ cit.\ where the term invertible is used for what we call here Poincaré). The $\infty$-category $\FPMp$ of Poincaré form parameters can then be identified with the full subcategory of $\Ringp$ spanned by those $(A,\M)$ such that $A$ is discrete and $\M$ is projective and Mackey discrete. It sits inside the intermediate full subcategory 
\[
\FPMp \subseteq \Ringp_{\geq 0} \subseteq \Ringp
\]
where $\Ringp_{\geq 0}$ is spanned by those $(A,\M)$ such that $A$ is connective and $\M$ is projective and Mackey connective.
In particular, since the Mackey t-structure is multiplicative we have that $\Ringp_{\geq 0}$ contains the unit $(\SS,\SS^{\uni})$ of $\Ringp$ and is closed under tensor products, so that it inherits from $\Ringp$ a symmetric monoidal structure.

\begin{proposition}%
\label{proposition:sym_mon_structure_on_cat_of_form_parameters}%
The full subcategory $\FPMp \subseteq \Ringp_{\geq 0}$ is a symmetric monoidal localisation of $\Ringp_{\geq 0}$ with reflection functor $\Ringp_{\geq 0} \to \FPMp$ given by $(A,\M) \mapsto (\pi_0(A),\pi^{\mack}_0(\M))$, where $\pi^{\mack}_0$ denotes the $\pi_0$ functor with respect to the Mackey t-structure. In particular, $\FPMp$ inherits a symmeric monoidal structure making its inclusion in $\Ringp_{\geq 0}$ lax symmetric monoidal. The unit of this symmetric monoidal structure is $(\ZZ,\uline{\ZZ}^{\burn})$, where $\uline{\ZZ}^{\burn}$ is the Burside $\N\ZZ$-module, whose underlying $(\ZZ \otimes_{\SS} \ZZ)$-module is $\ZZ$ and whose genuine fixed points object is the Burnside ring of $\ZZ$, see Example~\refthree{example:dotto-ogle}. 
\end{proposition}

Concretely, the induced tensor product in $\FPMp$ is given by 
\[
(R,\M)\otimes (R',\M) = (\pi_0(R \otimes_{\SS} R'),\pi^{\mack}_0(\M\otimes\M')) = (R\otimes_{\ZZ} R',\M\otimes^{\mack} \M')
\]
where $\otimes^{\mack}$ is the induced tensor product on the heart of the Mackey t-structure, also known as the box product of $\Ab$-valued Mackey functors, see Definition~\reffour{definition:box-mackey} for an explicit algebraic description. We note that the inclusion $\Catpi \subseteq \Catp$ as idempotent complete Poincaré $\infty$-categories is lax symmetric monoidal, and so we get that the composite functor
\[
\FPMp \to \Catp \quad\quad (R,\lambda) \mapsto (\Dperf(R),\QF^{\g\lambda})
\]
is lax symmetric monoidal.
In particular, $\GW$ and $\L$ (and similarly all the functors considered above) promote to lax symmetric monoidal functors on $\FPMp$, and for any $(R, \lambda)$ in $\FPMp$, its $\GW$ and $\L$-theories are canonically modules over $\GW^{\burn}(\ZZ)$ and $\L^{\burn}(\ZZ)$ respectively.

\begin{proof}[Proof of Proposition~\reffour{proposition:sym_mon_structure_on_cat_of_form_parameters}]
Let $(A,\M) \in \Ringp_{\geq 0}$, so that $\pi_0(A)$ inherits the structure of a ring. We then observe that the induced map $\N A \to \N\pi_0(A)$ in $\Sp^{\gCt}$ has a Mackey 1-connective fibre, and so induces an isomorphism on $\pi_0^{\mack}$. Now since the Mackey t-structure is multiplicative we may refine $\pi_0^{\mack}$ to a left adjoint to the inclusion of Mackey discrete algebra objects inside Mackey connective algebra objects. It hence follows that the algebra map $\N A \to \pi_0^{\mack}(\N A)$ factors uniquely through an algebra map $\N\pi_0(A) \to \pi^{\mack}_0(\N A)$, through which we may view $\pi_0^{\mack}(\M)$ as an $\N\pi_0(A)$-module, or, equivalently, as a module with genuine involution over $\pi_0(A)$. Now since the Mackey t-structure is pointwise under the Mackey functor description we have that the underlying $\pi_0(A)\otimes_{\SS} \pi_0(A)$-module of $\pi_0^{\mack}(\M)$ is given by $\pi_0(M)$, where $M$ is the $A \otimes A$-module underlying $\M$. Since $M$ is projective as an $A$-module we have that $\pi_0(M)$ is projective as a $\pi_0(A)$-module.
We now claim that $\pi_0(M)$ is invertible as a module with involution over $\pi_0(A)$, that is, that the induced map
\[
\pi_0(A) \to \map_{\pi_0(A)}(\pi_0(M),\pi_0(M))
\]
is an equivalence. %
For this, note that the left hand side is a-priori connective since $\pi_0(M)$ is projective. Refining $\pi_0$ to a left adjoint to the inclusion $\Mod_{\pi_0(A)}(\Ab) \to \Mod_{A}(\Sp)$ we may then replace the right hand side with $\map_{A}(M,\pi_0(M))$. We then note that $\map_A(M,-)$ commutes with truncation since $M$ is projective and connective, so that the desired claim follows from the fact that the original $M$ was invertible as a module with involution over $A$.
We now claim that the resulting map
\[
(A,\M) \to (\pi_0(A),\pi_0(\M))
\]
is Poincaré, that is, that the induced map $M\otimes_A \pi_0(A) \to \pi_0(M)$ is an equivalence. Indeed, the fact that the standard t-structure is multiplicative implies that this map is an isomorphism on $\pi_0$, where on the other hand the both sides are discrete since $M$ is projective as an $A$-module.

Finally, we now show that $(\pi_0(A),\pi_0(\M))$ is initial in $\FPMp_{(A,\M)/}$.  Unwinding all definitions and using $\pi_0$ and $\pi_0^{\mack}$ are left adjoints to the inclusion of discrete objects inside connective one, it will suffice to show that if $(R,\mathcal{N})$ is an object of $\FPMp$ equipped with a ring map $\pi_0(A) \to R$ and an $\N\pi_0(A)$-module map $\pi_0^{\mack}(\M) \to \mathcal{N}$ then the map $(\pi_0(A),\pi_0(\M)) \to (R,\mathcal{N})$ is Poincaré as soon as the composite $(A,\M) \to (\pi_0(A),\pi^{\mack}_0(\M)) \to (R,\mathcal{N})$ is Poincaré. In other words, the induced map $\pi_0(M)\otimes_{\pi_0(A)} R \to N$ is an equivalence as soon as the induced map $M \otimes_{A} R \to N$ is so, where $N$ is the underlying module with involution over $R$. Indeed, this is simply because the map $M\otimes_A \pi_0(A) \to \pi_0(M)$ is an equivalence, as we showed above.

We have thus shown that the inclusion $\FPMp \subseteq \Ringp_{\geq 0}$ is reflective with reflection functor $F(A,M) = (\pi_0(A),\pi^{\mack}_0(\M))$. To show that the resulting localisation is multiplicative we need to show the collection of $F$-equivalences are closed under taking tensor products with fixed objects. Now by construction the $F$-equivalence are exactly the maps $(A,\M) \to (B,\mathcal{N})$ such that the ring spectrum map $A \to B$ is a $\pi_0$-isomorphism and the $\N A$-module map $\M \to \mathcal{N}$ is a $\pi_0^{\mack}$-isomorphism. This collection of maps is closed under tensoring with fixed objects simply because the standard t-structure on $\Sp$ and the Mackey t-structure on $\Sp^{\g\Ct}$ are both multplicative.
\end{proof}

\subsection{Examples of multiplicative structures on Grothendieck-Witt- and L-spectra}
\label{subsection:examples}%

In this section we will collect some examples of naturally occurring symmetric monoidal Poincaré $\infty$-categories (and modules over them) to which the multiplicative invariants summarised in the cube of Corollary~\reffour{corollary:kgw-kl-monoidal} can be applied, yielding cubes of commutative ring spectra (and modules over such), and we will describe the multiplicative structure on several (previously additively computed) Grothendieck-Witt and L-spectra. We begin with categorical examples.

\begin{example}
\label{example:initial-terminal}%
If $(A,\sig)$ is an $\Einf$-ring spectrum with involution, the $\infty$-category $\Modp{A}$ of compact $A$-modules carries a symmetric monoidal structure and $\Dual_{\sig}$ is a symmetric monoidal duality on $\Modp{A}$. The associated Tate construction $A^{\tC}$ then inherits the structure of an $\Einf$-$A$-algebra spectrum via the twisted Tate-Frobenius and symmetric monoidal Poincaré refinements of $(\Dperf(A),\Dual_{\sig})$ are classified by $\Einf$-$A$-algebras $N$ (representing the linear part of the associated Poincaré refinement) equipped with an $\Einf$-$A$-algebra map to $A^{\tC}$. Among those, the initial such algebra is given by the unit map $A \to A^{\tC}$, and corresponds to the Tate Poincaré structure $\QF^{\tate}_{\sig}$, while the terminal such algebra is given by the identity map $A^{\tC} \to A^{\tC}$, and corresponds to the symmetric Poincaré structure $\QF^{\sym}_{\sig}$. For $A=\SS$ the Tate Poincaré structure on $\Modp{\SS} = \Spaf$ coincides with the universal one of Example~\reffour{example:universal}.
In particular, for $(A,\sig)$ as above, we have a (pullback) squares of $\Einf$-ring spectra 
\[
\begin{tikzcd}
\GW^{\tate}(A,\sig) \ar[r]\ar[d] & \L^{\tate}(A,\sig) \ar[d] &\KGW^{\tate}(A,\sig) \ar[r]\ar[d] & \KL^{\tate}(A,\sig) \ar[d] \\
\GW^{\sym}(A,\sig)\ar[r] & \L^{\sym}(A,\sig) & \KGW^{\sym}(A,\sig)\ar[r] & \KL^{\sym}(A,\sig)..
\end{tikzcd}
\]
In addition, if $B$ is an $\Eone$-$A$-algebra equipped with a module with involution $M$ whose $\Ct$-action is $A$-sesquilinear then $(\Modp{B},\QF^{\sym}_M)$ is a module over $(\Modp{A},\QF^{\sym}_{\sig})$ and any Poincaré refinement $(\Modp{B},\QF_M)$ of $\Dual_{M}$ is a module over $(\Modp{A},\QF^{\tate}_{\sig})$. Consequently, $\GW^\s(\Modp{B},\Dual_M))$ and $\L^\s(\Modp{B},\Dual_M))$ are modules over $\GW^s(A,\sig)$ and $\L^\s(A,\sig)$, respectively, and for any Poincar\'e refinement $\QF_M$ of $\Dual_M$, we find that $\GW(\Modp{B},\QF_M)$ and $\L(\Modp{B},\QF_M)$ are modules over $\GW^\tate(A,\sig)$ and $\L^\tate(A,\sig)$, respectively. The same holds for $\KGW$ and $\KL$ in place of $\GW$ and $\L$.
\end{example}

\begin{example}
\label{example:geq-zero}%
If $(A,\sig)$ is a connective $\Einf$-ring spectrum with involution, the Poincaré structure $\QF^{\geq 0}_{\sig}$ on $\Modp{A}$ is symmetric monoidal: its associated $\Einf$-$A$-algebra is the connective cover $\tau_{\geq 0} A^{\tC}$ of $A^{\tC}$, which is indeed an $\Einf$-$A$-algebra since the t-structure on $\Mod_A$ is compatible with its symmetric monoidal structure. In fact, $\QF^{\geq0}_{\sig}$ constitutes the terminal symmetric monoidal Poincaré refinement of $(\Modp{A},\Dual_{\sig})$ with connective linear part. If $B$ is a connective $\Eone$-$A$-algebra equipped with an $A$-sesquilinear module with involution $M$ then $(\Modp{B},\QF^{\geq m}_M)$ is a module over $(\Modp{A},\QF^{\geq 0}_{\sig})$ for all $m\in \ZZ \cup \{ \pm \infty \}$. Consequently, $\GW(\Modp{B},\QF^{\geq m}_M)$ and $\L(\Modp{B},\QF^{\geq m}_M)$ are modules over $\GW(\Modp{A},\QF^{\geq0}_\sig)$ and $\L(\Modp{A},\QF^{\geq0}_\sig)$, respectively. The same holds for $\KGW$ and $\KL$ in place of $\GW$ and $\L$.
\end{example}

\begin{example}
\label{example:ideal}%
If $(A,\sig)$ is an $\Einf$-ring spectrum with involution, then $(\Modp{A},\QF^\qdr_\sig)$ is not a symmetric monoidal Poincar\'e $\infty$-category, but it is a module over $(\Modp{A},\QF^\sym_\sig)$. In fact, the canonical functor $(\Modp{A},\QF^\qdr_\sig) \to (\Modp{A},\QF^\sym_\sig)$ is a Smith ideal, that is, a commutative algebra in $\Fun([1],\Catp)$ with respect to the Day convolution symmetric monoidal structure for the cartesian product on $[1]$, see \S\reffour{subsection:normal-L} for details. Consequently, $\GW^\qdr(A,\sig)$ and $\L^\qdr(A,\sig)$ are Smith ideals in $\GW^\s(A,\sig)$ and $\L^\s(A,\sig)$, respectively, whose quotients agree (by the fundamental fibre sequence) and are called symmetric normal $\L$-theory; 
see again \S\reffour{subsection:normal-L} for more on normal L-theory. If $A$ is connective and $m\geq 0$, then $(\Modp{A},\QF^{\geq m}_\sig) \to (\Modp{A},\QF^{\geq0}_\sig)$ is a Smith ideal, so $\GW^{\geq m}(A,\sig)$ and $\L^{\geq m}(A,\sig)$ are Smith ideals in $\GW^{\geq 0}(A,\sig)$ and $\L^{\geq0}(A,\sig)$, respectively.
The same holds for $\KGW$ and $\KL$ in place of $\GW$ and $\L$.
\end{example}

\begin{example}
\label{example:modules-over-genuine}%
If $(R,\sig)$ is a discrete commutative ring with involution then we may also consider the genuine Frobenius Poincar\'e structure $\QF^{\g\fr}_\sig$ of $\Dperf(R)$. It participates in the following sequence of symmetric monoidal refinements of $(\Dperf(R),\Dual_{\sig})$
\[
\QF^{\tate}_{\sig} \Rightarrow \QF^{\g\fr}_{\sig} \Rightarrow \QF^{\gs}_{\sig} \Rightarrow \QF^{\sym}_{\sig} ,
\]
where $\QF^{\tate}_{\sig}$ and $\QF^{\sym}_{\sig}$ are the initial and terminal refinements as in Example~\reffour{example:initial-terminal} and the genuine symmetric structure $\QF^{\gs}_{\sig} = \QF^{\geq 0}_{\sig}$ is the terminal refinement with connective linear part as in of Example~\reffour{example:geq-zero}.
The genuine Frobenius Poincaré structure $\QF^{\g\fr}_{\sig}$ on $\Dperf(R)$ is obtained from $\QF^{\tate}_{\sig}$ by applying the reflection functor of Proposition~\reffour{proposition:sym_mon_structure_on_cat_of_form_parameters}, and 
constitutes the initial symmetric monoidal %
Poincaré structure refining $\Dual_{\sig}$ whose associated module with involution $(R,N^{\g\fr} \to R^{\tC})$ is Mackey 0-truncated. %
The $\Einf$-$R$-algebra spectrum $N^{\g\fr}$ is connective with $\pi_n(N^{\g\fr}) = \pi_n(R^{\tC})$ for $n \geq 2$, $\pi_1(N^{\g\fr}) = 0$ and $\pi_0(N^{\g\fr}) = R$, so that $(R,N^{\g\fr})$ is also Mackey connective and $\QF^{\g\fr}_{\sigma}$ corresponds to the Frobenius form parameter of
Example~\reffour{example:form_parameter_over_Tate_Poincare}.
For $R=\ZZ$ the genuine Frobenius structure $\QF^{\g\fr}_{\ZZ}$ coincides with the Burnside Poincaré structure %
appearing in Example~\refthree{example:dotto-ogle}.

If $(R,\sig)$ is a discrete commutative ring with involution and $B$ is a discrete associative $R$-algebra equipped with an invertible module with involution $M$ such that the involution on $M$ is $R$-sesquilinear then $(\Dperf(B),\QF^{\gs}_M)$ is a module over $(\Dperf(B),\QF^{\gs}_{\sig})$, and if $\lambda$ is any form parameter over $B$ refining $M$ then $(\Dperf(B),\QF^{\g\lambda}_M)$ is a module over $(\Dperf(R),\QF^{\g\fr}_{\sig})$. In general, $(\Dperf(B),\QF^{\g\lambda}_M)$ is not a module over $(\Dperf(R),\QF^{\gs}_{\sig})$, for instance $(\Dperf(\ZZ),\QF^{\g\fr}_\ZZ)$ is not a module over $(\Dperf(\ZZ),\QF^\gs_\ZZ)$, see also Remark~\reffour{remark:criterion-genuine-symmetric}.

In particular, $\GW(\Dperf(B),\QF^{\g\lambda}_M)$ and $\L(\Dperf(B),\QF^{\g\lambda}_M)$ are modules over $\GW(\Dperf(R),\QF^{\g\fr}_\sig)$ and $\L(\Dperf(R),\QF^{\g\fr}_\sig)$, but in general not over $\GW(\Dperf(R),\QF^{\gs}_\sig)$ and $\L(\Dperf(R),\QF^{\gs}_\sig)$. For instance, $\GW^{\g\fr}(\ZZ)$ is not a module over $\GW^{\gs}(\ZZ)$, and similarly $\L^{\g\fr}(\ZZ)$ is not a module over $\L^{\gs}(\ZZ)$; indeed it it were a module, since $\L^{\g\fr}_0(\ZZ) \to \L^{\gs}_0(\ZZ)\cong \ZZ$ is surjective, we would obtain a splitting $\L^{\g\fr}(\ZZ) \simeq \L^\gs(\ZZ) \oplus M$ contradicting that $\L^{\g\fr}_1(\ZZ) = 0$ and $\L^{\gs}_1(\ZZ) = \ZZ/2$; a similar argument also shows the case of $\GW$.
\end{example}

Next, we aim to describe more explicitly the multiplicative structure on some GW- and L-spectra, focussing first on ring structures on homotopy groups.
\begin{example}
\label{example:universal}%
The most fundamental example of a symmetric monoidal Poincaré $\infty$-category is the unit $(\Spaf,\QF^{\uni})$ of $\Catp$ and we obtain $\Einf$-ring spectra $\GW^\uni(\SS)$ and $\L^\uni(\SS)$ over which all GW- and L-spectra are modules. By Proposition~\reftwo{proposition:pi0gwuniv}, we have
$\GW_0^\uni(\SS) = \ZZ[e,h]/(e^2-8e,he-8h,h^2-2h)$ and $\L_0^\uni(\SS) = \ZZ[e]/(e^2-8e)$.
\end{example}

\begin{example}
\label{example:L-ring-Z}%
Another fundamental commutative ring is $\ZZ$. Let us record here the multiplication on various forms of its L-groups, due to Mishchenko, Kervaire--Milnor, and Ranicki, see \cite{Ranickiyellowbook}*{Section 4.3.1} or Corollary~\refthree{corollary:compute-L-groups-global-field} and Proposition~\refthree{proposition:compute-quadratic-L-groups}. Namely, we have $\L^\s_*(\ZZ) = \ZZ[x^{\pm 1},e]/(2e,e^2)$ with $|x|=4$ and $|e|=1$ and 
\[
\L^\qdr_n(\ZZ) = \begin{cases} 8\ZZ & \text{ for } n\equiv 0 (4) \\ 0 & \text{ for } n \equiv 1,3 (4) \\ \ZZ/2 & \text{ for } n \equiv 2(4) \end{cases}.
\]
which is an $\L_*^\s(\ZZ)$-module in a unique fashion (the module multiplication by $e$ is trivial for degree reasons). To determine the genuine symmetric L-ring of $\ZZ$, recall from \paperthree that the ring map $\L^\gs_*(\ZZ) \to \L^\s_*(\ZZ)$ is an isomorphism for $*\geq -2$ and the $\L^\gs_*(\ZZ)$-module map $\L^\qdr_*(\ZZ) \to \L^\gs_*(\ZZ)$ is an isomorphism for $*\leq -3$. Since the $\L^\gs(\ZZ)$-module structure on $\L^\qdr(\ZZ)$ is restricted from its canonical $\L^\s(\ZZ)$-module structure, one can infer from the above the ring structure on $\L^\gs_*(\ZZ)$, see \cite{HLN}*{Remark 6.3} for explicit formulae.
\end{example}

In fact, the computation of the symmetric and quadratic L-groups of $\ZZ$ together with the canonical module structure explained above can be generalised as follows.

\begin{example}
\label{example:L-rings-dedekind}%
Let $D$ be a Dedekind domain whose fraction field $F$ is a global field of characteristic different from 2, e.g.\ that $D$ is a number ring. We have computed in Corollary~\refthree{corollary:compute-L-groups-global-field} that its symmetric L-groups are 4-periodic and given by
\[
\L^\s_n(D) = \begin{cases} \W^\s(D) & \text{ for } n \equiv 0 (4),\\
						\bigoplus_{\I} \ZZ/2 & \text{ for } n \equiv 1 (4),\\
						0 & \text{ for } n \equiv 2 (4), \\
						\Pic(D)/2 & \text{ for } n \equiv 3 (4), \end{cases}
\]
where $\I$ is the set of dyadic primes in $D$. In fact, for $\mathfrak{p} \in \I$ such a prime, the map $D \to D_\mathfrak{p}$ induces the associated projection on $\L_1^\s(-)$. There are now three types of basic products in the ring $\L_\ast^\s(D)$, described by $\W^\s(D)$-module structures on $\Pic(D)/2$ and $\bigoplus_\I \ZZ/2$, respectively, and the product
\[
\Pic(D)/2 \times \bigoplus_\I \ZZ/2 \to \W^\s(D).
\]
We claim that this product vanishes. To show this, we may postcompose with the injection $\W^s(D) \to \W^\s(F)$, so that the claim follows from the naturality of the product and the fact that $\Pic(F) = 0$. Next, we claim that the two $\W^\s(D)$-module structures on $\L_{1}^\s(D)$ and $\L_{-1}^\s(D)$ are obtained from their canonical $\FF_2$-module structures via the rank mod 2 homomorphism $\W^\s(D) \to \FF_2$. We first argue this in case of $\L_1^\s(D)$. From the above concrete description of this group, we may assume that $D=D_\mathfrak{p}$ is local for some dyadic prime $\mathfrak{p}$. Since the endomorphism ring of $\ZZ/2$ is $\FF_2$, the claimed module structure is equivalently described by a ring homomorphism $\W^\s(D_\mathfrak{p}) \to \FF_2$, of which the rank mod 2 homomorphism is the unique one. Indeed, $\W^\s(D_{\mathfrak{p}})$ is generated by units represented by the 1-dimensional forms $\langle a \rangle$ for $a \in D_{\mathfrak{p}}^\times$, see \cite{knebusch-rosenberg-ware2}*{Theorem 1.16}, so it follows that $\langle a \rangle$ is sent to $1$ under any ring homomorphism $\W^\s(D_{\mathfrak{p}}) \to \FF_2$. It follows that any such map is the rank mod 2 homomorphism as claimed.
Finally, we need to argue that the $\W^\s(D)$-module structure on $\Pic(D)/2$ is also obtained via the rank mod 2 homomorphism $\W^\s(D) \to \FF_2$. In Example~\reffour{example:free-dedekind} we have shown that the composite of maps of $\Einf$-ring spectra
\[
\L^\s(D) \to \K(D)^{\tC} \to \K_0(D)^{\tC}
\]
induces an isomorphism on $\pi_{-1}$. The result then follows from the well-known computation that $\K_0(D)$ is isomorphic to the trivial square zero extension $\ZZ \oplus \Pic(D)$. Indeed, recall that 
\[
\widehat{\mathrm{H}}^n(\Ct;\K_0(D)) = \begin{cases} \FF_2 \oplus \Pic(D)[2] & \text{ for } n \equiv 0 (2) \\ \Pic(D)/2 & \text{ for } n \equiv 1(2) \end{cases}
\]
so that the lax symmetric monoidality of Tate cohomology implies that the product
\[
\Pic(D)[2] \times \Pic(D)/2 \subseteq \widehat{\mathrm{H}}^0(\Ct;\K_0(D)) \times \Pic(D)/2 \to \Pic(D)/2
\]
vanishes as claimed.

Next, recall that the quadratic L-groups of $D$ are given as follows
\[
\L^\qdr_n(D) = \begin{cases} \W^\qdr(D) & \text{ for } n \equiv 0 (4),\\
						0 & \text{ for } n \equiv 1 (4), \\
						\bigoplus_{\I} \W^\qdr(\FF_\mathfrak{p}) & \text{ for } n \equiv 2 (4),\\
						A \rtimes \Pic(D)/2 & \text{ for } n \equiv 3 (4) \end{cases}
\]
where $\FF_{\mathfrak{p}} = D/\mathfrak{p}$ is the residue field at $\mathfrak{p}$, $A$ is the cokernel of the map $\W^\s(D) \oplus \W^\qdr(D^\cwedge_2) \to \W^\s(D^\cwedge_2)$ and $A \rtimes \Pic(D)/2$ is a certain extension of $\Pic(D)/2$ by $A$, see Proposition~\refthree{proposition:compute-quadratic-L-groups}. Moreover, the projection $\L_2^\qdr(D) \to \W^\qdr(\FF_{\mathfrak{p}})$ is induced by the projection $D \to \FF_\mathfrak{p}$.
Let us assume that $2$ is not invertible in $D$ (and hence the fraction field of $D$ is a number field), as else the map $\L^\qdr(D) \to \L^\s(D)$ is an equivalence and we have already computed all relevant products. 

The $\L^\s_*(D)$-module structure on $\L^\qdr_*(D)$ is now determined by various multiplications. Using similar arguments as above, one can show that these products are given as follows:

The $\W^\s(D)$-module structures on $\L_n^\qdr(D)$ are the evident ones for $n=0,1$, and is the restricted $\FF_2$-module structure on $\L^\qdr_2(D)$ obtained from the rank mod 2 homomorphism $\W^\s(D) \to \FF_2$. Moreover, the extension $A \rtimes \Pic(D)/2$ is one of $\W^\s(D)$-modules, where $\Pic(D)/2 = \L_{-1}^\s(D)$ is again a $\W^\s(D)$-module via the rank mod 2 map $\W^\s(D) \to \FF_2$ as above, and where the surjection $\L_0^\s(D^\cwedge_2) \to A$ is also one of $\W^\s(D)$-modules.
Next, we describe the products with $\L_1^\s(D)$ and $\L_{-1}^\s(D)=\Pic(D)/2$. We begin with $\L_1^\s(D)$ in which case there are products
\[
\L_1^\s(D) \times \L^\qdr_2(D) \to \L^\qdr_3(D) \quad \text{ and } \quad \L_1^\s(D) \times \L^\qdr_3(D) \to \W^\qdr(D).
\]
of which the latter vanishes (for this, one may use that the symmetrisation map $\W^\qdr(D) \to \W^\s(D)$ is injective, see the proof of Corollary~\refthree{corollary:finite-generation} and the previously established products on $\L^\s_*(D)$). 
For the first product, the composite with the surjection $\L^\qdr_3(D) \to \L^\s_3(D)=\Pic(D)/2$ vanishes, so it is determined by a map
\[
\L_1^\s(D) \times \L_2^\qdr(D) \to A
\]
about which we do not know more at this point.
Similarly, there are products
\[
\Pic(D)/2 \times \W^\qdr(D) \to \L_{3}^\qdr(D) \quad \text{ and } \quad \Pic(D)/2 \times \L^\qdr_3(D) \to \L^\qdr_2(D).
\]
of which the latter vanishes again. 
Moreover, the former product followed by the surjection $\L^\qdr_3(D) \to \L^\s_3(D)=\Pic(D)/2$ vanishes, and is hence determined by a product 
\[
\Pic(D)/2 \times \W^\qdr(D) \to A
\]
about which we again do not know more at this point.
Indeed, the composite under investigation is restricted form the tautological $\FF_2$-module structure of $\Pic(D)/2$ via the composite $\W^\qdr(D) \to \W^\s(D) \to \FF_2$ where the latter map is the rank mod 2 homomorphism. Now since we have assumed that $2$ is not invertible in $D$ there exists at least one dyadic prime $\mathfrak{p}$ and the  composite factors through the map $\W^\qdr(\FF_\mathfrak{p}) \to \W^\s(\FF_\mathfrak{p})$ which vanishes, see \cite{milnor-symmetric}*{Appendix 1}. 

Finally one can use the comparison results of \paperthree and \cite{Schlichting-genuine}, to find that the map $\L^\gs_*(D) \to \L^\qdr_*(D)$ is an isomorphism for $* \geq -2$ and the map $\L^\qdr_*(D) \to \L^\gs_*(D)$ is one for $* \leq -3$, so the above also determines (up to the ambiguity involving products with image in $A$) the ring structure on $\L^\gs_*(D)$.
\end{example}

Next, we will describe several GW- and L-spectra as $\Einf$-rings in terms of more familiar $\Einf$-rings. The first two examples are of the kind where the comparison maps to homotopy or Tate fixed points on K-theory are in fact (close to) equivalences:

\begin{example}
\label{example:perfect-fields-of-characteristic-2}%
For any perfect field of characteristic 2, the map of $\Einf$-ring spectra $\L^\s(F) \to \K(F)^{\tC}$ is an equivalence, see Remark~\refthree{remark:K-of-perfect-field}. Moreover, part of the argument in loc.\ cit.\ is that the map $\K(F) \to \ZZ$ is a 2-adic equivalence and hence induces an equivalence of $\Einf$-ring spectra $\K(F)^{\tC} \to \ZZ^{\tC}$. As a result, we obtain an $\Einf$-equivalence $\L^\s(F) \simeq \ZZ^{\tC}$, and deduce that the map of $\Einf$-ring spectra $\GW^\s(F) \to \ZZ^{\hC}$ is an equivalence after 2-adic completion. We also note that there exists a canonical equivalence of $\Einf$-ring spectra $\mathrm{THH}(\FF_2) \simeq \tau_{\geq0}\ZZ^{\tC}$ whose source is canonically an $\Einf$-$\FF_2$-algebra, see \cite{NS}*{Corollary IV.4.16}. In particular, we deduce that $\L^\s(A)$ is acquires to structure of an $\Einf$-$\FF_2$-algebra whenever $A$ is itself an $\Einf$-$\FF_2$-algebra.
\end{example}

\begin{example}
\label{ex:algebraically-closed-fields}%
For a quadratically closed field $K$ of characteristic different from 2, we argue that $\L^\s(K)$ is canonically equivalent to $\ku^{\tC}$, where the $\Ct$-action on $\ku$ is via complex conjugation. To see this, first note that $\L^\s_*(K) = \FF_2[u^{\pm1}]$ with $|u|=4$. Indeed, $\L^\s_*(K)$ is 4-periodic and vanishes in degrees not divisible by $4$, while $\L_0^\s(K) \cong \W_0^\s(K) \cong \FF_2$, see Remark~\refthree{remark:odd-vanishing}
and~\cite{milnor-symmetric}*{\S 3 Remark 3.4}. If $K$ has characteristic $0$ then there is an algebraically closed field $L$ which contains $K$ and $\CC$ as subfields. The maps $\L^\s(K) \to \L^\s(L) \leftarrow \L^\s(\CC)$ are then equivalences of $\Einf$-rings, because all homotopy rings are isomorphic to the graded ring $\FF_2[u^{\pm1}]$ with $|u|=4$, which has only the identity as endomorphism. If $K$ has positive characteristic different from $2$ we may choose an algebraic closure $\overline{K}$ of $K$ and the map $\L^\s(K) \to \L^\s(\overline{K})$ is an equivalence of $\Einf$-rings for the same reason as before. Now $\overline{K}$ is a perfect field of characteristic $p >2$, so we may consider the complete discrete valuation ring $D$ of $p$-typical Witt vectors of $\overline{K}$ and denote by $F$ an algebraic closure of the fraction field of $D$, an algebraically closed field of characteristic $0$. We claim that both maps $\L^\s(\overline{K}) \leftarrow \L^\s(D) \to \L^\s(F)$ are equivalences of $\Einf$-rings. To see that the left map is an equivalence, we note that $2$ is invertible in $D$ and $\overline{K}$, so that the map is equivalent to the map $\L^\qdr(D) \to \L^\qdr(\overline{K})$ which is an equivalence by Proposition~\refthree{proposition:hensel} since $D$ is $p$-adically complete and $D/p \to K$ is an isomorphism. In particular, $\L^\s_*(D)$ is also isomorphic to $\FF_2[u^{\pm1}]$ with $|u|=4$. The same argument as before then implies that also the map $\L^\s(D) \to \L^\s(F)$ is an equivalence. We have therefore now argued that any two quadratically closed fields of characteristic different from $2$ have canonically equivalent symmetric L-spectra. Note then that there are canonical maps of $\Einf$-rings $\L^\s(\CC) \to \K(\CC)^{\tC} \to \ku^{\tC}$. By \cite{Greenlees-reality}*{Corollary 3.3}, $\pi_*(\ku^{\tC})$ is again isomorphic to $\FF_2[u^{\pm1}]$ with $|u|=4$, so our favourite argument implies that the map $\L^\s(\CC) \to \ku^{\tC}$ is an equivalence of $\Einf$-rings.

As for $\GW$-theory, if $K$ is separably closed of characteristic different from 2 then it follows from Suslin's rigidity results that there is an 
equivalence of $\Einf$-rings $\K(K)^\cwedge_2 \simeq \ku^\cwedge_2$, and the above implies that $\GW^\s(K)^\cwedge_2 \to (\ku^\cwedge_2)^{\hC}$ is an equivalence. In particular, we obtain an equivalence of $\Einf$-rings $\tau_{\geq0} \GW^\s(K)^\cwedge_2 \simeq \ko^\cwedge_2$ for $K$ separably closed of characteristic different from $2$.
\end{example}

To move on, let us recall that for a commutative ring with involution $(R,\sigma)$ and $\sigma$-fixed unit $u \in R^\times$, there is associated the 1-dimensional $\sigma$-symmetric form $(x,y) \mapsto xu\sigma(y)$, giving rise to a unit $\langle u \rangle$ in the classical symmetric Grothendieck-Witt ring $\GW_0^\s(R,\sigma)$ of $(R,\sigma)$. Since multiplication by $u$ provides an isomorphism between the forms associated to $1$ and to $u\sigma(u)$, this construction gives rise to a map of commutative rings 
\[
\ZZ\Bigl[\frac{(R^\times)^\sigma}{\mathrm{N}(R^\times)}\Bigr] \to \GW_0^\s(R,\sigma),
\]
where $\mathrm{N}(R^\times)$
refers to the subgroup given by elements of the form $u\sigma(u)$. This map is instrumental to many of the classical descriptions of the symmetric Grothendieck-Witt and Witt rings of commutative rings. The following provides a lift of this algebraic construction to the level of the $\Einf$-ring spectrum $\GW^\s(R,\sigma)$. We will then use this to lift several known algebraic descriptions of 
symmetric Grothendieck-Witt and Witt rings
to a description of the $\Einf$-ring structure on suitable symmetric Grothendieck--Witt and L-spectra.

\begin{construction}
\label{construction:Bgl-to-gl}%
Let $(\C,\QF)$ be a symmetric monoidal Poincar\'e $\infty$-category. Since the map $\Poinc \to \GWspace$ is a lax symmetric monoidal transformation, we find that the map $\Poinc(\C,\QF) \to \GWspace(\C,\QF)$ is an $\Einf$-map with respect to tensor product on each side. Passing to $\otimes$-invertible components, we obtain a natural map of $\Einf$-groups $\Poinc(\C,\QF)^\times \to \gl_1(\GW(\C,\QF))$.
In the case where $(\C,\QF) = (\Modp{R},\QF^{\sym}_{\sig})$ for an $\Einf$-ring spectrum with involution $(R,\sig)$ %
this gives a map $\Pic(R)^{\hC} \to \gl_1(\GW^{\sym}(R,\sig))$,
which fits into a commutative diagram
\[
\begin{tikzcd}
	\Pic(R)^{\hC} \ar[r] \ar[d] & \gl_1(\GW^\s(R,\sig)) \ar[r] & \gl_1(\K(R,\sig)^{\hC}) \ar[d] \\
	\Pic(R) \ar[rr] && \gl_1(\K(R))
\end{tikzcd}
\]
whose lower horizontal map is well-known in algebraic K-theory.
Now suppose that $R$ is discrete. The inclusion of the identity component $\B R^{\times} = \mathrm{Bgl}_1(R) \subseteq \Pic(R)$
then induces a map
\[
(\B R^\times)^{\hC} \to \Pic(R,\sig)^{\hC} \to \gl_1(\GW^\s(R,\sig))
\]
whose first map is an equivalence if the canonical map $\ZZ \to \pi_0(\Pic(R))$ is an isomorphism.
Observe now that $\B R^\times$ is the space underlying the $\ZZ$-module spectrum $\Sigma \GEM R^\times$ and that the $\Ct$-action induced by $\sigma$ on it is automatically $\ZZ$-linear; in fact, it is induced by the $\Ct$-action $u \mapsto \sigma(u)^{-1}$ on the abelian group $R^\times$. It follows that also the fixed point spectrum $(\B R^\times)^{\hC}$ is canonically a $\ZZ$-module, and hence (non-canonically) splits as $\Einf$-group into a direct sum $\mathrm{H}^1(\Ct;R^\times)\oplus \B \mathrm{H}^0(\Ct;R^\times)$.
The choice of such a splitting then induces an $\Einf$-map 
\[
\frac{(R^\times)^\sigma}{\mathrm{N}(R^\times)} = \mathrm{H}^1(\Ct;R^\times) \to \Pic(R)^{\hC} \to \gl_1(\GW^\s(R,\sig)),
\]
and by adjunction, %
maps of $\Einf$-ring spectra
\[
\SS\Bigl[\frac{(R^\times)^\sigma}{\mathrm{N}(R^\times)}\Bigr] \to \GW^\s(R,\sig) \to \L^\s(R,\sig).
\]
On $\pi_0$, this map sends the class of a $\sig$-invariant unit $u \in R^\times$ to the 1-dimensional symmetric form $\langle u \rangle$. In particular, if $R$ is equipped with the identity involution, there exist maps of $\Einf$-ring spectra
\[
\SS\Bigl[\frac{R^\times}{(R^\times)^2}\Bigr] \to \GW^\s(R) \to \L^\s(R)
\]
which on $\pi_0$ send a unit $u \in R^\times$ to $\langle u \rangle$.
\end{construction}

\begin{example}
In \paperthree, we have utilised a topological version of Grothendieck-Witt theory $\GW_\topo ^\s(\RR)$ introduced in the context of Banach algebras with involution by Schlichting in \cite{schlichting-derived}*{\S 10} and which also features in the analysis of symmetric L-theory of real $C^*$-algebras in \cite{LNS}. By construction, it comes with $\Einf$-ring map $\GW^\s(\RR) \to \GW^\s_\topo(\RR)$ and by \cite{LNS}*{Remark 4.2} it also comes with an $\Einf$-ring map $\ko \to \GW^\s_\topo(\RR)$; informally this map is induced by sending a finite dimensional real vector space to its (essentially unique) positive definite symmetric bilinear form. As a consequence of this, Construction~\reffour{construction:Bgl-to-gl}, and the observation that $\RR^\times/(\RR^\times)^2 = \{\pm 1\} = \Ct$, there exists a map of $\Einf$-rings
\[
\ko[\Ct] \to \tau_{\geq0}\GW_\topo^\s(\RR)
\]
which, by the argument used in the proof of \cite{LNS}*{Corollary 3.6}, is an equivalence.
\end{example}

\begin{example}
Let $K$ be a complete discretely valued field with ring of integers $D$, a complete discrete valuation ring with residue field $k$. Assume that the characteristic of $k$ is different from $2$ and choose a uniformiser $\pi \in D$. 
Consider the map $\Ct \to K^\times/(K^\times)^2$ with image the subgroup generated by $\pi \in D$. Applying Construction~\reffour{construction:Bgl-to-gl}, we deduce the existence of a map of $\Einf$-rings $\L^\s(D)[\Ct] \to \L^\s(K)$ which on $\pi_0$ sends the generator of $\Ct$ to $\langle \pi \rangle$.
As in Example~\reffour{ex:algebraically-closed-fields}, the projection $D \to k$ induces an equivalence $\L^\s(D) \to \L^\s(k)$ since $D$ is $\pi$-adically complete and $D/\pi \to k$ is an isomorphism: Indeed, since $2$ is invertible in $D$ and $k$, it suffices to show that $\L^\qdr(D) \to \L^\qdr(k)$ is an equivalence, which follows from the $\pi$-completeness of $D$, see e.g.\ Proposition~\refthree{proposition:hensel}. In total, we obtain a map of $\Einf$-ring spectra
\[
\L^\s(k)[\Ct] \to \L^\s(K)
\]
which induces the Springer decomposition $\W^\s(k)[\Ct] \cong \W^\s(K)$ of symmetric Witt groups on $\pi_0$, see e.g.\ \cite{Lam}*{VI Corollary 1.7}. Since both homotopy rings are 4-periodic, with non-trivial homotopy only in degrees divisible by $4$, see again Remark~\refthree{remark:odd-vanishing}, the above map is an equivalence of $\Einf$-rings.
For instance, we obtain that for $p$ an odd prime, there is an equivalence of $\Einf$-ring spectra 
\[
\L^\s(\FF_p)[\Ct] \simeq \L^\s(\QQ_p).
\]
That the same result is not true for $p=2$ is well-known and already fails on the level of $\pi_0$; Here the source has 4 elements while the target has $32$-elements.
\end{example}

\begin{example}
\label{example:Gm}%
Let $R$ be a regular noetherian ring of finite Krull dimension. Consider the ring of Laurent polynomials $R[t^{\pm1}]$ as a commutative ring with trivial involution as in Remark~\reffour{remark:all-involutions}. Similarly as in the example above, consider the map 
\[
\Ct \to \Ct \times \Ct \cong \ZZ[t^{\pm1}]^\times/(\ZZ[t^{\pm1}]^\times)^2 \to R[t^{\pm1}]^\times/(R[t^{\pm1}]^\times)^2
\]
with image the subgroup generated by $t$. From Construction~\reffour{construction:Bgl-to-gl}, we obtain maps of $\Einf$-rings 
\[
\GW^\s(R)[\Ct] \to \GW^\s(R[t^{\pm1}]) \quad \text{and} \quad \L^\s(R)[\Ct] \to \L^\s(R[t^{\pm1}]).
\]
Under the the equivalences indicated in Remark~\reffour{remark:all-involutions}, one can show that the underlying map of spectra for the latter map is the identity while for $\GW$ it is induced by the boundary map in the Genauer sequence (which is not an equivalence unless $R = 0$).

In \S\reffour{subsection:mult-shaneson} below, we will discuss the multiplicative properties of the Ranicki-Shaneson decompositions for the L-theory of the group ring $R[\ZZ]$, with involution induced by the group inversion on $\ZZ$.
\end{example}

We deduce the following computation from Example~\reffour{example:Gm}.
\begin{example}
Let $R$ be a regular noetherian commutative ring of finite Krull dimension and denote by $\mathbb{P}^1(R)$ the projective line over $R$. In \cite{motives} it is shown that the square
\[
\begin{tikzcd}
	\Perf(\mathbb{P}^1(R)) \ar[r]\ar[d] & \Perf(\mathbb{A}^1(R)) \ar[d] \\
	\Perf(\mathbb{A}^1(R)) \ar[r] & \Perf(\mathbb{G}_m(R))
\end{tikzcd}
\]
canonically refines to a Poincar\'e-Verdier square using the symmetric Poincar\'e structures throughout. Moreover, it is shown that $\GW^\s$ and $\L^\s$ are $\mathbb{A}^1$-invariant in this situation, so that there are pullback squares
\[
\begin{tikzcd}
	\GW^\s(\mathbb{P}^1(R)) \ar[r] \ar[d] & \GW^\s(R) \ar[d] && \L^\s(\mathbb{P}^1(R)) \ar[r] \ar[d] & \L^\s(R) \ar[d] \\
	\GW^\s(R) \ar[r] & \GW^\s(\mathbb{G}_{m}(R)) && \L^\s(R) \ar[r] & \L^\s(\mathbb{G}_{m}(R))
\end{tikzcd}
\]
of $\Einf$-ring spectra. Through the canonical section $1 \colon R \to \mathbb{G}_m(R) \subseteq \mathbb{P}^1(R)$, these diagrams are in fact pullback diagrams of augmented $\Einf$-$\GW^\s(R)$- and $\L^\s(R)$-algebras, respectively. Example~\reffour{example:Gm} implies that $\L^\s(\mathbb{P}^1(R))$ as an $\Einf$-ring spectrum is given by $\L^\s(R) \otimes P$ where $P = \SS \times_{\SS[\Ct]} \SS$ is the loop object of $\SS[\Ct]$. If $R$ is a commutative $\FF_2$-algebra, recall from Example~\reffour{example:perfect-fields-of-characteristic-2} that $\L^\s(R)$ admits a canonical $\Einf$-$\FF_2$-algebra structure. In particular, as an $\Einf$-ring, we have $\L^\s(\mathbb{P}^1(R)) \simeq \L^\s(R) \otimes_{\FF_2}[ \FF_2 \times_{\FF_2[\Ct]} \FF_2]$ and $\FF_2[\Ct]$ is a trivial square-zero extension. Consequently, if $R$ is a commutative $\FF_2$-algebra, we deduce that $\L^\s(\mathbb{P}^1(R))$ is equivalent to the trivial square-zero extension $\L^\s(R) \oplus \Omega\L^\s(R)$ as $\Einf$-$\L^\s(R)$-algebra.

Moreover, by \cite{formality}*{Theorem B}, it follows that as $\Eone$-$\GW^\s(R)$-algebra, $\GW^\s(\mathbb{P}^1(R))$ is the trivial square-zero extension on $\GW(R,\Omega \QF^\s)$. Similarly, as $\Eone$-$\L^\s(S)$-algebra, $\L^\s(\mathbb{P}^1(R))$ is the trivial square-zero extension on $\Omega \L^\s(R)$, also when $R$ is not an $\FF_2$-algebra.

As a concrete example, consider the projective line $\mathbb{P}^1(D)$ over the ring of integers in a number field. Then we obtain an equivalence of augmented $\Eone$-rings $\L^\s(\mathbb{P}^1(D)) \simeq \L^\s(D) \oplus \Omega \L^\s(D)$, where the latter is the trivial square-zero extension. In particular, by Corollary~\refthree{corollary:compute-L-groups-global-field}, we find 
\[
\L^{\sym}_0(\mathbb{P}^1(D)) = \L^{\sym}_0(D) \oplus \L^{\sym}_1(D) = \W^\s(D) \oplus (\ZZ/2)^{\oplus n}
\]
where $n$ is the number of dyadic primes in $D$ and the multiplication on $(\ZZ/2)^n$ is trivial, and the $\W^\s(D)$-module structure on $(\ZZ/2)^{\oplus n}$ is via the rank mod 2 homomorphism $\W^\s(D) \to \FF_2$ as explained in Example~\reffour{example:L-rings-dedekind}.
\end{example}

\begin{remark}
As a consequence of the $\Einf$-ring map $\SS[R^\times/(R^\times)^2] \to \GW^\s(R)$ from Construction~\reffour{construction:Bgl-to-gl}, all the units $\langle u \rangle$ in $\GW_0^\s(R)$ are \emph{strict} units in the $\Einf$-ring $\GW^\s(R)$, i.e.\ are represented by $\Einf$-maps $\SS[\ZZ] \to \GW^\s(R)$, and the same therefore holds true for their images in $\L_0^\s(R)$. In fact, these strict units are strictly of order $2$, that is, their classifying maps factor as the composite of $\Einf$-ring maps $\SS[\ZZ] \to \SS[\Ct] \to \GW^\s(R) \to \L^\s(R)$. 

Strict units allow for more constructions than general units. For instance, for an $\Einf$-ring spectrum $A$ and $u \in \pi_0(A)$ a strict unit classified by an $\Einf$-map $\SS[\ZZ] \to A$, one may consider the base-change $A \otimes_{\SS[\ZZ]} \SS$ of $A$ along the augmentation $\SS[\ZZ] \to \SS$ induced by the group homomorphism $\ZZ \to \{1\}$ on spherical group rings. One can compute that the underlying $A$-module of this base-change is given by $A/(u-1)$, so that this $A$-module is canonically an $\Einf$-$A$-algebra. In fact, more generally, one can show that $A/(u-1)^k$ is an $\Einf$-$A$-algebra for all $k\geq 1$: Considering the universal case $\SS[\ZZ]$ with strict unit $t$ corresponding to the generator $1 \in \ZZ$, we find that the fibre of the augmentation map $\SS[\ZZ] \to \SS$ is a Smith ideal in $\SS[\ZZ]$. Taking (pointwise) $k$th powers of the Smith ideal gives a new Smith ideal, whose quotient can be shown to deserve the name $\SS[t]/(t-1)^k$. Then $A/(u-1)^k$ is given by $A \otimes_{\SS[\ZZ]} \SS[t]/(t-1)^k$.

Now, it follows from the above construction that for a discrete commutative ring $R$, the element $\langle -1 \rangle = -1 \in \L_0^\s(R)$ is a strict unit. In particular, for all $k\geq 1$, we deduce that $\L^\s(R)/2^k$ is canonically an $\Einf$-ring spectrum, in fact an $\Einf$-algebra over $\L^\s(\ZZ)/2^k$. Moreover, one can form the base-change $\L^\s(R) \otimes_{\SS[\Ct]} \SS$ of $\L^\s(R)$ along the augmentation of $\SS[\Ct]$ which is again an $\Einf$-$\L^\s(R)$-algebra whose underlying $\L^\s(R)$-module can be shown to be given by $\L^\s(R)_{\hC}$, where the $\Ct$-action is induced by multiplication by the strict unit of order 2 given by $-1$.

As a consequence of the $\Einf$-equivalences $\L^{\sym}(\ZZ)[1/2] = \L^{\sym}(\RR)[1/2] = \KO[1/2]$, the second arising from \cite{LN}*{Theorem 5.2}, $-1 \in \pi_0(\KO[1/2])$ is also given the structure of a strict unit. This is not surprising, though: Indeed, it follows from results of Carmeli that the spectrum of strict units in $\SS[1/2]$ is discrete with $\pi_0$ given by $\{\pm 1\}$, while the spectrum of strict units of $\SS$ is contractible \cite{Carmeli}*{Theorem B}. In particular, in general, $-1 \in \pi_0(A)$ does not admit the structure of a strict unit.
\end{remark}

\subsection{The periodicity elements of Karoubi and Ranicki}
\label{subsection:PeriodicityviaMultiplication}%

In this section we aim to explain the relation between the multiplicative structure on Grothendieck-Witt- and $\L$-spectra and their periodicity properties, and then use this to reconcile our formulations of these in \paperone and \papertwo with Ranicki's formulation of his periodicity isomorphisms for $\L$-groups and Karoubi's original formulation of his fundamental theorem. In particular, we will show that the equivalences
\[
\mathbb S^2 \otimes \L^\sym(R) \simeq \L^{-\sym}(R) \quad \text{and} \quad \mathbb S^1 \otimes \mathrm V^\sym(R) \simeq \mathrm U^{-\sym}(R),
\]
produced in \papertwo are induced by multiplication with certain periodicity classes $\rho$ and $\kappa$ along with appropriate generalisations to form parameters and $\Einf$-rings, where %
$\mathrm V$ and $\mathrm U$ are given by
\[
\mathrm V(\C,\QF) = \fib(\GW(\C,\QF) \xrightarrow{\fgt} \K(\C))  \quad \text{and} \quad \mathrm U(\C,\QF) = \fib(\K(\C) \xrightarrow{\hyp} \GW(\C,\QF)),
\]
see \S\reftwo{subsection:Genauer-Bott-Karoubi}.

To get started, %
recall that by Corollary~\refone{corollary:quadratic-genuine} any hermitian structure $\QF \colon \C\op \rightarrow \Sp$ canonically lifts to a functor $\tilde\QF\colon \C\op \to \Sp^{\g\Ct}$ to genuine $\Ct$-spectra along the genuine fixed point functor $(-)^\gC \colon \Sp^\gC \rightarrow \Sp$. The underlying $\Ct$-spectrum of $\tilde \QF(X)$ is given by $\Bil_\QF(X,X)$ and its geometric fixed points by $\Lambda_\QF(X)$. We can now tensor this lift valuewise with the genuine representation sphere $\SS^{m+n\sigma} \in \Sp^{\g\Ct}$ for $m,n \in \mathbb Z$ and apply genuine fixed points. Slightly abusing notation we will denote the resulting functor by $\QF\qshift{m+n\sigma}$. Inspection of the isotropy separation square of $\mathbb S^{[m+n\sigma]} \otimes \tilde\QF$ then immediately verifies that $\QF\qshift{m+n\sigma}$ is again a hermitian structure on $\C$ with
\[
\Bil_{\QF\qshift{m+n\sigma}} \simeq \SS^{m+n\sigma} \otimes \Bil_\QF \quad \text{and} \quad \Lambda_{\QF\qshift{m+n\sigma}} \simeq \mathbb S^m \otimes \Lambda_\QF,
\]
connected by the composite
\[
\SS^m \otimes \Lambda_\QF \longrightarrow \SS^m \otimes (\Bil_\QF \circ \Delta)^\tC \simeq (\SS^{m+n\sigma} \otimes \Bil_\QF \circ \Delta)^\tC,
\]
where the first map is the image of the glueing map for $\QF$ and the second induced by the inclusion $\SS^0 \rightarrow \SS^{n\sigma}$ as fixed points. In particular, the operation $(-)\qshift{m+n\sigma}$ preserves both symmetric and quadratic Poincaré structures.

To obtain the interaction of this construction with the tensor product of Poincaré $\infty$-categories, we note the following alternative description. Theorem~\refone{theorem:classification-genuine-modules} and Remark~\refone{remark:NormModules} imply that $\Funq(\Spaf) \simeq \Sp^{\g\Ct}$ as symmetric monoidal $\infty$-categories, see \S\refone{subsection:examples-monoidal}). In particular, we may identify hermitian structures on $\Spaf$ with genuine $\Ct$-spectra by regarding their classifying squares as %
isotropy separation ones. Consequently, the symmetric monoidal structure on $\Cath$ determines an action of $\Sp^{\g\Ct}$ on $\Funq(\C)$ for any stable $\infty$-category $\C$. Unwinding definitions one finds that acting with $\mathbb S\qshift{m+n\sigma}$ on some $\QF$ reproduces exactly $\QF\qshift{m+n\sigma}$. %
We therefore have that
\[
(\C,\QF\qshift{m+n\sigma}) \otimes (\D,\Psi\qshift{k+l\sigma}) = (\C \otimes \D,(\QF \otimes \Psi)^{(m+k)+(n+l)\sigma}).
\]
In particular, the operation $(-)\qshift{m+n\sigma}$ is invertible with inverse $(-)\qshift{-m-n\sigma}$.  %
Now for the periodicity elements:

\begin{construction}
We will first define two elements 
\[
\kappa \in \Omega^\infty\GW(\Spaf,(\QF^\uni)\qshift{-1-\sigma}) \quad \text{and} \quad \delta \in \Omega^{\infty-1}\GW(\Spa^\fin,(\QF^\uni)\qshift{-1}),
\]
multiplication with which is equivalent to the maps 
\[
\Omega_\ast \colon \GW(\C,\QF) \longrightarrow \GW(\C,\QF\qshift{-1-\sigma}) \quad \text{and} \quad \partial\colon \GW(\C,\QF) \longrightarrow \SS^1 \otimes \GW(\C,\QF\qshift{-1})
\]
for any Poincaré $\infty$-category $(\C;\QF)$; the former map here is induced by the equivalence $\Om\colon (\C,\QF) \rightarrow (\C,\QF\qshift{-1-\sigma})$ from Proposition~\refone{proposition:shiftofqingeneral}, and the latter is the boundary map of the Bott-Genauer sequence. The classes $\kappa$ and $\rho$ must then simply be given by applying these two maps to the unit of the multiplication, that is, to $[\mathbb S,q^\uni] \in \Omega^\infty\GW(\Spaf,\QF^\uni)$.
In the former case we thus set
\[
\kappa \coloneq \Omega_\ast[\mathbb S,q^\uni] = [\mathbb S^{-1},q^\uni] \in \Omega^\infty\GW(\Spaf,(\QF^\uni)\qshift{-1-\sigma}),
\]
and we will call $\kappa$ the \emph{Karoubi periodicity element}. In the latter case we put
\[
\delta \coloneq \partial[\mathbb S,q^\uni] \in \Omega^{\infty-1}\GW(\Spa^\fin,(\QF^\uni)\qshift{-1}),
\]
which we refer to as the \emph{Bott-Genauer element}. To give an explicit representative of it, recall that generally 
\[
\Omega^{\infty-1}\GW(\D,\Psi\qshift{-1}) \simeq |\Cob(\D,\Psi\qshift{-1})|,
\]
which receives a map from $\Poinc(\D,\Psi) \simeq \grpcr\Cob(\D,\Psi\qshift{-1})$. By Proposition \reftwo{proposition:boundarymapGB} this map agrees with the composite of $\partial$ and the tautological map $\Poinc(\C,\QF) \rightarrow \Omega^{\infty}\GW(\C,\QF)$, so that $\delta$ ends up simply represented by the object $(\mathbb S,q^\uni)$ of $\Cob(\Spaf,(\QF^\uni)\qshift{-1})$.

To see that multiplication with these classes does indeed agree with the indicated maps, note in the case of $\kappa$ that the Poincaré functor $\Omega \colon (\C,\QF) \rightarrow (\C,\QF\qshift{-1-\sigma})$ naturally identifies with 
\[
\id \otimes \Om\colon (\C,\QF) \otimes (\Spa^\fin,\QF^\uni) \longrightarrow (\C,\QF) \otimes (\Spa^\fin,(\QF^\uni)\qshift{-1-\sigma}).
\]
That the induced map $\Omega_\ast \colon \GW(\C,\QF) \longrightarrow \GW(\C,\QF\qshift{-1-\sigma})$ is given by multiplication with $\kappa$ is then witnessed by the lax symmetric monoidal structure of $\GW$, established in Corollary~\reffour{corollary:GW-universal-multiplicative}. More generally,
if $\F\colon \Catp \to \E$ is an arbitrary additive functor to some $\infty$-category $\E$ with finite limits then $\F$ is a $\GWspace$-module (for the induced $\Sps$-module structure on $\cP(\E) = \Fun(\E\op,\Sps)$ over $\Sps$, say) by Corollary~\reffour{corollary:GW-universal-multiplicative} and so for $(\C,\QF) \in \Catp$ we obtain a commutative square
\[
\begin{tikzcd}[column sep = 30pt]
\F(\C,\QF) \otimes \GWspace(\Spaf,\QF^{\uni}) \ar[d,"\mu"]\ar[r,"{(\id,\Om_*)}"] & \F(\C,\QF) \otimes \GWspace(\Spaf,(\QF^{\uni})\qshift{-1-\sig}) \ar[d,"\mu"] \\
\F(\C,\QF) \ar[r,"\Om_*"] & \F(\C,\QF\qshift{-1-\sig})
\end{tikzcd}
\]
in $\cP(\E)$, natural in $(\C,\QF)$, witnessing the fact that the lower horizontal map is induced by multiplication with $\kappa$. 
If, furthermore $\mathcal F$ is a module over some lax monoidal additive functor $\mathcal G \colon \Catp \rightarrow \Sps$, and $(\C,\QF)$ is a module over a monoidal Poincaré $\infty$-category $(\D,\Psi)$, then the action of $\GWspace(\Spaf,\QF\qshift{-1-\sigma})$ on $\F$ above factors through that of $\mathcal G(\D,\Psi\qshift{-1-\sigma})$, and so the action of $\kappa$ on $\F(\C,\QF)$ only depends on its image in $\mathcal G(\D,\Psi\qshift{-1-\sigma})$ under the unit map. We will continue to denote such images under unit maps by $\kappa$; for example, we obtain $\kappa \in \Lspace(\Dperf(\mathbb Z),\QF^{\mathrm g\sym})$ this way.

The analogous statements for $\delta$ are obtained almost verbatim by noting that according to Lemma \refone{lemma:pairing-projection} the metabolic Poincaré-Verdier sequence
\[
(\C,\QF\qshift{-1}) \to \Met(\C,\QF) \xrightarrow{\met} (\C,\QF)
\]
of Example~\reftwo{example:metabolicfseq} is naturally obtained from that of $(\Spaf,\QF^{\uni})$ by tensoring with $(\C,\QF)$. Arguing as above, we then get that for any additive $\F\colon \Catp \to \E$ to a stable $\E$, the boundary map $\partial\colon \F(\C,\QF) \rightarrow \mathbb S^1 \otimes \F(\C,\QF\qshift{-1})$ is induced by multiplication with $\delta$. Again we will regard $\delta$ as an element of $\Omega^{\infty-1}\mathcal G(\D,\Psi)$ for any monoidal $(\D,\Psi) \in \Catp$ and lax monoidal additive $\mathcal G \colon \Catp \rightarrow \Spa$ without changing notation. 

Finally, note that as the Poincaré functors $\Omega$ are equivalences, the Karoubi periodicity element $\kappa$ is invertible in the sense that there exists a unique element $\kappa^{-1} \in \mathcal G(\D,\Psi\qshift{1+\sig})$ such that $\kappa\kappa^{-1}=\kappa^{-1}\kappa=1$. The same is of course generally true for the images of $\delta$ if and only if $\mathcal G$ is bordism invariant.  
In particular, the image $\delta \in \Om^{\infty-1}\L(\Spaf,(\QF^\uni)\qshift{-1})$ is invertible with inverse an element $\delta^{-1} \in \Om^{\infty+1}\L(\Spaf,(\QF^\uni)\qshift{1})$. %
We may therefore also consider
\[
\rho \coloneq \kappa \cdot \delta^{-2} \in \Om^{\infty+2}\L(\Spaf,(\QF^\uni)\qshift{1-\sigma})
\]
which we will refer to as the \emph{Ranicki periodicity element}. Multiplication by $\rho$ provides an equivalence
\[
\rho \cdot - \colon \F(\C,\QF) \longrightarrow \SS^{-2}\otimes \F(\C,\QF\qshift{1-\sigma})
\]
for every bordism invariant additive functor $\mathcal F$ which, by construction, agrees with the composite equivalence
\[
\F(\C,\QF) \xrightarrow{\partial^{-1} \circ \partial^{-1}} \SS^{-2}\otimes \F(\C,\QF\qshift{2}) \xrightarrow{\Omega_*} \SS^{-2}\otimes \F(\C,\QF\qshift{1-\sigma}).
\]
To compare this element to Ranicki's incarnation of it, let us finally give an explicit representative for $\rho$. To this end, recall that there is a natural map $\Poinc(\C,\QF\qshift{-n}) \rightarrow \Omega^{\infty+n}\L(\C,\QF)$, that is traditionally used to identify the homotopy groups of $\L$-spectra as bordism groups: Using $\Omega^\infty\L(\C,\QF) = |\Poinc \rho(\C,\QF)|$ this originates from the equivalence $(\C,\QF\qshift{-n}) \hookrightarrow \rho_n(\C,\QF)$ onto the subcategory with trivial boundary through the general map $\fib(X_n \rightarrow X^{\partial \Delta^n}) \rightarrow \Omega^n|X|$ %
where $X^- \colon \mathrm s\Sps\op \rightarrow \Sps$ is the limit preserving extension of a simplicial space $X$. The element $\delta$ is then the image of $(\SS^{-1},q^\uni) \in \Poinc(\Spaf,(\QF^\uni)\qshift{-1-\sigma})$ under this map (for $n=2$): Since this Poincaré form represents $\kappa$ the claim boils down to the composite
\[
\Poinc(\C,\QF\qshift{-n}) \rightarrow \Omega^{\infty+n}\L(\C,\QF) \xrightarrow{\partial^{n}}  \Omega^{\infty}\L(\C,\QF)
\]
being the tautological map. For $n=1$ this follows immediately from the analogous (finer) statement for Grothendieck-Witt spaces in Proposition~\reftwo{proposition:boundarymapGB} and the general case reduces to this by considering the embedding of $\rho_n(\C,\QF)$ into the $n$-fold metabolic Poincaré $\infty$-category in Lemma~\reftwo{lemma:rho_vs_power_set}. 
\end{construction}

\begin{remark}\
\begin{enumerate}
\item For additive $\F$, we deduce from Proposition~\reftwo{proposition:Q-bordism} the formula 
\[
\F^\bord(\C,\QF) \simeq \colim [\F(\C,\QF) \xrightarrow{\delta\cdot -} \Sigma \F(\C,\QF\qshift{-1}) \xrightarrow{\delta\cdot -} \Sigma^2 \F(\C,\QF\qshift{-2}) \to \cdots ]
\]
and thus $\F^\bord = \F[1/\delta]$, where $\F[1/\delta]$ refers to the initial functor under $\F$ for which multiplication by $\delta$ is an equivalence, or equivalently, the initial functor under $\F$ which vanishes on $\Met(\C,\QF)$ for all Poincar\'e $\infty$-categories; that this initial functor agrees with $\F^\bord$ is a consequence of Lemma~\reftwo{lemma:criterion}. In particular, we have $\L = \GW[1/\delta]$; the latter term (defined as the above displayed filtered colimit) was used as a definition of L-theory in Schlichting's work, see \cite{schlichting-derived}*{Definition 7.1} (where the class we call $\delta$ corresponds to the class called $\eta$ in loc.\ cit.).
\item Recall from \S\reftwo{subsection:Mackey} that any additive $\F\colon \Catp \to \Spa$ has a canonically associated genuine refinement $\F^{\g\hyp} \colon \Catp \to \Spa^{\g\Ct}$ along the fixed point functor $(-)^{\g\Ct}$ such that the geometric fixed points of $\F^{\g\hyp}$ canonically identify with $\F^\bord$. The above formula $\F^\bord(\C,\QF) = \F(\C,\QF)[1/\delta]$ is then well-known in equivariant homotopy theory, where the geometric fixed points are obtained from the genuine fixed points by inverting the equivariant Euler class of the reduced regular representation: And indeed, unravelling the definitions, (and using that the reduced regular representation for $\Ct$ is precisely the sign representation) this equivariant Euler class for $\KR(\C,\QF)=\GW^{\g\hyp}(\C,\QF)$ is precisely the Bott-Genauer class $\delta$.
\end{enumerate}
\end{remark}

Let us now explain the relation of the Karoubi periodicity element $\kappa$ to his fundamental theorem. Recall that its general implementation in \papertwo is the composite
\[
\mathbb S^1 \otimes \mathrm V(\C,\QF) \simeq \GW(\C,\QF\qshift{1}) \simeq \mathrm U(\C,\QF\qshift{2})
\]
with both intermediate equivalences arising from the Bott-Genauer sequence, see Corollary~\reftwo{corollary:karoubi-fundamental-theorem}. Further composing with multiplication by $\kappa$ on the additive functor $\mathrm{U}$, we obtain an equivalence 
\[
\mathbb S^1 \otimes \mathrm V(\C,\QF) \simeq \mathrm U(\C,\QF\qshift{2}) \xrightarrow{\kappa \cdot -} \mathrm U(\C,\QF\qshift{1-\sigma}).
\]
which by construction agrees with the map used in \papertwo to establish our extensions of Karoubi's fundamental theorem, see the discussion after Corollary~\reftwo{corollary:karoubi-fundamental-theorem}.

In the special case of discrete rings we can unwind further:

\begin{example}
\label{example:perclasses}%
Let $R$ be a commutative ring, $A$ and $R$-algebra and $M$ an invertible $R$-module with involution over $A$. Recall that in this case for a Poincaré structure on $\QF_M$ on $\Dperf(A)$, whose bilinear part is classified by $M$, the bilinear part of $\QF\qshift{1-\sigma}$ is classified by $-M$, simply given by inserting a sign into the original action. 
\begin{enumerate}
\item In the simplest case of the symmetric and quadratic structures we consequently have $(\QF^\sym_M)\qshift{1-\sigma} = \QF^\sym_{-M}$ and periodicity elements 
\[
\kappa \in \Om^{\infty}\GW(\Dfree(R),(\QF^{-\sym})\qshift{-2}) \quad \text{and} \quad \rho \in \Om^{\infty+2}\L(\Dfree(R),\QF^{-\sym}),
\]
both represented by the complex $R[-1]$ with its tautological 2-shifted skew-symmetric form. Multiplication with these elements now gives equivalences
\[
\mathbb S^1 \otimes \mathrm V(\D^c(A),\QF^r_{M}) \xrightarrow{\kappa \cdot -} \mathrm U(\D^c(A),\QF^r_{-M}) \quad \text{and} \quad \L(\D^c(A),\QF^r_M) \xrightarrow{\rho \cdot -} \SS^{-2}\otimes \L(\D^c(A),\QF^r_{-M}),
\]
respectively, where $r \in \{\qdr,\sym\}$ and $c \subseteq \K_0(A)$ is some choice of control group. By the discussion above these agree with the ones established in Corollaries~\reftwo{corollary:karoubi-fundamental} and~\reftwo{corollary:Ranickipergeneral}, respectively, 
and multiplication with $\rho^2 \in \Om^{\infty+4}\L(\Dfree(R),\QF^\sym)$, represented by $R[-2]$ with its tautological 4-shifted symmetric form 
gives Ranicki's classical periodicity isomorphisms
\[
\L(\D^c(A),\QF^r_M) \longrightarrow \SS^{-4}\otimes \L(\D^c(A),\QF^r_{M}).
\]
\item As for the genuine symmetric structure, we have $(\QF^{\mathrm g\sym}_R)\qshift{1-\sigma} \simeq \QF^{-\gev}_R$ by Corollary~\refone{corollary:shift-genuine}, and so we may view the Karoubi and Ranicki periodicity elements as living in
\[
\kappa \in \Om^{\infty}\GW(\Dfree(R),(\QF^{-\gev})\qshift{-2}) \quad \text{and} \quad \rho \in \Om^{\infty+2}\L(\Dfree(R),\QF^{-\gev}).
\]
If $\lambda$ is an invertible form parameter of symmetric type on $A$ with underlying $R$-module with involution $M$, then $(\QF^{\mathrm g\lambda}_M)\qshift{1-\sigma} = \QF^{\mathrm g \check{\lambda}}_{-M}$, see Proposition~\refone{proposition:form-parameter-shift}, and we find that multiplication with the periodicity classes
gives the equivalences
\[
\mathbb S^1 \otimes \mathrm V(\D^c(A),\QF^{\mathrm g\lambda}_{M}) \longrightarrow \mathrm U(\D^c(A),\QF^{\mathrm g\check{\lambda}}_{-M}) \quad \text{and} \quad \mathbb S^{2} \otimes \L(\D^c(A),\QF^{\mathrm g\lambda}_{M}) \longrightarrow \L(\D^c(A),\QF^{\mathrm g\check{\lambda}}_{-M}),
\]
from Corollary~\reftwo{corollary:formparakarper} and%
~\reftwo{corollary:karoubi-periodicity-form-parameters}, the former of which is very much in the spirit of Karoubi's original formulation in \cite{karoubi-periodicity}. 

\end{enumerate}
\end{example}

\begin{remark}
To avoid confusion, let us mention the following misalignment of the example above with \cite{karoubi-periodicity}. Karoubi predicted the existence of an element $x \in \Omega^{\infty+2}\GW^{-\mathrm g\sym}(\mathbb Z)$ analogous to Ranicki's periodicity element $\rho \in \Omega^{\infty+2}\L^{-\mathrm g\sym}(\mathbb Z)$, multiplication with which ought to induce an equivalence
\[
\mathbb S^1 \otimes \mathrm V(\D^c(A),\QF^{\mathrm g\lambda}_{M})\longrightarrow\mathrm U(\D^c(A),\QF^{\mathrm g\check{\lambda}}_{-M}),
\]
in place of our $\kappa \in \Omega^\infty\GW(\Dfree(\mathbb Z),(\QF^{-\mathrm g\ev})\qshift{-2}) = \Omega^\infty\GW(\Dperf(\mathbb Z),(\QF^{-\mathrm g\sym})\qshift{-2})$. The only multiplication of the above shape induced by such a class $x$ seems to be that with $\partial^2(x) \in \Omega^\infty\GW(\Dperf(\mathbb Z),(\QF^{-\gs})\qshift{-2})$. However, any element of $\GW_0(\Dperf(\mathbb Z),(\QF^{-\gs})\qshift{-2})$ that lies in the image of $\partial$ vanishes in $\K_0(\mathbb Z)$. Thus multiplication with it kills all hyperbolic forms and it consequently cannot induce an equivalence as above for all $A$. 

In fact, in Proposition~\refthree{proposition:gwz-low-deg} we calculated $\GW_2^{-\mathrm g\sym}(\mathbb Z) \cong \mathbb Z$ and from the fundamental fibre sequence one reads off that the bordification map to $\L_2^{-\mathrm g\sym}(\mathbb Z) = \mathbb Z$ is injective with index $4$. Similarly one finds 
\[
\GW_0(\Dperf(\mathbb Z),(\QF^{-\mathrm g\sym})\qshift{-2}) \cong \mathbb Z \times_{\mathbb Z/2} \mathbb Z
\]
via signature and Euler characteristic. It follows that $\partial^2$ takes classes in $\GW_2^{-\mathrm g\sym}(\mathbb Z)$ to multiples of $(4,0)$, while $\kappa$ corresponds to $(1,-1)$, and also that Ranicki's periodicity class $[\rho] \in \L_2^{-\g\sym}(\mathbb Z)$ does not lift to $\GW_2^{-\g\sym}(\mathbb Z)$. 

In total we consider the above implementation of the equivalence in Karoubi's conjecture as the closest possible approximation to the original formulation.
\end{remark}

Let us now extend the discussion of Ranicki's periodicity isomorphism 
\[
\rho^2 \cdot -\colon  \L(\D^c(R),\QF^\sym_M) \longrightarrow \SS^{-4}\otimes \L(\D^c(R),\QF^\sym_{M}).
\]
and its relatives to $n\sigma$-oriented $\Einf$-rings $R$ for $n \geq 1$, in place of ordinary rings (where $n=2$). 
To this end, recall from Definition~\refone{def:nsigmaorpaperI} that an $n\sigma$-orientation on an $\Einf$-ring with involution $(R,\tau)$ is a
$\Ct$-equivariant $R$-linear equivalence $\eta\colon R \rightarrow \mathbb S^{n\sigma-n} \otimes R$ %
lifting the identity on underlying $R$-modules, where $\Ct$ acts on the left term via $\tau$ on the right term via the corresponding tensor of actions. %
Equivalently, the data of such an equivalence is a lift 
\[
\eta\in (\QF^{\sym}_{\tau})\qshift{n\sigma-n}(R) = \Omega^\infty (\mathbb S^{n\sigma-n} \otimes R)^\hC
\]
of the unit $1 \in \Omega^\infty R$ along the forgetful map $\Omega^\infty (\mathbb S^{n\sigma-n} \otimes R)^\hC \rightarrow \Omega^\infty R$. Note that if $R$ is discrete, then it indeed carries a unique $2\sigma$-orientation, since this forgetful map is then an equivalence. 

In fact, it is convenient to upgrade the notion of orientation to general Poincaré $\infty$-categories:

\begin{definition}
For a natural number $n \geq 1$, an $n\sigma$-orientation on a commutative algebra $(\D,\Phi)$ in $\Catp$ is an element $\eta \in \Omega^\infty\Phi\qshift{n\sigma-n}(\mathbb 1)$ lifting the unit of
\[
\Omega^\infty \Bil_{\Phi\qshift{n\sigma-n}}(\mathbb 1, \mathbb 1) \simeq \Omega^\infty \Bil_{\Phi}(\mathbb 1, \mathbb 1)
\]
under the forgetful map $\Phi\qshift{n\sigma-n} \Rightarrow \Bil_{\Phi\qshift{n\sigma-n}} \circ \Delta$. 
\end{definition}

An $n\sigma$-orientation on an $\Einf$-ring $R$ with involution $\tau$ is now by definition the same thing as an $n\sigma$-orientation on $(\Dperf(R),\QF^\sym_\tau)$.

\begin{construction}
An $n\sigma$-orientation $\eta$ on $(\D,\Phi)$ determines a Poincaré object $o_\eta = (\one,\eta) \in \Poinc(\D,\Phi\qshift{n\sigma-n})$, multiplication with which gives Poincaré functors
\[
(\D,\mathbb \Phi\Phi\qshift{-n\sigma-n}) \longrightarrow (\D,\Phi\qshift{-2n}) \quad \text{and} \quad (\D,\mathbb \Phi\Phi\qshift{-n\sigma+n}) \longrightarrow (\D,\Phi)
\]
covering the identity of $\D$, and thus
\[
o_{\eta}\kappa^n \in \Om^{\infty}\GW(\D,\Phi\qshift{-2n}) \quad \text{and} \quad o_{\eta}\rho^n \in \Om^{\infty+2n}\L(\D,\Phi),
\]
which we will refer to as the \emph{potential periodicity elements} of $\eta$. For any module $(\C,\QF)$ over $(\D,\Phi)$ they induce maps
\[
o_{\eta}\kappa^n \cdot - \colon \GW(\C,\QF) \longrightarrow \GW(\C,\QF\qshift{-2n}) \quad \text{and} \quad  o_{\eta}\rho^n \cdot - \colon \L(\C,\QF) \longrightarrow \SS^{-2n}\otimes \L(\C,\QF).
\]
\end{construction}

\begin{examples}\
\begin{enumerate}
\item
For an $\Einf$-ring spectrum with involution $(R,\tau)$ an $n\sig$-orientation $\eta$ on $(\Modf{R},\QF^{\sym}_{\tau})$ amounts to an $n\sigma$-orientation on $(R,\tau)$. In this case the Poincaré object $o_{\eta}$ is invertible (since for symmetric Poincaré $\infty$-categories this is detected on underlying objects) and the composite equivalence
\[
(\Modf{R},\QF^{\sym}_{\tau}) \xrightarrow{\Om^n} (\Modf{R},(\QF^{\sym}_{\tau})\qshift{-n-n\sig}) \xrightarrow{o_{\eta}} (\Modf{R},(\QF^{\sym}_{\tau})\qshift{-2n})
\]
coincides (up to a shift) with the equivalence of Corollary~\refone{corollary:periodicity-poincare} specialised to the case at hand, since $\eta\colon R \rightarrow \mathbb S^{n\sigma-n} \otimes R$ is $R$-linear.
In particular, the potential periodicity elements 
\[
o_{\eta}\kappa^n \in \Om^{\infty}\GW(\Modf{R},(\QF^\sym_\tau)\qshift{-2n}) \quad \text{and} \quad o_{\eta}\rho^n \in \Om^{\infty+2n}\L(\Modf{R},\QF_\tau^\sym)
\]
are invertible, 
and we find periodicity equivalences 
\[
o_\eta\rho^n \cdot -\colon  \L(\Mod^c(R),\QF^r_M) \longrightarrow \SS^{-2n}\otimes \L(\Mod^c(R),\QF^r_{M})
\]
for $r \in \{\mathrm q,\sym\}$ recovering the ones from Corollary~\reftwo{corollary:Ranickipergeneral}. The same is true for any $\Eone$-$R$-algebra $A$ with invertible $R$-module with involution.

A completely analogous statement also holds, whenever $\D$ is any symmetric monoidal $\infty$-category and $\Phi = \QF^{\sym}_{\Dual}$ is the symmetric Poincaré structure associated to some symmetric monoidal duality on $\D$ (that is, a duality yielding a commutative algebra object in $\Catps$), and $(\C,\QF)$ is a module over $(\D,\Psi)$ for which $\QF$ is either symmetric or quadratic. 
\item
If $(R,\tau)$ is a connective $\Einf$-ring spectrum with involution then an $n\sigma$-orientation uniquely lifts along the map $\QF_\tau^{\geq 0} \rightarrow \QF_\tau^\sym$, since the map
\[
\Omega^\infty (\QF_{\tau}^{\geq 0})\qshift{n\sigma-n}(R) \longrightarrow \Omega^\infty (\QF_{\tau}^{\sym})\qshift{n\sigma-n}(R)
\]
is then an equivalence. The potential periodicity elements 
\[
o_{\eta}\kappa^n \in \Om^{\infty}\GW(\Modf{R},(\QF^{\geq 0}_\tau)\qshift{-2n}) \quad \text{and} \quad o_{\eta}\rho^n \in \Om^{\infty+2n}\L(\Modf{R},\QF_\tau^{\geq 0})
\]
are then generally not invertible, the latter giving rise to the stabilisation maps
\[
o_{\eta}\rho^n \cdot - \colon \L(\Mod^c(R),\QF^{\geq k}_{\tau}) \to \SS^{-2n}\otimes \L(\Mod^c(R),\QF^{\geq k}_{\tau}).
\]
for arbitrary $k \in \mathbb Z \cup \{\pm \infty\}$. 

\item In the situation of the previous point, the polarisation map
\[
\Omega^\infty(\QF^{\geq n}_\tau)\qshift{n\sigma-n}(R) \longrightarrow \Omega^\infty(\QF^{\geq 0}_\tau)\qshift{n\sigma-n}(R)
\]
is an equivalence, yielding a lift $o^\mathrm g_\eta \in \Poinc(\Modf{R},(\QF^{\geq n}_\tau)\qshift{n\sigma-n})$ of $o_\eta$ and consequently a factorisation
\[
(\Modf{R},(\QF^{\geq 0}_\tau)\qshift{-n\sigma+n}) \xrightarrow{o_\eta^g \cdot -} (\Modf{R},\QF^{\geq n}_\tau) \xrightarrow{\mathrm{pol}} (\Modf{R},\QF^{\geq 0}_\tau),
\]
of the multiplication with $o_\eta$, whose first map is an equivalence by direct inspection of linear parts. Multiplication with $o_\eta^\g$ then recovers the equivalences 
\[
o^\mathrm g_{\eta}\rho^n \cdot - \colon \L(\Mod^c(R),\QF^{\geq k}_{\tau}) \to \SS^{-2n}\otimes \L(\Mod^c(R),\QF^{\geq k+n}_{\tau}).
\]
induced by the equivalences of Corollary~\refone{corollary:periodicity-poincare-truncated}, refining the stabilisation maps from the previous point.
\item 
The simplest example that is not quite of this nature 
is $\Dperf(\ZZ)$ equipped with the Poincar\'e structure $\QF^{\g\fr}_{\ZZ}$ associated to the Frobenius
form parameter of Example~\reffour{example:form_parameter_over_Tate_Poincare}.
In this case $(\QF^{\g\fr}_{\ZZ})\qshift{2\sigma-2}(\ZZ)$ evaluates to the (Eilenberg-Mac Lane spectrum of the) pullback of
\[
\ZZ^\hC \longrightarrow \ZZ^\tC \longleftarrow \ZZ[-2] \times_{\ZZ/2[-2]} \tau_{\geq -2}\ZZ^\tC.
\]
In particular, at the level of infinite loop spaces, the fibre over the unit of $\Bil_{\QF^{\g\fr}_{\ZZ}}(\ZZ, \ZZ)^{\hC} \simeq \ZZ^\hC$ is trivial, giving a unique $2\sigma$-orientation $\eta$ to $\QF^{\g\fr}_{\ZZ}$ lifting the $2\sigma$-orientation of $\QF^{\sym}_{\ZZ}$. The associated potential periodicity element $o_{\eta}\rho^2  \in \Om^{\infty}\L^{\g\fr}(\ZZ)$ then lifts the one in $\Om^{\infty+4}\L^{\sym}(\ZZ)$; in fact the map between these two $\Om^{\infty+4}\L$ spaces is an equivalence by Example~\refthree{example:dotto-ogle}, and so such a potential periodicity element is unique as well.
Since for every ring with form parameter $(R,M,\lambda)$, the spectrum $\L(\D^c(R),\QF^{\mathrm g\lambda}_M)$ is a module over $\L^{\g\fr}(\ZZ)$
by Proposition~\reffour{prop:tensormackeyvstensorfpm}, we learn that these spectra always carry stabilisation maps
\[
\L(\D^c(R),\QF^{\mathrm g\lambda}_M) \longrightarrow \SS^{-4}\otimes \L(\D^c(R),\QF^{\mathrm g\lambda}_M)
\]
induced by the multiplication with the potential periodicity element $o_{\eta}\rho^2 \in \Om^{\infty+4}\L^{\g\fr}(\ZZ)$ 
(and by its image in $\Om^{\infty+4}\L^{\mathrm g\sym}(\mathbb Z)$ when $\lambda$ is of symmetric type). 
\item 
Not every lax-symmetric monoidal Poincar\'e structure over an $n\sigma$-oriented $\Einf$-ring admits an $n\sigma$-orientation itself: The simplest example is the Tate Poincar\'e structure on $\ZZ$: In this case
$(\QF^\tate_{\ZZ})\qshift{2\sigma-2}(\ZZ)$ evaluates to the pullback of
\[
\ZZ^\hC \longrightarrow \ZZ^\tC \longleftarrow \ZZ[-2],
\]
which does not admit the required lift of the element $1 \in \Om^{\infty}(\ZZ^\hC)$. And indeed, the forthcoming calculation of $\L^{\tate}(\ZZ)$ in \cite{TateL} show that Ranicki's periodicity element in $\Om^{\infty+4}\L^{\g\fr}(\ZZ)$ does not admit a further lift to this spectrum.
\end{enumerate}
\end{examples}

\subsection{Multiplicativity of the circle decomposition theorem}
\label{subsection:mult-shaneson}%

Let $\F\colon \Catp \to \A$ be a functor valued in a cocomplete $\infty$-category $\A$ %
and recall from \S\reffour{subsection:shaneson} the assembly map
\[
\B\ZZ \otimes \F(\C,\QF) \to \F((\C,\QF)_{\B\ZZ})
\]
defined for any Poincaré $\infty$-category $(\C,\QF)$ as the colimit-exchange map. In this section we refine this map to a symmetric monoidal transformation in the case where $\A$ is symmetric monoidal and $\F$ is lax symmetric monoidal. We begin with the following general observation.

\begin{proposition}
Let $\A$ be a cocomplete symmetric monoidal $\infty$-category in which the tensor product preserves colimits in each variable, $\F\colon \Catp \to \A$ a lax symmetric monoidal functor, $\I$ a small symmetric monoidal $\infty$-category and $\varphi \colon \I \to \Catp$ a functor. Then the associated assembly map
\[
\colim_\I \F\circ \varphi \to \F(\colim_\I \varphi)
\]
is the component of a lax symmetric monoidal transformation among lax symmetric monoidal $\A$-valued functors in $\varphi$. In particular, if $\varphi$ is lax symmetric monoidal, then the map 
\[
\colim_\I \F\circ \varphi \to \F(\colim_\I \varphi)
\]
is one of commutative algebra objects in $\A$.
\end{proposition}
\begin{proof}
Recall from Remark~\reffour{remark:left-kan-monoidal} that under the assumptions of the proposition the colimit functor refines to a symmetric monoidal adjunction
\[
\colim_\I \colon \Fun(\I,\Catp)^{\Day} \adj \Catp \cocolon \const
\]
where the symmetric monoidal structure on the left is given by Day convolution, %
$\colim_\I$ is symmetric monoidal, %
and $\const$, the constant diagram functor, is canonically lax symmetric monoidal. In particular, the associated unit and counit maps are lax symmetric monoidal transformations among lax symmetric monoidal functors.
Since $\F$ is lax symmetric monoidal,
the associated commutative square 
\[
\begin{tikzcd}
\Catp \ar[r,"\F"]\ar[d,"\const"'] & \A \ar[d,"\const"] \\
\Fun(\I,\Catp)^{\Day} \ar[r,"\F_*"] & \Fun(\I,\A)^{\Day}
\end{tikzcd}
\]
of lax symmetric monoidal functors determines a (Beck-Chevalley type) natural transformation
\[
\colim_\I \circ \F_* \Rightarrow \colim_\I \circ\F_*\circ \const\circ\colim_\I = \colim_\I\circ \const\circ \F \circ\colim_\I \Rightarrow \F\circ\colim_\I
\]
which, being a composite of lax symmetric monoidal transformations, gives the desired lax symmetric monoidal transformation.
\end{proof}

Combining this with Corollary~\reffour{corollary:iterated-shaneson-splitting} we conclude:
\begin{corollary}
\label{corollary:iteratedassmonoidal}%
Let $\A$ be a cocomplete symmetric monoidal $\infty$-category in which the tensor product preserves colimits in each variable and let $\F\colon \Catp \to \A$ be a lax symmetric monoidal bordism invariant Karoubi-localising functor. %
Then for every $n \geq 0$ the assembly map
\[
\B\ZZ^n \otimes \F(-) \Rightarrow \F((-)_{\B\ZZ^n})
\]
is an equivalence of lax symmetric monoidal functors $\Catp \to \A$. %
\end{corollary}

\begin{example}
\label{example:group-rings}%
For $(R,\sig)$ a commutative ring with involution the symmetric monoidal Poincaré $\infty$-category $(\Dperf(R),\QF^{\sym}_{\sig})_{\B\ZZ^n}$ is equivalent to $(\Dperf(R[\ZZ^n]),\QF^{\sym}_{\bar{\sig}})$ where we view $R[\ZZ^n]$ as a ring with involution $\bar{\sig}$ extending $\sig$ and inducing the inversion on $\ZZ^n$. By Corollary~\reffour{corollary:iteratedassmonoidal} we then obtain 
a canonical equivalence of $\Einf$-rings
\[
\KL^{\sym}(R[\ZZ^n],\bar{\sig}) \simeq  \KL^{\sym}(R,\sig)[\B\ZZ^n].
\]
where for an $\Einf$-ring $A$ we denote by
$A[\B\ZZ^n]=\B\ZZ^n \otimes A$ the associated group-algebra.
\end{example}

For not necessarily bordism invariant Karoubi-localising lax symmetric monoidal functors $\F$, like $\KK$ or $\KGW$, the map 
\[
\B\ZZ^n \otimes \F(-) \Rightarrow \F((-)_{\B\ZZ^n})
\]
is still one of lax symmetric monoidal functors, and an equivalence in certain cases:

\begin{example}
\label{example:regular-noetherian}%
For $(R,\sig)$ a regular noetherian commutative ring with involution 
the map of $\Einf$-rings
\[
\KGW^{\sym}(R,\sig)[\B\ZZ^n] \to \KGW^{\sym}(R[\ZZ^n],\bar{\sig})
\]
induced by the assembly map, is an equivalence. Indeed, in these cases the inductively appearing complementary summands $\N\KK(\Dperf(R[\ZZ^i]))$ vanish. 
\end{example}

Let us now specialise to the case $n=1$, for which Examples \reffour{example:group-rings} and \reffour{example:regular-noetherian} motivate one to understand the $\Einf$-algebra structure on 
the group $\Einf$-$A$-algebra $A[\B\ZZ]$ associated to an $\Einf$-ring spectrum $A$ and the $\Einf$-group $\B\ZZ$. Being a group algebra, it comes equipped with an augmentation $A[\B\ZZ] \to A$ induced by the terminal map $\B\ZZ \to \ast$ whose augmentation ideal is given by $\Sigma A$. Hence $A[\B\ZZ]$ is equivalent to $A\oplus\Sigma A$ as an $A$-module over and under $A$, and one may ask if this equivalence is multiplicative when viewing $A \oplus \Sigma A$ as
the split square-zero extension.

\begin{lemma}
\label{lemma:sq0}%
Let $A$ be an $\Einf$-algebra. Refinements of the $A$-module-over-under-$A$ equivalence $A[\B\ZZ] \simeq A \oplus \Sigma A$ to one of augmented  
$\Ek{n}$-$A$-algebras naturally correspond to factorisations of the unit $\SS \to A$,
as a map of spectra, 
through the map $\SS \to \SS^{-n-1}\otimes\Sig^{\infty}K(\ZZ,n+1)$. In particular, for $n=\infty$, such refinements correspond to factorisations of the unit though the Eilenberg-MacLane spectrum $\GEM \ZZ = \colim_n\SS^{-n-1}\otimes\Sig^{\infty}K(\ZZ,n+1)$.
\end{lemma}
\begin{proof}
For an $A$-module $M$, we claim that $\fib(\gl_1(A\oplus M) \to \gl_1(A))$ is naturally equivalent to $\Omega^\infty M$: Indeed, as functors of $M$ both expressions are right adjoints $\Mod_A \to \Grp_{\Einf}$ among pointed presentable $\infty$-categories which agree after post-composition with the right adjoint $\Grp_{\Einf} \to \Sps_*$, and hence must agree before post-composition since $\Mod_A$ is stable and the induced functor $\Sp(\Grp_{\Einf}) \to \Sp(\Sps_*)$ is an equivalence.
Given this claim, the desired result now follows from the following calculation of the mapping space of augmented $\Ek{n}$-$A$-algebras:
\begin{align*}
\Map_{\Alg^{\mathrm{aug}}_{\Ek{n}}(\Mod_A)}(A[\B\ZZ],A\oplus \Sig A) &= \Map_{\Grp_{\Ek{n}}}(\B\ZZ,\fib[\gl_1(A\oplus \Sig A) \to \gl_1(A)]) \\
&= \Map_{\Grp_{\Ek{n}}}(\B\ZZ,\Om^{\infty}(\Sigma A))\\
&= \Map_{\Grp_{\Ek{n+1}}}(\ZZ,\Om^\infty A) \\
&= \Map_{\Sp}(\Om^{n+1}\Sig^{\infty}K(\ZZ,n+1),A).\qedhere
\end{align*}
\end{proof}

As a sample application of the above, we obtain:
\begin{corollary}
Let $(R,\sig)$ be a discrete commutative ring with involution. Then there is a preferred equivalence
\[
\KL^\s(R[\ZZ],\bar{\sig}) \simeq \KL^\s(R,\sig) \oplus \Sig \KL^\s(R,\sig)
\]
as augmented $\Eone$-$\KL^\s(R,\sig)$-algebras which in general does not extend to one of $\Ek{2}$-algebras. After localisation at 2, however, there is such an equivalence as $\Einf$-algebras.
\end{corollary}
\begin{proof}
Note that for $n=1$, the condition of Lemma~\reffour{lemma:sq0} is equivalent to equipping $A$ with a complex orientation and recall that $\L^\s(\ZZ)$ admits a preferred complex orientation, the Sullivan-Ranicki orientation. In particular, from the $\Einf$-map $\L^\s(\ZZ) \to \KL^\s(R,\sig)$, we deduce that $\KL^\s(R,\sig)$ is equipped with an induced complex orientation, so that the first claim follows from combining Example~\reffour{example:group-rings} with Lemma~\reffour{lemma:sq0} in case $n=1$. Then note that for $n=\infty$, the condition of Lemma~\reffour{lemma:sq0} holds whenever the underlying spectrum of $A$ admits a $\ZZ$-module structure. This is the case for $\KL^\s(R,\sig)_{(2)}$ as is again transported from the Sullivan-Ranicki orientation and the classical fact that $\MSO_{(2)}$ admits a $\ZZ$-module structure. Finally, to see the claim about $\Ek{2}$-structures, we show that $\L^\s(\ZZ)[\B\ZZ]$ is not $\Ek{2}$-equivalent to the trivial square-zero extension $\L^\s(\ZZ) \oplus \Sig \L^\s(\ZZ)$. By Lemma~\reffour{lemma:sq0}, it suffices to show that there is no map $K(\ZZ,3) \to \Omega^\infty \Sigma^3\L^\s(\ZZ)$ inducing an isomorphism on $\pi_3$. If such a map exists, we conclude from the equivalence $\L^\s(\ZZ)[1/2] \simeq \KO[1/2]$ that $K(\ZZ,3)$ is, after localising away from 2, a retract of $\Omega^\infty \Sigma^3 \KO[1/2]$ which in turn is a retract of $\Omega^\infty \Sigma^3 \ku[1/2] \simeq \SU[1/2]$. But the $\ZZ$-homology of $\SU[1/2]$ is torsion-free, while $H_8(K(\ZZ,3),\ZZ)$ contains non-trivial $3$-torsion.
\end{proof}

\begin{example}
\label{remark:L-group-sq0}%
As a consequence of the above, for a commutative ring $(R,\sig)$ with involution, $\KL^{\sym}_0(R,\sig)[\B\ZZ]$ is a square-zero extension of $\KL^{\sym}_0(R,\sig)$ by its module $\KL_{-1}^\s(R,\sig)$. For instance, for $R$ a Dedekind domain whose fraction field is a global field, we obtain that $\KL^\s_0(R[\ZZ],\bar{\id})$ is isomorphic to the split square-zero extension $\W^\s(R) \oplus \Pic(R)/2$, see Corollaries~\refthree{corollary:compute-L-groups-global-field} and~\refthree{corollary:compute-L-groups-global-Dedekind-ring-characteristic-2}.
\end{example}

\begin{remark}
Let $(R,\sig)$ be a regular noetherian ring with involution for which $2 \neq 0$ holds. Then $\KGW^\s(R,\sig) = \GW^\s(R,\sig)$ is not complex orientable: indeed, in $\pi_1$ of any complex orientable ring spectrum the image of the Hopf element $\eta \in \pi_1\SS$ vanishes, whereas 
the above assumptions %
imply that this image is non-zero in $\K_1(R)$, hence already in $\GW_1^\s(R,\sig)$.
In particular, Lemma~\reffour{lemma:sq0} implies that the $\Ek{0}$-equivalence $\GW^\s(R[\ZZ],\bar{\sig}) \simeq \GW^\s(R,\sig) \oplus \Sig\GW^\s(R,\sig)$ does not extend to one of $\Eone$-algebras. Nevertheless, in this situation,
$\GW_0^\s(R[\ZZ],\bar{\sig})$ is isomorphic to the square-zero extension $\GW_0^\s(R,\sig) \oplus \GW_{-1}^\s(R,\sig)$. Indeed, to see this, it suffices to argue that the augmentation ideal $\GW_{-1}^\s(R,\sig)$ has trivial multiplication. But the map $\GW_{-1}^\s(R,\sig) \to \L_{-1}^\s(R,\sig)$ is an isomorphism and is in addition compatible with its multiplication since $\GW \to \L$ is lax symmetric monoidal. Hence, the claim follows from Remark~\reffour{remark:L-group-sq0}. 
\end{remark}

\subsection{Multiplicative structures on normal L-spectra}
\label{subsection:normal-L}%

The purpose of the present section is to record that for any $\Einf$-Poincar\'e $\infty$-category $(\C,\QF)$ the associated map $\L(\C,\QF^\qdr) \rightarrow \L(\C,\QF)$ exhibits the source as a Smith ideal in the target, and consequently that its cofibre, the normal $\L$-spectrum inherits a multiplicative structure. To set this up, first observe that for any Poincar\'e $\infty$-category $(\C,\QF)$ we have a cube 
\[
\begin{tikzcd}
[column sep=0.5em, row sep=0.7em]
  & \L(\C,\QF^\qdr)\ar[rr]\ar[dd] && \KL(\C,\QF^\qdr)\ar[dd]\\
\GW(\C,\QF^\qdr)\ar[ru]\ar[rr, crossing over]\ar[dd] && \KGW(\C,\QF^\qdr)\ar[ur]\\
  & \L(\C,\QF)\ar[rr] && \KL(\C,\QF)\\
\GW(\C,\QF)\ar[rr] \ar[ru] && \KGW(\C,\QF)\ar[uu,<-, crossing over] \ar[ru]\\
\end{tikzcd}
\]
whose front, back and side faces are exact by Theorem \reffour{theorem:Karoubi-pullback-squares} and Corollary~\reftwo{corollary:tate-square-L}.
In particular, we find a square of equivalences
\[
\begin{tikzcd}
\cof(\GW(\C,\QF^\qdr) \rightarrow \GW(\C,\QF)) \ar[r] \ar[d] & \cof(\L(\C,\QF^\qdr) \rightarrow \L(\C,\QF)) \ar[d] \\
\cof(\KGW(\C,\QF^\qdr) \rightarrow \KGW(\C,\QF)) \ar[r] & \cof(\KL(\C,\QF^\qdr) \rightarrow \KL(\C,\QF))
\end{tikzcd}
\]
natural in $(\C,\QF) \in \Catp$.

\begin{definition}
We denote by $\NL \colon \Catp \rightarrow \Sp$, \emph{normal $\L$-theory}, the functor given by any of the four formulae in the square above.
\end{definition}

This notion was first introduced by Ranicki.
We will often use the abbreviation $\NL^r(R,M)$ for $\NL(\Dperf(R),\QF^r_M)$, where $M$ is an invertible module with involution over $R$.

Normal $\L$-spectra are closely connected to real topological cyclic homology for Poincar\'e $\infty$-categories and will be more thoroughly analysed in the forthcoming \cite{Harpaz-Nikolaus-Shah}, one of the main results of which equates them with the geometric fixed points of $\mathrm{TCR}(\C,\QF)$ under suitable assumptions.
Here, we will equip $\NL$ with a lax symmetric monoidal structure under $\KL \colon \Catp \rightarrow \Sp$ (and thus also under $\GW,\KGW$ and $\L$). 
To this end, we first recall that the pushout-tensor product 
\[
(e \rightarrow e') \otimes (f \rightarrow f') = [e\otimes f' \amalg_{e \otimes f} e' \otimes f \rightarrow e' \otimes f']
\]
on the arrow category of a symmetric monoidal $\infty$-category $\E$ with finite colimits is nothing but Day convolution 
with respect to the symmetric monoidal structure on $[1]$ given by the minimum (or equivalently, by the product). We will denote the corresponding symmetric monoidal $\infty$-category by 
\[
\Ar(\E)^\Sm := \Fun([1],\E)^\Day.
\]
Its unit is the unique arrow $\emptyset \to \one$ from the initial object to the unit of $\E$, and its algebras %
are precisely the Smith ideal inclusions.
There is also the pointwise symmetric monoidal structure on arrows, and it will be denoted by $\Ar(\E)^{\mathrm{pw}}$; it is similarly given by Day convolution with respect to the maximum (or equivalently, the coproduct).

\begin{proposition}
\label{proposition:qdr-lax-monoidal}%
The functor 
\[
\Catp \longrightarrow \Ar(\Catp)^{\Sm}, \quad (\C,\QF) \longmapsto [(\C,\QF^\qdr) \rightarrow (\C,\QF)]
\]
admits a canonical lax symmetric monoidal structure.
\end{proposition}  
\begin{proof}
For $(\C,\QF) \in \Catp$, the natural map 
\[
(\C,\QF^{\qdr}) \to (\C,\QF) \otimes (\Spaf,\QF^\qdr),
\]
adjoint to the equivalence on underlying stable $\infty$-categories with duality on both sides under the adjunction of Proposition~\refone{proposition:adjunction-for-poincare-categories},
is an equivalence by the explicit formula for the tensor product in $\Catp$.
In then follows that the object
\[
[(\Spaf,\QF^\qdr) \to (\Spaf,\QF^\uni)] \in \Ar(\Catp)
\]
is idempotent as an $\Ek{0}$-algebra in $\Ar(\Catp)^{\Sm}$. Hence, it essentially uniquely refines to an idempotent $\Einf$-algebra and tensoring with it determines a smashing localisation $L \colon \Ar(\Catp)^{\Sm} \to \Ar(\Catp)^{\Sm}$ which is canonically a lax symmetric monoidal functor. Furthermore, the left Kan extension along the symmetric monoidal functor $\{1\} \to [1]$ induces a symmetric monoidal functor $\Catp \to \Ar(\Catp)^{\Sm}$. The composite of these two functors is then given by
\[
\Catp \to \Ar(\Catp)^{\Sm} \xrightarrow{L} \Ar(\Catp)^{\Sm}, \quad (\C,\QF) \mapsto (0\to (\C,\QF)) \mapsto (\C,\QF^\qdr) \to (\C,\QF),
\]
and so the proof is complete.
\end{proof}

We thus find that for any lax symmetric monoidal $\F \colon \Catp \rightarrow \Sp$ the assignment
\[
(\C,\QF) \longmapsto [\F(\C,\QF^\qdr) \rightarrow \F(\C,\QF)]
\]
upgrades to a lax symmetric monoidal functor $\Catp \rightarrow \Ar(\Sp)^{\Sm}$, in fact naturally in $\F$. In particular, for any $(\C,\QF) \in \Alg_{\Einf}(\Catp)$ we find that $\F(\C,\QF^\qdr)$ is a Smith ideal in the $\Einf$-ring $\F(\C,\QF)$. Finally, we note that taking the cofibre of an arrow induces a symmetric monoidal equivalence
\[
\cof \colon \Ar(\Sp)^{\Sm} \to \Ar(\Sp)^{\mathrm{pw}},
\]
see e.g.\ the argument for \cite{almostpaper}*{Proposition 2.3}.
We obtain:
\begin{corollary}
\label{corollary:normallmult}%
The tautological transformation $\KL \Rightarrow \NL$ of functors $\Catp \rightarrow \Sp$ admits a canonical lax symmetric monoidal structure extending that on the source constructed in Proposition \reffour{corollary:kl-monoidal}. Moreover, $\NL$ is Karoubi-localising and bordism invariant.
\end{corollary}
\begin{proof}
The first part follows immediately from Proposition~\reffour{proposition:qdr-lax-monoidal}. Then $\NL$ is a module over $\KL$ and hence bordism invariant. Finally, the functor $(\C,\QF) \mapsto (\C,\QF^\qdr)$ preserves Poincar\'e-Karoubi squares, for instance since it is given by the functor $- \otimes (\Spaf,\QF^\qdr)$ as noted above, so we may appeal to Proposition~\reffour{prop:closure-properties-of-Rs}. It follows that $\NL$ is Karoubi-localising as it is a cofibre of two Karoubi-localising functors.
\end{proof}

Ranicki already observed (see, e.g., \cite{RanickiTSO}) that normal $\L$-spectra of commutative rings (based on the symmetric Poincar\'e structures) admit homotopy associative and commutative multiplications. An upgrade regarding $\Einf$-structures on normal L-spectra was envisioned by Lurie \cite{Lurie-L-theory}*{Lecture 37}.

\begin{example}
\label{example:F2}%
Perhaps the simplest non-trivial example of a normal $\L$-spectrum is $\NL^{\sym}(\FF_2)$, which is a 2-periodic ring spectrum all of whose homotopy groups are $\FF_2$. 
Since $\FF_2$ has cohomological dimension $0$, Corollaries~\refthree{corollary:genuine-is-quadratic} and \refthree{corollary:genuine-is-symmetric} together implies that
\[
\NL^{\geq m}(\FF_2) \simeq \tau_{\geq 2m-1}\NL^{\sym}(\FF_2) .
\]
In particular, for $m=0$ we get that the $(-1)$-connective cover map
\[
\NL^{\gs}(\FF_2) = \tau_{\geq -1}\NL^{\sym}(\FF_2) \to \NL^{\sym}(\FF_2)
\]
is an $\Einf$-ring spectrum map. %
The ring structure on $\NL^{\sym}_*(\FF_2)$ is thus
\[
\NL^{\sym}_*(\FF_2) = \FF_2[\beta,\beta^{-1},\eps]/\eps^2 = \L^{\sym}_*(\FF_2)[\eps]/\eps^2 \quad\text{and}\quad \NL^{\gs}_*(\FF_2) = \FF_2[\beta,\eps]/\eps^2 = \L^{\gs}_*(\FF_2)[\eps]/\eps^2
\]
with $|\beta| = 2$ and $|\eps| =-1$.
In particular, $\NL^{\sym}(\FF_2)$ has the same homotopy ring as $\L^{\sym}(\FF_2)^{S^1}$, though it is not a-priori clear if and why there would be any map relating the two.
This point is elucidated in forthcoming work \cite{Harpaz-Nikolaus-Shah}, where it is shown that there is always a canonical map $\L(\C,\QF)^{S^1} \to \NL(\C,\QF)$. In the case of $\FF_2$, this map induces equivalences of $\Einf$-rings
\[
\NL^{\sym}(\FF_2) \simeq \L^{\sym}(\FF_2)^{S^1} \quad\text{and}\quad \NL^{\gs}(\FF_2) \simeq \L^{\gs}(\FF_2)^{S^1} .
\]
Using the results of Heuts--Land \cite{formality}*{Theorem A \& C}, it follows that $\NL^\gs(\FF_2)$ and $\NL^\s(\FF_2)$ are split square zero extensions as augmented $\Eone$-algebras, but not
as augmented $\Ek{2}$-$\L^\gs(\FF_2)$-algebras. 
\end{example}

\begin{example}
\label{example:normal-L-of-Z}%
In the case of the integers, the normal symmetric $\L$-groups are given by
\[
\NL_i^{\sym}(\mathbb Z) \cong \begin{cases} \mathbb Z/8 & i \equiv 0 \ (4) \\ \mathbb Z/2 & i \equiv 1 \ (2) \\ 0 & i \equiv 2 \ (4)\end{cases}
\]
and the ring structure is $\NL^{\sym}_\ast(\mathbb Z) = \mathbb Z/8[t^{\pm 1},e,f]/(2e,2f,e^2,f^2,4-ef)$ where $|t|=4$, $|e| =1$ and $|f| = -1$, originally determined by Ranicki in unpublished work, see \cite{HLN}*{Appendix A} for a proof. 
Let us also recall from \cite{HLN}*{Corollary 6.2} that the map $\NL^\gs(\ZZ) \to \NL^\s(\ZZ)$ exhibits the source as the $(-1)$-connective cover of the target, just as is the case for $\FF_2$ recalled in Example~\reffour{example:F2}.
In particular, $\tau_{\geq-1}\NL^\s(\ZZ)$ admits an $\Einf$-ring structure, a phenomenon also originally observed in \cite{HLN}.

More generally, from Corollary~\refthree{corollary:genuine-is-quadratic} we find that $\NL^{\gs}(R,M)$ is $(-1)$-connective for every ordinary (in fact, connective $\Eone$-) ring $R$ and invertible module with involution $M$, and by Corollary~\refthree{corollary:genuine-is-symmetric} the resulting map
\[
\NL^{\gs}(R,M) \to \tau_{\geq -1}\NL^{\sym}(R,M),
\]
is an isomorphism in degrees $\geq d-1$ if $R$ is coherent of finite global dimension $d$. In a recent paper \cite{Schlichting-genuine}, Schlichting proved that if $R$ is a commutative Noetherian domain of finite Krull dimension either containing a field or smooth over a Dedekind domain and $M=R$ with the trivial $\Ct$-action, then the above map is in fact an equivalence.
In particular, for such $R$ %
we have that $\tau_{\geq -1}\NL^{\sym}(R)$ admits an $\Einf$-structure, yielding a non-trivial constraint on the $\Einf$-structure of $\NL^{\sym}(R)$. 
\end{example}

\begin{corollary}
\label{corollary:normalL8tors}%
Let $R$ be an $\Eone$-$\mathbb Z$-algebra
and $M$ an invertible module with involution over $R$. Then the normal $L$-spectra $\NL^{\geq m}(R,M)$ are naturally $\NL^{\gs}(\ZZ)$-modules for $m \in \ZZ \cup \{\infty,-\infty\}$. In particular, for every $i \in \ZZ$ and every $m,n \in \mathbb Z \cup \{\infty,- \infty\}$ with $m\geq n$, both kernel and cokernel of the natural maps
\[
\GW_i^{\geq m}(R,M) \longrightarrow \GW_i^{\geq n}(R,M)
\]
are $8$-torsion, and the same holds for any of the functors $\L,\KL,\KGW$ in place of $\GW$.
\end{corollary}
\begin{proof}
By Example~\reffour{example:ideal} and Corollary~\reffour{corollary:normallmult}, the normal L-spectra $\NL^{\geq m}(R,M)$ are module spectra over $\NL^{\gs}(\ZZ)$, and for $m\geq n$, the map $\NL^{\geq m}(R,M) \to \NL^{\geq n}(R,M)$ is one of $\NL^{\gs}(\ZZ)$-modules. Using then the equivalences
\begin{align*}
	\cof[\GW^{\geq m}(R,M) \to \GW^{\geq n}(R,M)] & \simeq \cof[\KGW^{\geq m}(R,M) \to \KGW^{\geq n}(R,M)] \\
		& \simeq \cof[\KL^{\geq m}(R,M) \to \KL^{\geq n}(R,M)] \\
		& \simeq \cof[\L^{\geq m}(R,M) \to \L^{\geq n}(R,M)] \\
		& \simeq \cof[\NL^{\geq m}(R,M) \to \NL^{\geq n}(R,M)]
\end{align*}
we find that all of these objects carry the structure of a $\NL^\gs(\ZZ)$-module. Hence the claim follows from the computation $\NL^\gs_0(\ZZ) \cong \ZZ/8$, see Example~\reffour{example:normal-L-of-Z}.
\end{proof}

Finally, we record the following invariance property of normal L-theory of $\ZZ$-algebras.
\begin{example}
Let $(A,\sig)$ be a $\ZZ$-algebra spectrum with involution in which 2 is invertible and $m\in \ZZ \cup \{\pm \infty \}$. Then $\NL^{\geq m}(A,\sig)$ vanishes since $\L^\qdr(A,\sig) \to \L^{\geq m}(A,\sig)$ is an equivalence. As a consequence, if $(R,\sig)$ is a discrete ring and $m\in \ZZ \cup \{\pm \infty \}$, the maps 
\[
\NL^{\geq m}(R,\sig) \to \NL^{\geq m}(R_{(2)},\sig) \to \NL^{\geq m}(R_2^\cwedge,\sig)
\]
are equivalences, that is, genuine normal L-spectra of $\ZZ$-algebras are insensitive to localising or (derived) completing at 2; this was implicitly already proven and used in \paperthree. Indeed, since $\NL$ is Karoubi-localising, this follows from the Poincar\'e-Karoubi square
\[
\begin{tikzcd}
	(\Dperf(R),\QF^{\geq m}_\sig) \ar[r] \ar[d] & (\Dperf(R_2^\cwedge),\QF^{\geq m}_\sig) \ar[d] \\
	(\Dperf(R[\tfrac{1}{2}],\QF^{\geq m}_\sig) \ar[r] & (\Dperf((R_2^\cwedge)[\tfrac{1}{2}],\QF^{\geq m}_{\sig})
\end{tikzcd}
\]
established in Proposition~\refthree{proposition:LC} and the previous claim. Note that the apparent assumption about the boundedness of the 2-torsion in $R$ in loc.\ cit.\ was only made to ensure that the derived completion $R^\cwedge_2$, for which the above Poincar\'e-Karoubi square always holds, is a discrete ring.
\end{example}

\subsection{The automorphism rings of Grothendieck-Witt- and L-theory}
\label{subsection:automorphism-GW-L}%

Finally, we use Proposition~\reftwo{proposition:pi0gwuniv} to determine the automorphisms of the Grothendieck-Witt and $\L$-theory functors. Yoneda's lemma, the universal property of the Grothendieck-Witt spectrum, Corollary~\reffour{corollary:GW-universal-multiplicative} and Proposition~\refone{proposition:corepresentability-of-poinc} provide equivalences
\[
\Nat(\GW,\GW) \simeq \Nat(\Poinc,\Omega^\infty\GW) \simeq \Nat(\Hom_{\Catp}((\Spa^\perf,\QF^\uni), -), \GWspace) \simeq \GWspace^\uni(\SS),
\]
with the $\Eone$-structure induced by composition on the left hand side which corresponds to the $\Eone$-structure on $\GWspace^\uni(\SS)$ induced by the underlying monoidal structure of $\GWspace$ and the symmetric monoidal structure of $(\Spa^\perf,\QF^\uni)$ as the unit of $\Catp$. Similarly,
$
\Nat(\L,\L) \simeq \Lspace^\uni(\SS),
$
as $\Eone$-spaces and bordification induces the map
\[
\Nat(\GW,\GW) \longrightarrow \Nat(\L,\L),
\]
which identifies with
\[
\GWspace^\uni(\SS) \xrightarrow{\bord} \Lspace^\uni(\SS).
\]
Applying Proposition~\reftwo{proposition:pi0gwuniv} which identifies this map on $\pi_0$ we find:
\begin{corollary}
\label{corollary:automorphisms-GW-L}%
The above identifications provide isomorphisms
\[
\pi_0\Nat(\GW,\GW) \cong \ZZ[e,h]/I \quad \text{and} \quad \pi_0\Nat(\L,\L) \cong \ZZ[e]/(e^2-8e)
\]
with $I = (e^2-8e,he-8h,h^2-2h)$ as in Proposition~\reftwo{proposition:pi0gwuniv}. In particular, we have 
\[
\pi_0\Aut(\GW) = \{\pm \id, \pm (\id-\hyp \circ \fgt)\} \cong (\Ct)^2 \quad \text{and} \quad \pi_0\Aut(\L) = \{\pm \id\} \cong \Ct.
\]
\end{corollary}

Finally, we note that this computation passes to the Karoubi versions of $\GW$ and $\L$:
\begin{corollary}
The canonical maps $\GW \to \KGW$ and $\L \rightarrow \KL$ induce equivalences
\[
\Nat(\GW,\GW) \simeq \Nat(\KGW,\KGW) \quad \text{and} \quad \Nat(\L,\L) \simeq \Nat(\KL,\KL)
\]
\end{corollary}
\begin{proof}
By the universal properties of the functors involved and Proposition~\refone{proposition:corepresentability-of-poinc}, the statement is equivalent to the maps 
\[
\Omega^\infty\GW^\uni(\SS) \longrightarrow \Om^\infty\KGW^\uni(\SS) \quad \text{and} \quad \Omega^\infty\L^\uni(\SS) \longrightarrow \Omega^\infty \KL^\uni(\SS)
\]
being equivalences, which is a consequence of Corollary~\reffour{corollary:regular-noetherian}.
\end{proof}

\begin{appendices}
\renewcommand{\appendixtocname}{{\bf Appendices}}
\addappheadtotoc

\section{Multiplicative aspects of localisation theory}
\label{section:appendix-localisations}%

Let $\C$ be an $\infty$-category. Recall that a full subcategory $\iota\colon \D \subseteq \C$ is called reflective if the inclusion functor $\iota$ admits a left adjoint $L\colon \C \to \D$. In this case we say that $L$ exhibits $\D$ as a left Bousfield localisation of $\C$. In particular, $L$ is a localisation functor: it exhibits $\D$ as the localisation of $\C$ by the collection of $L$-equivalences, that is, the maps $f\colon x \to y$ in $\C$ such that $L(f)$ is an equivalence in $\D$.
In this situation, there is well-known machinery to transport a symmetric monoidal structure on $\C$, or more generally an $\cO$-monoidal structure on $\C$ over some $\infty$-operad $\cO$, to $\D$, 
see e.g.\ \cite{HA}*{\S 2.2.1}. The goal of the present appendix is to treat the more general case where $\C$ is only required to be an $\infty$-operad (over $\cO$), and discuss left Bousfield localisations in this setting, as we introduce in \S\reffour{subsection:multiplicative-localisation} under the name of $\cO$-multiplicative localisations. A special type of multiplicative localisations is given by the smashing localisations associated to idempotent objects, which can be defined when $\C$ and $\cO$ admit compatible unit objects. We develop the relevant notions and study smashing localisations in this generality in \S\reffour{subsection:units-sym-mon-objects} and 
\S\reffour{subsection:smashing-localisations}.

The application we have in mind is the following. Suppose given a small symmetric monoidal $\infty$-category and a presentably symmetric monoidal $\infty$-category $\E$. Then the functor $\infty$-category $\Fun(\C,\E)$ inherits a symmetric monoidal structure, known as Day convolution. Given a collection $\cR$ of cone diagrams in $\C$, one can consider the full subcategory $\Fun^\cR(\C,\E) \subseteq \Fun(\C,\E)$ spanned by those functors taking all diagrams in $\cR$ to limit diagrams in $\E$. The adjoint functor theorem then implies that this inclusion has a left adjoint, and there are well-known criteria for this adjoint to be compatible with Day convolution (see, e.g., Lemma~\refone{lemma:day-localisation}). %
If we now relax the condition that $\C$ is small and only assume that $\C$ is, say, locally small, then one can still define the Day convolution $\infty$-operad $\Fun(\C,\E)^\Day$, which is generally not a symmetric monoidal $\infty$-category. In addition, given a collection of cone diagrams $\cR$, a left adjoint to the inclusion $\Fun^\cR(\C,\E) \subseteq \Fun(\C,\E)$ might only be defined on some intermediate subcategory. %
As a test case of interest, we have dedicated \S\reffour{subsection:day-convolution} to a discussion of multiplicative and smashing localisations in the setting of Day convolution $\infty$-operads. The results we obtain here play a crucial role in the body of the paper for establishing universal lax symmetric monoidal structures on Grothendieck-Witt and L-theory; we refer the reader to the beginning of \S\reffour{subsection:multiplicative-GW} for a discussion of exactly which results are used and how.

\subsection{Multiplicative Bousfield localisations}
\label{subsection:multiplicative-localisation}%

We begin with some basic definitions and notations. Following Lurie's conventions, for an $\infty$-operad $\cO^\otimes$, we write $\cO := \cO^\otimes_{\langle 1 \rangle}$ for its $\infty$-category of colours and refer to maps in $\cO^{\otimes}$ with target in $\cO$ as \defi{operations}. Given $p\colon\C^\otimes \rightarrow \cO^\otimes$ a map of $\infty$-operads, we say that $p$ exhibits $\C^\otimes$ as \defi{$\cO$-monoidal} if $p$ is a cocartesian fibration of $\infty$-operads. A functor $f^\otimes \colon \C^\otimes \to \D^\otimes$ over $\cO^\otimes$ between $\cO$-monoidal $\infty$-categories is called $\cO$-monoidal if it is a map of $\infty$-operads which preserves cocartesian arrows over $\cO^{\otimes}$. We also recall that the maps in any $\infty$-operad support a factorisation system into \defi{active} and \defi{inert} maps, see \cite{HA}*{\S 2.1} for the definition and a thorough discussion.

\begin{definition}
\label{definition:left-bousfield}%
We say that a map $L^{\otimes}$ as in the following diagram
\[
\begin{tikzcd}
\C^{\otimes}\ar[rr, "L^{\otimes}"']\ar[dr,"p"'] && \D^{\otimes}\ar[dl,"q"]\ar[ll,dotted, bend right = 20, "R^{\otimes}"',"\rotatebox{180}{$\perp$}"] \\
& \cO^{\otimes} &
\end{tikzcd}
\]
is 
an \defi{$\cO$-multiplicative} left Bousfield localisation if it admits a
right adjoint $R^{\otimes}$ relative to $\cO^{\otimes}$ (see~\cite{HA}*{Definition 7.3.2.2}) which is furthermore %
fully faithful and preserves inert maps. %

If in addition $L^{\otimes}$ is an $\cO$-monoidal functor among $\cO$-monoidal $\infty$-categories then we say that $L^{\otimes}$ is an \defi{$\cO$-monoidal} left Bousfield localisation.
\end{definition}

\begin{remark}
\label{remark:automatically-monoidal}%
If $L^{\otimes}\colon \C^{\otimes} \to \D^{\otimes}$ is an $\cO$-multiplicative left Bousfield localisation and $\C^{\otimes}$ is $\cO$-monoidal then $\D^{\otimes}$ is also $\cO$-monoidal by \cite{NS}*{Proposition A.14} and $L^{\otimes}$ preserves cocartesian arrows over $\cO^{\otimes}$ since it is a relative left adjoint and hence preserves relative colimits. In particular, an $\cO$-multiplicative left Bousfield localisation is $\cO$-monoidal as soon as its domain is $\cO$-monoidal.
\end{remark}

\begin{remark}
\label{remark:fibrewise-L}%
In the situation of Definition~\reffour{definition:left-bousfield},
the left adjoint $L^{\otimes}$ is also a fibrewise left adjoint to $R^{\otimes}\colon \D^{\otimes} \hrar \C^{\otimes}$ by~\cite{HA}*{Proposition 7.3.2.5}, and hence determines an adjunction $L_a\colon \C_a \adj \D_a \cocolon R_a$ for every colour $a \in \cO$. In addition, for every $A=(a_1,\dots,a_n) \in \cO^{\otimes}$ 
the induced adjunction on fibres over $A$ breaks as a pointwise adjunction of cartesian products
\[
\prod L_{a_i} \colon \prod \C_{a_i} \adj \prod \D_{a_i} \cocolon \prod R_{a_i},
\]
since $L^{\otimes}$ and $R^{\otimes}$ are maps of $\infty$-operads. 
\end{remark}

In the situation of Definition~\reffour{definition:left-bousfield},
the counit of the adjunction $\C^{\otimes} \adj \D^{\otimes}$ is an equivalence. It hence induces, for every $\infty$-operad $\E^{\otimes} \to \cO^{\otimes}$ over $\cO^{\otimes}$, an adjunction
\begin{equation}
\label{eq:adj-algberas}%
\Alg_{\C/\cO}(\E) \locadj \Alg_{\D/\cO}(\E) 
\end{equation}
on $\infty$-categories of algebras, whose left and right adjoints are given by restriction along $R^{\otimes}\colon \D^{\otimes} \hookrightarrow \C^{\otimes}$ and $L^{\otimes}\colon \C^{\otimes} \to \D^{\otimes}$, respectively, and whose counit is again an equivalence. The essential image of the resulting fully faithful functor 
\[
(L^{\otimes})^*\colon \Alg_{\D/\cO}(\E) \to \Alg_{\C/\cO}(\E)
\]
then consists of those $\infty$-operad maps $\Psi^{\otimes}\colon \C^{\otimes} \to \E^{\otimes}$ over $\cO$ which send the components of the unit of $\C^{\otimes} \locadj \D^{\otimes}$ to equivalences. By Remark~\reffour{remark:fibrewise-L}, this is the same as saying that $\Psi_a\colon \C_a \to \E_a$ sends the components of the unit of $\C_a \locadj \D_a$ to equivalences for every $a \in \cO$, or equivalently, that $\Psi_a$ inverts $L_a$-equivalences for every $a \in \cO$. We conclude:

\begin{corollary}
Let $L^{\otimes} \colon \C^{\otimes} \to \D^{\otimes}$ be an $\cO$-multiplicative left Bousfield localisation.
Then $L^{\otimes}$ exhibits $\D^{\otimes}$ as the initial recipient of an $\infty$-operad map over $\cO^{\otimes}$ from $\C^{\otimes}$ which inverts $L_a$-equivalences for every $a \in \cO$. %
\end{corollary}

\begin{lemma}
\label{lem:monoidal-bousfield}%
Let $L^{\otimes}\colon \C^{\otimes} \to \D^{\otimes}$ be an $\cO$-monoidal left Bousfield localisation and $\E^{\otimes} \to \cO^{\otimes}$ an $\cO$-monoidal $\infty$-category.
Then the adjunction~\eqreffour{eq:adj-algberas}
restricts to an equivalence 
\[
\Fun^{\cO-\otimes}_L(\C^{\otimes},\E^{\otimes}) \locadjL{\simeq} \Fun^{\cO-\otimes}(\D^{\otimes},\E^{\otimes})
\]
where on the left hand side we have the full subcategory 
spanned by those $\cO$-monoidal functors $\Psi^{\otimes}\colon \C^{\otimes} \to \E^{\otimes}$ 
which invert $L_a$-equivalences 
for every $a \in \cO$.
\end{lemma}
\begin{proof}
Note first that since $L^{\otimes}$ is an $\cO$-monoidal functor which inverts $L_a$-equivalences for every $a \in \cO$ the same holds for $L^{\otimes} \circ \Phi^{\otimes}$ for every $\cO$-monoidal functor $\Phi^\otimes \colon \D^{\otimes} \to \E^{\otimes}$.
Moreover, if $e\colon (x_1,...,x_n) \to y$ is a cocartesian arrow in $\D^{\otimes}$ over an arrow of the form $(a_1,...,a_n) \to b$ in $\cO^{\otimes}$, then its image in $\C^{\otimes}$ factors as a composite of a cocartesian arrow $(x_1,...,x_n) \to y'$ in $\C^{\otimes}$ over $\cO^{\otimes}$ followed by an $L_b$-equivalence $y'\to y$. It then follows that if $\Psi^{\otimes} \colon  \C^{\otimes} \to \E^{\otimes}$ inverts $L_a$-equivalences for every $a \in \cO$ then the composite $\Psi^{\otimes} \circ \iota^{\otimes}$ is an $\cO$-monoidal functor. This shows that the adjunction~\eqreffour{eq:adj-algberas}
indeed restricts as stated. At the same time, the counit of this adjunction is an equivalence already before restriction, and the unit is an equivalence after restriction since its components are equivalences when evaluated on any $\infty$-operad map that inverts $L_a$-equivalences. 
\end{proof}

We now come to the main result of this subsection, whose goal is to establish a recognition criterion for $\cO$-multiplicative left Bousfield localisations. Our result is a direct generalisation of \cite{HA}*{Proposition 2.2.1.9}, where the analogous criterion is established in the $\cO$-monoidal case.

\begin{definition}
\label{definition:multiplicative-localisation}%
Let $\C^{\otimes} \to \cO^{\otimes}$ be a map of $\infty$-operads. %
We will say that a collection of left Bousfield localisations $\{L_a\colon \C_a \to \D_a \subseteq \C_a\}_{a \in \cO}$ is \defi{$\cO$-multiplicative} if for every operation $\alp\colon (a_1,...,a_n) \to b$ in $\cO$, any choice of an $L_{a_i}$-equivalence $f_i\colon x_i \to y_i$ in $\C_{a_i}$ for each $i=1,...,n$, and any $d \in \D_b$, restriction along the $f_i$ induces an equivalence
\[
\Mul^{\alp}_{\C}(y_1,\dots,y_n;d) \xrightarrow{\simeq} \Mul^{\alp}_{\C}(x_1,\dots,x_n; d)
\]
on spaces of operations in $\C$ covering $\alp$. If $\cO = \Com$ then we will drop $\cO$ from the notation and just write \defi{multiplicative}.
\end{definition}

\begin{remarks}\
\label{remark:full-suboperad}%
\begin{enumerate}
\item If $\C^{\otimes} \to \cO^{\otimes}$ is an $\cO$-monoidal $\infty$-category, 
then the condition that a given collection of localisations $\{L_a\}_{a \in \cO}$ is $\cO$-multiplicative can be rephrased by saying that for every operation $\alp\colon (a_1,...,a_n) \to b$ and every $i = 1,...,n$, the induced map $\C_{a_1} \times \dots \times \C_{a_n} \to \C_{b}$ sends $L_{a_i}$-equivalences in the $i$th coordinate to $L_{b}$-equivalences. In particular, we recover in this case \cite{HA}*{Definition 2.2.1.6}.
\item
\label{item:full-suboperad}%
It follows from the definition that if $\C^{\otimes}$ is an $\infty$-operad over $\cO^{\otimes}$ equipped with an $\cO$-multiplicative collection of left Bousfield localisations $\{L_a\colon \C_a \to \D_a \subseteq \C_a\}_{a \in \cO}$ and $\E^{\otimes} \subseteq \C^{\otimes}$ is a full suboperad such that $\E_a$ contains $\D_a$ for every $a \in \cO$ then the collection of restricted localisations $\{L_a|_{\E_a}\colon \E_a \to \D_a \subseteq \E_a\}_{a \in \cO}$ is again $\cO$-multiplicative.
\item
\label{item:one-at-a-time}%
In the situation of Definition~\reffour{definition:multiplicative-localisation}, one can break the induced map on spaces of operations as a composite where each step features just one of the $f_i$, with the rest of the maps being identities. In establishing multiplicativity one may therefore assume that all the $f_i$'s are identities except for a single $f_{i_0}$. In addition, up to replacing $\alp$ with one of its permutations we may as assume that $i_0=1$.
\end{enumerate}
\end{remarks}

\begin{proposition}
\label{proposition:multiplicative-localisation}%
Let $p\colon\C^{\otimes} \to \cO^{\otimes}$ be map of $\infty$-operads and $\{L_a\colon \C_a \to \D_a\}_{a \in \cO}$ an $\cO$-multiplicative collection of left Bousfield localisations as in Definition~\reffour{definition:multiplicative-localisation}. 
Then the inclusion $\D^{\otimes} \subseteq \C^{\otimes}$ of the full suboperad spanned by the $\D_a$ admits a left adjoint
\[
\begin{tikzcd}
\C^{\otimes} \ar[rr,"L^{\otimes}"]\ar[dr,"p"'] && \D^{\otimes} \ar[dl,"q"] \\
& \cO^{\otimes} & 
\end{tikzcd}
\]
relative to $\cO^{\otimes}$, and such that $L^{\otimes}$ preserves inert arrows. %
In particular, $L^{\otimes}$ is an $\cO$-multiplicative left Bousfield localisation.
\end{proposition}
\begin{proof}
We show that the inclusion $\D^{\otimes} \subseteq \C^{\otimes}$ admits a left adjoint $L^\otimes$ relative to $\cO^{\otimes}$. Then $L^\otimes$ preserves colimits relative to $\cO^\otimes$, and hence in particular $\cO^\otimes$-cocartesian arrows. In particular, $L^\otimes$ is a map of $\infty$-operads.

Applying~\cite{HA}*{Proposition 7.3.2.1} it will suffice to show that for every $\uline{x} := (x_1,\dots,x_n) \in \C^{\otimes}$, the slice category $\D^{\otimes} \times_{\C^{\otimes}} \C^{\otimes}_{\uline{x}/}$ admits an initial object $f\colon \uline{x}\to \uline{d}$ satisfying the additional property that the image of $f$ in $\cO^{\otimes}$ is an equivalence. Let $p(\uline{x}) := \uline{a} = (a_1,\dots,a_n)$. Since each inclusion $\D_{a_i} \subseteq \C_{a_i}$ admits a left adjoint $L_{a_i}\colon \C_{a_i} \to \D_{a_i}$, each of the slice categories $\D_{a_i} \times_{\C_{a_i}} (\C_{a_i})_{x_i/}$ admits an initial object $f_i\colon x_i \to d_i$. Let $f\colon \uline{x}=(x_i,\dots,x_n) \to \uline{d}=(d_1,\dots,d_n)$ be the induced map, so that $\uline{d} \in \D^{\otimes}$ and $f$ maps to the identity on $\uline{a}$ by construction. We now claim that $f$ encodes an initial object of $\D^{\otimes} \times_{\C^{\otimes}} \C^{\otimes}_{\uline{x}/}$. 
For this, we need to show that for every $\uline{e}=(e_1,\dots,e_m) \in \D^{\otimes}$ the pre-composition map
\[
(-) \circ f\Map_{\C^{\otimes}}(\uline{d},\uline{e}) \to \Map_{\C^{\otimes}}(\uline{x},\uline{e})
\]
is an equivalence. Writing $b_i=f(e_i)$ and $\uline{b}=(b_1,\dots,b_m)$, the above map fits in a commutative diagram
\[
\begin{tikzcd}
\Map_{\C^{\otimes}}(\uline{d},\uline{e}) \ar[rr]\ar[dr] && \Map_{\C^{\otimes}}(\uline{x},\uline{e})\ar[dl] \\
&\Map_{\cO^{\otimes}}(\uline{a},\uline{b}) &
\end{tikzcd}
\]
It will then suffice to show that the top horizontal map induces an equivalence on vertical fibres for every map $\alp\colon \uline{a} \to \uline{b}$ in $\cO^{\otimes}$. Let $\rho\colon \langle n\rangle \to \langle m\rangle$ be the image of $\alp$ in $\Com$. Then $\alp$ determines a map $\alp_j\colon (a_i \;|\; \rho(i)=j) \to b_j$ for each $j=1,...,m$ and the mapping spaces in the top row decompose compatibly as products of mapping spaces covering the various $\alp_i$. We may hence assume without loss of generality that $m=1$. Then, we need to show that for any operation 
$\alp\colon \uline{a} \to b$ and any $e \in \D_b$, the induced map on spaces of operations in $\C^\otimes$ covering $\alp$
\[
\Mul^{\alp}_{\C}(d_1,\dots,d_n;e) \to \Mul^{\alp}_{\C}(x_1,\dots,x_n;e)
\]
is an equivalence, which is the case by the assumption that the collection of localisations $\{L_a\}_{a \in \cO}$ is $\cO$-multiplicative.
\end{proof}

When $\C^{\otimes} \to \cO^{\otimes}$ is an $\cO$-monoidal $\infty$-category Proposition~\reffour{proposition:multiplicative-localisation} recovers \cite{HA}*{Proposition 2.2.1.9}, see Remark \reffour{remark:automatically-monoidal}, as we record now (we also record the corresponding universal property of $\D^{\otimes}$, which is not mentioned in loc.\ cit.)

\begin{proposition}
\label{proposition:multiplicative-localisation-2}%
In the situation of Proposition~\reffour{proposition:multiplicative-localisation}, suppose in addition that $\C^{\otimes} \to \cO^{\otimes}$ is an $\cO$-monoidal $\infty$-category. Then $\D^{\otimes} \to \cO^{\otimes}$ is also $\cO$-monoidal as is the functor $L^{\otimes}$ provided by Proposition~\reffour{proposition:multiplicative-localisation}. In addition, this functor
enjoys the following universal property: for every $\cO$-monoidal $\infty$-category $\E^{\otimes} \to \cO^{\otimes}$, the functor
\[
(L^{\otimes})^*\colon \Fun^{\cO-\otimes}(\D^{\otimes},\E^{\otimes}) \to \Fun^{\cO-\otimes}(\C^{\otimes},\E^{\otimes})
\]
is fully faithful, and its essential image is spanned by those $\cO$-monoidal functors $\Psi^{\otimes}\colon \C^{\otimes} \to \E^{\otimes}$ such that $\Psi_a\colon \C_a \to \D_a$ sends $L_a$-equivalences to equivalences for every $a \in \cO$.
\end{proposition}
\begin{proof}
Combine Remark \reffour{remark:automatically-monoidal} and Lemma~\reffour{lem:monoidal-bousfield}.
\end{proof}

\begin{example}
\label{example:sm-localisation}%
Suppose $\cO = \Com$. Then an $\cO$-monoidal $\infty$-category amounts to a symmetric monoidal $\infty$-category $\C^{\otimes}$, and a left localisation $L\colon \C \to \D \subseteq \C$ is multiplicative if and only if for every pair of $L$-equivalences $f\colon x \to y$ and $g\colon z \to w$, the map $f \otimes g \colon x \otimes z \to y \otimes w$ is again an $L$-equivalence, equivalently, the collection of $L$-equivalences is closed under the operation $x \otimes (-)\colon \C \to \C$ for every $x \in \C$. 
When this holds, the symmetric monoidal structure of $\C$ descends to one on $\D$, and the resulting symmetric monoidal functor
\[
L^{\otimes} \colon \C^{\otimes} \to \D^{\otimes}
\]
is initial in symmetric monoidal functors out of $\C^{\otimes}$  which send $L$-equivalences to equivalence, as an instance of Proposition~\reffour{proposition:multiplicative-localisation}.
The induced symmetric monoidal structure on $\D$ is then given explicitly by the formula $x \otimes y = L(\iota(x) \otimes \iota(y))$.
\end{example}

\begin{example}
\label{example:mcom-loc}%
Suppose $\cO=\MCom$, the $\infty$-operad which classifies pairs $(R,M)$ of a commutative algebra object $R$ equipped with an $R$-module $M$, see~\cite{hinich-modules}*{\S 5.2}. Then $\MCom$-monoidal $\infty$-categories correspond to pairs $(\C^{\otimes},\M)$ where $\C$ is a symmetric monoidal $\infty$-category and $\M$ is an $\infty$-category equipped with an action of $\C$. A pair of left localisations $L_{\alg}\colon \C \to \D \subseteq \C$ and $L_{\md}\colon \M \to \cN \subseteq \M$ is then $\MCom$-multiplicative if $L_{\alg}$ is a multiplicative localisation of $\C$ as in Example~\reffour{example:sm-localisation} and the action of $\C$ on $\M$ satisfies that for $f\colon x \to y$ in $\C$ and $g\colon z \to w$ in $\M$, if $f$ is an $L_{\alg}$-equivalence and $g$ is an $L_{\md}$-equivalence then $f \otimes g \colon x \otimes z \to y \otimes w$ is an $L_{\md}$-equivalences. In that case we will also say that the localisation $L_{\md}$ is \emph{linear over $L_{\alg}$}. By Proposition~\reffour{proposition:multiplicative-localisation}, when this holds the symmetric monoidal structure of $\C$ descends to one on $\D$ and the $\C$-action on $\M$ descends to a compatible $\D$-action on $\cN$, such that the resulting map
\[
L^{\otimes}\colon (\C^{\otimes},\M) \to (\D^{\otimes},\cN)
\]
of $\MCom$-algebras is initial in $\MCom$-algebras maps out of $(\C,\M)$ sending $L_{\alg}$-equivalences and $L_{\md}$-equivalences to equivalences. The induced $\D$ action on $\cN$ is given explicitly by the formula $x \otimes y = L_{\md}(\iota_{\alg}(x) \otimes \iota_{\md}(y))$, where $\iota_{\alg}\colon \D \subseteq \C$ and $\iota_{\md}\colon \cN \subseteq \M$ are the corresponding full inclusions.

An interesting special case of the above is when $\D=\C$ and $L_{\alg} = \id$.
In this case we can ignore the algebra component and just consider the left localisation $L\colon \M \to \cN \subseteq \M$ on $\M$. To check multiplicativity it is then enough to check that the collection of $L$-equivalences in $\M$ is preserved by the operation $x \otimes (-)\colon \M \to \M$ for every $x \in \C$. In this case we will also say that $L$ is a \defi{$\C$-linear localisation}. Then $\cN$ inherits from $\M$ its $\C$-module structure and the resulting $\C$-module refinement $\M \to \cN$ of $L$ is initial in $\C$-module functors out of $\M$ sending $L$-equivalences to equivalences.
\end{example}

\subsection{Operadic units and idempotents}
\label{subsection:units-sym-mon-objects}%

Our goal in this subsection is to recall and explore the notions of units and idempotent objects in $\infty$-operads which will be used in \S\reffour{subsection:smashing-localisations} for the study of smashing localisations of $\infty$-operads (with units). The notion of units we use here is closely related to that introduced by Lurie in~\cite{HA}*{\S 3.2.1}, see Remark~\reffour{remarks:symmetric-monoidal-objects-and-units}\;\reffouritem{item:HA}
for the precise relation between the two.

We begin with the auxiliary notion of an operadically cocartesian arrow. To that end, recall that for an $\infty$-operad $\cO^\otimes$ there is a symmetric monoidal structure on $\cO^\otimes_\act$ which we refer to as concatenation, see \cite{HA}*{Corollary 2.2.4.5 \& Remark 2.2.4.6}.  
Recall from~\cite{HA}*{Definition 3.1.1.2} that for a map of $\infty$-operads $p\colon \C^{\otimes} \to \cO^{\otimes}$, a cone diagram $\ovl{\alp} \colon K^{\triangleright} \to \C^{\otimes}_{\act}$ is said to be a weak operadic $p$-colimit cone if the map 
\[
\C^{\act}_{\ovl{\alp}/} \to \C^{\act}_{\alp/} \times_{\cO^{\act}_{p\alp/}} \cO^{\act}_{p\ovl{\alp}/}
\]
is a trivial Kan fibration, where $\alp = \ovl{\alp}|_{K}$ and the superscript $\act$ means that we restrict attention to cones involving only active maps and having tips lying over $\<1\> \in \Com^{\otimes}$. We then say that $\ovl{\alp}$ is an operadic $p$-colimit cone if it is a weak operadic $p$-colimit cone and remains so after concatenation with arbitrary objects in $\C^{\otimes}$.
In what follows we will implicitly identify $\Delta^1$ with $(\Delta^0)^\triangleright$, so that morphisms can be viewed as cone diagrams.

\begin{lemma}
\label{lemma:operadically}%
Let $p\colon \C^{\otimes} \to \cO^{\otimes}$ be a map of $\infty$-operads and $\alp\colon x \to y$ a map in $\C^{\otimes}_\act$ with $y \in \C := \C^{\otimes}_{\<1\>}$. Then the following conditions are equivalent:
\begin{enumerate}
\item
\label{item:def1}%
$\alp$ encodes an operadic $p$-colimit cone in $\C^{\otimes}$.
\item
\label{item:def2}%
$\alp$ is $p$-cocartesian and remains $p$-cocartesian after concatenation with arbitrary $z \in \C^{\otimes}$.
\end{enumerate}
\end{lemma}

\begin{definition}
\label{definition:operadically-cocartesian}%
Let $p\colon \C^{\otimes} \to \cO^{\otimes}$ be a map of $\infty$-operads and $\alp$ a map in $\C^{\otimes}_\act$ with target in $\C$. We will say that $\alp$ is \defi{operadically $p$-cocartesian} if it satisfies the equivalent conditions of Lemma~\reffour{lemma:operadically}.
\end{definition}

\begin{proof}[Proof of Lemma~\reffour{lemma:operadically}]
To see $\reffouritem{item:def2} \Rightarrow \reffouritem{item:def1}$, it suffices by definition to show that every $p$-cocartesian active arrow in $\C^{\otimes}$
is a weak operadic $p$-colimit diagram.
But such an arrow is also $p_{\act}$-cocartesian, where $p_{\act}\colon \C_{\act}^\otimes \to \cO_{\act}^\otimes$, and so this follows directly from the definition of weak operadic $p$-colimits.
To get $\reffouritem{item:def1} \Rightarrow \reffouritem{item:def2}$ we first note that what we have just proven implies that active equivalences, in particular identities, are operadic $p$-colimit cones. Then the claim follows rom \cite{HA}*{Prop.\ 3.1.1.10} using that an arrow is $p$-cocartesian if and only if it encodes a $p$-colimit diagram.
\end{proof}

\begin{definition}
\label{definition:monoidal-and-unit}%
Let 
\[
\begin{tikzcd}
\C^{\otimes} \ar[dr,"p"']\ar[rr,"q"] && \cO^{\otimes} \ar[dl,"r"] \\
& \Com^{\otimes} &
\end{tikzcd}
\]
be a diagram of $\infty$-operads.
\begin{enumerate}
\item 
A map $u\colon \emptyset \to \one_\cO$ in $\cO^\otimes$ lying over the unique map $\<0\> \to \<1\>$ in $\Com^{\otimes}$ is said to exhibit $\one_\cO$ as a \defi{unit object} of $\cO$ if it is operadically $r$-cocartesian. We say that a map $\tilde{u}\colon \emptyset \to \one_\C$ in $\C^{\otimes}$ exhibits $\one_\C$ as a \defi{unit over $\one_\cO$} if it exhibits $\one_\C$ as a unit in $\C$ and $u=q(\tilde{u})$.
\item
An object $x \in \cO$ is said to be \defi{symmetric monoidal} if for every $y \in \cO$ there exists an operadically $r$-cocartesian map $(x,y) \to x\otimes y$ (using suggestive notation) lifting the unique active map $\mu_2\colon \<2\> \to \<1\>$.
Given a symmetric monoidal object $x \in \cO$, we will say that an object $\tilde{x}$ in $\C$ lifting $x$ is \defi{symmetric monoidal over $x$} if $\tilde{x}$ is symmetric monoidal %
and for every operadically $p$-cocartesian map %
$\alp\colon (\tilde{x},\tilde{y}) \to \tilde{x}\otimes\tilde{y}$ %
covering $\mu_2$, the map $q(\alp)$ is operadically $r$-cocartesian. 
\end{enumerate}
\end{definition}

\begin{examples}\
\begin{enumerate}
\item
In the $\infty$-operad $\Com$, the unique map $\<0\> \to \<1\>$ exhibits the $\<1\>$ as a unit object. 
\item
If $\D^{\otimes}$ is a symmetric monoidal $\infty$-category then the unit $\one_{\D}$ of $\D^{\otimes}$ is also canonically a unit object in the sense of Definition~\reffour{definition:monoidal-and-unit}.
\item
In the $\infty$-operad $\MCom$ the unique null operation $\emptyset \to a$ exhibits the algebra colour $a$ as a unit object. 
\end{enumerate}
\end{examples}

\begin{lemma}
\label{lemma:absolute-is-relative}%
Let $q\colon \C^{\otimes} \to \cO^{\otimes}$ be a map of $\infty$-operads such that $\cO^{\otimes}$ admits a unit object $u\colon \emptyset \to \one_{\cO}$ and $\C^{\otimes}$ admits a unit $\wtl{u}\colon \emptyset \to \one_{\cO}$
over $\one_{\cO}$. Then the following holds:
\begin{enumerate}
\item
The map $\wtl{u}$
is an operadically $q$-cocartesian lift of $u$.
\item
Let $a \in \cO^{\otimes}$ be a symmetric monoidal object and $\alp\colon (a,b) \to a\otimes b$ an operadically cocartesian lift of the map $\<2\> \to \<1\>$ in $\Com^{\otimes}$. Then for any symmetric monoidal object in $x \in \C$ over $a$ and any lift $y \in \C$ of $b$ there exists an operadically $q$-cocartesian arrow $(x,y) \to x\otimes y$ lifting $\alp$.
\end{enumerate}
\end{lemma}
\begin{proof}
To prove the first claim, note that by definition both $\wtl{u}$ and $u=q(\wtl{u})$ are operadically cocartesian over $\Com^{\otimes}$, that is, they are cocartesian over $\Com^{\otimes}$ and remain so after concatenation with arbitrary objects. 
We conclude from~\cite{HTT}*{Proposition 2.4.1.3} that $\wtl{u}$ is $q$-cocartesian and remains $q$-cocartesian after concatenation with arbitrary objects. In other words, $\wtl{u}$ is operadically $q$-cocartesian.

For the second claim, let us fix $\alp\colon (a,b) \to a\otimes b$, an operadically cocartesian lift of the map $\<2\> \to \<1\>$ such that $a$ is symmetric monoidal and $x,y \in \C$ lifts of $a,b$ such that $x$ is symmetric monoidal over $a$. Then by definition there exists an operadically cocartesian arrow $\beta\colon (x,y) \to x\otimes y$ lift of $\<2\> \to \<1\>$ in $\C$ such that $q(\beta)$ is operadically cocatesian in $\cO$ over $\Com^{\otimes}$. We then have that $\alp$ and $q(\beta)$ are two cocartesian lifts of $\<2\> \to \<1\>$ with the same domain, and so there exists an equivalence $\tau\colon \alp\simeq q(\beta)$ as objects of $\cO^{\otimes}_{(a,b)/} \times_{\Com^{\otimes}_{\<2\>/}} \{\<2\>\to \<1\>\}$. Using $\tau$ we may hence view $\beta$ as a lift of $\alp$. Arguing as above we again conclude from~\cite{HTT}*{Proposition 2.4.1.3} that $\beta$ is operadically $q$-cocartesian, as claimed.
\end{proof}

\begin{remarks}\
\label{remarks:symmetric-monoidal-objects-and-units}%
\begin{enumerate}
\item
If $\cO^{\otimes}$ is any $\infty$-operad and $u\colon \emptyset \to \one_{\cO}$  is a map in $\cO^{\otimes}$ lying above $\<0\> \to \<1\>$, then $u$ exhibits $\one_{\cO}$ as a unit object if and only if it exhibits it as a unit object over the unit of $\Com$. %
\item
\label{item:HA}%
If follows from Lemma~\reffour{lemma:absolute-is-relative} that a map $\tilde{u}\colon \emptyset \to \one_\C$ in $\C^{\otimes}$ exhibits $\one_\C$ as a unit over $\one_\cO$ if and only if it exhibits $\one_{\C}$ as an $\one_{\cO}$-unit in the sense of~\cite{HA}*{Definition 3.2.1.1}. To avoid confusion, let us note that loc.\ cit. also considers the unit property for objects of $\C$ lying over arbitrary objects in $\cO$, a notion which has no counterpart in Definition~\reffour{definition:monoidal-and-unit}.
\item
\label{item:unit-is-monoidal}%
Any unit object is symmetric monoidal, and if $\C^{\otimes} \to \cO^{\otimes}$ is a map of $\infty$-operad such that $\cO^{\otimes}$ has a unit $\one_{\cO}$ then any unit in $\C$ over $\one_{\cO}$ is also symmetric monoidal over $\one_{\cO}$. %
\item
\label{item:simplifies}%
In the $\infty$-operad $\Com^{\otimes}$ the unique active arrow $\mu_n\colon \<n\> \to \<1\>$ with $n \geq 2$ can be written as a composite of arrows obtained from the active arrow $\mu_2\colon \<2\> \to \<1\>$ via concatenation with objects. Given an $\infty$-operad $\C^{\otimes}$ it then follows by induction that if $x \in \C$ is a symmetric monoidal object then for every $y \in \C$ there exists an operadically cocartesian map of the form $\alp\colon (x,\dots,x,y) \to x\otimes ... \otimes x\otimes y$ lying over $\mu_n$. Similarly, if $x \in \C$ is symmetric monoidal over some symmetric monoidal object in $\cO$ then the operadically cocartesian maps $\alp$ just mentioned are sent to operadically cocartesian arrows in $\cO^{\otimes}$ (over $\Com^{\otimes}$). %
\item
\label{item:local-unit}%
Suppose that $\C^{\otimes}$ admits a unit $\emptyset \to \one_{\C}$ over a unit $\one_{\cO}$ in $\cO$ and let $L^{\otimes} \colon \C^{\otimes} \to \D^{\otimes}$ be an $\cO$-multiplicative left Bousfield localisation. %
Then the composite 
\[
\emptyset \xrightarrow{u} \one_{\C} \xrightarrow{\eps} L(\one_\C)
\]
exhibits its target as a unit object in $\D$ over $\one_{\cO}$. Indeed, since $L^{\otimes}$ is surjective on objects and preserves cocartesian maps over $\cO^{\otimes}$ that exist in $\C^{\otimes}$, it also preserves all operadically cocartesian maps that exist in $\C^{\otimes}$.
\end{enumerate}

\end{remarks}

The following proposition summarises the main facts we need about unit objects. It is obtained by combining~\cite{HA}*{Proposition 3.2.1.8, Lemma 3.2.1.10 \& Proposition 3.4.2.1}: 

\begin{proposition}
\label{proposition:units}%
Let $p\colon \D^\otimes \to \Com^{\otimes}$ be an $\infty$-operad.
Suppose that $u\colon \emptyset \to \one_{\D}$ exhibits $\one_{\D}$ as a unit object. Then the map 
$\mathbb{E}_0^{\otimes} \to \D^{\otimes}$ encoded by $u$ extends in an essentially unique manner to a commutative algebra object 
\[
\wtl{\one}_\D\colon \Com^{\otimes} \to \D^{\otimes}, \quad \<1\>  \mapsto \one_{\D},
\]
often referred to as trivial algebra object. In addition, this commutative algebra object satisfies the following properties:
\begin{enumerate}
\item
$\wtl{\one}_\D$ sends every map in $\Com^\otimes$ to a $p$-cocartesian map in $\D^{\otimes}$.
\item
With this algebra structure, $\one_{\D}$ is initial in $\CAlg(\D)$. 
\item
The forgetful functor $\Mod_{\one_{\D}}(\D)^{\otimes} \to \D^{\otimes}$ is an equivalence. In other words, any object in $\D$ is a $\one_\D$-module in any essentially unique manner.
\end{enumerate}
\end{proposition}

\begin{proposition}
\label{proposition:action-of-mon}%
In the situation of Proposition~\reffour{proposition:units}, if
\[
\mu \colon \MCom^\otimes \to \D^\otimes
\]
is an $\infty$-operad map encoding an algebra-module pair $(\wtl{\one}_{\D},d)$ in which the algebra object is trivial then 
$\mu$ sends every arrow in $\MCom^{\otimes}$ to an $p$-cocartesian arrow in $\D^{\otimes}$. %
\end{proposition}
\begin{proof}
Since $\mu|_{\Com^{\otimes}}$ encodes the trivial algebra object $\one_{\D}$, it sends the map from $\emptyset$ to the algebra colour in $\MCom$ to an operadically $p$-cocartesian arrow in $\D^{\otimes}$. As a consequence, $\mu$ sends every map that can be obtained by concatenating this one with objects in $\MCom^{\otimes}$ to a $p$-cocartesian arrow. By direct inspection, this includes all active arrows in $\MCom^{\otimes}$ which lie over (active) injective maps $\<n\> \to \<m\>$ in $\Com^{\otimes}$. It hence suffices to show that $\mu$ sends maps lying over active surjective maps $\<n\> \to \<m\>$ in $\Com^{\otimes}$ to $p$-cocartesian maps. But any such map in $\MCom^{\otimes}$ admits a section which is itself of the injective kind, and so the desired result follows from~\cite{HTT}*{Proposition 2.4.1.7}.
\end{proof}

Given a map of $\infty$-operads $q\colon \C^{\otimes} \to \cO^{\otimes}$ and a unit object $\one_{\cO} \in \cO$ with associated trivial algebra object $\tilde{\one}_{\cO}\colon \Com^{\otimes} \to \cO^{\otimes}$, let us denote by $\C^\otimes_{\mon/\one_\cO}$ the full suboperad of $\C^{\otimes} \times_{\cO^{\otimes}} \Com^{\otimes}$ spanned by the symmetric monoidal objects over $\one_\cO$.

\begin{proposition}
\label{proposition:main-point}%
Let $q\colon \C^{\otimes} \to \cO^{\otimes}$ be a map of $\infty$-operads such that $\cO^{\otimes}$ admits a unit object $\one_{\cO}\in \cO$ and $\C^{\otimes}$ admits a unit $\one_{\C}$ over $\one_{\cO}$. 
Suppose given a map $\mu\colon \MCom^{\otimes} \to \cO^{\otimes}$ encoding an algebra-module pair $(\wtl{\one}_{\cO},b)$ whose algebra is trivial, and form the commutative square
\[
\begin{tikzcd}
\M_b^{\otimes} \ar[r]\ar[d] & \C^{\otimes} \ar[d] \\
\MCom^{\otimes} \ar[r,"\mu"] & \cO^{\otimes}
\end{tikzcd}
\]
where $\M_b^{\otimes}$ is the full suboperad of $\MCom^{\otimes}\times_{\cO^{\otimes}} \C^{\otimes}$ spanned by the objects $(a, c)$ with $a \in \MCom$ the algebra colour and $c \in \C$ symmetric monoidal over $\one_{\cO}$ and all the objects of the form $(m,c)$ where $m \in \MCom$ is the module colour.
Then the left vertical map is a cocartesian fibrations of $\infty$-operads. In particular, $\C^{\otimes}_{\mon/\one_{\cO}}$ is a symmetric monoidal $\infty$-category and $\M_b^{\otimes}$ encodes an action of $\C^{\otimes}_{\mon/\one_{\cO}}$ on $\C_b$.
\end{proposition}

The proof of Proposition~\reffour{proposition:main-point} will make use of the following lemma:

\begin{lemma}
\label{lemma:sm-closed-under-tensor}%
Let $q\colon \C^{\otimes} \to \cO^{\otimes}$ be a map of $\infty$-operads.
Let $\bet\colon (x,y) \to x\otimes y$ be an operadically cocartesian lift of the map $\mu_2\colon \<2\> \to \<1\>$ such that $\alp=q(\bet)\colon (a,b) \to a\otimes b$ is an operadically cocartesian lift of $\mu_2$ in $\cO^{\otimes}$. Then the following holds:
\begin{enumerate}
\item
If $a$ and $b$ are symmetric monoidal then $a\otimes b$ is symmetric monoidal.
\item
If $a,b$ and $a\otimes b$ are symmetric monoidal and $x,y$ are symmetric monoidal over $a,b$, respectively, then $x\otimes y$ is symmetric monoidal over $a\otimes b$.
\end{enumerate}
\end{lemma}
\begin{proof}
The first claim can then be deduced from the second by applying it to the $\infty$-operad map $\cO^{\otimes} \to \Com^{\otimes}$.
We now prove the second claim. We need to check that operadically $p$-cocartesian lifts of the active map $\<2\> \to \<1\>$ exist for every source of the form $(z,w)$ with $w \in \C$, and that these lifts map to operadically $r$-cocartesian maps in $\cO$. For this, consider the commutative square
\[
\begin{tikzcd}
\<3\>\ar[r]\ar[d] & \<2\> \ar[d] \\
\<2\> \ar[r] & \<1\>
\end{tikzcd}
\]
of active arrows in $\Com^{\otimes}$ where the two maps $\<3\> \to \<2\>$ are the two distinct active maps sending $1$ to $1$ and $3$ to $2$. Since both $x$ and $y$ are symmetric monoidal over $a$ and $b$ respectively
we can lift the left vertical and bottom horizontal arrows to operadically $p$-cocartesian arrows 
\[
(x,y,w) \to (x,y \otimes w) \to (x \otimes y \otimes w)
\]
in $\C^{\otimes}$ whose image in $\cO^{\otimes}$ is operadically $r$-cocartesian, where as above we use the tensor notation suggestively. Let $\gam$ denote the composite of these two maps.
Since $\beta\colon (x,y) \to x\otimes y$ is operadically cocartesian its concatenation $(\beta,w)\colon (x,y,w) \to (x\otimes y,w)$ with $w$ is cocartesian and so we may then factor $\gam$ as a composite
\[
(x,y,w) \xrightarrow{(\beta,w)} (x\otimes y,w) \to (x \otimes y \otimes w) .
\]
Since $\gam$ and the first arrow are both operadically cocartesian it follows that the second arrow $(x\otimes y,w) \to (x \otimes y \otimes w)$ is an operadically cocartesian lift of $\<2\> \to \<1\>$ by the cancellation property for cocartesian arrows. Its image in $\cO$ is also an operadically cocartesian lift of $\<2\> \to \<1\>$ by the same cancellation property. %

\end{proof}

\begin{proof}[Proof of Proposition~\reffour{proposition:main-point}]
The morphisms in $\MCom^{\otimes}$ are generated under concatenation and composition by the unit $\emptyset \to a$, the multiplication map $(a,a) \to a$, and the action map $(a,m) \to m$. It will hence suffice to show that any of these three arrows admits operaically cocartesian lifts in $\M_b^{\otimes}$ for any lift of their domain. Since $\mu \colon \MCom^\otimes \to \D^\otimes$ sends every arrow to an operadically cocartesian arrow over $\Com$ by Proposition~\reffour{proposition:action-of-mon} the desired result now follows form Lemmas~\reffour{lemma:absolute-is-relative} and \reffour{lemma:sm-closed-under-tensor}.
\end{proof}

Proposition~\reffour{proposition:main-point} allows to define the notion of idempotent objects in the $\infty$-operadic context, which we will use in the subsequent section to define and study the notion of smashing localisations:
\begin{definition}
Let $\cO^{\otimes}$ be an $\infty$-operad with unit $\one_{\cO}$
and $q\colon \C^{\otimes} \to \cO^{\otimes}$ a map of $\infty$-operads with unit 
$\one_{\C}$
over $\one_\cO$. 
We will say that a map $\eps\colon \one_{\C} \to e$ in $\C$ exhibits $e$ as \defi{idempotent over $\cO$}, or \defi{$\cO$-idempotent}, if the following three conditions hold:
\begin{enumerate}
\item
The map $q(\eps)$ is an equivalence in $\cO$.
\item
The object $e$ is symmetric monoidal over $\one_{\cO}$. %
\item
The map
$\eps \otimes \id\colon e \to e \otimes e$ is an equivalence in the symmetric monoidal $\infty$-category $\C_{\mon/\one_{\C}}$.
\end{enumerate}
\end{definition}

\begin{remark}
\label{remark:action-of-idempotent}%
Any symmetric monoidal object $e$ in $\C$ over $\one_{\cO}$ determines an endo-functor $e \otimes (-)\colon \C_a \to \C_a$ for any $a \in \cO$ via the action of 
$\C_{\mon/\one_{\C}}$
on $\C_a$ of Proposition~\reffour{proposition:main-point}.
If $\eps\colon \one_{\C} \to e$ exhibits $e$ is idempotent over $\cO$ then the resulting endo-functor of $\C_a$ inherits a natural transformation from $\id_{\C_a}$ which satisfying the conditions of~\cite{HTT}*{Proposition 5.2.7.4}, and so $e \otimes (-)$ exhibits its image as a left Bousfield localisation of $\C_a$. The map $\epsilon\colon \one_{\C} \to e$ can then be recovered as the local approximation of $\one_{\C}$ in $\C_{\one_{\C}}$. Collections of left Bousfield localisations which arise in this manner are always $\cO$-multiplicative, and hence refine to $\cO$-multiplicative left Bousfield localisations of $\infty$-operads 
by Proposition~\reffour{proposition:multiplicative-localisation}. These special types of left Bousfield localisations are called \emph{smashing localisations}, which we discuss next.
\end{remark}

\subsection{Smashing localisations of \texorpdfstring{$\infty$-}{}operads}
\label{subsection:smashing-localisations}%

In this subsection we develop the notion of smashing localisations in the context of $\infty$-operads.
Let
\[
\begin{tikzcd}
\C^{\otimes} \ar[dr,"q"]\ar[ddr,"p"']\ar[rr,"L^{\otimes}"] && \ar[dl] \D^{\otimes}\\
& \cO^{\otimes} \ar[d,"r"] & \\
& \Com^{\otimes} &
\end{tikzcd}
\]
be an $\cO$-multplicative left Bousfield localisation of $\infty$-operads and assume that %
$\cO^{\otimes}$ admits a unit $\emptyset \to \one_{\cO}$ and 
$\C^{\otimes}$ admits a unit $\emptyset \to \one_{\C}$ over $\one_{\cO}$.

\begin{construction}
\label{construction:eta}%
Write $e \coloneq L(\one_{\C})$ and fix an object $x \in \C$ lying above a colour $a \in \cO$. 
Then the map $x \to (\one_\C,x)$ obtained by concatenating $\emptyset \to \one_{\C}$ with $x$ is $p$-cocartesian and maps to an $r$-cocartesian arrow in $\cO$. We may then factor the map $x \to L(x)$ as a composite in $\C^\otimes$ of the form appearing on the left
\[
\begin{tikzcd}[column sep = 10pt]
& (\one_\C,x)\ar[dr] &     &&  & (\one_{\cO},a)\ar[dr] &   &&  & \<2\> \ar[dr] & \\
x \ar[ur]\ar[rr] && L(x) &&  a \ar[ur]\ar[rr,equal] && a && \<1\> \ar[ur]\ar[rr,equal] && \<1\>
\end{tikzcd}
\]
lying over a triangle in $\cO^\otimes$ of the form appearing in the middle and over a triangle in $\Com^\otimes$ of the form appearing on the right, and where the arrow $(\one_{\cO},a) \to a$ is operadically $r$-cocartesian (by applying~\cite{HTT}*{Proposition 2.4.1.7} after concatenating with arbitrary objects in $\cO^{\otimes}$).
Since the map $\one_{\C} \to e$ is an $L$-equivalence, %
so is the map $(\one_{\C},x) \to (e,x)$, and so %
the map $(\one_{\C},x) \to L(x)$ factors through a map 
\[
\eta_x\colon (e,x) \to L(x)
\]
lying over the same operadically $r$-cocartesian arrow $(\one_{\cO},a) \to a$ (up to an equivalence). 
\end{construction}

\begin{definition}
\label{definition:smashing}%
Keeping the above notation and assumptions, we say that the $\cO$-multiplicative left Bousfield localisation $L^{\otimes}$ is \emph{$\cO$-smashing} if the map $\eta_x\colon (e,x) \to L(x)$ of Construction~\reffour{construction:eta} is operadically $p$-cocartesian for every $x \in \C$.
\end{definition}

The following is the main result of this subsection, which can be viewed as a recognition criterion as the one of Proposition~\reffour{proposition:multiplicative-localisation}, this time for smashing localisations.

\begin{proposition}
\label{proposition:equivalent-1}%
Let $\cO^{\otimes}$ be an $\infty$-operad with unit $\one_{\cO}$
and $q\colon \C^{\otimes} \to \cO^{\otimes}$ a $\infty$-operad over $\cO^{\otimes}$ which admits a unit $\emptyset \to \one_{\C}$ over $\one_\cO$. Let $\D^{\otimes} \subseteq \C^{\otimes}$ be a full suboperad.
Then the following are equivalent
\begin{enumerate}
\item
\label{item:smashing}%
The inclusion $\D^{\otimes} \subseteq \C^{\otimes}$ admits an $\cO$-smashing left Bousfield localisation $L^{\otimes}\colon \C^{\otimes} \to \D^{\otimes}$. %
\item
\label{item:idem-1}%
There exists an $\cO$-idempotent object $\eps\colon \one_{\C} \to e$ in $\C$ such that for every $a \in \cO$ the full subcategory $\D_a \subseteq \C_a$ is exactly the image of $e \otimes (-)\colon \C_a \to \C_a$. 
\item
\label{item:idem-2}%
Each inclusion $\D_a \subseteq \C_a$ admits a left adjoint $L_a\colon \C_a \to \D_a$
and there exists an $\cO$-idempotent object $\eps\colon\one_{\C} \to e$ in $\C$ such that $e \otimes (-) \simeq L_a(-)$ as endo-functors under $\id_{\C_a}$ for every $a \in \cO$. In particular, %
$e \simeq L_{\one_{\cO}}(\one_{\C})$ under $\one_{\C}$. %
\item
\label{item:closed}%
Each inclusion $\D_a \subseteq \C_a$ admits a left adjoint $L_a\colon \C_a \to \D_a$
such that the collection $\{L_a\}_{a \in \cO}$ is $\cO$-multiplicative, the object $e := L_{\one_{\cO}}(\one_{\C})$ is symmetric monoidal over $\one_{\cO}$,
and for any $a \in \cO$, $x \in \C_{\mon/\one_\cO}$ and $y \in \C_a$, if either $x$ belongs to $\D_{\one_{\cO}}$ or $y$ belongs to $\D_{a}$, then $x \otimes y \in \C_a$ belongs to $\D_a$. 
\end{enumerate}
\end{proposition}
\begin{proof}
We begin with proving \reffouritem{item:smashing} $\Rightarrow$ \reffouritem{item:idem-1}, and so we 
suppose that the inclusion $\D^{\otimes} \subseteq \C^{\otimes}$ admits an $\cO$-smashing left Bousfield localisation $L^{\otimes}\colon \C^{\otimes} \to \D^{\otimes}$.
Then by definition for every $x \in \C$ there exists an operadically cocartesian map $(e,x) \to L_{q(x)}(x)$ lifting $\<2\> \to \<1\>$ and which maps to an operadically cocartesian map in $\cO^{\otimes}$ over $\Com^{\otimes}$, so that $e$ is symmetric monoidal over $\one_{\cO}$
by Remark~\reffour{remarks:symmetric-monoidal-objects-and-units}\;\reffouritem{item:simplifies}. In addition, taking $x = e$ we obtain a triangle
\[
\begin{tikzcd}
& (e,e) \ar[dr] & \\
e \ar[rr,"\simeq"]\ar[ur] && L_{\one_{\cO}}(e)
\end{tikzcd}
\]
in which the right diagonal map is operadically cocartesian over $\Com^{\otimes}$ and the bottom horizontal map is an equivalence since $e$ belongs to $\D_{\one_{\cO}}$, and so the map $\one_{\C} \to e$ exhibits $e$ as idempotent over $\cO$. Finally, for every $a \in \cO$ and $x \in \C_a$ the $q$-cocartesian map $\eta_x \colon (e,x) \to L_a(x)$ induces an equivalence $e \otimes x \simeq L_a(x)$, and hence $\D_a$ coincides with the image of $e \otimes (-)\colon \C_a \to \C_a$, proving \reffouritem{item:idem-1}.
The implication \reffouritem{item:idem-1} $\Rightarrow$ \reffouritem{item:idem-2} follows from Remark~\reffour{remark:action-of-idempotent}. %
To prove \reffouritem{item:idem-2} $\Rightarrow$ \reffouritem{item:closed}, we assume that each inclusion $\D_a \subseteq \C_a$ admits a left adjoint $L_a\colon \C_a \to \D_a$ and that we have an $\cO$-idempotent $\one_{\C} \to e$ %
such that $\sig\colon e \otimes (-) \simeq L_a(-)$ as endo-functors under $\id_{\C_a}$ for every $a \in \cO$. Let us first show that $\{L_a\}_{a \in \cO}$ is $\cO$-multiplicative. For this, we need to show that for every operation $\alp\colon (a_1,...,a_n) \to b$ in $\cO$, any choice of an $L_{a_i}$-equivalence $f_i\colon x_i \to y_i$ in $\C_{a_i}$ for each $i=1,...,n$, and any $d \in \D_b$, the restriction map
\[
\Mul^{\alp}_{\C}(y_1,\dots,y_n;d) \to \Mul^{\alp}_{\C}(x_1,\dots,x_n; d)
\]
on spaces of operations in $\C^{\otimes}$ covering $\alp$, is an equivalence. By Remark~\reffour{remark:full-suboperad}~\reffouritem{item:one-at-a-time} we may assume without loss of generality that all the $f_i$ are identities (or just equivalences) except for $f_1$. 
Now consider the commutative diagram
\[
\begin{tikzcd}
\Mul^{\alp}_{\C}(y_1,\dots,y_n;d) \ar[r]\ar[d] & \Mul^{\alp}_{\C}(x_1,\dots,x_n; d) \ar[d] \\
\Mul^{\alp}_{\C}(e \otimes y_1,y_2,...,y_n;e \otimes d) \ar[r,"\simeq"]\ar[d] & \Mul^{\alp}_{\C}(e\otimes x_1,x_2,...,x_n; e\otimes d) \ar[d] \\
\Mul^{\alp}_{\C}(y_1,\dots,y_n;e \otimes d) \ar[r] & \Mul^{\alp}_{\C}(x_1,\dots,x_n; e\otimes d)  \\
\end{tikzcd}
\]
where the maps from the top row to the middle row are induced by tensoring with $e$ (well-defined since $e$ is a symmetric monoidal object over $\one_{\cO}$),
and the maps from the middle to the bottom row are induced by precomposition with $y_1 \to e \otimes y_1$ and $x_1 \to e \otimes x_1$, respectively. 
Then the middle horizontal map is an equivalence by the assumption that each $f_i\colon x_i \to y_i$ is an equivalence except for $f_1$, which is an $L_{a_1}$-equivalence (and given that we have an equivalence $L_{a_1}(-) \simeq e \otimes (-)$). At the same time, the composite vertical maps coincide with the ones induced by post-composition with $d \to e \otimes d$ by naturality of concatenation and of factorisation through cocartesian lifts, and is hence an equivalence by the assumption that $d$ belongs to $\D_a$. The entire diagram is hence equivalent to a retract diagram in the arrow category of spaces, and since the middle row is an equivalence we conclude that the top (and bottom) horizontal arrows are equivalences, as desired.

Moving forward, note that since $e \simeq L_{\one_{\cO}}(\one_\C)$ is assumed $\cO$-idempotent, it is in particular symmetric monoidal over $\one_{\cO}$.
To establish~\reffouritem{item:closed} it is hence left to show that for any $a \in \cO$, $x \in \C_{\mon/\one_\cO}$, and $y \in \C_a$, if either $x$ belongs $\D_{\one_{\cO}}$ or $y$ belongs to $\D_{a}$ then $x \otimes y \in \C_a$ belongs to $\D_a$. Since $L_a(-) \simeq e \otimes (-)$ this amounts to showing that the induced map 
\[
\eps\otimes x \otimes y\colon x \otimes y \to e \otimes x \otimes y
\]
is an equivalence. Now if $x$ belongs to $\D_{\one_{\cO}}$ then $\eps \otimes x \colon x \to e \otimes x$ is an equivalence and hence $\eps \otimes x \otimes y$ is an equivalence as well, and if $y$ belongs to $\D_a$ then $\eps \otimes y \colon y \to e \otimes y$ is an equivalence and hence $\eps \otimes x \otimes y \simeq x \otimes \eps \otimes y$ is an equivalence as well. 

Finally, we prove \reffouritem{item:closed} $\Rightarrow$ \reffouritem{item:smashing}. By Proposition~\reffour{proposition:multiplicative-localisation} we have that the collection $\{L_a\}_{a \in \cO}$ refines to a left Bousfield localisation $L^{\otimes}\colon \C^{\otimes} \to \D^{\otimes}$ of $\infty$-operads over $\cO^{\otimes}$. It is left to show that $L^{\otimes}$ is $\cO$-smashing. %
Let us set $e := L(\one_{\C})$ and for $a \in \cO$ and $x \in \C_a$ consider the commutative triangle
\[
\begin{tikzcd}
& (e,x) \ar[dr,"\eta_x"] & \\
x \ar[rr]\ar[ur] && L(x) 
\end{tikzcd}
\]
resulting from Construction~\reffour{construction:eta}. We need to show that $\eta_x$ is operadically cocartesian over $\Com^{\otimes}$. Let us first point out that since $e$ is assumed symmetric monoidal over $\one_{\cO}$
we may factor $\eta_x$ as a composite 
\[
(e,x) \to e \otimes x \xrightarrow{\eta_x'} L(x)
\]
where the first map is an operadically cocartesian lift of $\<2\> \to \<1\>$ and maps to an operadically cocartesian 
lift of this map in $\cO$, while $\eta_x'$ is a map that lies in the fibre $\C_a$ over $a \in \cO$. It will hence suffice to show that $\eta_x'$ is an equivalence. 
By assumption, we have that $e \otimes x$ belongs to $\D_a$. Since $L(x)$ also belongs to $\D_a$ %
it is enough to show that $\eta_x'$ is a $\D_a$-local equivalence, that is, that for every $d \in \D_a$ the induced map 
\[
\Map_{\C_a}(L(x),d) \to \Map_{\C_a}(e \otimes x,d) \simeq \Mul_{\C}(e,x; d)
\]
is an equivalence. By the multiplicativity of $\{L_a\}_{a \in \cO}$ and the unitality of $\one_{\C}$, the restriction maps
\[
\Mul_{\C}(e,x; d) \to \Mul_{\C}(\one_{\C},x; d) \to \Map_{\C_a}(x,d)
\]
are both equivalences. Finally, the total composite 
\[
\Map_{\C_a}(L(x),d) \to \Map_{\C_a}(x,d)
\]
coincides with restriction along $x \to L(x) = L_a(x)$, and is hence an equivalence for every $d \in \D_a$, completing the proof of the proposition.
\end{proof}

\begin{corollary}
\label{corollary:action-smashing}%
Let $\cO^{\otimes}$ be an $\infty$-operad with unit $\one_{\cO}$, $q\colon \C^{\otimes} \to \cO^{\otimes}$ a map of $\infty$-operads with unit $\one_{\C}$ over $\one_\cO$, and $L^{\otimes} \colon \C^{\otimes} \to \D^{\otimes}$ an $\cO$-smashing left Bousfield localisation. 
Write $e := L(\one_{\C})$ for the image of $\one_{\C}$, so that $e$ is a unit object in $\D$ by Remark~\reffour{remarks:symmetric-monoidal-objects-and-units}~\reffouritem{item:local-unit}.
Then in the square
\[
\begin{tikzcd}
\Mod_e(\D)^{\otimes} \ar[r]\ar[d,"\simeq"'] & \Mod_e(\C)^{\otimes} \ar[d]\\
\D^{\otimes} \ar[r] & \C^{\otimes},
\end{tikzcd}
\]
where the left vertical map is an equivalence by Proposition~\reffour{proposition:units},
the top horizontal map is also an equivalence, yielding an equivalence $\D^{\otimes} \simeq \Mod_e(\C)^{\otimes}$ of $\infty$-operads over $\C^{\otimes}$.
\end{corollary}
\begin{proof}
Since $\D^{\otimes} \to \C^{\otimes}$ is a full suboperad inclusion, the top horizontal map in the displayed diagram is again a full suboperad inclusion. It will hence suffice to show that it is essentially surjective. In other words, we need to show that if $x \in \C^{\otimes}$ refines to an $e$-module object then $x$ lies in $\D$. Now, by Proposition~\reffour{proposition:equivalent-1}, the map $\eps\colon \one_{\C} \to e$ exhibits $e$ as idempotent over $\cO$, and so $e$ is in particular symmetric monoidal over $\one_{\cO}$.
We may consequently view $e$ as a commutative algebra object in the symmetric monoidal $\infty$-category $\C_{\mon/\one_\cO}$. 
Since $e$ lies over the trivial algebra object $\one_{\cO}$ in $\cO^{\otimes}$, the object $a$ admits an essentially unique $\one_{\cO}$-module structure, and so any $e$-module structure on $x$ will induce this structure on $a := q(x)$ upon applying $q$. 
In particular, any $e$-module structure on $x$ uniquely refines to an $e$-module structure on $x$ as an object of $\C_a$, where we view $\C_a$ as a module over the symmetric monoidal $\infty$-category $\C_{\mon/\one_\cO}$ (see Proposition~\reffour{proposition:action-of-mon}) and $e$ as a commutative algebra object there. In particular, if $x$ is an $e$-module in $\C_a$ then the unit and action maps
\[
x \to e \otimes x \to x
\]
exhibit $x$ as a retract of $e \otimes x$. Applying Proposition~\reffour{proposition:equivalent-1} again we have that $e \otimes x$ belongs to $\D_a$, and hence $x$ belongs to $\D_a$ as well (since reflective subcategories are closed under retracts).  
\end{proof}

Finally, let us also summarise the conclusions of the above discussion in the case where $\cO^{\otimes} = \Com^\otimes$ and $\C^{\otimes}$ is a symmetric monoidal $\infty$-category.

\begin{corollary}
\label{corollary:smashing-symmetric-monoidal}%
Let $L^{\otimes} \colon \C^{\otimes} \to \D^{\otimes}$ be a left Bousfield localisation of symmetric monoidal $\infty$-categories with fully faithful (and lax symmetric monoidal) right adjoint $\iota^{\otimes}\colon \D^{\otimes} \hrar \C^{\otimes}$.
Then the following conditions are equivalent:
\begin{enumerate}
\item
\label{item:counit-tensor}%
$L$ is smashing. %
\item
\label{item:left-adjoint-is-tensoring}%
There exists an idempotent $e \in \C_{\one/}$ such that $\iota L$ is equivalent, as a functor under $\id_{\C}$, to $x \mapsto e \otimes x$.
\item
\label{item:closed-under-tensoring}%
The image of $\iota\colon \D \hrar \C$ is closed under tensoring with arbitrary objects in $\C$, that is, it is an ideal.
\item
\label{item:smashing-as-algebra}%
There exists an idempotent algebra object $e \in \CAlg(\C^{\otimes})$ and a symmetric monoidal equivalence $\D^{\otimes} \simeq \Mod_e(\C)^{\otimes}$ under which $\iota$ identifies with the forgetful functor $\Mod_e(\C) \to \C$.
\end{enumerate}
\end{corollary}
\begin{proof}
The equivalence of the first three condition follows from Proposition~\reffour{proposition:equivalent-1}, while Corollary~\reffour{corollary:action-smashing} shows that~\reffouritem{item:smashing-as-algebra} follows from~\reffouritem{item:counit-tensor}. To finish the proof it will hence suffices to note that~\reffouritem{item:smashing-as-algebra} implies~\reffouritem{item:closed-under-tensoring}, since the forgetful functor $\Mod_e(\C) \to \C$ refines to one of $\C$-modules. %
\end{proof}

\begin{remark}
\label{remark:smashing-closed-colimits}%
In the situation of Corollary~\reffour{corollary:smashing-symmetric-monoidal}, if $\C$ admits small colimits and its tensor product preserves small colimits in each variable then $\iota(\D) \subseteq \C$ is closed under colimits by the characterisation in \reffouritem{item:smashing-as-algebra}. Let us emphasise that the latter condition is strictly weaker than being smashing, see Example~\reffour{example:cant-remove} below for a multiplicative localisation closed under colimits which is not smashing. Nonetheless, since in any multiplicative localisation $\iota(\D)$ is closed under tensoring with dualisable objects, if $\C$ is generated under colimits by dualisable objects and $\iota(\D) \subseteq \C$ is closed under colimits then the localisation is smashing by the characterisation in~\reffouritem{item:closed-under-tensoring}. This is the case, for example, when $\C=\Spa$, whose smashing localisations have been the subject of much study.
\end{remark}

\subsection{Multiplicative localisations of Day convolution}
\label{subsection:day-convolution}%
\label{subsec:large-indexing}%

In this subsection we specialise the general discussion above to localisations of Day convolution $\infty$-operads obtained by requiring a certain collection of cone diagrams to be sent to limits. Our main goal is to establish Propositions~\reffour{proposition:multiplicative-day} and \reffour{proposition:smashing-day} below, which give sufficient criterion for such localisations to be multiplicative and smashing, respectively.

Let $\cO^{\otimes}$ be an $\infty$-operad and $\C^{\otimes} \to \cO^{\otimes}$ an $\cO$-monoidal $\infty$-category. Then 
the functor
\[
(-) \times_{\cO^{\otimes}} \C^{\otimes} \colon (\Op)_{/\cO} \to (\Op)_{/\cO}
\]
admits a right adjoint
\[
\Fun_{\cO}(\C,-)^{\Day}\colon (\Op)_{/\cO} \to (\Op)_{/\cO}
\]
given by the formation of Day convolution, see~\cite{HA}*{Construction 2.2.6.7}; an alternative proof is given in~\cite{hinich-enriched}*{Proposition 2.8.3}, which also gives a simple explicit model in the case where both $\C^{\otimes}$ and $\D^{\otimes}$ are $\cO$-monoidal. The latter context is in fact sufficient for our purposes, and we shall henceforth fix the assumption that $\D^{\otimes}$ (in addition to $\C^{\otimes}$) is $\cO$-monoidal.
The universal property of Day convolution then gives %
\[
\Alg_{\cP/\cO}(\Fun_{\cO}(\C,\D)^{\Day}) = \Alg_{\cP \times_{\cO} \C/\cO}(\D) =  \Alg_{\cP \times_{\cO} \C/\cP}(\cP \times_{\cO}\D) ,
\]
that is, the $\infty$-category of $\cP$-algebras in $\Fun_{\cO}(\C,\D)^{\Day}$ is equivalent to the $\infty$-category of lax $\cP$-monoidal functors $\C \to \D$ for any $\infty$-operad $\cP^{\otimes} \to \cO^{\otimes}$ over $\cO^{\otimes}$.

For every $a \in \cO$, the $\infty$-category of colours of $\Fun_{\cO}(\C,\D)^{\Day}$ lying over $a$ is given by the functor category $\Fun(\C_a,\D_a)$. The spaces of operations in $\Fun_{\cO}(\C,\D)^{\Day}$ can then be described explicitly as follows. Given an operation $\alp\colon (a_1,...,a_n) \to b$ in $\cO$ with associated $\cO$-monoidal operation $\alp_{\D}\colon \D_{a_1} \times \dots \times \D_{a_n} \to \D_b$, functors $f_i\colon \C_{a_i} \to \D_{a_i}$ for $i=1,...,n$ and a functor $g \colon \C_b \to \D_b$, the corresponding space of operations over $\alp$ is given by the space of natural transformations
\[
\Mul^{\alp}_{\Fun_{\cO}(\C,\D)}(f_1,...,f_n;g) = \Nat(f_1 \boxtimes \dots \boxtimes f_n,g|_{\C_{a_1} \otimes \dots \otimes \C_{a_n}} ),
\]
between the composite 
\[
f_1 \boxtimes \dots \boxtimes f_n\colon \C_{a_1} \times \dots \times \C_{a_n} \to \D_{a_1} \times \dots \times \D_{a_n} \xrightarrow{\alp_{\D}} \D_b
\]
and the restriction of $g$ along the $\alp$-product
$\alp_{\C}\colon \C_{a_1} \times \dots \times \C_{a_n} \to \C_{b}$ in $\C^{\otimes}$. %

\begin{definition}
\label{definition:cocomplete}%
Let $\D^{\otimes} \to \cO^{\otimes}$ be an $\cO$-monoidal $\infty$-category.  
\begin{enumerate}
\item
We say that $\D^{\otimes}$ is \defi{$\cO$-cocomplete} if $\D_a$ admits small colimits for every $a \in \cO$ and $\alp_{\D}\colon \D_{a_1} \times ... \times \D_{a_n} \to \D_b$ is entry-wise colimit preserving 
for every arrow $\alp\colon (a_1,...,a_n) \to b$ in $\cO^{\otimes}$ with target in $\cO$.
\item
We say that $\D^{\otimes}$ is \defi{$\cO$-closed} if for every $\alp\colon (a_1,...,a_n) \to b$ in $\cO^{\otimes}$ with target in $\cO$ and every $(d_2,\dots,d_n) \in \D_{a_2} \times \dots \times \D_{a_n}$, the operation $\alp_{\D}(-,d_2,...,d_n)\colon \D_{a_1} \to \D_b$ admits a right adjoint $R^{d_2,...,d_n}_{\alp}\colon \D_b \to \D_{a_1}$;
\item We say that $\D$ is \defi{presentably} $\cO$-monoidal if it is $\cO$-cocomplete and each $\D_a$ is presentable.
\end{enumerate} 
\end{definition}

By the adjoint functor theorem, any presentably $\cO$-monoidal $\infty$-category is $\cO$-closed.

The following proposition is a variant of~\cite{HA}*{Proposition 2.2.6.16}. We include the proof for the convenience of the reader.

\begin{proposition}
\label{proposition:day-cocartesian}%
Let $\D^{\otimes} \to \cO^{\otimes}$ be an $\cO$-cocomplete $\cO$-monoidal $\infty$-category and $\C$ a small $\cO$-monoidal $\infty$-category. Then the Day convolution $\infty$-operad $\Fun_{\cO}(\C,\D)^{\Day} \to \cO^{\otimes}$ is an $\cO$-monoidal $\infty$-category.
\end{proposition}

In fact, in the situation of Proposition~\reffour{proposition:day-cocartesian} one can show that $\Fun_{\cO}(\C,\D)^{\Day}$ is again $\cO$-cocomplete, and is further presentably $\cO$-monoidal if $\D$ is so. %

\begin{proof}[Proof of Proposition~\reffour{proposition:day-cocartesian}]
An arrow $\alp\colon (a_1,...,a_n) \to b$ in $\cO^{\otimes}$ with target in $\cO$ admits a locally cocartesian lift in $\Fun_{\cO}(\C,\D)^{\Day}$ with source $(f_1,...,f_n)$ if and only if 
\[
f_1 \boxtimes \dots \boxtimes f_n\colon \C_{a_1} \times ... \times \C_{a_n} \to \D_b
\]
admits a left Kan extension along $\alp_{\C}\colon \C_{a_1} \times ... \times \C_{a_n} \to \C_b$. This holds when $\C$ is small and each $\D_a$ admits small colimits. 
Theses locally cocartesian lifts then remain locally cocartesian after concatenation with fixed objects (since locally cocartesian arrows are closed under concatenation in any $\infty$-operad). The condition that all $\cO$-monoidal operation functors preserves small colimits in each entry then implies that all these locally cocartesian arrows are actually cocartesian.
Since cocartesian arrows are closed under composition we get that $\Fun_{\cO}(\C,\D)^{\Day} \to \cO^{\otimes} \to \cO^{\otimes}$ admits a sufficient supply of cocartesian lifts, as desired.
\end{proof}

\begin{examples}\
\begin{enumerate}
\item
Suppose $\cO=\Com$, $\C$ is a small symmetric monoidal $\infty$-category and $\D$ is a $\Com$-cocomplete symmetric monoidal $\infty$-category, that is, $\D$ admits small colimits and the monoidal product preserves small colimits in each entry. 
Then $\Fun_{\Com}(\C,\D)^{\Day}$ is the symmetric monoidal $\infty$-category whose underlying $\infty$-category is the $\infty$-category of functors $\C \to \D$ and whose tensor product $(f,g) \mapsto f \otimes_{\Day} g$ is given by a left Kan extension of the form 
\[
\begin{tikzcd}
\C \times \C \ar[d]\ar[r, "f \times g"] & \D \times \D \ar[r, "\otimes"] & \D \\
\C \ar[urr, dotted] && \ .
\end{tikzcd}
\]
This product is what is classically known as \emph{Day convolution}, see~\cite{day1970} for the original treatment of the idea in the setting of ordinary categories. 
symmetric monoidal $\infty$-category.
\item
Suppose that $\cO = \MCom$ (see Example~\reffour{example:mcom-loc}), $\C = (\C_a,\C_m)$ is small $\MCom$-monoidal $\infty$-category and $\D = (\cB,\E)$ is an $\MCom$-monoidal $\infty$-category %
such that $\cB$ is a $\Com$-cocomplete
symmetric monoidal $\infty$-category and $\E$ is a $\cB$-module in cocomplete $\infty$-categories.
Then the $\MCom$-monoidal $\infty$-category $\Fun^{\MCom}(\C,\D)^{\Day}$ encodes the data of the pair $(\Fun(\C_a,\cB),\Fun(\C_m,\E))$ where $\Fun(\C_a,\cB)$ is the %
symmetric monoidal functor category with its Day convolution as above, and $\Fun(\C_m,\E)$ is a $\Fun(\C_a,\cB)$-module %
with action of $f$ on $g$ given similarly by a left Kan extension of the form 
\[
\begin{tikzcd}
\C_a \times \C_m \ar[d]\ar[r, "f \times g"] & \cB \times \E \ar[r, "\otimes"] & \E \\
\C_m \ar[urr, dotted] && \ .
\end{tikzcd}
\]
\end{enumerate}
\end{examples}

\begin{remark}
\label{remark:left-kan-monoidal}%
Suppose given an $\cO$-monoidal functor $\rho^{\otimes}\colon \cB^{\otimes} \to \C^{\otimes}$ among small $\cO$-monoidal $\infty$-categories and let $\D^{\otimes} \to \cO^{\otimes}$ be an $\cO$-cocomplete $\cO$-monoidal $\infty$-category, so that $\Fun_{\cO}(\cB,\D)^{\Day}$ and $\Fun_{\cO}(\C,\D)^{\Day}$ are $\cO$-monoidal by Proposition~\reffour{proposition:day-cocartesian}. Then $\rho$ induces a (solidly indicated) restriction functor 
\[
\begin{tikzcd}
\Fun_{\cO}(\cB,\D)^{\Day} \ar[rr,bend left =15,dotted,"\rho_!"]\ar[dr]&& \ar[ll,"\rho^*","\perp"'] \Fun_{\cO}(\C,\D)^{\Day}\ar[dl] \\
& \cO^{\otimes} &
\end{tikzcd}
\]
which is generally only a map of $\infty$-operads over $\cO$ (that is, a lax $\cO$-monoidal functor), where for every $a \in \cO$, the induced map $\Fun(\C_a,\D_a) \to \Fun(\cB_a,\D_a)$ is given by restriction along $\rho$. Since each $\D_a$ is assumed cocomplete each of these restriction functors admits a left adjoint given by left Kan extension,  
and these left Kan extensions assemble to form a relative left adjoint $\rho_!$ over $\cO$ by~\cite{HA}*{Corollary 7.3.2.12}, which is automatically %
an $\cO$-monoidal functor (since relative left adjoints preserve relative colimits). 
We may summarise the situation as follows: Left Kan extension along symmetric monoidal functors among small symmetric monoidal $\infty$-categories induces an $\cO$-monoidal functor on the level of Day convolution $\infty$-operads
with $\cO$-cocomplete targets. 
\end{remark}

To start the passage to arbitrary locally small
source category $\C$, note that if $\C_0 \subseteq \C$ is some small full subcategory and $f\colon \C_0 \to \D$ is a functor with $\D$ cocomplete, then $f$ admits a left Kan extension to all of $\C$ by the pointwise formulae. We refer to functors which are obtained in the manner as \defi{small}, and write $\Fun^{\sm}(\C,\D) \subseteq \Fun(\C,\D)$ for the full subcategory spanned by the small functors. These are all generated under small colimits from functors of the following form:

\begin{notation}
\label{notation:representable}%
Given $c \in \C$ and $d \in \D$, we will denote by $c_!d \colon \C \to \D$ the functor obtained by left Kan extending the functor $\{c\} \to \D$ sending $c$ to $d \in \D$ along the inclusion $\{c\} \subseteq \C$. 
\end{notation}

Explicitly, $c_!d$ is given on some object $c' \in \C$ by the formula $d_{\Hom_\C(c,c')} \simeq \colim_{\Hom_\C(c,c')} \const_d$.

\begin{lemma}
\label{lemma:canonical-colimit}%
Let $f\colon \C \to \D$ be a functor from a locally small $\infty$-category $\C$ to an $\infty$-category with colimits over (small) spaces, and let $\vphi\colon \Twar(\C) \to \Fun(\C,\D)$ be the functor given by the formula 
\[
[x \to y] \mapsto y_!f(x) .
\]
Then the canonical cone diagram $\ovl{\vphi}\colon \Twar(\C)^{\triangleright} \to \Fun(\C,\D)$ on $\vphi$ with cone point $f$ is a colimit diagram.
\end{lemma}

We note that in the situation of Lemma~\reffour{lemma:canonical-colimit} the $\infty$-category $\Twar(\C)$ is potentially not small, and hence the colimit of $\vphi$ does not exist a-priori even if we ask $\D$ to be cocomplete, only a-posteriori. 

\begin{remark}
\label{remark:canonical-colimit}%
If $f$ is left Kan extended from some small full subcategory $\C_0 \subseteq \C$ then we may apply this lemma to $f|_{\C_0}$ to deduce that $f$ is a small colimit of functors of the form $c_!d$.
\end{remark}

\begin{proof}[Proof of Lemma~\reffour{lemma:canonical-colimit}]
Since colimits in functor $\infty$-categories that exists pointwise are also computed pointwise it will suffice to show that for every $z \in \C$ the corresponding cone diagram
\[
\ev_{z} \circ \ovl{\vphi}\colon \Twar(\C)^{\triangleright} \to \D
\]
with cone point $f(z)$ is a colimit diagram. Explicitly, this diagram sends $[x \to y] \in \Twar(\C)$ to $f(x)_{\Hom_{\C}(y,z)} \in \D$. In particular, we see that it is the left Kan extension along the projection 
\[
\Twar(\C_{/z}) \simeq \Twar(\C) \times_{\C\op} (\C_{/z})\op \to \Twar(\C)
\]
of the functor $[x \to y \to z] \mapsto f(x)$. It will hence suffice to show that the corresponding cone diagram
\[
\Twar(\C_{/z})^{\triangleright} \to \D
\]
with cone point $f(z)$ is a colimit diagram. We now observe that both the functor $[x \to y \to z] \mapsto f(x)$ and the natural transformation from it to $f(z)$ do not depend on $y$, and so the last cone diagram factors as
\[
\Twar(\C_{/z})^{\triangleright} \to \C_{/z}^{\triangleright} \to \D
\]
where the first functor sends $[x \to y \to z]$ to $[x \to z]$. Now $\Twar(\C_{/z}) \to \C_{/z}$ is a cartesian fibration with contractible fibres and is hence cofinal, see e.g.\ \cite{HHLN-two-var-fil}*{Lemma 5.5}. It will hence suffice to show that the cone diagram
\[
\C_{/z}^{\triangleright} \to \D
\]
is a colimit diagram. Indeed, this is because $\C_{/z}$ has a final object and the above diagram sends the arrow from the final object to the cone point to the identity on $f(z)$. 
\end{proof}

Now unlike the situation in Proposition~\reffour{proposition:day-cocartesian}, when $\C$ is only locally small the Day $\infty$-operad $\Fun_{\cO}(\C,\D)^{\Day}$ is generally not an $\cO$-monoidal $\infty$-category, even if $\D$ is $\cO$-cocomplete. In the following proposition and subsequent corollary we summarise a few properties which nonetheless persist to the locally small setting:

\begin{proposition}
\label{proposition:persist-locally-small}%
Let $\C$ be a locally small $\cO$-monoidal $\infty$-category and $\D$ an $\cO$-cocomplete $\cO$-monoidal $\infty$-category. 
Let $\alp\colon (a_1,\dots,a_n) \to b$ be an operation in $\cO$ and $f_i\colon \C_{a_i} \to \D_{a_i}$ be functors. Then $\alp$ admits an operadically cocartesian lift to $\Fun_{\cO}(\C,\D)$ starting from $(f_1,...,f_n)$ in either of the following cases:
\begin{enumerate}
\item
\label{item:all-small}%
The $f_i$ are all small. In particular, if $n=0$.
\item
The $f_i$ are small for all $i$ except at most one $i_0 \in \{1,...,n\}$ and $\C^{\otimes}$ is $\cO$-closed.
\end{enumerate}
In addition, in case~\reffouritem{item:all-small} the target of the resulting operadically cocartesian lift is again small.
\end{proposition}

\begin{corollary}
\label{corollary:persist-locally-small}%
Let $\C$ be a locally small $\cO$-monoidal $\infty$-category and $\D$ an $\cO$-cocomplete $\cO$-monoidal $\infty$-category. Then the following holds:
\begin{enumerate}
\item
The full suboperad $\Fun_{\cO}^{\sm}(\C,\D)^{\Day} \subseteq \Fun_{\cO}(\C,\D)^{\Day}$ spanned by $\Fun^{\sm}(\C,\D_a) \subseteq \Fun(\C,\D_a)$ for $a \in \cO$ is an $\cO$-monoidal $\infty$-category.
\item
If $\cO$ admits a unit object $\one_{\cO}$ then $\Fun_{\cO}(\C,\D)^{\Day}$ admits a unit object lying above $\one_{\cO}$.
\item
\label{item:small-monoidal}%
If $\C$ is $\cO$-closed and $a \in \cO$ is a symmetric monoidal object then any small $f\colon \C_a \to \D_a$ is a symmetric monoidal object in $\Fun_{\cO}(\C,\D)^{\Day}$ lying over $a$.
\end{enumerate}
\end{corollary}

\begin{proof}[Proof of Proposition~\reffour{proposition:persist-locally-small}]
Showing the existence of the desired operadically cocartesian lift amounts to showing that there exists a dashed left Kan extension 
\[
\begin{tikzcd}[column sep = 50pt]
\C_{a_1} \times \dots \times \C_{a_n} \ar[r,"f_1 \times ... \times f_n"]\ar[d,"\alp_{\C}"'] & \D_{a_1} \times \dots \times \D_{a_n} \ar[d,"\alp_{\D}"] \\
\C_b \ar[r,dashed] & \D_b
\end{tikzcd}
\]
of the top right composite along the left vertical arrow, which furthermore remains a left Kan extension after composing with all functors of the form $\beta(-,d_1,...,d_r)\colon \D_b \to \D_{b'}$ for operations $\beta\colon (b,b_1,...,b_r) \to b'$ and objects $d_i \in \D_{b_i}$. 
Let now $i_0 \in \{1,...,n\}$ be such that $f_i$ is small for $i \neq i_0$ and that either $f_{i_0}$ is also small or $\C$ is $\cO$-closed. The assumption that $\D$ is $\cO$-cocomplete implies that $\D_{a_i}$ and $\D_b$ admit small colimits and $\alp_D$ %
preserves small colimits in each argument, so that %
the top right composite depends on $f_1,...,f_n$ in a manner that preserves small colimits in each argument. Since $\D_b$ admits small colimits and the functors of the form $\beta(-,d_1,...,d_r)$ preserve small colimits by assumption the collection of $(f_1,....,f_n)$ for which the desired statement hold is closed under small colimits in each $f_i$ separately. Now by Lemma~\reffour{lemma:canonical-colimit} and Remark~\reffour{remark:canonical-colimit} small functors are all obtained from the basic ones of Notation~\reffour{notation:representable} via small colimits, and so we may assume without loss of generality that $f_{i} = (c_i)_!d_i$ for some $c_i \in \C,d_i \in \D$ for every $i \neq i_0$.
Using again that $\alp_{\D}$ commutes with small colimits in each argument we then get that 
\[
(c_1)_!d_1 \boxtimes \dots \boxtimes f_{i_0} \boxtimes \dots \boxtimes (c_n)_!d_n \colon \C \times \dots \times \C \to \D_{a_1} \times \dots \times \D_{a_n} \xrightarrow{\alp_{\D}} \D_b
\]
is left Kan extended from $\{c_1\} \times \dots \times \C \times \dots \times \{c_n\}$, and the same hold after post-composing with functors of the form $\beta(-,d_1,...,d_r)\colon \D_b \to \D_{b'}$.
It will hence suffice to show that the left Kan extension of $\alp_{\D}(d_1,...,f_{i_0},...,d_n)\colon \C_{a_{i_0}} \to \D_b$
along the functor $\alp_{\C}(c_1,...,-,...,c_n)\colon \C_{a_{i_0}} \to \C_b$
exists and preserved by post-composing with all functors of the form $\beta(-,d_1,...,d_r)$. If $\C$ is $\cO$-closed then the functor $\alp_{\C}(c_1,...,-,...,c_n)$ admits a right adjoint $R_{\alp}^{c_1,...,c_n}\colon \C_b \to \C_{a_{i_0}}$, and so left Kan extension along $\alp_{\C}(c_1,...,-,...,c_n)$ is given by pre-composition with $R_{\alp}^{c_1,...,c_n}$, an operation which commutes with post-composition with any functor. If $\C$ is not $\cO$-closed, then by assumption $f_{i_0}$ is also small, and so the argument above could also be applied at $i_0$ to reduced to the case where $f_{i_0} = (c_{i_0})_!d_{i_0}$ for some $c_{i_0} \in \C_{a_{i_0}}$ and $d_{i_0} \in \D_{a_{i_0}}$, in which case it will suffice to show that the functor $\ast \to \D_b$ picking the object 
$\alp_{\D}(d_1,...,d_n) \in \D_b$ admits a left Kan extension to $\C_b$ along the functor $\ast \to \C_b$ picking the object $\alp_{\C}(c_1,...,c_n)$, and that this left Kan extension is preserved by post-composition with any $\beta(-,d_1,...,d_r)$. Indeed, since $\C_b$ is locally small such left Kan extension exists as a pointwise one with small comma $\infty$-categories, and so the desired result follows. In addition, in this case we also get that the desired left Kan extension is itself small by definition, so that the last statement of Proposition~\reffour{proposition:persist-locally-small} follows from the fact that the collection of small functors is closed under small colimits.
\end{proof}

\begin{remark}
\label{remark:units}%
Unwinding the proof of Proposition~\reffour{proposition:persist-locally-small}
shows that if $\alp\colon \emptyset\to b$ is some nullary operation in $\cO$ with corresponding unit objects $\one^{\alp}_{\C} \in \C_b$ and $\one^{\alp}_{\D} \in \D_b$, then the corresponding unit object in $\Fun(\C_b,\D_b)$ is given by $(\one^{\alp}_{\C})_!\one^{\alp}_{\D}$ using Notation~\reffour{notation:representable}. In particular, if $\alp$ exhibits $b$ as a unit object of $\cO$ then the induced map $\emptyset \to (\one^{\alp}_{\C})_!\one^{\alp}_{\D}$ exhibits its target as a unit object in $\Fun_{\cO}(\C,\D)^{\Day}$ over $b$,
see Definition~\reffour{definition:monoidal-and-unit}. In this case, by~\cite{HA}*{Proposition 3.2.1.8 and Lemma 3.2.1.10} the objects $b$ and $(\one^{\alp}_{\C})_!\one^{\alp}_{\D}$ carry essentially unique commutative algebra structures, inducing symmetric monoidal structures on $\C_b$ and $\D_{b}$ together with a compatible lax symmetric monoidal structure on the functor $(\one^{\alp}_{\C})_!\one^{\alp}_{\D} \colon \C_b \to \D_b$ (by the universal property of Day convolution).
\end{remark}

\begin{definition}
\label{definition:almost-presentable}%
\label{definition:presentable}%
Let $\D^{\otimes} \to \cO^{\otimes}$ be an $\cO$-monoidal $\infty$-category. We say that $\D$ is \defi{almost closed} if it is $\cO$-cocomplete and every $\D_a$ is generated under small colimits by a 
(not necessarily small) full subcategory $\D'_a \subseteq \D_a$ 
satisfying the following property: for every operation $\alp\colon (a_1,...,a_n) \to b$ in $\cO$ and every $(d_2,\dots,d_n) \in \D'_{a_2}\times \dots \times \D'_{a_n}$, the operation $\alp_{\D}(-,d_2,...,d_n)\colon \D_{a_1} \to \D_b$ admits a right adjoint $R^{d_2,...,d_n}_{\alp}\colon \D_b \to \D_{a_1}$.
\end{definition}

\begin{remark}
Every closed $\cO$-monoidal $\infty$-category is almost closed. In particular, every presentably $\cO$-monoidal $\infty$-category is almost closed.
\end{remark}

\begin{example}
\label{example:mcom-almost-presentable}%
Suppose that $\cO^{\otimes} = \MCom^{\otimes}$, so that the data of a cocartesian fibration $\D^{\otimes} \to \MCom^{\otimes}$ corresponds to a pair $(\C,\M)$ where $\C$ is a symmetric monoidal $\infty$-category and $\M$ is a $\C$-module. Now if $\M$ is any cocomplete $\infty$-category, then $\M$ is canonically a module over $\Sps$ in a small colimit compatible manner (since $\Sps$ is the unit for the Lurie tensor product on cocomplete $\infty$-categories) in which case the pair $(\Sps,\M)$ is an $\MCom$-cocomplete $\MCom$-monoidal $\infty$-category. We claim that such an $(\Sps,\M)$ is almost closed. 
To see this, note that for every $m \in \M$ the functor $(-) \otimes m\colon \Sps \to \M$ admits a right adjoint $\Map_\M(m,-)$, and $\Sps$ is generated under colimits by its unit $\ast \in \Sps$, whose action $\ast \otimes (-) = \id_\M$ admits a right adjoint. Similarly, if $\M$ is a cocomplete additive $\infty$-category then $\M$ is canonically a module over $\Grp_{\Einf}$ in a manner compatible with small colimits and the pair $(\Grp_{\Einf},\M)$ constitutes an almost closed %
$\MCom$-monoidal $\infty$-category,
and if $\M$ is a cocomplete stable $\infty$-category then $\M$ is canonically a module over $\Spa$ in a manner compatible with small colimits and the pair $(\Spa,\M)$ constitutes an almost closed 
$\MCom$-monoidal $\infty$-category.
\end{example}

\begin{notation}
\label{notation:cone-diagrams}%
Let $\C$ be a locally small $\infty$-category equipped with a collection of cone diagrams $\cR = \{p_{\alp}\colon K_{\alp}^{\triangleleft} \to \C\}_{\alp}$. For any $\infty$-category %
$\E$ we write
$\Fun^{\cR}(\C,\E) \subseteq \Fun(\C,\E)$ for the full subcategory spanned by the $\cR$-local functors, that is, those functors $\C \to \E$ which send every cone diagram in $\cR$ to a limit cone. We will say that a natural transformation $\vphi \Rightarrow \vphi'$ is an $\cR$-local equivalence if it induces an equivalence on spaces of natural transformations to any $\cR$-local functor $\psi$. We will say that $\vphi \Rightarrow \vphi'$ is an $\cR$-local approximation if it is an $\cR$-local equivalence and $\vphi'$ is $\cR$-local.
We will say that a functor $\vphi\colon \C \to \E$ is \defi{$\cR$-localisable}  if it admits an $\cR$-local approximation $\vphi \Rightarrow \vphi'$. 
In particular, any $\cR$-local functor is $\cR$-localisable with the identity as an $\cR$-local approximation. We then write $\Fun^{\cR-\mathrm{able}}(\C,\E) \subseteq \Fun(\C,\E)$ for the full subcategory spanned by the $\cR$-localisable functors. We note that the inclusion $\Fun^{\cR}(\C,\E) \subseteq \Fun^{\cR-\mathrm{able}}(\C,\E)$ is reflective by design.

Now let $\cO$ be an $\infty$-operad, $\C^{\otimes},\D^{\otimes}$ two $\cO$-monoidal $\infty$-categories and let $\cR = \{\cR_a\}_{a \in \cO}$ be a collection of cone diagrams $\cR_a = \{p_{\alp}\colon K_{\alp}^{\triangleleft} \to \C_a\}$ in $\C_a$ for every $a \in \cO$. We then write $\Fun^{\cR}_{\cO}(\C,\D)^{\Day} \subseteq \Fun_{\cO}(\C,\D)^{\Day}$ for the full suboperad spanned by the $\cR_a$-local functors $f\colon \C_a \to \D_a$ for every $a$, and similarly $\Fun_{\cO}^{\cR-\mathrm{able}}(\C,\D)^{\Day} \subseteq \Fun_{\cO}(\C,\D)^{\Day}$ for the full suboperad spanned by the $\cR_a$-localisable functors for every $a \in \cO$. 

In this situation, we will say that $\cR$ is \defi{closed under $\cO$-tensoring} with objects in $\C$ if for every operation $\alp\colon (a_1,...,a_n) \to b$ in $\cO$ and objects $(c_2,\dots,c_n) \in \C_{a_2}\times \dots \times \C_{a_n}$, the functor $\alp(-,c_2,...,c_n)\colon \C_{a_1} \to \C_b$ sends diagrams in $\cR_{a_1}$ to diagrams in $\cR_b$. If $\C$ is in addition $\cO$-closed, we will say that $\cR$ is \defi{closed under $\cO$-cotensoring} with objects in $\C$ if for every operation $\alp\colon (a_1,...,a_n) \to b$ in $\cO$ and objects $c_i \in \C_{a_i}$ for $i=2,...n$, the right adjoint $R^{c_2,...,c_n}_{\alp}\colon \C_{b} \to \C_{a_1}$ sends diagrams in $\cR_{b}$ to diagrams in $\cR_{a_1}$.

\end{notation}

We now arrive at the first main result of this subsection.

\begin{proposition}
\label{proposition:multiplicative-day}%
Let $\C$ be an $\cO$-monoidal $\infty$-category
and let $\cR = \{\cR_a\}_a$ be a collection of cone diagrams $\cR_a = \{p_{\alp}\colon K_{\alpha}^{\triangleleft} \to \C_a\}_{\alp}$ for every $a \in \cO$. Let $\D$ be an almost closed $\cO$-monoidal $\infty$-category such that $\D_a$ admits $K_{\alp}$-indexed limits for every index $\alp$ and every $a \in \cO$ and suppose that $\cR$ is closed under $\cO$-tensoring with objects in $\C$.
Then the collection of left Bousfield localisations 
\[
\Fun_{\cO}^{\cR-\mathrm{able}}(\C,\D)_a = \Fun^{\cR-\mathrm{able}}(\C_a,\D_a) \locadjL{L_{a}} \Fun^{\cR}(\C_a,\D_a) = \Fun_{\cO}^{\cR}(\C,\D)_a
\]
for $a \in \cO$ is $\cO$-multiplicative in the sense of Definition~\reffour{definition:multiplicative-localisation} and so by Proposition~\reffour{proposition:multiplicative-localisation} assembles to form an $\cO$-multiplicative left Bousfield localisation of $\infty$-operads 
\[
\Fun_{\cO}^{\cR-\mathrm{able}}(\C,\D)^{\Day} \locadjL{L^{\otimes}} \Fun_{\cO}^{\cR}(\C,\D)^{\Day}.
\]
\end{proposition}

\begin{remark}
\label{remark:small-localisable}%
In the situation of Proposition~\reffour{proposition:multiplicative-day}, if $\D^{\otimes}$ is presentably $\cO$-monoidal and $\C$ is small 
then every functor $\C_a \to \D_a$ is $\cR_a$-localisable by the adjoint functor theorem and so, when $\cR$ is closed under $\cO$-tensoring with objects in $\C$, %
Proposition~\reffour{proposition:multiplicative-day} gives an $\cO$-monoidal left Bousfield localisation 
\[
\Fun_{\cO}(\C,\D)^{\Day} \locadjL{L^{\otimes}} \Fun_{\cO}^{\cR}(\C,\D)^\Day.
\]
of $\cO$-monoidal $\infty$-categories.
More generally, if $\C$ is accessible (and $\D^{\otimes}$ still presentably $\cO$-monoidal) and for every sufficiently large regular cardinal $\kappa$, every diagram in $\cR_a$ is a $\kappa$-filtered colimit of diagrams in $\cR_a$ involving only $\kappa$-compact objects, then the adjoint functor theorem implies that every small functor $f\colon \C_a \to \D_a$ admits an $\cR_a$-local approximation which is again small, and so the $\cO$-multiplicative left Bousfield localisation of Proposition~\reffour{proposition:multiplicative-day} restricts to an $\cO$-monoidal left Bousfield localisation 
\[
\Fun_{\cO}^{\sm}(\C,\D)^{\Day} \locadjL{L^{\otimes}_{\sm}} \Fun_{\cO}^{\cR,\sm}(\C,\D)^\Day
\]
on the level of small functors and small $\cR$-local functors, respectively, see Corollary~\reffour{corollary:persist-locally-small}. In this case, the collection of small functors also coincides with the collection of accessible ones.
\end{remark}

\begin{proof}[Proof of Proposition~\reffour{proposition:multiplicative-day}]
By Remark~\reffour{remark:full-suboperad}~\reffouritem{item:one-at-a-time} we need to show that for every operation $\alp\colon (a_1,...,a_n) \to b$ in $\cO$, any $\cR_{a_1}$-local equivalence $\sig\colon f_1 \Rightarrow f_1'$, any collection of functors $g_i\colon \C_{a_i} \to \D_{a_i}$ with $i=2,...,n$, and any $\cR_b$-local functor $h\colon \C_b \to \D_{b}$, the induced map
\[
\Nat(f_1' \boxtimes g_2 \boxtimes ... \boxtimes g_n,h|_{\C_{a_1} \times ... \times \C_{a_n}}) \to \Nat(f_1 \boxtimes g_2 \boxtimes ... \boxtimes g_n,h|_{\C_{a_1} \times ... \times \C_{a_n}})
\]
is an equivalence. We now observe that $\boxtimes$ preserves arbitrary colimits (not only small) in each entry separately by our assumption that $\D$ is almost closed. %
Writing each $g_i$ as a (possibly large) colimit of representables using Lemma~\reffour{lemma:canonical-colimit} we may reduce to the case where for each $i=2,\dots,n$ we have $g_i = (c_i)_!d_i$ for some $c_i \in \C$ and $d_i \in \D_{a_i}$. In addition, by the assumption that $\D^{\otimes} \to \cO^{\otimes}$ is almost closed
we may further assume that $d_i \in \D'_{a_i}$ (using the notation of Definition~\reffour{definition:almost-presentable}) for $i=2,...,n$, so that $\alp_{\D}(-,d_2,...,d_n)\colon \D_{a_1} \to \D_b$ admits a right adjoint $R^{d_2,...,d_n}_{\alp}\colon \D_b \to \D_{a_1}$. We then have
\begin{align*}
\Nat(f_1 \boxtimes (c_2)_!d_2 \boxtimes \dots \boxtimes (c_n)_!d_n, h|_{\C_{a_1}\times ... \times \C_{a_n}}) =& 
\Nat(\alp_{\D}(f_1,d_2,...,d_n), h(\alp_{\C}(-,c_2,...,c_n)))\\
=&\Nat(f_1,R^{d_2,...,d_n}_{\alp}(h(\alp_{\C}(-,c_2,...,c_n))))
\end{align*}
and similarly for $f_1'$.
Since $R^{d_2,...,d_n}_{\alp}$ preserves limits (being a right adjoint) it suffices to check that $h(\alp_{\C}(-,c_2,...,c_n))$
is $\cR_{a_1}$-local. This follows from the assumption that $\alp_{\C}(-,c_2,...,c_n)$ sends diagrams in $\cR_{a_1}$ to diagrams in $\cR_{b}$, and that $h$ is $\cR_b$-local.
\end{proof}

Combining Proposition~\reffour{proposition:multiplicative-day} with Remark~\reffour{remarks:symmetric-monoidal-objects-and-units}~\reffouritem{item:local-unit}
and Proposition~\reffour{proposition:units} %
we therefore conclude:

\begin{corollary}
\label{corollary:multiplicative-day}%
In the situation of Proposition~\reffour{proposition:multiplicative-day}, suppose in addition that $\cO$ admits a unit object $\emptyset \to \one_{\cO}$ (see Definition~\reffour{definition:monoidal-and-unit}), so that $\cB := \D_{\one_{\cO}}$ inherits a symmetric monoidal structure, and $\Fun_{\cO}(\C,\D)^{\Day}$ admits a unit $\one_{\Day} := (\one_\C)_!\one_{\cB}\colon \C_{\one_{\cO}} \to \D_{\one_{\cO}}$ over $\one_{\cO}$ (see Remark~\reffour{remark:units}). Assume that $\one_{\Day}$ is $\cR_{\one_{\cO}}$-localisable with $\cR_{\one_{\cO}}$-local approximation $\one_{\Day} \Rightarrow \one_\Day^{\cR}$. Then the composite $u\colon \emptyset \to \one_{\Day} \to \one_\Day^{\cR}$ exhibits $\one_\Day^{\cR}$ as a unit object in $\Fun_{\cO}^{\cR}(\C,\D)^{\Day}$ over $\one_{\cO}$
in the sense of Definition~\reffour{definition:monoidal-and-unit}. In particular, 
the functor 
\[
\one_\Day^{\cR}\colon \C_{\one_{\cO}} \to \D_{\one_{\cO}}
\]
carries an essentially unique lax symmetric monoidal structure with unit $u$ under which it is the initial $\cR$-local lax symmetric monoidal functor, and $\one_{\Day} \Rightarrow \one_\Day^{\cR}$ refines to a lax symmetric monoidal transformation in an essentially unique manner. In addition, for every $a \in \cO$, any $\cR_a$-local functor $\C_a \to \D_a$ carries an essentially unique $\one_\Day^\cR$-module structure. 
\end{corollary}

We now give a criterion for when the above multiplicative localisation is smashing.

\begin{proposition}
\label{proposition:smashing-day}%
Let $\cO^{\otimes}$ be an $\infty$-operad with a unit object $\emptyset \to \one_{\cO}$. 
Let $\C^{\otimes} \to \cO^{\otimes}$ be an $\cO$-closed $\cO$-monoidal $\infty$-category and let $\cR = \{\cR_a\}_a$ be a collection of cone diagrams $\cR_a = \{p_{\alp}\colon K_{\alpha}^{\triangleleft} \to \C_a\}_{\alp}$ for every $a \in \cO$. Let $\D^{\otimes} \to \cO^{\otimes}$ be an almost closed $\cO$-monoidal $\infty$-category such that $\D_a$ admits $K_{\alp}$-indexed limits for every index $\alp$ and every $a \in \cO$. In particular, $\C$ and $\D$ admits unit objects $\one_{\C} \in \C_{\one_{\cO}}$ and $\one_{\D} \in \D_{\one_{\cO}}$ over $\one_{\cO}$. Assume the following conditions are satisfied:
\begin{enumerate}
\item
$\cR$ is closed under $\cO$-tensoring and $\cO$-cotensoring with objects in $\C$ (see Notation~\reffour{notation:cone-diagrams}),
\item
\label{item:Rlocalclosedcolim}%
for every $a \in \cO$ the collection of $\cR_a$-local functors $\C_a \to \D_a$ is closed under all (potentially large) colimits that exist in $\Fun(\C_a,\D_a)$, 
and
\item
\label{item:OmonoidalpresKlimits}%
all the $\cO$-monoidal operations $\alp_{\D}\colon \D_{a_1} \times ... \times \D_{a_n} \to \D_b$ preserve, for each $i=1,...,n$, and each indexing $\infty$-category $K_{\alp}$ in $\cR_{a_i}$, $K_{\alp}$-indexed limits in the $i$'th argument.
\end{enumerate}
Then the following holds:
\begin{enumerate}
\item
\label{item:ideal}%
Let $a \in \cO$ and consider functors $F\colon \C_{\one_{\cO}} \to \D_{\one_{\cO}}$ and $G\colon \C_a \to \D_a$ at least one of whom is small, so that the Day convolution 
$F\otimes_{\Day} G\colon \C_a \to \D_a$ (covering the unique action of $\one_{\cO}$ on $a$) exists by Corollary~\reffour{corollary:persist-locally-small}\;\reffouritem{item:small-monoidal}. If at least one of $F,G$ is $\cR$-local then $F\otimes_{\Day} G$ is $\cR$-local.
\item
\label{item:equivalent}%
Let $\eta\colon\one_{\Day} = (\one_\C)_!\one_{\D} \Rightarrow F$ be a natural transformation of functors $\C_{\one_{\cO}} \to \D_{\one_{\cO}}$ such that $F$ is small and $\cR$-local. Then the following are equivalent:
\begin{enumerate}[label=\alph*),ref=\alph*)]
\item
\label{item:approximation}%
$\eta$ exhibits $F$ as an $\cR$-local approximation of $\one_{\Day}$.
\item
\label{item:localising}%
For every $a \in \cO$ and every $\cR$-local functor $G\colon \C_a \to \D_a$ the map $G \to F\otimes_{\Day} G$ is an equivalence.
\item
\label{item:F-smashing}%
The map $\eta$ exhibits $F$ as $\cO$-idempotent and for every $a \in \cO$ the collection of $\cR$-local functors $\C_a \to \D_a$ coincides with the image of $\F\otimes_{\Day}(-)\colon \Fun(\C_a,\D_a) \to \Fun(\C_a,\D_a)$.
\end{enumerate}
In particular, when these equivalent conditions hold all functors are $\cR$-localisable and
the resulting $\cO$-multiplicative left Bousfield localisation 
\[
\Fun_{\cO}(\C,\D)^{\Day} \locadjL{L^{\otimes}} \Fun_{\cO}^{\cR}(\C,\D)
\]
is $\cO$-smashing in the sense of Definition~\reffour{definition:smashing} with associated idempotent $F$.
\end{enumerate}
\end{proposition}

Under the same assumptions, the $\cO$-smashing localisation of Proposition~\reffour{proposition:smashing-day} restricts to an $\cO$-smashing localisation of $\cO$-monoidal $\infty$-categories 
\[
\Fun_{\cO}^{\sm}(\C,\D)^{\Day} \locadjL{L^{\otimes}_{\sm}} \Fun_{\cO}^{\cR,\sm}(\C,\D)
\]
on the level of small functors and small $\cR$-local functors, respectively, see Corollary~\reffour{corollary:persist-locally-small}.

We record the above special case of Proposition~\reffour{proposition:smashing-day}, which is frequently encountered:

\begin{corollary}
\label{corollary:smashing-day-small}%
Let $\C^{\otimes}$ be a small symmetric monoidal $\infty$-category, $\D^{\otimes}$ a presentably symmetric monoidal $\infty$-category and $\cR = \{p_{\alp}\colon K_{\alp}^{\triangleleft} \to \C\}_{\alp}$ a collection of cone diagrams in $\C$. Assume that the following conditions hold:
\begin{enumerate}
\item
the symmetric monoidal structure on $\C$ is closed,
\item
\label{item:tensor-and-cotensor}%
the collection $\cR$ is closed under both tensoring and cotensoring by objects in $\C$,
\item
\label{item:closed-colimits}%
the collection of $\cR$-local functors $\C \to \D$ is closed under all small colimits, 
and
\item
\label{item:closed-tensor}%
for $d \in \D$, the functor $d \otimes (-)\colon \D \to \D$ preserves $K_{\alp}$-indexed limits for every $K_{\alp}$ appearing in $\R$.
\end{enumerate}
Then the left Bousfield localisation $L\colon \Fun(\C,\D) \to \Fun^{\cR}(\C,\D)$ refines to a 
smashing left Bousfield localisation
\[
\Fun(\C,\D)^{\Day} \locadjL{L^{\otimes}} \Fun^{\cR}(\C,\D)^{\Day}
\]
of presentably symmetric monoidal $\infty$-categories.
\end{corollary}

\begin{example}
\label{example:excisive}%
Let $\C$ be a small stable $\infty$-category equipped with a closed symmetric monoidal structure and $\D$ a stable presentably symmetric monoidal $\infty$-category. Then tensoring and cotensoring with any object in $\C$ is automatically an exact operation and hence preserves for example, the collection of strongly cocartesian $n$-cubes. It then follows from Corollary~\reffour{corollary:smashing-day-small} that the full subcategory $\Exc^n(\C,\D)$ spanned by $n$-excisive functors is a smashing localisation of $\Fun(\C,\D)$, and similarly for the full subcategory $\Exc^n_*(\C,\D)$ of reduced $n$-excisive functors. For example, the $\infty$-category $\Funq(\C)$ of hermitian structures on $\C$ is a smashing localisation of $\Fun(\C\op,\Sp)$, a case of interest in the context of the present paper series.

For completeness, let us note that it is enough to assume above that the target $\D$ is closed and cocomplete, since the existence of an explicit $n$-excisive approximation is well-known in this case and does not require presentability, see~\cite{HA}*{Theorem 6.1.1.10}. We also point out that $\Exc^n(\C,\D)$ and $\Exc^n_*(\C,\D)$ are typically not smashing localisations of $\Fun(\C,\D)$ when $\C$ is not stable, see Example~\reffour{example:cant-remove} below.
\end{example}

\begin{remark}
\label{remark:rigid}%
Let $\C^{\otimes}$ be small symmetric monoidal $\infty$-category, $\D^{\otimes}$ a presentably symmetric monoidal $\infty$-category and $\cR$ a collection of cone diagrams in $\C$ which is closed under tensoring with objects in $\C$, so that $\Fun^{\cR}(\C,\D)$ is a multiplicative localisation of $\Fun(\C,\D)$ by Remark~\reffour{remark:small-localisable}. Now, if, for example, every object in $\C$ is dualisable then the monoidal structure on $\C$ is closed and $\cR$ is closed also under cotensoring, and similarly if $\D$ is generated under colimits by dualisable objects then Condition~\reffouritem{item:closed-colimits} implies Condition~\reffouritem{item:closed-tensor} from Corollary~\reffour{corollary:smashing-day-small}. In particular, when these dualisability assumptions hold for both $\C$ and $\D$, the conclusion of Corollary~\reffour{corollary:smashing-day-small} holds as soon as $\Fun^{\cR}(\C,\D)$ is closed under colimits. This is consistent with Remark~\reffour{remark:smashing-closed-colimits}: Indeed, in this case $\Fun^{\cR}(\C,\D)$ is itself generated under colimits by dualisable objects, and so 
its smashing localisations are exactly its multiplicative localisations which are closed under colimits. %
\end{remark}

\begin{example}
\label{example:cant-remove}%
In the situation of Corollary~\reffour{corollary:smashing-day-small}, the condition that $\cR$ is closed under cotensoring cannot be removed. For example, let $k$ be a field
and take
$\C = \Dperf_{\geq 0}(k)$ to be the symmetric monoidal $\infty$-category of connective perfect $k$-complexes, 
$\cR$ the collection of pushout squares in $\C$ together with the singleton diagram on the zero object, and $\D=\Spa$. Note that $\C$ is closed with internal mapping objects given by $\uline{\map}_{\geq 0}(P,Q) = \tau_{\geq 0}\hom_{\Dperf(k)}(P,Q)$.
Then all the conditions of Corollary~\reffour{corollary:smashing-day-small} are met except~\reffouritem{item:tensor-and-cotensor}, namely, $\cR$ is closed under tensoring but not under cotensoring. For example, cotensoring againt $k[1]$ %
doesn't preserve the cofibre sequence $k \to 0 \to k[1]$.
We claim that the resulting multiplicative localisation 
\[
\Fun(\Dperf(k)_{\geq 0},\Spa) \locadj \Fun^{\cR}(\Dperf(k)_{\geq 0},\Spa) = \Exc^1_*(\Dperf(k)_{\geq 0},\Spa)
\]
is not smashing. Indeed, for every $P \in \Dperf(k)_{\geq 0}$ the left adjoint sends $P_!\SS$ to $\map_{\Dperf(k)}(P,-)$, and in particular sends the unit $k \in \Dperf(k)_{\geq 0}$ to the underlying spectrum functor 
\[
U\colon \Dperf(k)_{\geq 0} \to \Sp .
\]
It hence suffices to show that the natural transformation
\[
U \otimes_{\Day} P_!\SS = U(\uline{\map}_{\geq 0}(P,-)) \Rightarrow \map_{\Dperf(k)}(P,-)
\]
is generally not an equivalence. %
Indeed, taking for example $P=k[1]$ and evaluating this transformation at $k$, the right hand side is $\map_{\Dperf(k)}(k[1],k) = k[-1]$ while the left hand side is $\tau_{\geq 0}k[-1] = 0$. %
\end{example}

\begin{proof}[Proof of Proposition~\reffour{proposition:smashing-day}]
We begin with proving~\reffouritem{item:ideal}. Assume first that $F$ is $\cR$-local.
Since the collection of $\cR_a$-local functors $\C_a \to \D_a$ is assumed to be closed under (potentially large) colimits and Day convolution preserves all colimits that exist in each entry since $\D$ is almost closed, it will suffice by Lemma~\reffour{lemma:canonical-colimit} to prove the claim for $G  = c_!d$ for some $c \in \C_a$ and $d \in \D_a$. As described in the proof of Proposition~\reffour{proposition:persist-locally-small},
the Day convolution $F \otimes_{\Day} c_!d$ is then given by $F(\uline{\Hom}_{\C_a}(c,-)) \otimes d$, which is $\cR_a$-local since $F$ is $\cR_{\one_{\cO}}$-local and we assumed that $\uline{\Hom}_{\C_a}(c,-)\colon \C_a \to \C_{\one_{\cO}}$ sends diagrams in $\cR_a$ to diagrams in $\cR_{\one_{\cO}}$ and that $d \otimes(-)\colon \D_{\one_{\cO}} \to \D_a$ preserves $K_{\alp}$-indexed limits for every $K_{\alp}$ in $\cR_{\one_{\cO}}$. The argument if $G$ is $\cR$-local is similar, where one now considers the case where $F=c_!d$ for some $c \in \C_{\one_{\cO}}$ and $d \in \D_{\one_{\cO}}$.

We now proceed to prove~\reffouritem{item:equivalent}. To see~\reffouritem{item:approximation}$\Rightarrow$\reffouritem{item:localising} note that if $\eta$ exhibits $F$ as an $\cR$-local approximation of $\one_{\Day}$ then by Proposition~\reffour{proposition:multiplicative-day} the map $G \Rightarrow F \otimes_{\Day} G$ is an $\cR_a$-local equivalence for every $G\colon \C_a \to \D_a$, and hence an equivalence for every $\cR$-local such $G$. Next, if \reffouritem{item:localising} holds then the induced map $F \to F\otimes_{\Day} F$ is an equivalence since $F$ is $\cR$-local, and so $\eta$ exhibits $F$ is as $\cO$-idempotent. In addition, \reffouritem{item:localising} implies that the the image of $\F\otimes_{\Day}(-)\colon \Fun(\C_a,\D_a) \to \Fun(\C_a,\D_a)$ contains all $\cR$-local functors and the inverse inclusion follows from~\reffouritem{item:ideal}. We thus have \reffouritem{item:localising}$\Rightarrow$\reffouritem{item:F-smashing}. Finally the implication \reffouritem{item:F-smashing}$\Rightarrow$\reffouritem{item:approximation} follows from Remark~\reffour{remark:action-of-idempotent}.

The last part of~\reffouritem{item:equivalent} now follows from Proposition~\reffour{proposition:equivalent-1}.
\end{proof}

Let us finish with a few comments on the conditions on the target $\infty$-category in Proposition~\reffour{proposition:smashing-day}.

\begin{examples}
\label{someexampleofRlocalfunctors}%
Let $\cO^{\otimes},\C^{\otimes}$ and $\cR$ be as in Proposition~\reffour{proposition:smashing-day}.
\begin{enumerate}
\item If each diagram in each $\cR_a$ has source $(\Delta^0)^\triangleleft = \Delta^1$, then each $\cR_a$ is simply a collection of arrows in $\C_a$ and a functor on $\C_a$ %
is $\cR_a$-local if and only if it inverts all these maps. Then any $\cO$-monoidal $\D^{\otimes}$ satisfies assumptions \reffouritem{item:Rlocalclosedcolim} and \reffouritem{item:OmonoidalpresKlimits} of Proposition \reffour{proposition:smashing-day}.
\item If each diagram in each $\cR_a$ has source $\emptyset^\triangleleft = \Delta^0$, then $\cR_a$ is simply a collection of objects in $\C_a$ and a functor on $\C_a$ is $\cR_a$-local if and only it takes all these objects to the terminal object. %
Then any pointedly $\cO$-monoidal $\D$ satisfies assumptions \reffouritem{item:Rlocalclosedcolim} and \reffouritem{item:OmonoidalpresKlimits}. %
\item If all sources of morphisms in each $\cR_a$ are cones on finite sets, then any semi-additivly $\cO$-monoidal $\D$ satisfies assumptions \reffouritem{item:Rlocalclosedcolim} and \reffouritem{item:OmonoidalpresKlimits}. %
Similarly, if the sources are all finite $\infty$-categories then the same holds for any stably $\cO$-monoidal $\D^{\otimes}$.
\item As a slightly more subtle example, observe that if $\cR_a$ consists of square diagrams whose vertical maps admit sections, then any pre-stable $\D$ can serve as the target, since squares in $\D$ which admit vertical sections are pullbacks if and only if they are pushouts. 

A relevant example for our paper series is $\C^{\otimes} = (\Catx)^{\otimes} \to \Com^{\otimes}$, $\cR$ the collection of split Verdier squares and $\D = \Grp_{\Einf}$, in which case $\one_\Day^{\cR} = \Kspace$.
\end{enumerate}
\end{examples}

\section{Compact generation of Poincaré categories}
\label{section:presentability}%

Our goal in this appendix is to prove that $\Catp$ and $\Cath$ are presentable and compactly generated by an explicit pair of generators. The analogous statement for $\Catx$ was established in \cite{BGT}*{Corollary 4.25}.

\subsection{A criterion for compact generation}
\label{subsec:criterion}%

Our proof of compact generation relies on the following rather well-known criterion, analogous to \cite{MakkaiPitts}*{Lemma~1.7.ii} in the setting of ordinary categories, the proof of which we give for lack of reference. 
Recall that a set of objects $S\subseteq\C$ is said to be \emph{jointly conservative} if the functor $\{\hom_\C(\z,-)\}_{\z\in S}\colon \C\to \prod_{\z\in S}\Sps$ is conservative,
and that it generates the category $\C$ if the smallest full subcategory of $\C$ closed under colimits and containing the objects of $S$ is all of $\C$.

\begin{proposition}
\label{proposition:criterion-compact-generation}%
Let $\kappa$ be a regular cardinal, let $\C$ be a cocomplete locally small category and let $S\subseteq \C$ be a jointly conservative set of $\kappa$-compact objects. Then $S$ is a set of $\kappa$-compact generators. In particular $\C$ is $\kappa$-compactly generated. Furthermore, the set of $\kappa$-compact objects in $\C$ coincides with the smallest full subcategory of $\C$ containing $S$ and closed under $\kappa$-small colimits and retracts (if $\kappa$ is uncountable then retracts are automatic).
\end{proposition}
\begin{proof}
Let $S_{\kappa} \subseteq \C$ be the smallest full subcategory of $\C$ containing $S$ and closed under $\kappa$-small colimits. We claim that $S_{\kappa}$ is small. For this, observe that we can produce $S_{\kappa}$ as follows.
Let $S_0=S$ and, for all ordinals $\alpha\ge 0$, let $S_{\alpha+1}$ be the full subcategory spanned by representatives of $\kappa$-small colimits in $\C$ of diagrams in $S_\alpha$. Finally, if $\lambda$ is a limit ordinal let $S_\lambda$ be the union of $S_\alpha$ for all $\alpha<\lambda$. By transfinite induction, $S_\alpha$ has a set of objects for every $\alpha$ (since the $\kappa$-small diagrams in $S_{\alpha}$ form a set by the assumption that $\C$ is locally small). Then 
\[
S_\kappa = \bigcup_{\alpha<\kappa} S_\alpha\,.
\]
also still has only a set of objects. Moreover $S_\kappa$ is closed under $\kappa$-small colimits: every $\kappa$-small diagram in $S_\kappa$ factors through some $S_\alpha$ and so has a colimit in $S_{\alpha+1}$. 

Now, since $\kappa$-compact objects are closed under $\kappa$-small colimits, $S_\kappa$ consists of $\kappa$-compact objects, and is still jointly conservative since it contains $S$. 
Since $\C$ is cocomplete, and in particular admits $\kappa$-filtered colimits, we have an induced functor $\Ind_{\kappa}(S_{\kappa}) \to \C$. We claim that this functor is an equivalence. By~\cite{HTT}*{Proposition 5.3.5.11} it will suffice to show that $S_{\kappa}$ generates $\C$ under $\kappa$-filtered colimits.
Let us fix $\x\in\C$. We will show that the canonical arrow
\[
\colim_{\z\in S_{\kappa}\times_\C\C_{/\x}}\z\to \x
\]
is an equivalence, where we note that the indexing $\infty$-category is $\kappa$-filtered since it admits $\kappa$-small colimits. Since $S_{\kappa}$ is jointly conservative, it will suffice to show that for every $\y\in S_{\kappa}$ the arrow
\[
\Map_\C\left(\y,\colim_{\z\in S_{\kappa}\times_\C\C_{/\x}}\z\right)\to \Map_\C(\y,\x)
\]
is an equivalence. Now, since $S_{\kappa}$ has $\kappa$-small colimits, the category $S_{\kappa}\times_\C\C_{/\x}$ is $\kappa$-filtered, and since $\y$ is $\kappa$-compact, we need to show that
\[
\colim_{\z\in S_{\kappa}\times_\C\C_{/\x}}\Map_\C(\y,\z)\to \Map_\C(\y,\x)
\]
is an equivalence. Recall that $\Map_\C(\y,-)$ classifies the left fibration $\C_{\y/}\to \C$. Thanks to \cite{HTT}*{Corollary 3.3.4.3}, this is equivalent to showing that the ``composition'' map
\[
\C_{\y/}\times_\C S_{\kappa}\times_\C\C_{/\x}\to \{\y\}\times_\C \C_{/\x}
\]
is a weak equivalence. Since $y$ is an object of $S_{\kappa}$, we have an equivalence 
\[
(S_{\kappa})_{\y/} \times_\C \C_{/\x} \simeq \C_{\y/}\times_\C S_{\kappa}\times_\C\C_{/\x}
\]
and the object $\id_{\y}$ in $(S_{\kappa})_{\y/}$ induces a functor 
\[
\{\y\} \times_\C \C_{/\x} \to (S_{\kappa})_{\y/} \times_\C \C_{/\x}
\]
which is a left adjoint (and a section) of the above composition map. In particular, it is a weak equivalence as needed. Finally, by~\cite{HTT}*{Lemma 5.4.2.4} every $\kappa$-compact object of $\C$ is a retract of an object in $S_{\kappa}$.
\end{proof}

\begin{remark}
\label{remark:sufficiently-large}%
Proposition~\reffive{proposition:criterion-compact-generation} implies that every $\kappa$-presentable $\infty$-category is also $\kappa'$-presentable for every $\kappa'\ge \kappa$, since it implies that a set of $\kappa$-compact generators is also a set of $\kappa'$-compact generators. 
This is not true for general $\kappa$-accessible $\infty$-categories (cf.\ \cite{Adamek-Rosicky}*{Example~2.13.(8)}). 

We note also that by the last part of Proposition~\reffive{proposition:criterion-compact-generation} if $\C$ is $\kappa$-presentable then for every $\kappa' > \kappa$ the collection of $\kappa'$-compact objects is generated under $\kappa'$-small colimits by the $\kappa$-compact ones. %
This also means that if $\C$ and $\D$ are two presentable $\infty$-categories and $f\colon \C \to \D$ is an accessible functor then there exists a $\kappa_0$ such that for all $\kappa \geq \kappa_0$ we have that $\C$ and $\D$ are $\kappa$-presentable and $f$ is both $\kappa$-accessible and preserves $\kappa$-compact objects. %
\end{remark}

\subsection{Compact generation for Poincaré and hermitian $\infty$-categories}
\label{subsec:compact-gen}%

We now come to the main result of this appendix:
\begin{proposition}
\label{proposition:catp-cath-are-compact}%
The $\infty$-categories $\Cath$ and $\Catp$ are presentable and compactly generated. In fact,
the $\infty$-category $\Cath$ is generated under colimits by the two compact hermitian $\infty$-categories
\[
(\Spaf,0) \quad\text{and}\quad(\Spaf,\QF^{\uni}),
\]
and the $\infty$-category $\Catp$ is generated under colimits by the two compact Poincaré $\infty$-categories
\[
\Hyp(\Spaf)  \quad \text{and} \quad \Met(\Spaf,\QF^{\uni}).
\]
\end{proposition}

\begin{corollary}
Finite limits in $\Catp$ and $\Cath$ commute with filtered colimits. 
\end{corollary}

\begin{remark}
The full subcategory $\Catpi \subseteq \Catp$ spanned by the idempotent complete Poincaré $\infty$-categories is reflective and closed under filtered colimits. It then follows that $\Catpi$ is compactly generated as well. The same holds for the analogous full subcategory of $\Cath$.
\end{remark}

\begin{corollary}
The collections of Poincaré-Verdier sequences, squares, inclusions and projections are all closed under filtered colimits (the latter actually under all colimits), and the same goes for the Poincaré-Karoubi variants of these notions.
\end{corollary}

The proof of Proposition~\reffive{proposition:catp-cath-are-compact} will require the following lemma:

\begin{lemma}
\label{lemma:hyp-U-preserve-compact}%
In the sequence of forgetful functors
\[
\Catp \to \Cath \to \Catx \to \Cat
\]
the first two preserve all colimits and even admit colimit preserving right adjoints, 
while the last one just preserves filtered colimits. In particular, the first two forgetful functors, 
as well as the left adjoints
\[
\begin{tikzcd}
\Cat \ar[r,"{\Spaf_{(-)}}"] & \Catx \ar[rr,bend right, "\Hyp"] \ar[r,"{(-,0)}"] & \Cath \ar[r,"\Pairings(-\qshift{1})"] & \Catp \ ,
\end{tikzcd}
\]
all preserve $\kappa$-compact objects for every regular cardinal $\kappa$. Here the left most left adjoint is the functor $\I \mapsto \Spaf_\I$, sending $\I$ to the stable $\infty$-category obtained by tensoring $\Spaf$ over $\I$, see Construction~\refone{construction:tensor} and Remark~\refone{remark:universal-underlying}. 
\end{lemma}
\begin{proof}
The statement that $\Catx \to \Cat$ preserves filtered colimits is proven in~\cite{HA}*{Proposition 1.1.4.6}. 
Concerning $p\colon \Cath \to \Catx$, since it is a cocartesian fibration its right adjoint is given by its objectwise terminal section, that is, by the functor sending $\C \in \Catx$ to $(\C,0)$. This is also the left adjoint, since $p$ is also a cartesian fibration and $0$ is also initial in each fibre. We conclude that the right adjoint of $p$ preserves colimits. The right adjoint of $q\colon \Catp \to \Cath$ is described in Proposition~\refone{proposition:first-adj-hermitian}. It is also very close to the left adjoint of $q$ described in Proposition~\refone{proposition:second-adj-hermitian}: the two adjoints differ by post-composing with the shift equivalence $\qshift{1}\colon \Catp \xrightarrow{\simeq} \Catp$. It hence follows that the right adjoint of $q$ preserves colimits. Finally, the conclusion for preservation of $\kappa$-compact objects is because any functor which preserves filtered colimits also preserves $\kappa$-filtered colimits for every $\kappa$, and hence its left adjoint preserves $\kappa$-compact objects for every $\kappa$.
\end{proof}

\begin{proof}[Proof of Proposition~\reffive{proposition:catp-cath-are-compact}]
Since $\Cath$ and $\Catp$ are cocomplete by Propositions~\refone{proposition:cath-has-limits-and-colimits} and \refone{proposition:Catp-cocomplete},
it suffices to establish the criterion of the Proposition~\reffive{proposition:criterion-compact-generation}. 
Consulting Proposition~\refone{proposition:second-adj-hermitian}, Example~\refone{examples:typical-pairings}~\refoneitem{item:zero-pairing-hyp} and Remark~\refone{remark:second-counit}, we observe that the proposed generators for $\Catp$ are exactly the images of the proposed generators for $\Cath$ under the left adjoint 
$\Pairings(-\qshift{1})\colon \Cath \to \Catp$. 
Since this left adjoint preserves compact objects by Lemma~\reffive{lemma:hyp-U-preserve-compact} and the forgetful functor $\Catp \to \Cath$ is conservative it will actually suffice to verify the criterion of Proposition~\reffive{proposition:criterion-compact-generation} in the case of $\Cath$.

Now $(\Spaf,0)$ is the image under the composed left adjoint
\[
\Cat \to \Catx \to \Cath
\]
of the object $\Del^0 \in \Cat$. In particular, $(\Spaf,0)$ corepresents the functor $(\C,\QF) \mapsto \core(\C)$. Since $\Del^0$ is compact in $\Cat$ we conclude from Lemma~\reffive{lemma:hyp-U-preserve-compact} that $(\Spaf,0)$ is compact. Similarly, by Proposition~\refone{proposition:corepresentability-of-poinc} $(\Spaf,\QF^{\uni})$ corepresents (in $\Cath$) the functor $(\C,\QF) \mapsto \spsforms(\C,\QF)$ which preserves filtered colimits by Proposition \refone{proposition:forms-poinc-compact}, so that $(\Spaf,\QF^{\uni})$ is also compact in $\Cath$.

We will now show that the proposed generators are jointly conservative. 
Suppose that $(f,\eta)\colon(\C,\QF)\to (\D,\QFD)$ is a hermitian functor inducing equivalences 
\[
\core\C \rightarrow \core\D \quad \text{and} \quad \spsforms(\C,\QF)\rightarrow \spsforms(\D,\QFD).
\]
Since the map on groupoid cores is an equivalence, $f$ is essentially surjective. Moreover, for all $\x,\y\in\C$, the space $\hom_\C(\x,\y)$ retract-embeds in $\hom_{\core\C}(x\oplus y,x \oplus y)$ as unipotent upper triangular matrices, so $f$ is also fully faithful and thus an equivalence.

To finish the proof, it is left to show that the map
\[
\eta_{\x} \colon \QF(\x)\to \QFD(f(\x))
\]
is an equivalence for every $\x\in \C$. By Lemma~\refone{lemma:coconnective-quadratic-functors},
it suffices to show this after applying $\Omega^\infty$. %
But this map is the one induced on vertical fibres (over the point $\x$) by the commutative square
\[
\begin{tikzcd}
	\spsforms(\C,\QF)\ar[r]\ar[d] & \spsforms(\D,\QFD)\ar[d]\\
 	\core\C\ar[r] & \core\D,
 \end{tikzcd}
\]
which is indeed an equivalence since the horizontal arrows are so by assumption.
\end{proof}

\begin{corollary}
\label{corollary:omega-monoidal}%
The symmetric monoidal structures on $\Catp$ and $\Cath$ are $\omega$-presentable, that is, the unit is compact and the collection of compact objects is closed under tensor products. In fact, for any $\kappa \geq \omega$ the collection of $\kappa$-compact objects in these $\infty$-categories is closed under tensor products.
\end{corollary}
\begin{proof}
We begin with $\Cath$. The compactness of the unit $(\Spaf,\QF^{\uni})$ is part of Proposition~\reffive{proposition:catp-cath-are-compact}.
In addition, the same proposition implies that $(\Spaf,0)$ and $(\Spaf,\QF^{\uni})$
generate all $\kappa$-compact objects under $\kappa$-small colimits and retracts. It will hence suffice to show that the tensor product of any two of $(\Spaf,0)$ and $(\Spaf,\QF^\uni)$ is compact, which follows from the fact that the latter one is the unit and the former one is idempotent.

For $\Catp$, the compactness of the unit is established in Proposition~\refone{proposition:forms-poinc-compact}, and we are again left with showing that the tensor product of any two of $\Hyp(\Spaf)$ and $\Met(\Spaf,\QF^{\uni})$ is again compact in $\Catp$. Since these are obtained from $(\Spaf,0)$ and $(\Spaf,\QF^\uni)$ by applying the left adjoint to the forgetful functor $\Catp \to \Cath$, and since both the forgetful functor and its left adjoint preserve compact objects by Lemma~\reffive{lemma:hyp-U-preserve-compact}, the desired result follows from the projection formula enjoyed by this adjunction, see Lemma~\refone{lemma:pairing-projection}.
\end{proof}

\section{Another description of the visible signature}
\label{section:appendix-visible-revisited}%

Recall that a Poincaré complex has a visible signature, originally due to Weiss-Williams, and described in the language of Poincaré $\infty$-categories in \S\refone{subsection:parametrised-spectra}. The purpose of this appendix is to give a second perspective on this visible signature. 
knowledge of either or the two mentioned sources, except for the very last part where we compare it to the construction from loc.\ cit.

Let us give an outline of the approach. Recall that the $\infty$-category of finite parametrised spectra over a space $X$, in which the underlying object of the visible signature lies, models the colimit of the constant $X$-indexed diagram with value $\Spaf$. This is an instance of a general phenomenon: Given a diagram 
\[
F\colon X \to \Catx
\]
of stable $\infty$-categories, the collection of Ind-right adjoints $i_x^R\colon \colim_X F \to \Ind(F(x))$ of the colimit structure maps $i_x\colon F(x) \to \colim_X F$ induce a fully-faithful embedding $e_{\Ind}\colon \colim_X F \to \lim_X \Ind\circ F$. Now when $X$ is finite the essential image of $e_{\Ind}$ contains $\lim_X F$. 
The resulting functor $\lim_X F \to \colim_X F$ is an instance of a rather general limit to colimit norm map, which we call the cocartesian conorm (although we refrain here from describing it in the most general setting). In the case of parametrised spectra, this corresponds to the inclusion of pointwise finite parametrised spectra inside of finite parametrised spectra. Dually, instead of taking the Ind-right adjoints $i_x^R$ one can take the Pro-left adjoints $i_x^L\colon \colim_X F \to \Pro(F(x))$, yielding an embedding $e_{\Pro}\colon \colim_X F \to \lim_X \Pro\circ F$, which again contains the image of $\lim_X F$ when $X$ is finite. This gives a dual variant $\lim_X F \to \colim_X F$ of the conorm map, which we call the cartesian conorm. Returning to the case where $F$ is constant with value $\Spaf$, it turns out that when $X$ is a Poincaré complex the cartesian and cocartesian conorm maps differ by tensoring with the dualising complex $\nu\colon X \to \Pic(\Spaf)$, and in particular share the same image. This common image is then invariant 
under the duality associated to the Poincaré structure of $\colim_{x \in X} (\Spaf,\QF^{\uni}\otimes \nu_x)$ (which is where the visible signature lives as a Poincaré object). In addition, restricted along the cartesian norm map, the resulting Poincaré structure on $\lim_X \Spaf$ identifies it with the limit $\lim_X (\Spaf,\QF^{\uni})$ of the constant diagram with value $(\Spaf,\QF^{\uni})$. In particular, this full Poincaré subcategory contains a canonical Poincaré object, corresponding to the constant diagram with value to the universal Poincaré object $(\SS,q_{\uni})$. This is our alternative construction of the visible signature.

\subsection{The (co)cartesian conorm}

Let $X$ be a finite space and $F\colon X\to \Catx$ a functor. The cocartesian conorm is an exact functor
\[
\ccconorm\colon \lim_X F \to \colim_X F
\]
defined as follows: The natural transformations of $X$-indexed diagrams %
\[
\lim_X F \xrightarrow{p_x} F(x) \xrightarrow{i_x} \colim_X F
\]
compose to a natural transformation between constant diagrams, which we may interpret as a functor 
\[
X \to \Hom_{\Catx}(\lim_X F, \colim_X F) \subset \Funx(\lim_X F, \colim_X F).
\]
The cocartesian conorm map is then, by definition, obtained by taking the colimit of this $X$-indexed family.

We emphasise that the colimit $\colim_X F$ is to be taken in $\Catx$, which (unlike the limit) does not agree with the colimit in $\Cat$. An explicit description of the colimit may be obtain as follows. For each $x \in X$ the associated exact functor $i_x\colon F(x) \to \colim_X F$ admits an Ind-right adjoint $i_x^R\colon \colim_X F \to \Ind(F(x))$, which together determine an exact functor
\[
e_{\Ind}\colon \colim_X F \to \lim_{X\op} \Ind\circ F = \lim_X\Ind \circ F ,
\]
where we have used the canonical equivalence $X \simeq X\op$.
The functor $e_{\Ind}$ is fully-faithful and exhibits its target as the Ind-completion of its source. To see this, note that $\Ind$ is colimit preserving when viewed as a functor from $\Catx$ to $\Pr^L_{\mathrm{st}}$ (as follows from the universal mapping property $\Fun^L(\Ind(\C),\D) = \Funx(\C,\D)$, see~\cite{HTT}*{Proposition 5.3.6.2, Example 5.3.6.8}). It follows that the collection of functors $\Ind(i_x)$ exhibit $\Ind(\colim_X F)$ as the colimit of $\Ind\circ F$ in $\Pr^L_{\mathrm{st}}$, and hence the collection of their right adjoints $i_x^R$ exhibit $\Ind(\colim_X F)$ as the limit of $\Ind \circ F$, but viewed as a diagram $X=X\op \to \Pr^R_{\mathrm{st}}$. This limit is then calculated in large categories %
(see~\cite{HTT}*{Theorem 5.5.3.18}), or in large stable ones for that matter,
and so the claim follows. This leads to a concrete description of $\colim_X F$ as the stable subcategory of 
$\lim_{X} \Ind\circ F$
spanned by the images of the composites 
\[
F(x) \to \Ind F(x) \xrightarrow{x_!} \lim_{X} \Ind\circ F
\]
where $x_!$ is left adjoint to the projection $p_x\colon \lim_{X} \Ind\circ F \to \Ind F(x)$ from the limit. %
Under this description we can describe the cocartesian conorm as follows:

\begin{lemma}
\label{lem:cocartesian_norm}%
For any finite space $X$ and $F\colon X\to \Catx$ the cocartesian conorm and the functor $e_{\Ind}$ above fit into a commutative triangle 
 \begin{equation}
\begin{tikzcd}
 & \lim_X F \ar[ld, "\ccconorm"'] \ar[d, "\mathrm{incl}"] \\
 \colim_X F \ar[r, "e_{\Ind}"'] & \lim_X \Ind\circ F.
\end{tikzcd}
\end{equation}
with vertical map the canonical inclusion, and all (co)limits computed in (possibly large) stable $\infty$-categories. In particular, the cocartesian conorm is fully faithful.
\end{lemma}
\begin{proof}
By abuse of notation let us write $p_x$ for both projections $\lim_X F \to F(x)$ and $\lim_X \Ind\circ F \to \Ind(F(x))$ and let $x_!$ be the left adjoint of the latter as above. %
We then calculate
\[
e_{\Ind}(\ccconorm(\vphi)) = e_{\Ind}(\colim_{x \in X}i_xp_x(\vphi)) = \colim_{x \in X}e_{\Ind}i_xp_x(\vphi) = \colim_{x \in X} x_!p_x(\vphi) \xrightarrow{\simeq} \vphi ,
\]
where the last equivalence, induced by the collection of counit maps $x_!p_x(\vphi) \to \vphi$, holds more generally for any $\vphi \in \lim_X \Ind\circ F$, as can be verified by mapping into an arbitrary test object $\psi \in \lim_X \Ind\circ F$.
\end{proof}

The \defi{cartesian conorm} 
\[
\cconorm\colon \lim_X F\to \colim_X F
\]
is the opposite of the cocartesian conorm map applied to $\mathrm{op}\circ F$, after identifying 
\[
\lim_X \mathrm{op}\circ F \simeq (\lim_X F)\op \quad \mathrm{and} \quad \colim_X \mathrm{op}\circ F \simeq (\colim_X F)\op.
\]
Equivalently, it is given by replacing in the formula the colimit by a limit, that is, $\cconorm(\vphi) = \lim_x i_xp_x(\vphi)$. %
A dual version of Lemma~\reffour{lem:cocartesian_norm} can then be obtained by replacing $e_{\Ind}$ with the dual construction
\[
e_{\Pro}\colon \colim_X F \to \lim_{X} \Pro\circ F
\]
induced by the $X$-indexed family of Pro-left adjoints $i_x^L\colon \colim_X F \to \Pro(F(x))$.
The functor $e_{\Pro}$ is then fully faithful and exhibits the target as the $\Pro$-completion of the source. As in Lemma \reffour{lem:cocartesian_norm} above, the cartesian conorm and the embedding $e_{\Pro}$ then fit in a commutative triangle
\begin{equation}
\label{eq:cartesian_conorm_and_e_Pro}%
\begin{tikzcd}
 & \lim_X F \ar[ld, "\cconorm"'] \ar[d, "\mathrm{incl}"] \\
 \colim_X F \ar[r, "e_{\Pro}"'] & \lim_X \Pro\circ F.
\end{tikzcd}
\end{equation}
with vertical map the canonical inclusion. In particular, the cartesian conorm is fully faithful.

The following result explains how to interpret the cartesian conorm in the (arguably more common) embedding $e_{\Ind}$. The formula for this uses the ``fibrewise'' tensoring of $\lim_X G$, for $G\colon X\to \PrL_{\mathrm{st}}$, over $\Spa^X$, the category of parametrised spectra over $X$, 
\[
\otimes \colon (\lim_X G) \otimes \Spa^X \to \lim_X G
\]
defined through the $X$-family of exact functors
\[
(\lim_X G)\otimes \Spa^X \xrightarrow{x^*\otimes x^*} G(x)\otimes \Spa \to G(x)
\]
where the last map is the usual tensoring of the presentable stable $\infty$-category $G(x)$ over spectra.

For what follows, recall that the dualising complex $\nu_X$ of a space $X$ is the universal parametrised spectrum equipped with a natural transformation 
\[
T\colon \colim_{x\in X} \xi\otimes \nu_X \to \lim_{x\in X} \xi
\]
for $\xi\in \Spa^X$. This transformation is %
an equivalence whenever $X$ is compact (in particular, finite). 

\begin{lemma}
\label{lem:e_Ind_of_cartesian_conorm}%
For any finite space $X$ and $F\colon X\to \Catx$ the %
composite functor 
\[
\lim_X F \xrightarrow{\cconorm} \colim_X F \xrightarrow{e_{\Ind}} \lim_X \Ind\circ F
\]
is canonically equivalent to tensoring with the dualising complex of $X$.
\end{lemma}

\begin{remark}
Applying this lemma to the constant functor $F=\const_{\Spa^f}\colon X\to \Cat$, and the object $\const_\SS$ in $\lim_X F = (\Spa^f)^X$,  we find that the dualising complex $\nu_X$ is explicitly given by the formula 
\[
\nu_X \simeq \lim_{x\in X} x_! \SS
\]
which corresponds to the definition of the dualising complex originally given by Klein in \cite{klein-dualizing}.
\end{remark}

For the proof of Lemma~\reffour{lem:e_Ind_of_cartesian_conorm}, we first observe that the transformation $T\colon \colim_{x\in X} \xi\otimes \nu_X \to \lim_{x\in X} \xi
$ %
extends to a similar transformation valid for arbitrary functors $\xi\colon X\to \C$ taking values in compactly generated stable $\infty$-categories (which is again an equivalence whenever $X$ is compact), which can be seen by identifying $\C\simeq \Ind(\C^\omega) \simeq \Funx((\C^\omega)\op, \Spa)$, and using that limits and colimits in $\Funx((\C^\omega)\op, \Spa)$ are computed pointwise.

\begin{proof}[Proof of Lemma \reffour{lem:e_Ind_of_cartesian_conorm}]
The composite functor $e_{\Ind} \circ \cconorm$ is given by 
\[
e_{\Ind}\big(\lim_{x \in X}i_xp_x(\vphi)\big) = \lim_{x\in X} x_!(p_x(\varphi))= \colim_{x\in X} (x_!(p_x(\varphi)) \otimes p_x(\nu_X)) = \colim_{x\in X} x_! p_x(\varphi\otimes \nu_X) = \varphi\otimes \nu_X
\]
where we used the equivalence $T$ in the compactly generated stable $\infty$-category $\lim_X \Ind\circ F$.
\end{proof}

\medskip

\subsection{The visible signature via the conorm map}

We now consider the analogue of the conorm map construction for diagrams of Poincaré $\infty$-categories. 
Suppose given a diagram $F\colon X \to \Catp$ where $X$ is a finite space. As above we may consider the resulting $X$-indexed family of Poincaré functors
\[
i_xp_x\colon \lim_X F \to \colim_X F .
\]
We then encounter the following immediate obstacle: neither of the exact functors $\cconorm = \lim_x i_xp_x$ and $\ccconorm = \colim_xi_xp_x$ is duality preserving; in fact, since each of the individual $i_xp_x$ is duality preserving we have that $\cconorm$ and $\ccconorm$ are the duals of one another in the Poincaré $\infty$-category $\Funx(\lim_X F,\colim_X F)$.
To solve this point, we first restrict to the case where the finite space is a \emph{Poincaré complex}, that is, its dualising complex $\nu_X$ takes values in invertible spectra. This in itself does not remove the above obstacle. However, it allows us to consider a twisted version of the conorm map, in the form a Poincaré functor
\[
\lim_X F \to \colim_X(F\otimes \nu_X) = \colim_{x \in X} F(x) \otimes \nu_x
\]
where for a hermitian $\infty$-category $(\C,\QF)$ and a spectrum $E$ we write $(\C,\QF)\otimes E := (\C,\QF\otimes E)$, which is Poincaré if $(\C,\QF)$ is so and $E$ is invertible. The exact functors $i_xp_x$ are then not duality preserving, but rather become $\nu_x\otimes i_xp_x$ upon conjugating with the dualities on both sides. Conjugating $\lim_x i_xp_x$ by the duality then gives $\colim_x \nu_x\otimes i_xp_x \simeq \lim_x i_xp_x$, so that the cartesian conorm %
has a chance of being promoted to a Poincaré functor. This is indeed the case:

\begin{proposition}
\label{proposition:hermitian_conorm}%
Let $X$ be a finite space with dualising complex $\nu_X$, and $F\colon X\to \Cath$ be a functor. 
The \emph{cartesian} conorm naturally upgrades to a hermitian functor 
 
\[
\cconorm\colon \lim_X F \to \colim_X (F\otimes \nu_X)
\]
which is in fact a hermitian embedding (i.e., identifies the source with a hermitian subcategory of the target).   
Furthermore, if $X$ is a Poincar\'e complex and $F$ takes values in $\Catp$, then this hermitian functor is a Poincar\'e functor between Poincar\'e categories.
\end{proposition}

Before giving a proof of this lemma, we employ it for the promised description of the visible signature. For a Poincar\'e complex $X$ with dualising complex $\nu_X$, we consider the constant functor $X\to \Catp$ with value $(\C, \QF)$. Its limit and colimit are the cotensor $(\C, \QF)^X$ and $(\C, \QF)_X$ respectively. We also write 
\[
(\C, \QF)_{(X, \nu_X)}:= \colim_X ((\C, \QF)\otimes \nu_X).
\]
According to Proposition~\reffour{proposition:hermitian_conorm}, the cartesian conorm thus refines to a Poincar\'e functor 
\[
\cconorm\colon (\C, \QF)^X \to (\C, \QF)_{(X, \nu_X)}.
\]
Applying $\Poinc(-)$ gives us a map of spaces
\[
\cconorm \colon \Hom_{\Sps}(X, \Poinc(\C, \QF)) \simeq \Poinc((\C, \QF)^X) \to \Poinc((\C, \QF)_{(X, \nu_X)}).
\]

\begin{definition}
\label{definition:visible-signature}%
The \defi{visible signature} of a Poincar\'e complex $X$ is the element 
\[
\vsign(X) \in \Poinc((\Spa^f, \QF^{\uni})_{(X, \nu_X)})
\]
obtained by applying this construction to the constant map $X\to \Poinc(\Spa^f, \QF^{\uni})$ with value the universal Poincar\'e object $\SS^u\in \Poinc(\Spa^f, \QF^{\uni})$. 
\end{definition}

We come to the proof of Proposition~\reffour{proposition:hermitian_conorm}. We first recall the $\Pro$-completion of a hermitian category $(\C, \QF)$, denoted $\Pro(\C, \QF)$: This is the stable $\infty$-category $\Pro(\C)$ equipped with the canonical colimit-preserving extension $\hat\QF\colon \Pro(\C)\op\simeq \Ind(\C\op)\to \Spa$, i.e., the hermitian structure left Kan extended from the inclusion $i\colon \C\to \Pro(\C)$. This canonically upgrades $i$ to a hermitian functor (in fact, hermitian embedding), and evidently the whole construction can be made functorial in $(\C, \QF)\in\Cath$. 

\begin{lemma}
\label{lemma:e_Pro_hermitian}%
Let $X$ be any space with dualising complex $\nu_X$, and $F\colon X\to \Cath$. The embedding $e_{\Pro}$ naturally refines to a hermitian functor
\[
e_{\Pro}\colon \colim_X (F\otimes \nu_X) \to \lim_X \Pro\circ F.
\]
If $X$ is compact, this hermitian functor exhibits the target as the $\Pro$-completion of the source hermitian category, and so is, in particular, a hermitian embedding.
\end{lemma}

\begin{proof}
 The hermitian structure on the colimit $\colim_X (F\otimes \nu_X)$ is formed, according to \S\refone{subsection:limits}, by left Kan extending each constituent $\QF_x(-)\otimes \nu_X(x)$ along the opposite of  the inclusion $i_x\colon F(x)\to \colim_X F$, and taking the colimit over these. Thus, its further left Kan extension along the $\Pro$-completion $e_{\Pro}$ (of stable $\infty$-categories) is given by the formula 
\[
(\lim_X \Pro\circ F)\op\to \Spa, \quad \varphi\mapsto \colim_{x\in X} \hat \QF_x(x^* \varphi)\otimes \nu_X(x)
\]
(since the left Kan extension along (the opposite of) $e_{\Pro}\circ i_x = x_*$ is given by pullback along its left adjoint $x^*$). 
On the other hand, again by \S\refone{subsection:limits}, the hermitian structure on $\lim_X \Pro\circ F$ is given by 
\[
(\lim_X \Pro\circ F)\op\to \Spa, \quad \varphi\mapsto \lim_{x\in X} \hat  \QF_x(x^* \varphi).
\]
The claimed refinement is given by the the transformation 
\[
\colim_{x\in X} \hat  \QF(x^*\varphi) \otimes \nu_X(x) \to \lim_{x\in X} \hat  \QF(x^*\varphi)
\]
issued from the defining property of the dualising complex. If $X$ is compact, this transformation is an equivalence. 
\end{proof}

\begin{proof}[Proof of Proposition~\reffour{proposition:hermitian_conorm}]
We consider the diagram in $\Cath$ 
 \begin{equation}
\label{eq:cocartesian_conorm_and_e_Ind}%
\begin{tikzcd}
 & \lim_X F \ar[ld, "\cconorm"', dotted] \ar[d, "\mathrm{incl}"] \\
 \colim_X (F\otimes \nu_X) \ar[r, "e_{\Pro}"'] & \lim_X \Pro\circ F.
\end{tikzcd}
\end{equation}
where the vertical map is the limit of the canonical hermitian embeddings (which is again a hermitian embedding) and the lower line is issued from Lemma~\reffour{lemma:e_Pro_hermitian}. We obtain a unique dotted lift from~\eqreffour{eq:cartesian_conorm_and_e_Pro},
which is again a hermitian embedding. 

It will suffice to show that for any $\varphi\in \lim_X F$, the spectral presheaf
\[
\Bil_{\colim_X (F\otimes \nu_X)}(-, \cconorm(\varphi))
\]
is represented by some object of the form $\cconorm(\psi)$ (where $\psi$ will then necessarily be the dual of $\varphi$ in $\lim_X F$). By the hermitian embedding $e_{\Pro}$ from Lemma~\reffour{lemma:e_Pro_hermitian}, this spectral presheaf is restricted from the spectral presheaf 
\[
\Bil_{\lim_X \Pro\circ F}(-, e_{\Pro}\cconorm \varphi) \simeq \Bil_{\lim_X \Pro\circ F}(-, \varphi) \simeq \lim_x \Bil_{\Pro F(x)}(p_x(-),p_x(\vphi)),
\]
and the latter is represented in $\lim_X \Pro\circ F$ by the pointwise dual $\Dual_{\lim F} \vphi$.
\end{proof}

\medskip

\begin{remark}
We mention the following two variants of this construction:
\begin{enumerate}
\item A version of the (co-)cartesian conorm still exists for spaces $X$ that are merely compact %
if the target category $\colim_X (F\otimes \nu_X)$ is replaced by its idempotent completion (that is, we form the colimit in $\Catxi$). Indeed, any idempotent complete stable $\infty$-category admits limit indexed over compact spaces, and the rest of the proof works verbatim. Consequently, a compact Poincar\'e complex admits a visible signature 
\[
\vsign(X)\in \Poinc(((\Spa^f, \QF^{\uni})_{(X, \nu_X)})^\natural).
\]
Under $e_{\Ind}$, this corresponds to the category of compact parametrised spectra (as opposed to finite ones). 
 \item For a finite space $X$, our argument really shows that the cartesian conorm is a Poincar\'e functor provided each $F(x)$ and each $F(x)\otimes \nu_X(x)$ are Poincar\'e categories.
For instance, if $F$ is constant with value $(\Mod^\omega(R), \QF^\tate)$
for an $\Einf$-ring spectrum $R$, this amounts to saying that each $\nu_X(x)\otimes R$ is invertible in $\Mod{R}$. We deduce that for such an ``$R$-Poincar\'e complex'', we can construct an ``$R$-visible signature''
\[
\vsign(X; R)\in \Poinc((\Mod^\omega(R), \QF^\tate)_{(X, \nu_X)})
\]
by applying the above construction to the constant map $X\to \Poinc(\Mod^\omega(R), \QF^\tate)$ with value the tautological Poincar\'e object $R$. 
\end{enumerate}
\end{remark}

\medskip

We conclude with a comparison of the present approach with the one pursued in \S\refone{subsection:parametrised-spectra}. In loc.\ cit.\ we denote by $\Spaf_X \subseteq \Sp^X$ the full subcategory spanned by the finite parametrised spectra. %
From the present point of view, it is exactly the essential image under $e_{\Ind}$ %
of the colimit of the constant $X$-diagram on $\Spaf$. 
Note the slight discrepancy in conventions: for a space $X$ we generally denote by $\C_X$ the abstract colimit of the constant diagram on $\C$, while in the case of $\Spaf_X$ we implicitly implement the embedding $e_{\Ind}$ and view it as a full subcategory of $\Sp^X$.  
For given spherical fibration $\xi\in \Pic(\SS)^X$, we then defined the \emph{$\xi$-twisted visible hermitian structure} $\QF^\vis_\xi$ on $\Spaf_X$, with linear and bilinear parts given as
\[
 \Lin^\vis_{\xi} := \Hom_{\Spa^X} (-, \xi) \quad \mathrm{and} \quad \Bil^\vis_{\xi} := \Hom_{\Spa^{X\times X}}(-\boxtimes -, \Delta_!(\xi)) 
\]
and obvious connecting transformation $\Lin^\vis_{\xi}\to (\Bil^\vis_{\xi})^\tC$. We write $
\Dual_{\xi}$ for the duality of $\QF^{\vis}_{\xi}$.
This construction is covariantly functorial in $(X,\xi)$ via left Kan extension of parametrised spectra, %
and we show in Proposition~\refone{proposition:visible-is-colimit-full} that the collection of Poincaré functors $x_!\colon (\Spaf,\QF^{\uni})\otimes \xi_x = (\Spaf_{\{x\}},\QF^{\uni}_{\xi_x}) \to (\Spaf,\QF^{\vis}_\xi)$ induces an equivalence
\[
\colim_{x \in X} (\Spaf,\QF^{\uni})\otimes \xi_x \xrightarrow{\simeq} (\Spaf,\QF^{\vis}_\xi)
\]
in $\Catp$.

Now suppose that $X$ is a Poincaré complex with dualising complex $\nu_X$. Note that the Ind-right adjoint of $x_!\colon \Spaf \to \Spaf_X$ is given by evaluation at $x$, and since $x_!$ is duality preserving, for instance, for the dualities associated to the trivial spherical fibrations $\SS_X$ and $\SS$ on $X$ and $*$, respectively, the Pro-left adjoint of $x_!$ is given by $\vphi \mapsto \hat{\Dual}_{\SS}x^*(\Dual_{\SS_X}\vphi)$, where $\hat{\Dual}_\SS\colon \Ind(\Spaf)\op = \Pro((\Spaf)\op) \to \Pro(\Spaf)$ is the equivalence induced by $\Dual_{\SS}\colon (\Spaf)\op \to \Spaf$. Proposition~\reffour{proposition:hermitian_conorm} then implies that the collection of Pro-left adjoints $\hat{\Dual}_{\SS}x^*(\Dual_{\SS_X}-)\colon \Spaf \to \Pro(\Spaf)$ induces a hermitian embedding
\[
(\Spaf,\QF^{\vis}_{\xi}) \to \lim_X \Pro(\Spaf,\QF^{\vis}_{\xi\otimes\nu_x^{-1}}) ,
\]
given on the level of underlying objects by
\begin{equation}
\label{equation:ind-to-pro}%
\vphi \mapsto \hat{\Dual}^{\fib}\Dual_{\SS_X}(\vphi) 
\end{equation}
where $\hat{\Dual}^{\fib} = \lim_X \hat{\Dual}_{\SS}\colon (\Ind(\Spaf)^X)\op \to \Pro(\Spaf)^X$ is given by applying $\hat{\Dual}_\SS$ fibrewise.
In particular, we obtain an equivalence
\[
\QF^{\vis}_{\xi}(\Dual_{\SS_X}\vphi) \simeq \lim_{x \in X} \QF^{\uni}(\Dual_{\SS_X}\vphi(x))\otimes \xi_x \otimes \nu_x^{-1} \simeq \colim_{x \in X}\QF^{\uni}(\Dual_{\SS}\vphi(x))\otimes \xi_x,
\]
a fact that was also observed in Remark~\refone{remark:useful-formula} (which, as written, also applies when $X$ is not finite, but reduces to the this identification when $X$ is finite). The construction of the visible signature in \S\refone{subsection:parametrised-spectra} then proceeds by noting that $\nu_X = \Dual_{\SS_X}\SS_X$, and so by the above identification we have
\[
\QF^{\vis}_{\xi}(\nu) = \QF^{\vis}_{\xi}(\Dual_{\SS_X}\SS_X) = \colim_{X}\QF^{\uni}(\SS)\otimes \xi ,
\]
so that the universal Poincaré form $\SS \to \QF^{\uni}(\SS)$ induces a map
\[
\colim_X \xi \to \QF^{\vis}_{\xi}(\nu) .
\]
Through this map for $\xi=\nu_X$, the Pontryagin-Thom collapse class $\SS\to \colim_X \nu_X$ provides us with the desired Poincaré form $q_{\mathrm{sig}}$ on $\nu_X$ for the $\nu_X$-twisted visible functor $\QF^\vis_{\nu_X}$. The Poincaré object $(\nu_X,q_{\mathrm{sig}}) \in \Poinc(\Spaf_X,\QF^{\vis}_{\nu_X})$ is the visible signature as constructed in loc.\ cit.

\begin{proposition}
\label{proposition:visible-signature-revisited}%
The visible signature $\vsign(X)$ of Definition~\reffour{definition:visible-signature} coincides with the visible signature $(\nu_X,q_{\mathrm{sig}})$ as defined in \S\refone{subsection:parametrised-spectra}.
In particular, under the embedding $e_{\Pro}$ the Poincaré object $(\nu_X,q_{\mathrm{sig}})$ maps to the constant diagram $\SS_X\colon X \to \Pro(\Spaf)$ with value $\SS$, with the Poincaré form $q \in (\QF^{\uni})^X(\SS_X) = \lim_X \QF^{\uni}(\SS)$ corresponding to the composite $\SS \to \QF^{\uni}(\SS) \to \lim_X \QF^{\uni}(\SS)$, where the first map is the universal Poincaré form.
\end{proposition}

\begin{proof}
By construction, $e_{\Pro}(\vsign(X))$ corresponds to the constant diagram $\SS_X\colon X \to \Pro(\Spaf)$ with value $\SS$ with the constant universal Poincaré form $\ovl{q}^{\uni} \in \lim_X \QF^{\uni}(\SS)$. Since $e_{\Pro}$ is a hermitian embedding by Lemma~\reffour{lemma:e_Pro_hermitian} it will suffice to show that this description also holds for the image of $(\nu_X,q_{\mathrm{sig}})$.
Now at the level of objects we have
\[
e_{\Pro}(\nu_X) = \SS_X
\]
by~\eqreffour{equation:ind-to-pro},
and the hermitian part of $e_{\Pro}$ at $\nu_X$ then implements the equivalence
\[
\QF^{\vis}_{\nu_X}(\nu_X) \simeq \colim_X\QF^{\uni}(\SS)\otimes \nu_X \xrightarrow{\simeq} \lim_X\QF^{\uni}(\SS) .
\]
It then suffices to recall that the constant class $\SS\to \lim_X \SS$ corresponds under this identification to the Pontryagin-Thom collapse class $\SS\to \colim_X \nu_X$ (by definition of the latter).
\end{proof}

\section{Form parameters revisited}
\label{section:appendix-form-parameters}%

The aim of this appendix is to give an explicit description of the symmetric monoidal structure on the category $\FPMp$ of Poincar\'e form parameters obtained in Proposition \reffour{proposition:sym_mon_structure_on_cat_of_form_parameters}, and of the resulting algebras and modules in this category, in terms of the classical definition of form parameters of Bak~\cite{bak-form-parameter} and  Schlichting~\cite{SchlichtinghigherI}.

\subsection{The box product of Mackey functors}

Recall that the $\infty$-category $\Sp^{\gCt}$ admits a $t$-structure which is pointwise if we identify $\Sp^{\gCt}$ with the category of spectral $\Ct$-Mackey functors, and that the heart of this t-structure is the category of  $\Ct$-Mackey functors valued in abelian groups. The symmetric monoidal structure on $\Sp^{\gCt}_{\geq 0}$ given by Day-convolution then localises onto the heart to a symmetric monoidal structure again given by the Day-convolution of $\Ct$-Mackey functors valued in abelian groups, which is usually called the box-product. 

Recall that a $\Ct$-Mackey functor can be explicitly described as a pair of abelian groups $A$ and $B$, an involution $\sigma$ on $A$, and group homomorphisms
\[
A_{\Ct}\xrightarrow{\tau}B\xrightarrow{\rho}A^{\Ct}
\]
whose composite is the map induced by $\id_A+\sigma$. The map $\rho$ is called the restriction and $\tau$ the transfer. We call $A$ the underlying abelian group and $B$ the fixed-point datum of the Mackey functor. The Mackey functors form the objects of a category that we denote by $\Mack$, where the morphisms are the pairs of group homomorphisms which commute with the involutions, transfers and restrictions.  Under this description of  $\Ct$-Mackey functors, the Day-convolution product takes the following form, see e.g. \cite{Bouc}*{Proposition 1.5.1}.

\begin{definition}
\label{definition:box-mackey}%
The box-product of two $\Ct$-Mackey functors is given by
\[
\big(A_{\Ct}\xrightarrow{\tau}B\xrightarrow{\rho}A^{\Ct}\big)
\boxutens
\big(A'_{\Ct}\xrightarrow{\tau'}B'\xrightarrow{\rho'}(A')^{\Ct}\big)
=\big((A\utens A')_{\Ct}\xrightarrow{}(A\utens A'\oplus B\utens B')/\!\!\sim\to(A\utens A')^{\Ct}\big)
\]
where the quotient is by the subgroup generated by:
\begin{enumerate}
\item $a\otimes a'-\sigma(a)\otimes \sigma'(a')$ for all $a\in A, a'\in A'$,
\item  $b\otimes \tau(a')-\rho(b)\otimes a'$, for all $a'\in A', b\in B$,
\item $\tau(a)\otimes b'-a\otimes \rho(b')$, for all $a\in A, b'\in B'$.
\end{enumerate}
The abelian group $A\utens A'$ has the involution $\sigma\utens\sigma'$,
the restriction map sends $[a\otimes a']$ to $a\otimes a'+\sigma(a)\otimes \sigma'(a')$ and $[b\otimes b']$ to $\rho(b)\otimes \rho(b')$, and the transfer sends $a\otimes a'$ to the equivalence class $[a\otimes a']$.
The unit of the resulting symmetric monoidal structure is the Burnside Mackey functor
\[
\mathbb{A}=\big(\mathbb{Z}\xrightarrow{\tau}\mathbb{A}(\Ct)\xrightarrow{\rho}\mathbb{Z}\big)
\]
where $\mathbb{Z}$ has the trivial involution, and $\mathbb{A}(\Ct)$ is the Burnside ring of $\Ct$ defined as the free abelian group on the (two) isomorphism classes of transitive $\Ct$-sets. The restriction sends such an isomorphism class to its number of elements, and the transfer sends $1$ to the class of the free transitive $\Ct$-set $\Ct$.
\end{definition}

\begin{remark}
\label{remark:green}%
The category of monoids in $\Mack$ with respect to the box-product is equivalent to the category of the so called Green functors, see e.g. \cite{Bouc}*{\S2.3}. A Green functor consists of a Mackey functor as above where $A$ and $B$ are equipped with a ring structure, such that $\sigma$ and $\rho$ are ring homomorphisms, and such that the transfer $\tau$ satisfies the Frobenius reciprocity identity
\[
\tau(\rho(c)\cdot a\cdot \rho(b))=c\cdot \tau(a)\cdot b
\]
for all $a\in A$ and $b,c\in B$. %
Similarly, a (left) module over such a Green functor consists of a Mackey functor 
\[
\big(N_{\Ct}\xrightarrow{\tau}P\xrightarrow{\rho}N^{\Ct}\big)
\]
with an $A$-module structure on $N$ and a $B$-module structure on $P$ such that
\begin{enumerate}
\item
$\sigma(a\cdot n)=\sigma(a)\cdot\sigma(n)$.
\item
$\rho(b\cdot p)=\rho(b)\cdot\rho(p)$
\item
$\tau(\rho(b)\cdot n)=b\cdot\tau(n)$
\item
$\tau(b\cdot \rho(n))=\tau(b)\cdot n$
\end{enumerate}
for all $a\in A$, $b\in B$, $n\in N$ and $p\in P$.
\end{remark}

\begin{example}
\label{example:mackey}%
Given a ring $R$, there is a monoid $\Nmack(R)$ in $\Mack$ defined by
\[
\Nmack(R)=\big((R\utens R)_{\Ct}\xrightarrow{\tau} W(R)\xrightarrow{\rho}(R\utens R)^{\Ct}\big)
\]
where the involution on $R\utens R$ flips the tensor factors, and $W(R)$ is the set $(R\times (R\utens R)_{\Ct})$ with the ring structure defined by
\begin{align*}
(r,x)+(r',x')&=(r+r',x+x'-[r\otimes r'])
\\
(r,[a\otimes b])\cdot(r',[a'\otimes b'])&=(r\cdot r',[r\otimes r\cdot  a'\otimes b']+[a\otimes b\cdot r'\otimes r']+[a\otimes b\cdot a'\otimes b'+a\otimes b\cdot b'\otimes a']).
\end{align*}
The transfer and restriction are defined by $\tau(x)=(0,x)$ and $\rho(r,[a\otimes b])=r\otimes r+a\otimes b+b\otimes a$. In the terminology of  \cites{DKNP1,DKNP2}, $W(R)$ is the ring of ``$2$-truncated $2$-typical Witt vectors of $\ZZ$ with coefficients in $R$''. Moreover, there is an isomorphism of monoids
\[
\Nmack(R)\cong \pi^{\mack}_0NR
\]
in $\Mack$, see \cite{DKNP2}*{Theorem C}.

In fact, $\Nmack$ arises from a symmetric monoidal functor $\Nmack\colon\Ab\to\Mack$. As such, there are canonical isomorphisms of monoids $\Nmack(R\utens R')\cong \Nmack(R)\boxutens \Nmack(R')$ (see \cite{Read}*{Proposition 4.27}) and $\Nmack(\mathbb{Z})\cong\mathbb{A}$, where the last one is induced by the isomorphism
\[
\mathbb{A}(\Ct)\stackrel{\cong}{\longrightarrow} W(\mathbb{Z})
\]
that takes the one-elements $\Ct$-set $1$ to $(1,0)$, and the free transitive  $\Ct$-set $\Ct$ to $(0,1)$.
\end{example}

\subsection{Rings and modules in form parameters}

Recall from Proposition \reffour{proposition:sym_mon_structure_on_cat_of_form_parameters} and its proof that in our framework a form parameter consists of a ring $R$ and a $\Ct$-Mackey functor $M$ which is a module over $\Nmack(R) = \pi^{\mack}_0NR$. %
Under the characterisation of $\pi^{\mack}_0NR$ and that of modules in $\Mack$ from Remark \reffour{remark:green}, one can unravel this definition
to the following, classical definition from \cite{SchlichtinghigherI} (see also Remark \refone{remark:form-parameters-and-modules-with-genuine-involution}).

\begin{definition}
A form parameter over a discrete ring $R$ consists of a $\Ct$-Mackey functor
\[
M_{\Ct}\xrightarrow{\tau}Q\xrightarrow{\rho}M^{\Ct}
\]
together with the structure of an $R\utens R$-module with involution on $(M,\sigma)$, an action on $Q$ by the multiplicative monoid of $R$ through group homomorphisms, such that $\rho$ and $\tau$ are $R$-equivariant and
\[
(r+s)x=rx+sx+\tau((r\otimes s)\rho(x))
\]
for all $r,s\in R$ and $x\in Q$. We denote such a form parameter by $\fpm$, and typically display it by its underlying Mackey functor as in the diagram above.
We say that a form parameter is Poincar\'e if $M$ is projective as an $R$-module and the map
\[
R\longrightarrow \hom_R(M,M)
\]
of Definition \refone{remark:invertible-criterion} is an isomorphism (so that $M$ is invertible).
\end{definition}

Here, the multiplicative action of $R$ on $Q$ corresponds to the action of the elements of the form $(r,0)\in W(R)$ under the $W(R)$-module structure on $Q$, and by the definition of the sum in $W(R)$ we find the relation
\begin{align*}
(r+s,0)\cdot x&=((r,0)+(s,0)+(0,[r\otimes s])))\cdot x=((r,0)+(s,0)+(0,[r\otimes s])))\cdot x
\\&=(r,0)\cdot x+(s,0)\cdot x+\tau([r\otimes s])\cdot x=(r,0)\cdot x+(s,0)\cdot x+\tau(r\otimes s\rho(x)).
\end{align*}

\begin{definition}
We let $\FPM$ be the category of form parameters, whose objects are pairs $(R,\fpm)$ of a ring $R$ and a form parameter $\fpm$ over $R$. A morphism from  $(R,\fpm)$ to  $(R',\fpm')$ consists of a ring homomorphism $f\colon R\to R'$ and a morphism of form parameters $g\colon \lambda\to f^\ast \fpm'$ over $R$, which is defined as a pair of morphisms $M\to M'$ and $Q\to Q'$ which commute with $\sigma$, $\rho$ and $\tau$, and with the respective actions of $R\utens R$ and $R$. 

We let $\FPMp$ be the subcategory of Poincar\'e form parameters, where the morphisms are those such that the underlying map of modules induces an isomorphism  $R'\utens M\cong M'$ (cf. Lemma~\refone{lemma:Poincresscalars}).
\end{definition}

By the discussion at the last part of \S\reffour{subsection:multiplicative-GW}, Poincaré form parameters on $R$ correspond to Poincaré structures $\QF$ on $\Dperf(R)$ such that $\Dual_{\QF}$ preserves $\Proj(R)$ and $\QF$ takes discrete values on $\Proj(R)$. 
The symmetric monoidal structure on $\FPMp$ from Proposition \reffour{proposition:sym_mon_structure_on_cat_of_form_parameters} is by construction that of algebras over $\pi^{\mack}_0N$. We thus obtain the following mild extension of the tensor product of form rings of  \cite{SchlichtinghigherI}*{\S 4} (see also Example \reffour{example:formrings} below).

\begin{proposition}
\label{prop:tensormackeyvstensorfpm}%
 The box product of Mackey functor lifts to a symmetric monoidal structure on the category $\FPM$ of form parameters, whose tensor product sends an $N^{\mack}(R)$-module $X$ and an $N^{\mack}(S)$-module $Y$ in $\Mack$ to $X\boxutens Y$, which we regard as a $N^{\mack}(R\utens S)$-module via the monoidal structure of $N^{\mack}$.

This symmetric monoidal structure restricts to that of Proposition \reffour{proposition:sym_mon_structure_on_cat_of_form_parameters} on the subcategory $\FPMp$ of Poincar\'e form parameters.

\qed
\end{proposition}

We can now conveniently unpack the definition of a monoid and of a module in  $\FPM$ as follows.

\begin{lemma}%
\label{lemma:monoidsFPM}%
A monoid in $\FPM$ consists of a commutative ring $R$, a Green functor $
X=\big(M_{\Ct}\xrightarrow{\tau}Q\xrightarrow{\rho}M^{\Ct}\big)
$  (i.e. a monoid in $\Mack$)
where $M$ has the structure of an $R\utens R$-algebra with involution, and a multiplicative unital map $n\colon R\to Q$ which satisfies 
\begin{enumerate}
\item $n(0)=0$,
\item $\rho n(r)=r\otimes r\cdot 1$
\item $n(r+s)=n(r)+n(s)+\tau(r\otimes s\cdot 1)$
\end{enumerate}
for all $r,s\in R$, where $1\in M$ is the multiplicative unit. Then $X$ is a form parameter with the action $r\cdot q=n(r)q$ for $r\in R$ and $q\in Q$. 

A monoid in $\FPMp$ is a monoid in $\FPM$ whose underlying form parameter is Poincar\'e and such that the canonical map of $R$-modules $R\to M$ induced by the unit of $M$ is an isomorphism.
\end{lemma}

\begin{proof}
By the definition of the tensor product of form parameters, if $X$ is a form parameter over $R$ equipped with the structure of a monoid, then in particular $R$ is a monoid in rings, that is a commutative ring. Then the remaining structure of the monoid is that of an $\Nmack(R)$-algebra, that is the structure of a Green functor on $X$ and a morphism of Green functors $\Nmack(R)\to X$. On underlying rings this map gives $M$ the structure of an $R\utens R$-algebra with involution. It remains to see that the ring homomorphim $g\colon W(R)\to Q$ from the map of Green functors is equivalent to a map $f$ with the properties above. 
Since the elements $(0,x)\in W(R)$ are in the image of the transfer, $g$ is determined by its values on the elements $g(r,0)$. This is the multiplicative map $n\colon R\to Q$ given by $n(r)=g(r,0)$, and the last identity follows from the additivity of $g$, since
\[
n(r+s)=g(r+s,0)=g((r,0)+(s,0)+(0,r\otimes s))=n(r)+n(s)+\tau(r\otimes s\cdot 1)
\]
for all $r,s\in R$. The second identity holds by the compatibility with the restriction maps, and $n$ preserves $0$ since $g$ does.
The condition that $R\cong M$ for monoids in $\FPMp$ is exactly the condition for the unit map $\mathbb{A}\to X$ to lie in $\FPMp$.
\end{proof}

\begin{example}%
\label{example:formrings}%
A commutative form ring in the sense of Schlichting \cite{SchlichtinghigherI}*{Definitions 4.1, 4.6} is exactly a commutative monoid in $\FPMp$, where one defines $n\colon R\to Q$ by $n(r)=r\cdot 1$, where $1$ is the ring unit of $Q$.
\end{example}

\begin{example} 
Recall that a Tambara functor for the group $\Ct$ is a Mackey functor
\[
A_{\Ct}\xrightarrow{\tau}B\xrightarrow{\rho}A^{\Ct}
\]
together with the structure of commutative rings on $A$ and $B$ making this into a Green functor, and a multiplicative unital map $n\colon A\to B$ which satisfies exactly the conditions of Lemma \reffour{lemma:monoidsFPM}. Thus, such a Tambara functor defines a commutative monoid in $\FPM$ with underlying ring $A$.
\end{example}

\begin{lemma}%
\label{lemma:modulesFPM}%
Let $X$ be a monoid in $\FPM$ as in Lemma \reffour{lemma:monoidsFPM}.
A module over $X$ consists of an $R$-algebra $S$ and a module 
$
Y=\big(N_{\Ct}\xrightarrow{\tau}P\xrightarrow{\rho}N^{\Ct}\big)
$
over the Green functor $X$ as in Remark \reffour{remark:green} which is additionally a form parameter over $S$, where $N$ is a module over $M$ in $R\utens R$-modules, and which satisfies the following compatibility condition with the form parameter structure:
\[
(n(r)q)(s\cdot p)=(rs)\cdot (qp)
\]
for all $r\in R, s\in S, q\in Q$ and $p\in P$, where $r\cdot p$ denotes the action of the form parameter $Y$.

If $X$ is a monoid in $\FPMp$, a module in $\FPMp$ over $X$ is a module in $\FPM$ over $X$ whose underlying form parameter is in $\FPMp$. 
\end{lemma}
\begin{proof}
We start by observing the the module action is given by a morphism of form parameters
\[
(R\utens S,X\boxutens Y)=(R,X)\utens (S,Y)\longrightarrow (S,Y).
\]
The map in the first component givens $S$ the structure of an $R$-algebra. The map in the second component is a map $X\boxutens Y\to Y$ of $\Nmack(R\utens S)$-modules, where $Y$ is regarded as a $\Nmack(R\utens S)$-module via the action of $R$ on $S$. This last linearity gives the compatibility of the module structure with the quadratic action, again using that the action of $W(R\otimes S)$ is determined by that of the elements of the form $(r\otimes s,0)$.

If $X$ is a monoid in $\FPMp$, so that we may assume that $M=R$, the condition for the action map of $X$ on $Y$ to lie in $\FPMp$ is that it induces an isomorphism
\[
S\otimes_{R\utens S}(R\utens N)\to N
\]
which is automatic since $N$ is a projective $S$-module. 
\end{proof}

\begin{example}
The tensor unit $\mathbb{A}$ is clearly a module over $\Nmack(\mathbb{Z})$ (it is isomorphic to $\Nmack(\mathbb{Z})$ by the monoidality of $\Nmack$). The multiplicative action of $\mathbb{Z}$ on $\mathbb{A}(\Ct)$ equipping it with the structure of $\Nmack(\mathbb{Z})$-module is given by
\[
a\cdot x=(a+\frac{a^2-a}{2}\Ct) x
\]
where $x$ is the isomorphism class of a $\Ct$-set and $a\in\mathbb{Z}$. This formula comes from the orbit decomposition of the set $k\times k$ with $k^2$-elements, $k$ a natural number, equipped with the $\Ct$-action that flips the product factors.

Since $\mathbb{A}$ is the tensor unit it is also a monoid in $\FPMp$. The ring structure on $\mathbb{A}(\Ct)$ is the usual one of the Burnside ring defined from the cartesian product of $\Ct$-sets, and the map $n\colon \mathbb{Z}\to \mathbb{A}(\Ct)$ is
\[
n(a)=a+\frac{a^2-a}{2}\Ct.
\]
Moreove, again because $\mathbb{A}$ is the tensor unit, any form parameter is uniquely a module over $\mathbb{A}$.
\end{example}

\begin{example}
\label{example:fixeptmonoid}%
Let $R$ be a ring and $M$ an $R\utens R$-module with involution $\sigma$. The fixed point Mackey functor
\[
M^{\sym}=\big(M_{\Ct}\xrightarrow{\id+\sigma}M^{\Ct}\xrightarrow{\id}M^{\Ct}\big)
\]
is a form parameter over $R$, called the \defi{symmetric} form parameter, where the multiplicative action of $R$ on $M^{\Ct}$ is defined by $r\cdot m=(r\otimes r) m$, for all $r\in R$ and $m\in M^{\Ct}$. The corresponding Poincaré structure on $\Dperf(R)$ is the genuine symmetric structure.

If $R$ is moreover commutative and $M$ is a commutative $R\utens R$-algebra with involution, then $M^{\sym}$ is a commutative monoid in $\FPM$ where the multiplication on $M^{\Ct}$ is restricted from that of $M$, and the map $n\colon R\to M^{\Ct}$ is defined by
\[
n(r)=r\otimes r\cdot 1,
\]
where $1$ is the unit of $M$. In particular if $R$ is a commutative ring with involution then $R^{\sym}$ is a commutative monoid in $\FPMp$.
\end{example}

\begin{example}
Let $S$ be a ring and $N$ an $S\utens S$-module with involution $\sigma$. Then the orbit Mackey functor
\[
N^{\qdr}=\big(N_{\Ct}\xrightarrow{\id}N_{\Ct}\xrightarrow{\id+\sigma}N^{\Ct}\big)
\]
is a form parameter over $S$, called the \defi{quadratic} form parameter, where the multiplicative action of $S$ on $N_{\Ct}$ is $s\cdot [n]=[(s\otimes s) n]$, for all $s\in S$ and $n\in N^{\Ct}$. The corresponding Poincaré structure on $\Dperf(R)$ is the genuine quadratic structure.

If $R$ is a commutative ring and $M$ is a commutative $R\utens R$-algebra with involution, $S$ is an $R$-algebra,  and $N$ an $M$-module with involution in $R\utens R$-modules, then $N^{\sym}$ and $N^{\qdr}$ are modules over $M^{\sym}$ in $\FPM$. %
The relation for a module from Lemma \reffour{lemma:modulesFPM} for the quadratic case holds since
\[
(n(r)m)(s\cdot [x])=
(r\otimes r)m[s\otimes s x]=[(rs\otimes rs) mx]=(rs)\cdot (m[x])
\]
for $m\in M^{\Ct}, x\in N, r\in R$ and $s\in S$, and similarly for the symmetric case.
\end{example}

\begin{example}%
\label{example:form_parameter_over_Tate_Poincare}%
Let $k$ be a commutative ring. There is a commutative monoid in $\FPMp$
\[
k^{\fr}=\big(k\xrightarrow{\tau}W_{\langle 2\rangle}(k)\xrightarrow{\rho}k\big)
\]
where $W_{\langle 2\rangle}(k)$ is the ring of $2$-truncated $2$-typical Witt vectors of $k$, so that the underlying set of $W_{\langle 2\rangle}(k)$ is $k\times k$. The transfer $\tau$ is the Verschiebung which sends $a$ to $(0,a)$, $\rho$ is the Frobenius that takes $a$ to $a^2+2b$, and $n\colon k\to W_{\langle 2\rangle}(k)$ is the Teichmüller lift that sends $a\in k$ to $(a,0)$. 

We call this form parameter the \defi{Frobenius form parameter}. The corresponding Poincaré structure is the \defi{genuine Frobenius structure} $\QF^{\g\fr}_k$, which coincides with the truncation $\tau_{\leq 0}\QF^{\tate}_k$ of the Tate Poincaré structure 
on $\Dperf(k)$, associated to the module with genuine involution $(k,k,k\to k^{\tC})$ whose structure map is the Tate-valued Frobenius. 
In the case of $k=\ZZ$ the ring $W_{\langle 2\rangle}(\ZZ)$ is isomorphic to the Burside ring and the Frobenius form parameter recovers the Burnside form parameter of Example~\refthree{example:dotto-ogle}.

This construction can also be adapted to the case where $k$ is a commutative ring with involution $\tau$ by replacing $W_{\langle 2\rangle}(k)$ with the ring $W_{\langle 2\rangle}(k,\tau)$ whose underlying set of elements is $k \times k_{\Ct}$ and ring structure defined similarly to the one on $W(k)$ in Example~\reffour{example:mackey}.
\end{example}

Let $k$ be a commutative ring. By Example \reffour{example:fixeptmonoid} the symmetric form parameter $k^{\sym}=(k\xrightarrow{2}k\xrightarrow{\id}k)$ is canonically a monoid in $\FPMp$.

\begin{lemma}
\label{lemma:moduleskg}%
A module over $k^{\sym}$ in $\FPMp$ consists of a form parameter
$
Y=\big(N_{\Ct}\xrightarrow{\tau}P\xrightarrow{\rho}N^{\Ct}\big)
$
over a ring $S$ together with a $k$-algebra structure on $S$ and a $k$-module structure on $P$, such that:
\begin{enumerate}
\item the $S\otimes S$-module structure on $N$ descends to a $S\otimes_k S$-module structure (so that $N$ is in particular a $k$-module and $\sigma$ is $k$-linear),
\item  $\tau$ and $\rho$ are $k$-linear,
\item the action of $S$ on $P$ of the form parameter is $k$-quadratic, in the sense that
\[
(as)\cdot p=a^2(s\cdot p)
 \]
for all $a\in k, s\in S$ and $p\in P$ where $a^2$ on the right acts via the $k$-module structure on $P$.
\end{enumerate}
\end{lemma}
 
\begin{proof}
From the characterisation of modules from Lemma \reffour{lemma:modulesFPM} we see that a $k^{\sym}$-module in $\FPM$ consists of a form parameter $Y$ as above over $S$ where $S$ is a $k$-algebra, $N$ is a $S\otimes_k S$-module, which is a module over the Green functor $k^{\sym}$ and where the action of $S$ on $P$ of the form parameter is $k$-quadratic. Moreover, by Remark \reffour{remark:green} the structure of a module over the Green functor $k^{\sym}$ is the same as a $k$-action on $P$ where $\tau$ and $\rho$ are $k$-linear (for $\tau$ this is because the restriction map of $k^{\sym}$ is surjective) and such that $\tau\rho =2$. We therefore only need to verify this last identity. By combining the formula for the action of the form parameter on a sum with the fact taht this action is $k$-quadratic we have that
 
\[
4(s\cdot p)=(2s)\cdot p=2(s\cdot p)+\tau((s\otimes s)\rho (p)),
\]
 which for $s=1$ gives $\tau(\rho (p))=2p$ for every $p\in P$.
\end{proof}

\begin{remark}
\label{remark:criterion-genuine-symmetric}%
By Lemma \reffour{lemma:modulesFPM}, a $\mathbb{Z}^{\sym}$-module structure on a form parameter $
Y=\big(N_{\Ct}\xrightarrow{\tau}P\xrightarrow{\rho}N^{\Ct}\big)
$
over a ring $S$ exists (and it is in this case unique) if and only if
\[
(as)\cdot p=a^2(s\cdot p)
\]
for all $a\in\mathbb{Z}, s\in S$ and $p\in P$. Moreover, this last property is in fact equivalent to the condition that $\tau\rho=2$. We saw in the proof of Lemma \reffour{lemma:moduleskg} that it is sufficient. Conversely, assume that $\tau\rho=2$. Since $\rho$ is $S$-equivariant we have that
 
\[
2(s\cdot p)=\tau\rho(s\cdot p)=\tau((s\otimes s) \rho (p)),
\]
 and by the property of the action of the form parameter on a sum that
\[
(2s)\cdot p=2(s\cdot p)+\tau((s\otimes s) \rho (p))=4(s\cdot p)=(2)^2(s\cdot p).
\]
By applying this argument inductively on $a$ we can prove that $(as)\cdot p=a^2(s\cdot p)$ for every $a\geq 2$. The formula for $a=0,1$ is obvious. Since
 
\[
0=(as+(-as))\cdot p=(as)\cdot p+(-as)\cdot p+\tau(as\otimes (-as)\rho(p))=(as)\cdot p+(-as)\cdot p-2((as)\cdot p)
\]
 we see that $(-as)\cdot p=(as)\cdot p$ and the formula also holds for negative $a\in\mathbb{Z}$.
\end{remark}

\begin{lemma}
A module over $k^{\g\fr}$ in $\FPM$ consists of a form parameter
$
Y=\big(N_{\Ct}\xrightarrow{\tau}P\xrightarrow{\rho}N^{\Ct}\big)
$
over a ring $S$ together with a $k$-algebra structure on $S$, such that the $S\otimes S$-module structure on $N$ descends to a $S\otimes_k S$-module structure.
\end{lemma}

 \begin{proof}
 An action of $W_{\langle 2\rangle}(k)$ on $P$ giving $Y$ the structure of a module is  uniquely determined by the $S$-action of the form parameter. Indeed by setting $s=1$ and $(a,b)=(1,0)$ in the condition
\[
((c,0)(a,b))(s\cdot p)=(cs)\cdot ((a,b)p)
\]
we find that $(a,0)p=a\cdot p$, and from the Frobenius reciprocity condition that
\[
(0,b)p=\tau(b\rho(p)),
\]
for every $a,b\in k$ and $p\in P$. One can directly verify that by defining the action of $W_{\langle 2\rangle}(k)$ on $P$ this way, all the conditions of a module are satisfied.
 \end{proof}

\subsection{Form parameters of quadratic and symmetric types}

The next definition captures the two flavours of form parameter originally considered by Bak in \cite{bak-form-parameter}.
 
 \begin{definition}
We call a form parameter $\big(N_{\Ct}\xrightarrow{\tau}P\xrightarrow{\rho}N^{\Ct}\big)$ over a ring $S$ of symmetric type if $\rho$ is injective, and of quadratic type if $\tau$ is surjective. We denote the respective full subcategories of $\FPM$ by $\FPM_\s$ and $\FPM_\qdr$.
 \end{definition}

\begin{example}\
\begin{enumerate}
\item
Any form parameter of symmetric type is uniquely a module over $\mathbb{Z}^{\sym}$. Indeed, since $\rho$ is injective one can check that the last condition of Lemma \reffour{lemma:moduleskg} holds after applying $\rho$, which is
\[
\rho\tau\rho=(\id+\sigma)\rho=\rho+\sigma\rho=2\rho=\rho 2.
\]
For example the form parameter $M^{\sym}$ associated to an $R$-module with involution is symmetric, and thus a module over $\mathbb{Z}^{\sym}$.
\item
Any form parameter of quadratic type is also uniquely a module over $\mathbb{Z}^{\sym}$. Since $\tau$ is surjective one can choose, for every $p\in P$, an element $n\in N_{\Ct}$ such that $\tau(n)=p$. Thus
\[
\tau\rho(p)=\tau\rho\tau(n)=\tau(n+\sigma(n))=\tau(n)+\tau\sigma(n)=2\tau(n)=2p.
\]
For example the form parameter $M^{\qdr}$ associated to an $R$-module with involution is quadratic, and thus a module over $\mathbb{Z}^{\sym}$. 
\end{enumerate}
\end{example}

We conclude this section with a comparison with the tensor product of form parameters of symmetric type introduced by Bak. Given any module with involution $M$, we write $-M$ for the same module with the involution changed by a sign.

\begin{proposition}%
\label{proposition:symmetric_form_parameters}%
\
\begin{enumerate}
 \item 
The subcategory $\FPM_{\s}$ is a multiplicative Bousfield localisation of $\FPM$, where the left adjoint $s$ of the inclusion is given by
\[
s\big(M_{\Ct}\xrightarrow{\tau}Q\xrightarrow{\rho}M^{\Ct}\big)=\big(M_{\Ct}\xrightarrow{\id+\sigma}Im(\rho)\xrightarrow{}M^{\Ct}\big)
\]
where both form parameters are over the same ring $R$.
In particular $\FPM_{\s}$ inherits a symmetric monoidal structure with unit $\s(\mathbb{A})= \ZZ^{\sym}$. 
\item The action of $\FPM$ onto itself restricts to an action on $\FPM_{\qdr}$, which in turn extends (uniquely) to an action of $\FPM_{\s}$ on $\FPM_{\qdr}$.
\item Tensoring with the form parameter $(-\mathbb{Z})^{\sym}$ gives a functor $(-\mathbb{Z})^{\sym}\boxutens(-)\colon \FPM\to\FPM_{\qdr}$ which moreover identifies with the dual form parameter
\[
 (-\mathbb{Z})^{\sym}\boxutens\fpm=\big((-M)_{\Ct}\xrightarrow{\pr}M/Im(\rho)\xrightarrow{\id-\sigma}(-M)^{\Ct}\big)=\check{\fpm}
\]
from \cite{SchlichtinghigherI} (see also Proposition \refone{proposition:form-parameter-shift}).
\end{enumerate}
\end{proposition}

\begin{proof}
For the first part, it is easy to see that $s$ is left-adjoint to the inclusion. We need to show that tensoring with a form parameter $X=(M_{\Ct}\xrightarrow{\tau}Q\xrightarrow{\rho}M^{\Ct})$ over $R$ preserves $s$-equivalences. These are the morphisms of form parameters $f\colon Y\to Y'$ whose underlying map of rings is an isomorphism, and whose underlying map of modules with involution $g\colon N\to N'$ is an isomorphism and it restricts to a surjection $Im(\rho)\to Im(\rho')$. One can easily see that these are preserved by $X\boxutens (-)$ using that, by the definition of the box product, the image of the restriction map of $X\boxutens Y'$ is generated by the elements of the form $m\otimes n'+\sigma(m)\otimes\sigma(n')$ and $\rho(q)\otimes\rho(p')$ for $m\in M, n'\in N', q\in Q$ and $p'\in P'$. The statement about the unit is clear, since the restriction map $\rho\colon\mathbb{A}(\Ct)\to\mathbb{Z}$ of $\mathbb{A}$ is surjective.

For the second part, let us show that if $X$ is any form parameter and $Y$ is one of quadratic type, then $X\boxutens Y$ is quadratic. By the definition of the transfer of the box product, we need to show that any element of the form $q\otimes p$ for $q\in Q$ and $p\in P$ is equivalent to one of the form $m\otimes n$ for some $m\in M$ and $n\in N$, under the defining equivalence relation of the box product. Since the transfer $\tau$ of $Y$ is assumed to be surjective, there is $n\in N$ with $\tau(n)=p$ and therefore
\[
q\otimes p=q\otimes \tau(n)\sim \rho(q)\otimes n
\]
as we wanted. For the last part, it is sufficient to show that $s$ preserves form parameters of quadratic type, that is that if the transfer $\tau$ of $Y$ is surjective then $\id+\sigma$ surjects onto the image of $\rho$. Indeed, given $p\in P$ we can write $p=\tau(n)$ for some $n\in N$, and therefore
\[
\rho(p)=\rho \tau(n)=n+\sigma(n).
\]

For the final claim, since $(-\mathbb{Z})^{\Ct}=(-\mathbb{Z})_{\Ct}=0$ we find that one of the summands in the definition of the box product $(-\mathbb{Z})^{\sym}\boxutens\fpm$ vanishes, and in the other summand $\mathbb{Z}\otimes M=M$ we have the relation $a\rho(q)\simeq \tau(a)q=0$ for $a\in\mathbb{Z}$ and $q\in Q$. From this one can easily see that $(-\mathbb{Z})^{\sym}\boxutens\fpm\cong \check{\fpm}$.
\end{proof}

\begin{remark}[Comparison to form parameters in the sense of Bak]
Given a ring $R$ and an $(R\otimes_{\ZZ} R)$-module with involution $(M,\sigma)$, 
forms parameters of symmetric type over $R$ extending $(M,\sigma)$ 
are classified by subgroups 
\[
\Lam \subseteq \widehat{\mathrm{H}}^0(\Ct;M) = \mathrm{coker}[\mathrm{N}\colon M_{\Ct} \to M^{\Ct}]
\]
which are stable under the operation $[m] \mapsto [(r\otimes r)(m)]$ for $m \in M^{\Ct}$. Here, the correspondence between these two types of data is given by $\Lambda =  \coker[M_{\Ct}\to Q]$ and $Q = M^{\Ct} \times_{\widehat{\mathrm{H}}^0(\Ct;M)} \Lam$.
 
 Such a triple $(R, M, \Lambda)$ %
 constitutes a form parameter in the sense of Bak (where Bak only considered the case where $R$ is a ring with involution and $M=\pm R$). Under the classification, the reflection $\s\colon \FPM \to \FPM_s$ (and in fact the reflection $\s\circ \tau_{\leq 0}$) are described by sending a module with genuine involution $(M,b\colon N\to M^\tC)$ over $R$ to $(R, M, \im(\pi_0(b))$. Thus, the tensor product in $\FPM_\s$ is described by Bak's tensor product on form parameters,
\[
(R, B, \Lambda)\utens (R', B', \Lambda') = (R\utens R', B\utens B', \im(\Lambda\utens \Lambda'\to \widehat{\mathrm{H}}^0(\Ct, B\utens B')).
\]

We also find that Bak-style form parameter $(R, B, \Lambda)$ classify quadratic-type form parameters just as well, where the duality however needs to be twisted by a sign. The passage from the symmetric-type parameter to the corresponding quadratic-type parameter corresponds to the passage to the ``dual form parameter'' from \S\refone{subsection:discrete-rings}. Under these correspondences, the action $\FPM_\s\times \FPM_\qdr\to \FPM_\qdr$ is again described by Bak's tensor product of form parameters. This explains the tensor product written down by Karoubi \cite{karoubi-periodicity}*{4.3} in our more categorical terms. 
\end{remark}

\end{appendices}

\end{document}